\documentclass[a4]{amsbook}
\usepackage[english]{babel}
\usepackage{amsfonts}
\usepackage{epsfig}
\usepackage{graphicx}
\usepackage{enumerate}
\title{A Fundamental Domain for $V_{3}$}
\author{Mary Rees}
\address{Department of Mathematical Sciences,\\ 
University of Liverpool,\\
Peach St.,\\
Liverpool L69 7ZL\\
U.K..}
\email{maryrees@liv.ac.uk}
\urladdr{http://www.liv.ac.uk/~maryrees/maryrees.homepage.html}
\def\Box{\hbox{$\sqcap \unskip \kern -6.5pt 
\sqcup$}}
\def \Amalg{\mathbin{\raise .5pt%
	\hbox{$\scriptstyle \amalg$}}}
	
\begin{document}
	\frontmatter
		\begin{abstract}
	We describe a fundamental domain for the punctured Riemann surface $V_{3,m}$ which parametrises (up to M\"obius conjugacy) the set of  quadratic rational maps with numbered critical points, such that the first critical point has period three, and such that the second critical point is not mapped in $m$ iterates or less to the periodic orbit of the first. This gives, in turn, a description, up to topological conjugacy, of all dynamics in all type III hyperbolic components in $V_{3}$, and gives indications of a topological model for $V_{3}$, together with the hyperbolic components contained in it.
	\end{abstract}
 \maketitle

    \tableofcontents 
\newtheorem{theorem}[section]{Theorem}
\newtheorem{lemma}[section]{Lemma}
\newtheorem{corollary}[section]{Corollary}
\newenvironment{ulemma}{\par\noindent\textbf{Lemma}\,\,\em}{\rm}
\newenvironment{utheorem}{\par\noindent\textbf{Theorem}\,\,\em}{\rm}
\newenvironment{maintheorem1}{\par\noindent\textbf{Main Theorem 
(first version)}\,\,\em}{\rm}
\newenvironment{maintheorem1.5}{\par\noindent\textbf{Main Theorem 
( vague version)}\,\,\em}{\rm}
\newenvironment{maintheorem2}{\par\noindent\textbf{Main Theorem 
(final version, easy cases)}\,\,\em}{\rm}
\newenvironment{maintheorem3}{\par\noindent\textbf{Main Theorem 
(final version, hard case)}\,\,\em}{\rm}
\numberwithin{section}{chapter}
\numberwithin{equation}{section}

\mainmatter

\chapter{Introduction}\label{1}

In this paper, we give a complete description of half the 
hyperbolic components in a  certain parameter space of quadratic 
rational maps. The parameter space is known as $V_{3}$, and consists 
of all maps
$$h_{a}:z\mapsto {(z-a)(z-1)\over z^{2}},\ \ a\in \mathbb C,\ a\neq 
0.$$
This map has two critical points,
$$c_{1}=0,\ c_{2}=c_{2}(a)={2a\over a+1},$$
and $c_{1}$ is periodic of period $3$, with orbit
$$0\mapsto \infty \mapsto 1\mapsto 0.$$
There is thus one free critical point, $c_{2}$. Describing the 
hyperbolic 
components in $V_{3}$ --- and in fact only  half, or in numerical 
terms, two thirds of them -- 
might seem rather a modest project, but has, in fact, been ongoing 
for some twenty years, and has generated an extensive theory, which 
will be summarised, as far as it concerns us here, in Section 
\ref{4}. 
It is essentially because of technical difficulties in the 
development 
of the theory that only half the hyperbolic components in $V_{3}$ are 
described. A conjectural picture for the remaining half is  
not too hard to obtain. Remarks are made on this at various points, 
but no explicit description is given.

One can obviously ask why $V_{3}$ has been chosen. The space $V_{n}$ 
is the quadratic rational maps with one critical point of period $n$, 
quotiented by M\"obius conjugacy. 
A quadratic rational map with one fixed critical point is M\" 
obius conjugate to a quadratic polynomial. Any quadratic polynomial 
is 
affinely conjugate to one of the form 
$$f_{c}:z\mapsto z^{2}+c$$
for some 
$c\in {\mathbb C}$, and this, thus, is the space $V_{1}$.
This parameter space must be one of the most studied in all dynamics. 
There are several reasons for this. One is that
dynamics within this parameter space is rich and varied. Another is 
that there is a pretty full description of the variation of dynamics 
within this parameter space, at least up to topological conjuacy, and 
at least conjecturally. Another is that this conjectural desciption 
of 
the dynamics, which is very detailed, but not fully complete, and 
which includes a conjecture on the topology of the 
Mandelbrot set, is  still not proved, although there has been very interesting progress recently.  The missing piece has analogues 
in the theory of Kleinian groups --- where the corresponding problem 
has been solved \cite{Min, B-C-M, R6, Bow1, Bow2, Bow3, BBES}  --- not to mention other, less 
well-understood, parameter spaces of holomorphic maps. So there are 
attractively simple questions one can ask, even in this case, which 
have very detailed, but not complete, answers, and one hopes that 
this case will shed some light on other parameter spaces of 
holomorphic 
maps.

Actually, the parameter space of quadratic polynomials differs from 
most other parameter spaces in one very fundamental respect. There is 
a natural ``base'' map in the space of quadratic polynomials, namely, 
the map $f_{0}:z\mapsto z^{2}$, and for any hyperbolic map $f$ (and 
hopefully some nonhyperbolic also) in the connectedness 
locus --- also called the Mandelbrot set --- there is an essentially 
unique description of the dynamics of $f$ in terms of $f_{0}$. This 
is because the complement of the Mandelbrot set is simply connected, and 
the Mandelbrot set itself has a natural tree-like structure. There 
appears to be essentially one natural path from $f_{0}$ to $f$ within 
the Mandelbrot set.  
This is very far from being the case in other parameter spaces of 
quadratic rational maps, although the structure of $V_{2}$ is quite 
simple, as 
we shall indicate --- but not prove --- in Section \ref{2}.

For $n\geq 3$, $V_{n}$ identifies with the family of quadratic 
rational maps 
$$f_{c,d}:z\mapsto 1+{c\over z}+{d\over z^{2}}$$
($d\neq 0$) for which the critical point $0$ is constrained to have 
period $n$. 
If $n\geq 3$, then the maps $f_{c,d}\in V_{n}$ are in one-to-one 
correspondence with M\" obius conjugacy classes of quadratic rational 
maps $f$ with named critical point $c_{1}(f)$ of period $n$, where we 
use only conjugacies which map $c_{1}(f)$ to $0$. These parameter 
spaces were the main objects of study in \cite{R1, R2, R3}. The main point  of 
\cite{R1, R2} was that it is possible to describe dynamics 
of hyperbolic maps $f_{c,d}$ in $V_{n}$ in terms of a path from some 
fixed base 
map $g_{0}$ - and, of course, in terms of that base map $g_{0}$. In 
\cite{R1} we used a polynomial (up to M\" obius conjugacy)  as base, 
and 
the resulting theorem was called the {\em{Polynomial-and-Path 
Theorem}}. It was shown that a path in parameter space from $g_{0}$ 
to $f_{c,d}$ gave rise to a 
path in the dynamical plane of $g_{0}$. So the theorem showed how to 
define $f_{c,d}$ in terms of $g_{0}$ and the path in the dynamical 
plane of $g_{0}$. The main problems with such a result were, firstly, 
identifying which paths in the dynamical plane of $g_{0}$ were 
associated with paths in $V_{n}$, and therefore, with hyperbolic 
maps 
in $V_{n}$, and secondly, to
determine when two paths gave rise to a description of the same map 
in $V_{n}$. Some progress was made in \cite{R2} in restricting the 
set of path in the dynamical plane of $g_{0}$ to a smaller set, which 
still described all hyperbolic maps in $V_{n}$. Further progress 
along 
these lines looked difficult.

In \cite{R3} there was a fundamental change of policy. The underlying 
idea in \cite{R1, R2} had been to attempt a conjectural 
description of parameter space $V_{n}$, with its hyperbolic 
components, 
as a quotient of a subset of some sort of univeral cover of the 
dynamical plane 
of some suitable base map $g_{0}$. This ``universal cover'' was an 
artificial construct. Basically, the idea was to remove all points in 
the full orbit of the period $n$ critical point of $g_{0}$, but it 
was not clear what was the right way to do this. 

In 1990 I decided to 
stop guessing on this, and just attempt to describe finitely many - 
but arbitrarily finitely many - hyperbolic components at a time, by 
looking at the universal cover of the complement, in the dynamical 
plane of a base map $g_{0}$, of just finitely many points in the full 
orbit of the period $n$ critical point. There was at least no problem 
about what this universal cover was. Immediately, a whole lot of 
structure became clear. Nevertheless, the proofs of the structure 
took a little over 10 years to write down correctly. Moreover, this 
was just a start. This was basically a topological result, describing 
the 
topology of parameter space minus finitely many hyperbolic components 
in terms of the complement of finitely many points in the dynamical 
plane of $g_{0}$. So there was no  comprehensive view of the whole of 
parameter space. In fact, technical difficultes meant that about half 
of all hyperbolic components - those of type IV - were omitted from 
consideration.
 
Hyperbolic components of quadratic rational maps come in four types. 
There is only one {\em{type I}} component, for which both critical 
points are attracted to the same attractive fixed point. This 
hyperbolic component intersects $V_{n}$ only for $n=1$, and contains 
$z\mapsto z^{2}+c$ for all large $c$. So it can be omitted from 
consideration if we are considering $V_{n}$ only for  any $n\geq 2$. 
A 
{\em{type II component of period $n$}} is one in which both critical 
points are in 
periodic components of the attractive basin of an attractive 
periodic point of 
period $n$. In this case, the critical points must be in different 
components of the attractive periodic basin. For any fixed period 
$n$, there are only finitely many type II hyperolic components for 
attractive periodic basins of period $n$. A {\em{type III component 
of preperiod $m$}} 
is one in which both critical points are in the full attractive basin 
of a single attractive periodic point, with one critical point in a 
periodic component (there must be at least one such), and the other 
in a nonperiodic component whose $m+1$'st forward iterate is 
periodic. 

A {\em{type IV}} hyperbolic component is one in 
which the critical points are in periodic components of attractive 
periodic points in distinct orbits.  The intersection of a  
type IV component with $V_{n}$
has {\em{period $m$}} if the period of the immediate attractive basin 
containing $c_{2}$ is $m$. There are finitely many type III (or type 
IV) 
components 
of preperiod (or period) $m$ intersecting $V_{n}$, for any integer 
$m$, 
but infinitely many 
altogether.  We shall see that the number of type III components in 
$V_{3}$ of preperiod $m$ can be computed exactly, and is 
$(1+{2\over 21}).2^{m+1}+O(1)$ . This is also the number 
of type IV 
components of period $m+1$ in $V_{3}$. 

The description of the hyperbolic components in $V_{1}$, the space of 
quadratic polynomials, is part of a description which conjecturally 
gives much more information, not just combinatorial, but topological 
and possibly geometric. Exactly one hyperbolic component in $V_{1}$ --- the type I component mentioned above ---
is unbounded, and was shown by Douady and Hubbard \cite{D-H1} in the 
1980's to be 
simply connected, once $\infty $ is added. The complement is a 
closed set known as the {\em{Mandelbrot set.}}
Hyperbolic components in $V_{1}$ are in 1-1 
correspondence with subsets in the {\em{combinatorial Mandelbrot 
set}}. The Mandelbrot set $M$ is conjecturally homeomorphic to the 
combinatorial Mandelbrot set $M_{c}$. There is a natural monotone 
map $\Phi :M\to M_{c}$, which is rather easily shown to be 
continuous, in the same way as a monotone map defined on a subset of 
$\mathbb R$, with dense image in $\mathbb R$, can be 
shown to extend continuously.  The question, therefore, is whether 
$\Phi $ is injective. Injectivity at parabolic points is relatively 
easy, though far from trivial. Much deeper results derive from use of 
the {\em{Yoccoz parapuzzle}}. It is known, for example, that $\Phi $ 
is injective at all points of $M$ that are not infinitely 
renormalisable. (This is equivalent to the more usual statement that 
the Mandelbrot set is locally connected at all points that are not 
infnitely renormalisable.)  Recently, some  15 years after Yoccoz 
first demonstrated his methods, these techniques have been used 
by Kahn and Lyubich and collaborators \cite{L1, L2, K-L3, Ka, K-L1, K-L2}
to prove similar results about 
unicritical polynomials, and to extend the results to some infinitely renormalisable quadratic polynomials.

Given the history of these techniques, it is 
over-optimistic to expect imminent application in other examples of 
parameter spaces, but there are certain features of the Yoccoz 
parapuzzle, and the Branner-Hubbard parapuzzle (\cite{B-H1, B-H2})
for the complement of 
the connectedness locus for cubic polynomials, which it is certainly 
worth noticing and trying to emulate. These parapuzzles are 
decompositions 
of parameter spaces, arise from Markov partitions of maps in the 
parameter space, Markov partitions which persist at least locally 
in parameter space. Such parapuzzles arising from Markov 
partitions have also been used by Kiwi \cite{K1} to describe other 
familes of 
cubic polynomials and to describe polynomials and rational maps on a 
certain nonArchimedean field. In the case of quadratic 
polynomials, Yoccoz uses a different Markov partition in each limb. 
In 
summary, persistent Markov partitions can give rise to parapuzzles 
and can be used to extract topological information about parameter 
spaces. So it seems interesting that the description of hyperbolic 
components in $V_{3}$ does involve the use of Markov partitions, 
three simple finite partitions and one countably infinite 
partition. They 
do give a topological model for the parameter space $(V_{3},H)$, 
where $H$ is the union of type III hyperbolic components, together 
with dynamical information about the components.  The model 
essentially identifies the parameter space locally with a subset of a 
fixed dynamical plane. But there is a difference from other models of 
parameter spaces, which may be significant. The identification of 
parameter space with dynamical plane is not one-to-one, but 
unboundedly-finite-to-one, on the part of parameter space associated 
with the countably infinite Markov partition mentioned above.  The 
topological 
picture is purely conjectural, and in fact no conjecture is being 
formalised at this stage, but it is a start.
 Countable Markov partitions have not 
arisen as parapuzzles before, to my knowledge, and they are not 
likely to be easy to use. Nevertheless, the way this one arises,
and its nature, is perhaps the most striking aspect of the current 
work. See Section \ref{9} for some very preliminary remarks on how 
this 
might develop. 

The organisation of this paper is as follows. In Section \ref{2} we 
collect together essentially elementary information about $V_{3}$, and recall some of the key concepts to be used, such as Thurston equivalence, captures and matings. Most of this material has, at 
least, been around for some considerable time.  This section ends with a first statement of the main theorem. In Section \ref{3}, we describe some less trivial equivalences between captures. This section ends 
with some straightforward counting, which might seem rather basic, 
but counting provides a very important check. In Section \ref{4} we 
summarise theory from earlier papers as it will be needed, in 
particular the Resident's View of \cite{R3}. In Section \ref{5} we consider the 
simple question of how fundamental domains are found in general, and 
specialise to our context. We give a second, rather vague, statement of the main theorem. In Section \ref{6}, we give the final statement of the main 
theorem in the three ``easy'' cases, and prove these. Section \ref{7} is 
largely given up to examples, showing how to calculate the fundamental domain for the first few preperiods in the ``hard'' region, and illustrating factors that have to be taken into account as the preperiod increases. The final statement of the ``hard'' case of the main theorem is given at the end of the section. The proof is given in Section \ref{8}. In Section \ref{9} we give 
some 
questions arising from this work.

I am greatly indebted to Vladlen Timorin, who read an earlier version of this paper, who made very perceptive suggestions on how to improve the presentation of the work, both on general strategy and on the level of detail. I have tried to respond to his suggestions as best I can. His suggestions and comments were concentrated on the first half of the work. During revision, quite substantial errors were found in the ``hard'' case of the theorem, necessitating a restatement, and, of course, new details of proof. This version of the paper is therefore quite significantly different from the 2005 version.

I should like to thank the referee of this paper for some very helpful comments, undoubtedly produced after a  lot of hard work.
\chapter{The space $V_{3}$}\label{2}

\section{About $V_{n}$}\label{2.1}

The map 
$$f_{c,d}(z)=1+{c\over z}+{d\over z^{2}}$$
has critical points 
$$0,\ \ {-2d\over c},$$
and critical values 
$$f_{c,d}(0)=\infty ,\ \ f_{c,d}(-2d/c)=1-{c^{2}\over 4d}.$$
Also $f_{c,d}(\infty )=1$. So if 
$n\geq 3$, the set of maps $f_{c,d}$ for which $0$ is of period $n$ 
is an algebraic curve minus finitely many punctures, defined by a 
polynomial 
in $(c,d)$. For 
example, $V_{3}$ identifies with the plane $1+c+d=0$ ($d\neq 0$). It 
is shown in Stimson's thesis \cite{Sti}, and also in Chapter 7 of 
\cite{R3}, that $V_{n}$ is nonsingular, that is, the only 
singularities on the algebraic curve containing it occur at points 
with $d=0$ or with both $c$, $d\to \infty $. See also \cite{M1, M2}.

We shall sometimes refer to $0$ as the {\em{first critical point}} 
$c_{1}$, and its image $\infty $ under $f_{c,d}$ as the {\em{first 
critical value}} $v_{1}$, to $-2d/c$ as the {\em{second critical 
point}} $c_{2}$ or, more precisely, $c_{2}(f_{c,d})$, and the image 
$1-(c^{2}/4d)$ under $f_{c,d}$ as the {\em{second critical value}}  
$v_{2}$ or $v_{2}(f_{c,d})$ For $m\geq 0$, we define $V_{n,m}$ to be 
the set of $f_{c,d}$ for 
which the $m$'th image of $v_{2}$ is not in the periodic orbit of 
$0=c_{1}$, that is, $f_{c,d}^{m+1}(-2d/c)\neq f_{c,d}^{j}(0)$ for any 
$j\geq 0$, equivalently, for any $0\leq j<n$. This is the complement 
in $V_{n}$ of 
finitely many points. If $n=3$ and $m=0$, then this identifies with 
the set of $d$ in  $ {\mathbb C}\setminus \{ 0,\pm 1\} $ with 
$1+c+d=0$. The set $V_{n,0}$ is the complement in $V_{n}$ of the 
centres of the type II components. For $m\geq 1$, the set $V_{n,m}$ 
is 
the complement in $V_{n,0}$ of the centres  of 
type III components of preperiod $\leq m$. So $V_{n,m+1}\subset 
V_{n,m}$ for 
all $m$. For all $n\geq 3$, the universal cover of $V_{n,m}$, for 
any $m\geq 0$, is conformally the unit disc $D$. We write $P_{n,m}$ 
for 
the set of punctures of $V_{n,m}$, that is, the union of the sets of 
singularities 
of $V_{n}$, of centres of type II components, and of centres of type 
III components of preperiod $\leq m$. In the case $n=3$, we have $0$, 
$\infty $, $\pm 1\in P_{3,m}$ for all $m\geq 0$. The other punctures 
are 
the centres of the type III hyperbolic components of preperiod $\leq 
m$.
The space $V_{n,m}$ identifies naturally with subspace of  a larger 
space 
$B_{n,m}$. The structure of $B_{n,m}$ was the main subject of  \cite{R5}: it is a countable union of geometric pieces connected by handles. One of the geometric pieces is $V_{n,m}$. For a degree two branched covering $g$, we number the two 
critical values $v_{1}=v_{1}(g)$ and $v_{2}=v_{2}(g)$. Suppose that 
$v_{1}$ is periodic, of period $n$. Write
$$Z_{m}(g)=g^{-m}(\{ g^{j}(v_{1}):j\geq 0\} ),\ \ 
Y_{m}(g)=Z_{m}(g)\cup \{ v_{2}\} .$$
Then $B_{n,m}$ is defined to be the space of all triples 
$[g,Z_{m}(g),v_2(g)]$, such that $v_1(g)$ is of period $m$ and $v_2(g)\notin Z_m(g)$, where $[.]$ denotes the quotient by M\"obius 
conjugation preserving the numbering of critical values. Hence, $[g,Z_m(g),v_2(g)]=[h,Z_m(h),v_2(h)]$ if and only if there is a M\"obius transformation $\sigma $ such that $\sigma \circ g\circ \sigma ^{-1}=h$, and $\sigma (v_i(g))=v_i(h)$ for $i=1$, $2$. This automatically implies that $\sigma (Z_m(g))=Z_m(h)$.
 We write $N_{n,m}$ for the set of 
topological (Freudenthal) ends of $B_{n,m}$. Then $P_{n,m}$ identifies with the set of ends of 
$V_{n,m}$, and with  a subset of $N_{n,m}$. Overall, it turns out that
many of the ends of $B_{n,m}$ identify naturally with spaces of 
critically 
finite branched coverings with $v_{2}$ of preperiod $\leq m$ 
\cite{R3}.

For the remainder of this paper, we retrict to the study of the case 
$n=3$, with some reference to the case of $n=1$, that is, the space 
$V_{1}$ of quadratic polynomials, and a few  references to the case 
$n=2$. 
Note that $V_{3,m}$ is a punctured sphere for all $m\geq 0$. Also, 
we can dispense with M\"obius conjugation in the definition of 
$B_{3,m}$ if we take all branched coverings $g$ in $B_{3,m}$ to have 
$c_{1}=0$, $v_{1}=\infty $ and $g(v_{1})=1$, as is true for 
$V_{3,m}$, 
which identifies with a subset of $B_{3,m}$. In 
Chapter 1 of \cite{R3} it is shown that $(V_{3,0},P_{3,0})$ 
is homotopy equivalent, under inclusion, to $(B_{3,0},N_{3,0})$. But 
for all $m>0$ inclusion of  $(V_{3,m},P_{3,m})$ in 
$(B_{3,m},N_{3,m})$ is injective, but not surjective, on $\pi _{1}$. 
That is a long story.

\section{Polynomials and Type II components in $V_{3}$}\label{2.2}
    
We start with a picture of $V_{3}$, 
in which polynomials (up to M\"obius conjugacy) and type II 
hyperbolic components are marked. 
\begin{center}
\includegraphics[width=6cm]{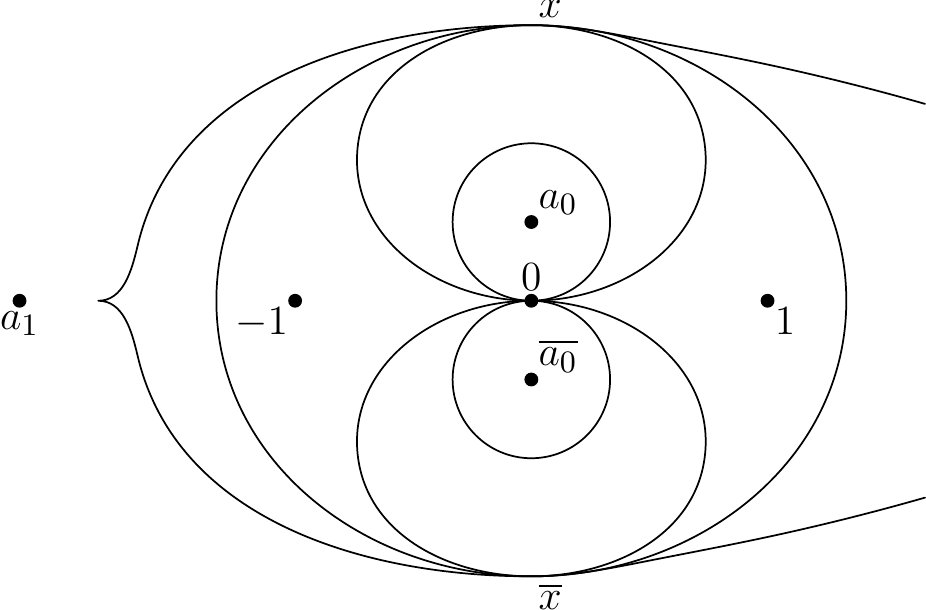}
\end{center}
\begin{center}$V_3$
\end{center}


 We write
$$h_{a}(z)={(z-1)(z-a)\over z^{2}}.$$
The critical points are $0=c_{1}$ and 
$$c_{2}=c_{2}(a)=\frac{2a}{a+1}.$$
The critical values are $\infty =v_{1}$ and 
$$v_{2}=v_{2}(a)=-\frac{(a-1)^{2}}{4a}.$$
There are three polynomials up to M\" obius conjugacy in $V_{3}$, 
represented by the conjugate parameter values $a=a_{0}$ or  
$\overline{a_{0}}$ 
and the real value $a=a_{1}<-1$. There are two type II hyperbolic 
components, containing $\pm 1$ respectively, and $h_{\pm 1}$ is the 
unique critically finite map in its hyperbolic component. For 
$h_{1}$, the critical points are $0$ and $1$ and the critical values 
are $\infty $ and $0$. For $h_{-1}$, the critical points are $0$ and  
$\infty $, and the critical values are $\infty $ and $1$. The hyperbolic 
components $H_{\pm 1}$ are drawn, identifying them with subsets of 
$\mathbb C$ and, as the very rough sketch 
indicates, their boundaries meet in three accessible points: $0$, 
$x$ and $\overline{x}$, where $x\in \mathbb C$ satisfies 
${\rm{Im}}(x)>0$. The points $a_{0}$, $\overline{a_{0}}$ and $a_{1}$ 
are then contained in distinct components of the complement of 
$H_{1}\cup H_{-1}$. The points $a_{0}$ and $\overline{a_{0}}$ are 
contained in bounded components of $\mathbb C\setminus (H_{1}\cup 
H_{-1})$, and hence their 
hyperbolic components $H_{a_{0}}$ and $H_{\overline{a_{0}}}$ are also 
contained in these bounded sets. The point $a_{1}$ and its hyperbolic 
component are contained in the unique unbounded component of  
$\mathbb C\setminus (H_{1}\cup 
H_{-1})$. This hyperbolic component is also drawn, and, as indicated, 
it also has  as accessible boundary points $x$ and $\overline{x}$. Proof of 
accessibility of  $x$ and  $\overline{x}$ is very 
similar to case 4b) of the Theorem of \cite{R4}. There is a natural 
parametrisation of $H_{1}$ (or $H_{-1}$) by the unit disc, using the 
B\"ottcher coordinate of the second critical value. Write $\varphi 
_{a}$ for the map such that 
$$\varphi _{a}(z^{2})=h_{a}^{3}\circ 
\varphi _{a}(z),\ \ \varphi _{a}(0)=1,\ \ \varphi _{a}'(0)=1.$$ 
Then 
$\vert \varphi _{a}^{-1}(v_{2}(a))\vert <1$ and 
$$\Phi _{1}:a\mapsto 
\varphi _{a}^{-1}(v_{2}(a))$$
is a degree three branched covering map 
from $H_{1}$ to the open unit disc, with $0$ as the unique branch 
point. The result of 4b) of \cite{R4}  says that each of the three components of 
$\Phi _{1}^{-1}(\{ r\zeta :r\in [0,1)\} )$ has a unique limit point as 
$r\to 1$, if $\zeta $ is a root of unity, $\zeta \neq 1$. The limit 
point is a parabolic parameter value, also in the closure of a type 
IV hyperbolic component. The Thurston equivalence class of the centre 
of this type IV hyperbolic component can be computed from the centre 
of the type II hyperbolic component --- in the current case, the 
hypebolic component is $H_{1}$, with centre $h_{1}$  --- and the 
choice of $\zeta $ and component of  $\Phi _{1}^{-1}(\{ r\zeta :r\in 
[0,1)\} )$. The result also holds for $\zeta =1$, and a particular  
component of  $\Phi _{1}^{-1}(\{ r :r\in [0,1)\} )$ if the associated 
critically finite branched covering is known to be Thurston 
equivalent to a rational map. (This is proved in exactly the same way 
as the result in \cite{R4}.) In the case of $H_{1}$, the branched 
coverings associated to two of the components  of  $\Phi _{1}^{-1}(\{ 
r :r\in [0,1)\} )$ are  indeed equivalent to rational maps. The 
branched covering associated to the third component is not equivalent 
to a rational map, but in that case the component has $a=0$ as a 
unique accumulation point.  It is not possible to be more precise 
without pre-empting the definition of mating in Section \ref{2.6}, and 
the subsequent results.

Now we fix some notation which will be used later. 
For any $a$ such that $h_{a}$ is hyperbolic, we write $H_{a}$ for the 
intersection with $V_{3}$ of the hyperbolic component of $h_{a}$. We 
write $V_{3}(a)$  
 for the component 
of $V_{3}\setminus (H_{1}\cup H_{-1})$ containing $a$, $a=a_{0}$ or  
$\overline{a_{0}}$. 
respectively. Write $V_{3}(a_{1},-)$ for the bounded component of 
$V_{3}\setminus (H_{1}\cup H_{-1}\cup H_{a_{1}})$ --- which 
intersects only the negative real axis --- and 
$V_{3}(a_{1},+)$ for the unbounded component of 
$V_{3}\setminus (H_{1}\cup H_{-1}\cup H_{a_{1}})$ --- which 
intersects only the positive real axis. We write 
$V_{3,m}(a)=V_{3}(a)\cap V_{3,m}$, 
$P_{3,m}(a)=P_{3,m}\cap V_{3}(a)$, and so on.

 The polynomial represented by $a_{0}$ is 
often known as the rabbit polyonomial, because of the shape of its 
Julia set. The Julia set for $\overline{a_{0}}$ is the conjugate of 
that for $a_{0}$. The rabbit polynomial for $a=a_{0}$ has the 
critical attractive basin rotated anticlockwise to the  critical 
value 
attractive basin, the next one round. To check this, it suffices to 
show that for $a$ near $0$ in the upper half plane in the 
hyperbolic component of $a_{0}$, the multiplier at the fixed point in 
the common boundary of the attractive basins is near $e^{2\pi i/3}$. 
Now this has to be a repelling fixed point, and, since $a$ is in the 
hyperbolic component of a polynomial, there also has to be 
an attractive fixed point, and the multiplier there has to be near 
$e^{-2\pi i/3}$. Now this is easily checked. The fixed points of 
$h_{a}$ are the 
roots of 
$$F(z,a)=1-{a+1\over z}+{a\over z^{2}}-z=0,$$
which for $a$ near $0$ are near the roots of 
$$F(z,0)=1-{1\over z}-z=0,$$
that is
$$z-1-z^{2}=0$$
or $z=e^{\pm \pi i/3}$. The multiplier at a fixed point is 
$${a+1\over z^{2}}-{2a\over z^{3}}.$$
Fix $\zeta $ with $F(\zeta ,0)=0$. Thus, $\zeta =e^{\pm i\pi /3}$.  
Let $z(a)$ be continuous in $a$ near $a=0$ with $z(0)=\zeta $ and 
$F(z(a),a)=0$.Then the Taylor expansion of $z(a)$ in $a$ near $0$ is  
$z=\zeta (1+ca)+o(a)$ with $c\zeta =-(\partial F/\partial a)(\zeta 
,0)/(\partial F/\partial z)(\zeta ,0)$. So
$$c={\overline{\zeta }-\zeta \over \vert 1+\zeta \vert 
^{2}}=-2i{{\rm{Im}}(\zeta )\over \vert 1+\zeta \vert 
^{2}},$$
 which is purely imaginary, and the derivative is
$$-\zeta (a+1)(1-2ca)+2a+o(a)=-\zeta  (1+a(1+2\zeta ^{2})-2ca).$$
So the modulus is $<1$ if and only if 
$${\rm{Re}}(1+a(1+2\zeta ^{2})-2ca)<0.$$
Since $1+2\zeta ^{2}$ is purely imaginary and equal to 
$2{\rm{Im}}(\zeta )$, the modulus is $<1$ if and only if 
$$2{\rm{Im}}(\zeta ){\rm{Im}}(a)+4{{\rm{Im}}(\zeta )\over \vert 
1+\zeta \vert 
^{2}}{\rm{Im}}(a)>0.$$
If the multiplier at $\zeta $ is near $e^{-2\pi i/3}$, then since the 
derivative is approximately $\zeta ^{-2}$, we must have $\zeta 
=e^{i\pi /3}$ and ${\rm{Im}}(a)>0$, as claimed.

 The immediate disc regions containing $a_{0}$ and 
$\overline{a_{0}}$ indicate the sets $H_{a_{0}}$ and  
$H_{\overline{a_{0}}}$.
The next largest  regions in parameter space are the regions 
$V_{3,0}(a_{0})$ and $V_{3}(\overline{a_{0}})$. It is known that all 
hyperbolic components in 
these regions are matings or captures with $h_{a_{0}}$ or  
$h_{\overline{a_{0}}}$. 
These include a mating with the polynomial 
represented by $a_{1}$ in each case. These hyperbolic components are 
adjacent to the points $x$ and $\overline{x}$ respectively. The points 
$x$ and $\overline{x}$ are also common boundary points of the 
sets $H_{1}$ and $H_{-1}$, and in the 
boundary of the set $H_{a_{1}}$, which is unbounded. 
There are two other regions of $V_{3}$ 
which contain infinitely many hyperbolic components. 
One is the region $V_{3,0}(a_{1},-)$ bounded by $H_{-1}$ and 
$H_{a_{1}}$, 
and the other is 
the region $V_{3,0}(a_{1},+)$  bounded by $H_{1}$ and $H_{a_{1}}$. 

\section{Lamination maps}\label{2.3}
Invariant laminations were introduced by Thurston \cite{T} to 
describe 
the dynamics of polynomials with locally connected Julia sets. The 
theory is most developed for quadratic polynomials, but a number of 
people have worked on developing the theory for higher degree 
polynomials, e.g. \cite{K2}, \cite{K3}, \cite{K4}. See also  
\cite{R1} (and 
\cite{R3}) for a slightly more detailed summary in the quadratic case 
than is given here. 
The {\em{leaves}} of a lamination $L$  
are straight line segments in $\{ z:\vert z\vert \leq 1\} $. 
Invariance of a lamination means that if 
there is a leaf with endpoints $z_{1}$ and $z_{2}$ then there is also a 
leaf with endpoints $z_{1}^{2}$ and $z_{2}^{2}$, a leaf with endpoints 
$-z_{1}$ and $-z_{2}$, and a leaf with endpoints $w_{1}$ and $w_{2}$ where 
$w_{1}^{2}=z_{1}$ and $w_{2}^{2}=z_{2}$.  A leaf 
with endpoints $e^{2\pi ia_{1}}$ and $e^{2\pi ia_{2}}$ for $0\leq 
a_{1}<a_{2}<1$,
is then said to have {\em{length}} 
${\rm{min}}(a_{2}-a_{1},a_{1}+1-a_{2})$. {\em{Gaps}} of the 
lamination are components of $\{ z:\vert z\vert <1\} \setminus (\cup 
L)$. 
If the longest leaf of $L$ has length $<{1\over 2}$  then there are 
exactly two with the same image which is called the {\em{minor 
leaf}}.  A lamination is {\em{clean}} if finite-sided gaps are never 
adjacent. Minor leaves of clean laminations are either equal or have 
disjoint interiors. For a clean lamination, the lamination 
equivalence 
relation $\sim _{L}$ is closed, where the nontrivial equivalence 
classes are the 
leaves of $L$ and closures of finite-sided gaps of $L$. The quotient 
space $\overline{\mathbb C}/\sim _{L}$ is a topological sphere.  A 
partial order on minor leaves is therefore defined by: $\mu 
_{1}\leq \mu _{2}$ if $\mu _{1}$ separates $\mu _{2}$ from $0$. For 
any minor leaf $\mu $, the set $\{ \mu ':\mu '\leq \mu \} $ is 
totally 
ordered, and has a unique minimal element.
If the gap containing $0$ has infinitely many sides, then it is periodic, and  one 
can define a {\em{lamination map}} $s=s_{L}$ which maps $L$ to $L$, 
$s(z)=z^{2}$ for $\vert z\vert \geq 1$, and maps gaps to gaps. Such a 
lamination is uniquely determined by its minor leaf. The 
laminations with this property are precisely those with minor leaf 
with endpoints $e^{2\pi ia_{1}}$ and $e^{2\pi ia_{2}}$, where $a_{1}$ and  
$a_{2}$ are odd 
denominator rationals, of the same period under the map $x\mapsto 
2x{\rm{\ mod\ }}1$. We can then choose $s$ so that $0$ is critical 
and 
periodic under $s$, of the same period as the gap containing $0$, and 
as the endpoints of the minor leaf, which bounds the gap containing 
$s(0)$. Since a lamination map $s_{L}$ preserves $L$, it descends to 
a 
map $[s_{L}]:\overline{\mathbb C}/\sim _{L}\to \overline{\mathbb 
C}/\sim _{L}$.
Any quadratic polynomial $f_{c}:z\mapsto z^{2}+c$, for which 
$0$ is periodic, is {\em{Thurston equivalent}} (to be defined in 
the next section) to exactly one such lamination map $s$, and 
topologically conjugate to the 
quotient lamination 
map $[s]$, and conversely.

Let $p$ be any odd denominator rational. We write $L_{p}$ for the 
invariant lamination with minor leaf with 
endpoint at $e^{2\pi ip}$, and $\mu _{p}$ for this minor leaf.  Thus,
$L_{p}=L_{q}$ if and only if either $p=q$ or $e^{2\pi ip}$ and
$e^{2\pi iq}$ are at opposite ends of the minor leaf of 
$L_{p}=L_{q}$. 
We write $s_{p}$ for $s_{L_{p}}$.  The maps $h_{a}$ for $a=a_{0}$,
$\overline{a_{0}}$ and $a_{1}$ are conjugate to quotient lamination maps
$[s_{1/7}]$, $[s_{6/7}]$ and $[s_{3/7}]$ respectively.  These will
feature large in our description of the type III hyperbolic components
in $V_{3}$.

\section{Thurston equivalence}\label{2.4}

Thurston's notion of equivalence of critically finite branched 
coverings, and Thurston's Theorem which characterises those 
critically 
finite branched coverings which are equivalent to rational 
maps, underpins this work.  A topological branched covering 
$f:\overline{\mathbb C}\to \overline{\mathbb C}$ is {\em{critically 
finite}} if
$$X(f)=\{ f^{n}(c):n>0,\ c{\rm{\ critical}}\} $$
 is a finite set. In the current work the definition of $f$ also 
 includes a fixed identification of $X(f)$ with a finite set $X_{0}$ 
 with dynamics 
 such that the identification preserves dynamics.  Two critically 
finite maps $f_{0}$ and $f_{1}$ are 
{\em{(Thurston) equivalent}} if there is a homotopy $f_{t}$ from 
$f_{0}$ to $f_{1}$ such that $\# (X(f_{t}))$ is constant in $t$, and 
thus the finite set $X(f_{t})$ varies isotopically with $t$, and the 
isotopy between $X(f_{0})$ and $X(f_{1})$ preserves identification 
with $X_{0}$. If $\# 
(X(f_{t}))\geq 3$, this is equivalent to the existence of 
homeomorphisms $\varphi $ and  $\psi :\overline{\mathbb C}\to 
\overline{\mathbb C}$ with $\varphi $ and $\psi $ isotopic via an 
isotopy mapping $X(f_{0})$ to $X(f_{1})$, preserving the 
identifications of $X(f_{0})$ and $X(f_{1})$ with $X_{0}$, and such 
that 
$$\varphi \circ f_{0}\circ \psi ^{-1} =f_{1}.$$
Thurston's theorem gives a characterisation of those critically 
finite  
branched coverings which are equivalent to rational maps. For this 
characterisation the reader is referred, for example, to 1.8 and 2.4 
of \cite{R3}, or to \cite{D-H2}. But it is worth noting that, if 
$f_{0}$ 
and  $f_{1}$ are equivalent to a hyperbolic rational map, then Thurston's 
Theorem 
says that the rational map is unique up to M\"obius conjugacy, and 
its Thurston equivalence class in the space of branched coverings is 
simply connected if $\#(X(f_{0}))\geq 3$. This depends on our 
definition of Thurston equivalence as being via isotopies preserving 
identification with $X_{0}$. It follows that $\varphi $ 
and $\psi $ are unique up to isotopy constant on $X(f_{0})$, and 
$\psi $ is actually unique up to isotopy constant on 
$f_{0}^{-1}(X(f_{0}))$. In the situation described here, we write 
$$f_{0}\simeq _{\varphi }f_{1},$$
which is slightly dubious notation as $\simeq _{\psi }$ is not an 
equivalence relation, although Thurston equivalence certainly is. We 
shall sometimes write
$$(f_{0},X(f_{0}))\simeq _{\varphi }(f_{1},X(f_{1})),$$
meaning that conjugacy by $\varphi $ is followed by an isotopy 
constant on $X(f_{1})$. 
Note that  since $\varphi $ and $\psi $ are isotopic via an isotopy 
constant on $X(f_{0})$, there is a homeomorphism $\psi _{2}:\mathbb 
C\to \mathbb C$ which maps $f_{0}^{-1}(X(f_{0}))$ to 
$f_{1}^{-1}(X(f_{1}))$ defined by the relation
$$f_{0}\circ \psi _{2}=\psi \circ f_{0}.$$
We thus have 
$$f_{0}\simeq _{\psi }f_{1},$$
or more precisely
$$(f_{0},f_{0}^{-1}(X(f_{0})))\simeq _{\psi 
}(f_{1},f_{1}^{-1}(X(f_{1}))).$$
Similarly we have
$$(f_{0},f_{0}^{-2}(X(f_{0})))\simeq _{\psi 
_{2}}(f_{1},f_{1}^{-2}(X(f_{1}))).$$
We can similarly define $\psi _{n}$ for all $n$ with
$$(f_{0},f_{0}^{-n}(X(f_{0})))\simeq _{\psi 
_{n}}(f_{1},f_{1}^{-n}(X(f_{1}))).$$
See also 1.7 of \cite{R3} and 1.3 of \cite{R1}.

\section{Some notation}\label{2.5}

Throughout this paper we shall use the notation $\sigma _{\beta }$ to 
denote a certain type of homeomorphism associated to a path $\beta $ 
in $\overline{\mathbb C}$. If $\beta :[0,1]\to \overline{\mathbb C}$ 
is an arc, then $\sigma _{\beta }$ is the identity outside a small 
disc neighbourhood of $\beta $, and $\sigma (\beta (0))=\beta (1)$. 
More generally, if $*$ denotes the usual multiplication on paths, and 
$\beta =\beta _{1}*\cdots *\beta _{r}$ where each 
$\beta _{i}$ is an arc, then
$$\sigma _{\beta }=\sigma _{\beta _{r}}\circ \cdots \circ \sigma 
_{\beta _{1}}.$$
This does not of course define $\sigma _{\beta }$ pointwise, but if 
we 
take a sufficiently small neighbourhood of $\beta $, then $\sigma 
_{\beta }$ will be defined up to isotopy constant on any finite set 
which intersects $\beta $ at most in the endpoints $\beta (0)$, 
$\beta 
(1)$. In this paper $\beta $ will always (as it turns out) be either 
an arc, or a simple closed loop, and we shall describe maps up to 
Thurston equivalence by compositions of the form $\sigma _{\beta 
}\circ f$ for various $\beta $ and various $f$. 

Throughout this paper, if $\beta :[a,b]\to \mathbb C$ is a path 
then $\overline{\beta }:[a,b]\to \mathbb C$ is the reverse path 
defined by $\overline{\beta }(t)=\beta (a+b-t)$.

\section{Captures and Matings}\label{2.6}

In this paper, a {\em{capture}} is a critically finite branched 
covering of a particular type, up to Thurston equivalence. This 
appears to be what is meant by ``capture'' in Wittner's 1988 thesis 
\cite{W}. Some authors take the word to mean any type III component, 
but the use made of the term in this paper coincides with the use 
made in \cite{R1, R2, R3}, the reason being that 
captures in this sense are the most natural (type III) analogue of 
the matings which are also discussed by Wittner. These were 
introduced by  Douady and Hubbard and have been very important 
examples since their introduction. 
Presumably because of the connection with the concept of mating, 
Wittner actually suggests the concept of capture is already known, 
but does 
not give any earlier written reference, and from recent conversation 
it seems likely that his thesis gives the first written reference to captures. 

 Let 
$s=s_{L}$ be a lamination map which is Thurston equivalent to a 
polynomial preserving the 
corresponding invariant lamination $L$. The second critical point 
$c_{2}(s)$ is fixed, coincides with $v_{2}(s)$, and is $\infty $ in 
the standard model. Fix any gap $G$ of $L$ which 
contains some point $x$ in the backward orbit of $0$ under $s$. Then 
$x$ is 
eventually periodic under $s$.  Let $\beta $ be an arc from 
$c_{2}(s)=v_{2}(s)$ to $x$ which crosses the unit circle 
exactly once, at a point of $\partial G$, such that, apart from this 
one point, all of $\beta $ is in the union of $G$ and $\{ z:\vert 
z\vert >1\} $. If $x$ is periodic under $s$ then let $\zeta $ be 
the arc in $s^{-1}(\beta )$ from $c_{2}(s)$ to the periodic 
point in $s^{-1}(x)$. Then the associated capture, a critically 
finite branched covering defined up to Thurston equivalence, is 
$$\sigma _{\beta }\circ s\ \ {\rm{or}}\ \ \sigma _{\zeta }^{-1}\circ \sigma 
_{\beta }\circ s$$
depending on whether $x$ is strictly preperiodic or periodic under 
$s$. These can be distinguished as {\em{type III captures}} 
and {\em{type II captures}}.

The similar definition of a mating $s_{p}\Amalg s_{q}$ is defined for 
any odd denominator rationals $p$ and $q$, and satisfies
$$s_{p}\Amalg s_{q}(z)=\begin{array} {l}s_{p}(z){\rm{\ if\ }}\vert 
z\vert \leq 1,\cr (s_{q}(z^{-1}))^{-1}{\rm{\ if\ }}\vert z\vert \geq 
1.\cr \end{array}$$
A rational map which is Thurston equivalent to a mating $s_{p}\Amalg 
s_{q}$ must be the centre of a type IV hyperbolic component. So in 
this paper we concentrate on the captures.

 There is a simple characterisation of the
captures which are Thurston equivalent to rational maps, which was 
proved by Tan Lei: see \cite{TL}. Actually that work concentrates on 
matings, but the modifications needed to consider captures are 
straightforward. The map $\sigma _{\beta }\circ s$ 
(or $\sigma _{\zeta }^{-1}\circ \sigma _{\beta }\circ s$ in the case 
of type II captures) is 
Thurston-equivalent to a rational map if and only if the 
$S^{1}$-crossing point of $\beta $ is not strictly inside the smaller 
region of the circle bounded by $\mu _{0}$, where $\mu _{0}$ is the 
minimal minor leaf with $\mu _{0}\leq \mu _{L}$, and, if $\beta $ 
crosses the circle precisely at an endpoint of $\mu _{0}$, then the 
second 
endpoint of 
$\beta $ is not $v_{1}$. This last case, when the $S^{1}$-crossing 
point is an endpoint of $\mu _{0}$, perhaps needs a bit more 
adapting from Tan Lei's 
results than the others, but it is in fact a very special case, 
implying that $\mu 
_{0}=\mu _{L}$ and that the associated capture is a type II capture 
of 
a very special type. So these few cases can be dealt with separately 
and quite easily. A more general result was proved in the 
Non-Rational Lamination Map Theorem of \cite{R2}.

  In any case, the dynamics of a rational map equivalent 
to a capture defined by 
$s=s_{L}$ and $\beta $ is easily described in terms of 
the dynamics of $s$, and its 
Julia set is easily described in 
terms of $L$. The description is proved in a 
more general setting in the Lamination Map Conjugacy Theorem of 
\cite{R1}.  We shall say that the hyperbolic 
component of a rational map which is a capture, up to Thurston 
equivalence is a {\em{capture hyperbolic component}}, or simply a 
capture.

\section {Captures in $V_{3}$}\label{2.7}

As a consequence of Tan Lei's theorem, or a slight generalisation of 
it,  in the family $V_{3}$,  there are three families of type III 
captures. The first is
$$\sigma _{\beta }\circ s_{1/7},$$
where $\beta $  is any path with single $S^{1}$-crossing into  a 
nonperiodic gap  $G$ in the full orbit of the critical gap of 
$s_{1/7}$, where $G$ is in the larger region of the disc bounded by 
the minor leaf of $L_{1/7}$, which has endpoints at $e^{2\pi i(1/7)}$ and $e^{2\pi 
i(2/7)}$. The second family is similar:
$$\sigma _{\beta }\circ s_{6/7},$$
where the single $S^{1}$-crossing by $\beta $ is in the boundary of a preperiodic gap in the full orbit of the critical gap of $L_{6/7}$, in the larger region of the disc bounded by the minor leaf of $L_{6/7}$. The third family is
$$\sigma _{\beta }\circ s_{3/7},$$
where the single $S^{1}$-crossing by $\beta $ is  in the boundary of a preperiodic gap in the full orbit of the critical gap of $L_{3/7}$, in the larger region of the disc bounded by the leaf of $L_{3/7}$ with endpoints at $e^{2\pi i(1/3)}$ and $e^{2\pi i(2/3)}$. This 
time the boundaries of all gaps are disjoint, so one can describe 
$\beta $ (up to homotopy keeping endpoints fixed) simply by its 
$S^{1}$-crossing. The $S^{1}$-crossing must be a nonperiodic point in 
the backward orbit under $z\mapsto z^{2}$ of $e^{2\pi i(2/7)}$ or 
$e^{2\pi i(5/7)}$, since  every leaf in the boundary of a gap is in 
the backward orbit of the leaf with these endpoints in the boundary 
of the critical gap. If $\beta $ has just one  $S^{1}$-crossing at 
$e^{2\pi 
ip}$, we shall simply write $\beta =\beta _{p}$ and  $\sigma _{\beta 
}=\sigma _{\beta _{p}}=\sigma _{p}$.  Strictly 
speaking, such maps are only equal up to a suitable homotopy, but 
this 
uniquely defines $\sigma _{p}\circ s_{3/7}$ up to Thurston 
equivalence. It should be pointed out that this notation could be 
confusing, because it is somewhat at odds with the notation used for 
matings. The mating $s_{q}\Amalg s_{p}$ preserves the lamination 
$L_{q}\cup L_{p}^{-1}$, and the minor leaf of $L_{p}^{-1}$ has an 
endpoint at $e^{-2\pi ip}$, not at $e^{2\pi ip}$. But the notation 
does seem rather natural and hopefully will not cause much confusion, 
since matings are somewhat in the background of this paper.

We shall use a similar notation for type II captures. If $e^{2\pi ip}$ is periodic under $z\mapsto z^2$, and $\beta _p$ is a capture path for $s$, where $s=s_{3/7}$ or $s_{1/7}$ or $s_{6/7}$, then we write $\zeta _p$ for the path with periodic endpoint which is mapped homeomorphically to $\beta _p$ by $s$, and define $\sigma _p=\sigma _{\zeta _p}^{-1}\circ \sigma _{\beta _p}$. Then $\sigma _p\circ s$ is a type II capture.

\section{Trivial Thurston equivalences between captures in 
$V_{3}$}\label{2.8}

It is in general a nontrivial question as to whether two captures, 
$\sigma _{\beta }\circ s$ and $\sigma _{\beta '}\circ s$ are 
Thurston equivalent. But there are some rather trivial equivalences 
within the families above which can be seen immediately. Suppose 
that $\beta _{1}$ and $\beta _{2}$ are two paths giving captures in 
the first family above, with $\beta _{1}$ and $\beta _{2}$ having the 
same endpoint $x$. Then $\beta _{1}*\overline{\beta _{2}}$ is a 
closed 
path which bounds a disc $D$ disjoint from the forward orbit of $x$. 
It follows that $\sigma _{\beta _{1}}$ and $\sigma _{\beta _{2}}$ are 
isotopic via an isotopy constant on $\{ s_{1/7}^{n}x:n\geq 0\} \cup 
\{ \infty \} $, and hence
$$\sigma _{\beta _{1}}\circ s_{1/7}\simeq \sigma _{\beta _{2}}\circ 
s.$$
Thus, type III captures  $\sigma _{\beta }\circ s_{1/7}$ in the first 
family 
are uniquely determined, up to Thurston equivalence, by the second 
endpoint of $\beta $. In fact, the same is true for type II captures, 
by the same argument, because if we take two captures $\sigma _{\zeta 
_{1}}^{-1}\circ \sigma _{\beta _{1}}\circ s$ and $\sigma _{\zeta 
_{2}}^{-1}\circ \sigma _{\beta _{2}}\circ s$, then the 
paths $\zeta _{1}$ and $\zeta _{2}$ are 
also homotopic via an isotopy constant on $\{ s_{1/7}^{n}x:n\geq 0\} 
\cup 
\{ \infty \} $, noting that in this case $x$ is periodic and 
the common endpoint of $\zeta 
_{1}$ and $\zeta _{2}$ is in the periodic orbit of $x$. Exactly 
similar statements hold for the second family of type III captures 
$\sigma _{\beta }\circ s_{6/7}$, and for type II captures $\sigma 
_{\zeta 
}^{-1}\circ 
\sigma _{\beta }\circ s_{6/7}$. 

It is not too hard to show directly that the Thurston equivalence 
classes in these two families are in one-to-one correspondence with 
the  
endpoints $x(\beta )$ of the corresponding paths $\beta $. The idea 
is basically that if one takes the tree formed by joining points in 
the forward orbit of $x(\beta )$ to the centres of adjacent triangle 
gaps, then the ambient isotopy class of this tree, with some vertices 
marked by a point in the forward orbit of $x(\beta )$, is fairly 
easily 
shown to be an invariant of the Thurston equivalence class, 
essentially because any Thurston equivalence between $\sigma _{\beta 
} 
\circ s_{q}$,  for $q={1\over 7}$ or ${6\over 7}$, and the corresponding 
polynomial $g$, must isotope the tree to a tree joining up the full 
orbit 
of the finite critical point, which is entirely in the Fatou set, 
apart 
from common boundary points between Fatou components. We omit 
the details because the main theorems of this paper subsume this, and 
use different methods. 

Somewhat more restricted statements are true for the third family of 
type III captures $\sigma _{\beta }\circ s_{3/7}$, and type II 
captures $\sigma _{\zeta }^{-1}\circ \sigma _{\beta }\circ s_{3/7}$. 
It is no longer generally true that $\beta _{1}$ and $\beta _{2}$ are 
homotopic via an isotopy constant on $\{ s_{3/7}^{n}x:n\geq 0\} \cup 
\{ \infty \} $, if $\beta _{1}$ and $\beta _{2}$ have the same 
endpoint $x$. But it is true if, in addition:
\begin{itemize}
    \item the $S^{1}$-crossing 
points 
of $\beta _{1}$ and $\beta _{2}$ are both in the upper half-plane;
\item or 
both in the lower half-plane;
\item or both in the clockwise unit circle arc from $e^{2\pi 
i(1/7)}$ to $e^{2\pi 
i(6/7)}$;
\item  or both in the anticlockwise circle arc from $e^{2\pi 
i(1/7)}$ to $e^{2\pi i(2/7)}$;
\item  or both in 
the anticlockwise circle arc from 
$e^{2\pi i(5/7)}$ to $e^{2\pi i(6/7)}$.
\end{itemize}
These all follow from a basic 
property of 
invariant laminations established by Thurston \cite{T}, that a leaf 
$\ell '$ in the forward orbit of a leaf $\ell $ can only be shorter 
than $\ell $ if $\ell $ is longer than the minor leaf. The minor leaf 
of $L_{3/7}$ joins $e^{2\pi i(3/7)}$ and $e^{2\pi i(4/7)}$ and any 
leaf in the region bounded by $\beta _{1}*\overline{\beta _{2}}$, 
under our hypotheses, has length less than this. So the disc bounded 
by $\beta _{1}*\overline{\beta _{2}}$ is again disjoint from the set 
$\{ s_{3/7}^{n}x:n\geq 0\} \cup 
\{ \infty \} $, and it is again true that $\sigma _{\beta _{1}}$ and 
$\sigma _{\beta _{2}}$ are 
isotopic via an isotopy constant on$\{ s_{3/7}^{n}x:n\geq 0\} \cup 
\{ \infty \} )$ and similarly for type II captures. 
 So, this time, if we restrict to lamination maps which are Thurston equivalent to rational maps,
there is at most one Thurston equivalence classes of type III 
captures $\sigma _{\beta }\circ s_{3/7}$, for each choice of endpoint 
of $\beta $, except when $\beta $ ends at a gap $G$ of $L_{3/7}$, in 
the smaller region of the disc bounded by the leaf with endpoints at 
$e^{2\pi i(2/7)}$ and $e^{2\pi i(5/7)}$ and the leaf with endpoints at 
$e^{2\pi i(1/3)}$ and $e^{2\pi i(2/3)}$, and $\partial G$ intersects the 
unit circle in both the upper and lower half plane. In this case, 
there are still at most two different Thurston equivalence classes, 
depending on whether the $S^{1}$-crossing of $\beta $ is in the lower 
or upper half-plane.  Similar statements hold for 
type II captures.

From the statements above it follows that 
any type II capture $\sigma _{\zeta }^{-1}\circ \sigma _{\beta }\circ 
s_{3/7}$ must be Thurston equivalent to one where $\beta $ has the 
its only $S^{1}$-crossing at $e^{2\pi i(2/7)}$, $e^{2\pi i(5/7)}$ or 
$e^{2\pi i(1/7)}$. We shall see later, in \ref{3.3}, that there is a 
simple Thurston 
equivalence between the type II captures corresponding to  
$S^{1}$-crossings at $e^{2\pi i(2/7)}$ and $e^{2\pi i(5/7)}$, where the 
homeomorphism $\varphi $ realising the equivalence fixes $0$ and 
$\infty $.

The number of gaps of $s_{3/7}$ in the full orbit of the critical gap 
and of preperiod $n$ is $2^{n}$. The number of these in the large 
region of the disc bounded by the leaves with endpoints $e^{\pm 2\pi 
i(1/3)}$ and $e^{\pm 2\pi i(2/7)}$ is ${1\over 21}2^{n}(1+o(1))$. A 
gap whose boundary intersects the circle in both upper and lower 
half-planes must have its entire forward orbit in between the leaves 
with endpoints at $e^{\pm 2\pi i(3/7)}$ and $e^{\pm 2\pi i(1/7)}$  
The number of such gaps whose boundaries intersect the unit circle in 
both the upper and lower half planes can be shown to be
 $${1\over 21}.\left({3\over 2}\right) ^{n}+O(1).$$
This is $o(2^{n})$. So, asymptotically, we can bound the number of 
captures 
up to Thurston equivalence by the number of endpoints of the paths 
$\beta $.

Regions of parameter space such as $V_{3,0}(a_{0})$ and
$V_{3,0}(\overline{a_{0}})$, that is, parts of quadratic rational
parameter space with a natural association to some polynomial with
star-like Julia set, were studied by the Hubbard school in the early
1990's.  In particular, J. Luo \cite{Luo} wrote his thesis on this
subject, concentrating on the space $V_{2}$.  He claimed a number of
results, including the following: all type III and type IV hyperbolic
components in $V_{2}$ have centres which are Thurston equivalent to
matings or captures.  This is true, as is an analogue for
$V_{3,0}(a_{0})$ and $V_{3,0}(\overline{a_{0}})$.  All type III
hyperbolic components in $V_{3,0}(a_{0})$ or $V_{3,0}(\overline{a_{0}})$
have centres which are Thurston equivalent to captures $\sigma _{\beta
}\circ s_{1/7}$ or $\sigma _{\beta }\circ s_{6/7}$, and that all such
captures are in this region, and similarly for type IV hyperbolic
components in this region and matings.   Luo
formulated and stated much stronger results: an analogue, for $V_{2}$,
of the Yoccoz puzzle and parapuzzle for quadratic polynomials 
\cite{H}, and
associated results about local connectivity of Julia sets of
nonrenormalisable polynomials, and uniqueness of quadratic rational
maps for any specified nonrenormalisable combinatorics.  The proofs in
the thesis were not complete, although it seems very likely that they
were completely known to the Hubbard school.  But these, and similar
results, have recently been reproved and written by Aspenberg and
Yampolsky \cite{Asp-Yam}.  Another, different exposition is given by
Timorin, \cite{Tim}, which concentrates on the boundary of the type II
hyperbolic component in $V_{2}$, where, rather surprisingly, only
nonrenormalisable combinatorics occur.  

As regards the simple result
about all components in $V_{3,0}(a_{0})$ and
$V_{3,0}(\overline{a_{0}})$ being represented by captures and matings,
a simple analytical argument is available.  (The argument in the case
of $V_{2}$, as explained by Timorin \cite{Tim}, is even simpler.) 
Basically, the
idea is to consider the point $z_{1}(h)$ for all $h$ in the region
bounded by  $H_{1}$ and $H_{-1}$,  where $z_{1}(h)$ is the common boundary
point of the three periodic attractive basins.  It can be shown that
$z_{1}(h)$ exists throughout this region, essentially because the set
where it exists is an open neighbourhood of any point where $z_{1}(h)$
is a repelling fixed point.  So $z_{1}(h)$ only disappears when one
passes through a parabolic fixed point or if $h\mapsto z_{1}(h)$ is a
multivalued map, remembering that the set of fixed points of $h$ is a
multivalued map in general.  But $h\mapsto z_{1}(h)$ is
single-valued in the region under consideration.  Simple computer
graphics do show clearly the rabbit-like connnections between the
different capture hyperbolic components in this part of parameter
space, the part between $H_{1}$and $H_{-1}$.  They also show the
trunctated Mandelbrot sets, each with a period three limb removed, one
in the upper half-plane between $H_{\pm 1}$ and the other in the lower
half-plane.

Part of the statement of the main theorem (\ref{2.10}, \ref{6.1}) is 
that all 
centres of type III 
hyperbolic components in $V_{3,0}(a_{1},-)$
are Thurston equivalent to captures $\sigma _{\beta 
}\circ s_{3/7}$, where $\beta $ has $S^{1}$-crossing at $e^{2\pi ip}$ 
with $p\in (-{1\over 7},{1\over 7})$, and in fact one can restrict to 
$p\in (0,{1\over 7})$. I do not know a simple analytical argument, 
like the one above, to prove this, and, so far as I know, this result 
is new.

\section{Symbolic Dynamics}\label{2.9}

Let $p$ be an odd denominator rational. 
It is natural to make a Markov partition for the lamination map 
$s_{p}$ using preperiodic leaves of $L_{p}$ as boundaries of sets of 
the partition.  In this paper,
we are especially interested in doing this for $p={1\over 7}$,
${3\over 7}$, ${6\over 7}$.  For the case $p={3\over 7}$, we look at
the partition of $\{ z:\vert z\vert \leq 1\} $ by the leaves in
$L_{3/7}$ with endpoints at $e^{2\pi ia}$, $e^{2\pi ib}$ where $(a,b)$
is one of the following: $$({2\over 7},{5\over 7}),\ ({3\over
7},{4\over 7})\ ,({1\over 7},{6\over 7}),\ ({5\over 14},{9\over 14}),\
({3\over 14},{11\over
14}),\ ({1\over 14},{13\over 14}).$$
In the case of $p={1\over 7}$ we take a similar partition where 
$(a,b)$ is one of the following:
$$({1\over 7},{2\over 7}),\ ({2\over 7},{4\over 7})\ ,({4\over 
7},{1\over 7}),\ ({9\over 14},{11\over 14}),\ ({11\over 14},{1\over 
14}),\ ({9\over 14},{1\over 14}).$$
The partition for $s_{6/7}$ is obtained from this partition by 
conjugation. Here are the Markov partitions for $s_{3/7}$ and 
$s_{1/7}$, with the different sets of the partition labelled. We use 
$C$ for the ``central'' or ``critical'' region, which contains the
critical gap of $L_{p}$.  
\begin{center}
\includegraphics[width=7cm]{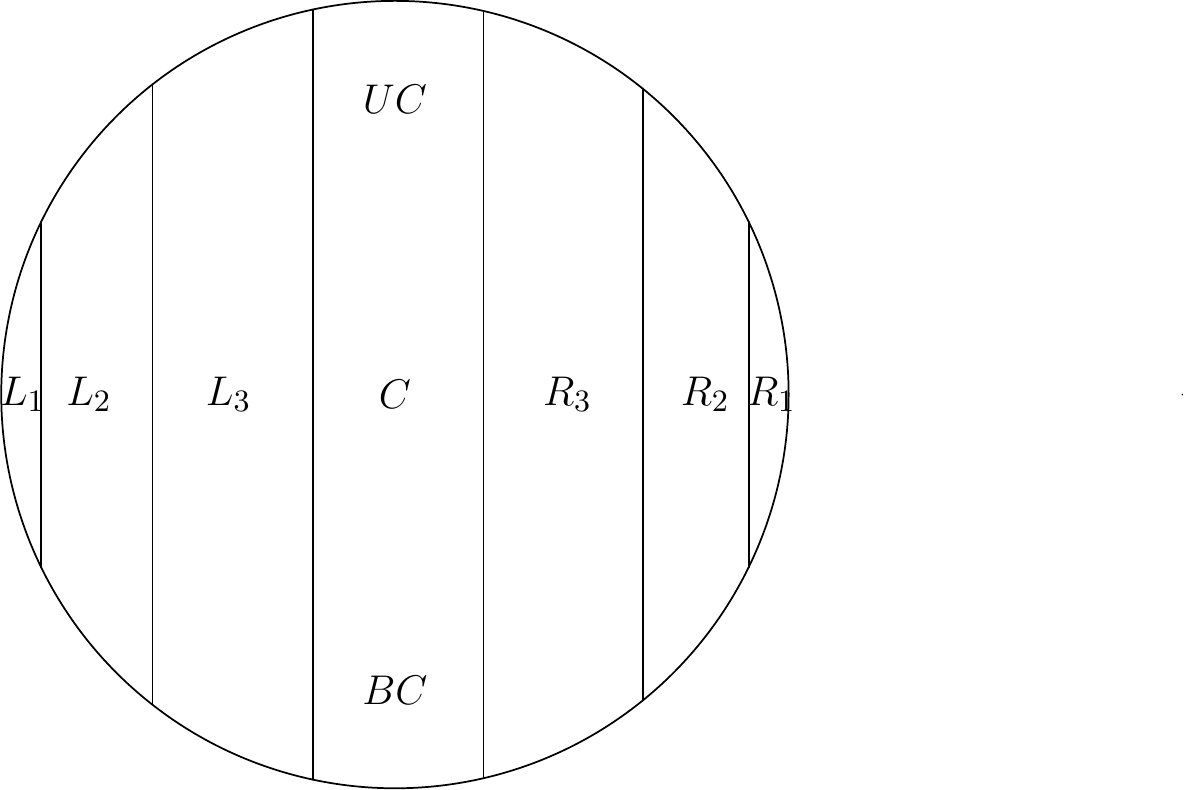}\includegraphics[width=4.7cm]{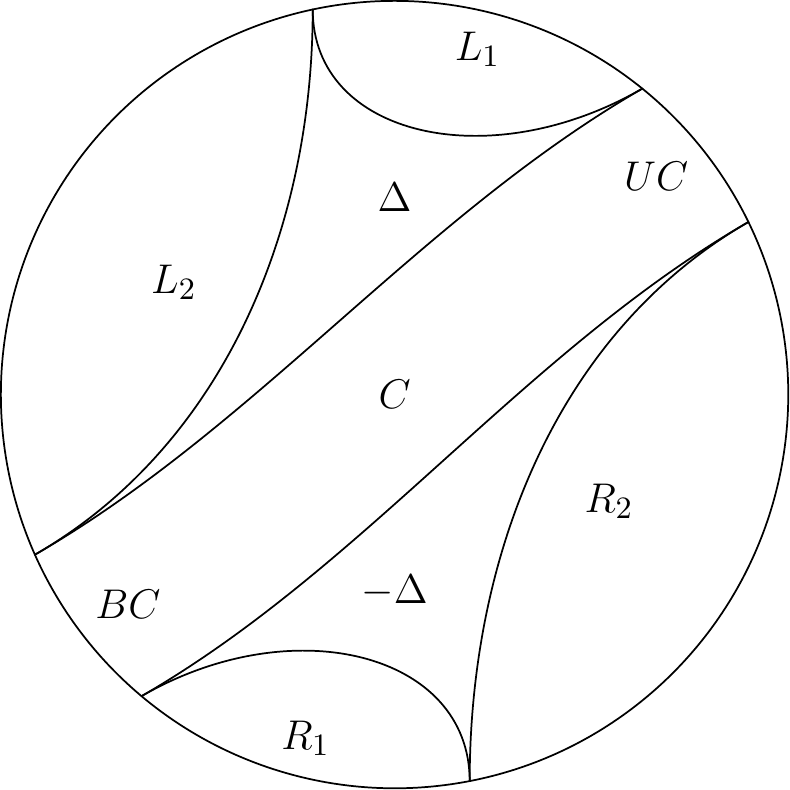}
\end{center}
\begin{center}Symbolic dynamics for $s_{3/7}$ and $s_{1/7}$
\end{center}

We denote by $BC$ and $UC$ with
intersection with the upper and lower half-planes respectively of
the complement in $C$ of the critical gap.  Thus, the intersection of
each of $BC$ and $UC$ with the open unit disc has countably many
components.  The letters $B$ and $U$ are for ``below'' and ``upper''. 
The sets $BC$ and $UC$ are not part of the Markov partition, but
nevertheless we can, and shall, use them in the symbolic dynamics. For
$s_{1/7}$ we also have triangular regions labelled $\Delta $ and
$-\Delta $.  The labels $L_{j}$ are for ``left'', and $R_{j}$ for
``right''.  This obviously makes more sense for $s_{3/7}$, when the
sets certainly do run from left to right, than for $s_{1/7}$, when the
``left'' regions are on average only slightly left of the ``right''
regions, and more above than to the left of the ``right'' regions.

In the case of $s_{3/7}$, a more common partition might be that
generated by the leaf $\mu _{1/3}$ with endpoints $e^{\pm 2\pi 
i(1/3)}$.  This
translates to the Yoccoz partition for the polynomial which is
Thurston equivalent to $s_{3/7}$.  An important feature of the Yoccoz
partition is that, although the partition varies with the polynomial,
all the partitions for polynomials in a fixed limb have the same
generator.  Our choice of partition for $s_{3/7}$ is made for similar
reasons.  The same generator gives a partition which is valid in part
of parameter space.
The part of parameter space which motivates the
choice of partition is $V_{3,0}(a_{1})$, rather than a limb of the
Mandelbrot set.  The leaf $\mu _{3/7}$ with endpoints $e^{\pm 2\pi 
i(3/7)}$
generates a smaller number of partitions than $\mu _{1/3}$, and is 
more
closely adapted to $V_{3,0}(a_{1})$.

The regions in the Markov partition for $s_{3/7}$
are mapped as follows by $s_{3/7}$:
$$L_{1},\ R_{1}\to R_{1}\cup R_{2},$$
$$L_{2},\ R_{2}\to R_{3}\cup C,$$
$$L_{3},\ R_{3}\to L_{2}\cup L_{3},$$
$$C\to L_{1}.$$
The regions in the Markov partition for $s_{1/7}$  are mapped as 
follows by 
$s_{1/7}$:
$$L_{1},\ R_{1}\to L_{2},$$
$$L_{2},\ R_{2}\to R_{1}\cup R_{2}\cup -\Delta \cup C,$$
$$\Delta ,\  -\Delta \to \Delta .$$

We now consider the case of $L_{3/7}$ and $s_{3/7}$.
Infinite words in the letters $L_{i}$, $R_{j}$ and $BC$, $UC$, $C$ 
label points 
of the circle, leaves of $L_{3/7}$ and gaps of $L_{3/7}$. The labels 
determine 
gaps uniquely. Points and leaves also have unique labels except when 
they map forward to boundaries of sets in the Markov partition.
 Any  
gap  of $L_{3/7}$ is labelled by an infinite word which ends with the 
infinite word $(CL_{1}R_{2})^{\infty }$.   The preperiod  of $w$
is the length of $w'$ 
where $w=w'w''$ and either $w''=(CL_{1}R_{2})^{\infty }$ with $w'$ 
ending with $L_{2}$, or $w''=(R_{2}CL_{1})^{\infty }$ with $w'$ 
ending with $R_{1}$.  Occurrences of $BC$ or $UC$ are allowed in $w'$,
always followed by $L_{1}$ and preceded by $L_{2}$ or $R_{2}$.  In
fact, in $w'$, we never use the letter $C$, but always $BC$ or $UC$,
whichever is appropriate. A gap $G$ of preperiod $m$ contains a point 
$x$ of
$Z_{m}(s_{3/7})$, and is therefore labelled by the same word
$w(G)=w(x)$ as $x$.  Points of $Z_{m}(s)$ are uniquely determined by
the prefix of $w(x)$ of length $m$, which ends in $R_{1}$ or $L_{2}$. 
In fact, we shall usually use the prefix of $w(x)$ of length $m+2$ if
the length $m$ prefix ends in $R_{1}$, and of length $m+1$ if the
length $m$ prefix ends in $L_{2}$, so that points in $Z_{m}(s)$, and the gaps containing them, are
represented by nonempty words ending in $C$.  Thus, if $x$ is of
preperiod $\geq 1$, $w(x)$ ends in $R_{1}R_{2}C$ or $L_{2}C$.  The
gaps between the leaves with endpoints $(1/3,2/3)$ and $(2/7,5/7)$
have words starting $L_{3}^{k}L_{2}$ where $k$ is odd.  The gaps
between the leaves with endpoints $(2/7,5/7)$ and $(9/28,19/28)$ have
words starting $L_{3}L_{2}$.  Gaps which cross the real axis have
words $w'w''$, where $w'$ is the preperiodic part as above, with no
occurrence of $L_{1}$ in $w'$, equivalently no occurrence of $C$ or
$R_{1}$ or $R_{2}$ in $w'$.

In future, for a finite word $w$, we shall use the notation $D(w)$ to 
denote 
the subset of the disc of points $z$ labelled by $w$, that is with 
$s^{i-1}(z)\in X_{i}$ if $X_{i}$ is the $i$'th letter of $w$. If $w$ 
does not end  with a string containing only the letters $BC$, $UC$, 
$L_{1}$ and $R_{2}$ then $D(w)$ is bounded by one or 
two leaves of $L_{3/7}$. If $w$ does end with such a string, then 
$D(w)$ 
is a countable union of components of complement of a gap.

\section{The Main Theorem (first version)}\label{2.10}

We are now ready to state the first version of our main theorem. 
 The description of 
$V_{3,m}$ starts from the four subsets $V_{3,m}(a_{0})$, 
$V_{3,m}(\overline{a_{0}})$ and $V_{3,m}(a_{1},\pm )$ defined in 
\ref{2.2}, together with the corresponding sets $P_{m}(a_{0})$, and 
so on. The cases of $V_{3,m}(a_{0})$ and 
$V_{3,m}(\overline{a_{0}})$ are described in item 1 below and 
$V_{3,m}(a_{1},-)$ in item 2. These descriptions are quite simple: 
all the type III hyperbolic components in these regions are 
captures, essentially described in exactly one way. The description 
for $V_{3,m}(a_{1},+)$ is not so simple. The later versions of this 
theorem will be given in \ref{5.7}, \ref{6.1} and \ref{7.8}. It is not claimed 
that the 
results about $V_{3,m}(a_{0})$ and $V_{3,m}(\overline{a_{0}})$ are 
new, but the statement in the case of $V_{3,m}(a_{1},-)$ probably is.

We define 
$$q_{p}=\frac{1}{3}-2^{-p}\frac{1}{21},$$
so that $q_{0}=\frac{2}{7}$, $q_{1}=\frac{9}{28}$, and so on. In terms of the symbolic dynamics introduced in \ref{2.9}, the code of the leaf of $L_{3/7}$ with endpoints $e^{2\pi iq_p}$ is $L_3^{2p+1}(L_2R_3L_3)^{\infty }$.

\begin{maintheorem1}
There are injective maps 
$$a\mapsto \beta (a):P_{m}(a_{0})\to \pi _{1}(\overline{\mathbb 
C}\setminus Z_{m}(s_{1/7}),Z_{m}(s_{1/7}),v_{2}),$$
$$a\mapsto \beta (a):P_{m}(\overline{a_{0}})\to \pi 
_{1}(\overline{\mathbb C}\setminus 
Z_{m}(s_{6/7}),Z_{m}(s_{6/7}),v_{2}),$$
$$a\mapsto \beta (a):P_{m}(a_{1})\to \pi _{1}(\overline{\mathbb 
C}\setminus Z_{m}(s_{3/7}),Z_{m}(s_{3/7}),v_{2}),$$
such that $h_{a}$ is Thurston equivalent to $\sigma _{\beta (a)}\circ 
s$, for $s=s_{1/7}$, $s_{6/7}$, $s_{3/7}$ respectively.
Moreover, we have 
the following additional information.
\begin{itemize}
\item[1.]  For $a\in P_{3,m}(a_{0})$, the path $\beta (a)$ is a 
capture path, intersecting $S^{1}$ exactly once, in the boundary of 
the 
gap of $L_{1/7}$ containing the endpoint of $\beta $, and the 
possible 
endpoints of such paths $\beta $ are all points of $Z_{m}(s)$  
in the larger region of the unit disc bounded by the leaf with 
endpoints at $e^{2\pi i(1/7)}$, $e^{2\pi i(2/7)}$. Thus $h_{a}$ is 
Thurston equivalent to the capture $\sigma _{\beta }\circ s_{1/7}$.
A similar statement 
holds with $a_{0}$ replaced by $\overline{a_{0}}$, and with $1/7$ and 
$2/7$ 
replaced by $6/7$ and $5/7$. 

\item[2.] For $a\in P_{3,m}(a_{1},-)$, the path $\beta (a)$  is again a 
capture path intersecting $S^{1}$ exactly once, in the boundary of the 
gap of $L_{3/7}$ containing the endpoint of $\beta $, and the 
possible 
endpoints of such paths $\beta $ are all points of $Z_{m}(s)$  
in the smaller region of the unit disc bounded by the leaf with 
endpoints at $e^{2\pi i(1/7)}$ and $e^{2\pi i(6/7)}$.

\item  In cases 1 and 2, the path $\beta (a)$ is completely determined by its second endpoint. There is thus a natural one-to-one correspondence between the set $P_{m}(a_{0})$, $P_{m}(\overline{a_{0}})$, $P_{m}(a_{1},-)$ and the set of points of $Z_{m}(s)$ in the corresponding region of the dynamical plane of $s$, for $s=s_{1/7}$, $s_{6/7}$ and $s_{3/7}$ respectively.

\item[3.] For $a\in P_{m}(a_{1},+)$, the set of paths $\beta (a)$ includes the capture paths crossing $S^{1}$ at the points $e^{\pm 2\pi iq_{p}}$ into the gap of preperiod $2p$. There is a bijection between  $P_{3,m}(a_{1},+)$ and a set
$$\cup _{p=0}^{\infty }(U^{p}\cap Z_{m})\times \{ p\} .$$
Although the bijection does not send $a$ to the endpoint of $\beta (a)$, the image of $a$ under the bijection does determine the path $\beta (a)$ algorithmically. The set $U^{p}$ is defined as follows:
$$U^{0}=\begin{array}{l}D(L_{3}L_{2})\cup D(BC)\\
\cup \cup _{k=1}^{\infty }(D(L_3^{3k+1}L_{2})\cup D(L_{3}^{3k-1}L_{2}BC))\\ \setminus (\cup _{k\geq 1}D(L_{3}L_{2}(UCL_{1}R_{2})^{k}BC)),\end{array}$$
and if $p\geq 1$, and $S_{1,p,k}$ and $S_{2,p,k}$ denote the local inverses corresponding to the words $L_{3}^{2p+1+k(2p+3)}$ and $L_{3}^{2p-1+k(2p+3)}$, then
$$U^{p}=\begin{array}{l}\cup _{k=0}^{\infty }(S_{1,p,k}D(L_{2})\cup \cup _{k=0}^{\infty }(\cup _{n=0}^{\infty }(S_{2,p,k}D(u_{n})\cup(\cup _{1\leq t\leq p}S_{2,p,k}D(v_{t,n})))\\
\cup \cup _{n=0}^{\infty }S_{1,p,k}D((L_{2}R_{3})^{2}L_{3}^{2p-1}u_{n})\\
\setminus (\cup _{k=0}^{\infty }\cup _{n=0}^{\infty }(S_{1,p,k}D(u_{n})\cup S_{2,p,k}D((L_{2}R_{3})^{2}L_{3}^{2p-1}u_{n})\cup \cup _{t=1}^{p}S_{1,p,k}D(v_{t,n})))),\end{array}$$
where  $v_{t,n}=L_{2}(UCL_{1}R_{2})^{n}R_{3}L_{2}R_{3}L_{3}^{2t-1}L_{2}$ and $u_{n}=L_{2}(UCL_{1}R_{2})^{n}BC$.

 \end{itemize}\end{maintheorem1}

\chapter{Captures and Counting}\label{3}

\section{Some nontrivial equivalences between 
captures}\label{3.1}

There is a more nontrivial Thurston equivalence between captures 
$\sigma _{\beta }\circ s_{3/7}$ and a subset of the captures $\sigma 
_{\beta '}\circ s_{1/7}$ which has, nevertheless, been known since 
the 1980's. Its mating analogue, known as ``shared mating'' is 
somewhat better known and plays a role in Adam Epstein's 
(unpublished) proof of 
noncontinuity of mating. These nontrivial equivalences arise from a 
rather simple fact: the existence of two isotopically distinct 
circles $\gamma _{1}$ and $\gamma _{2}$ such that $\gamma 
_{j}=(s_{1/7}\Amalg s_{3/7})^{-1}(\gamma _{j})$, for both $j=1$, $2$, 
modulo isotopy fixing $X(s_{1/7}\Amalg s_{3/7})$. 
\begin{center}
\includegraphics[width=6cm]{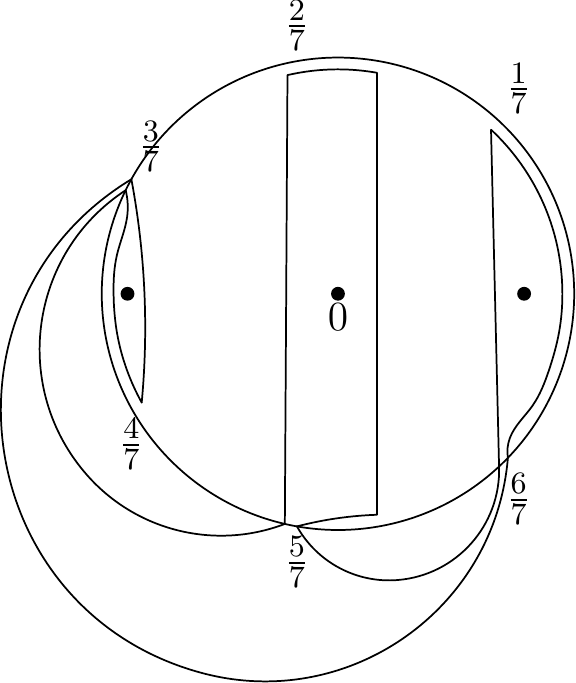}

The two invariant circles
\end{center}

An
isotopically invariant disc for a critically periodic branched 
covering $f$
gives
rise to a description of $f$ as a mating up to Thurston equivalence. 
This is because we can then make a new branched covering, keeping $f$
in the invariant disc and putting a single fixed
critical point in its complement.  It can then be shown quite easily that
there is no Thurston obstruction for this new branched covering, which
is therefore, by Thurston's theorem, equivalent to a polynomial. 
Applying this argument to the invariant disc and its complement, which is also an invariant disc, we see
that $f$ is Thurston equivalent to a mating of the two polynomials.

Note that the two circles $\gamma _1$ and $\gamma _2$ intersect in an approximate triangle, but that the points of $X(s_{1/7}\Amalg s_{3/7})$ on the same
side of $\gamma _{j}$ as the triangle are different for $\gamma _{1}$,
$\gamma _{2}$.  We shall call this approximate triangle $ET$.  It
follows that
$$s_{3/7}\Amalg s_{1/7}\simeq _{\varphi }s_{1/7}\Amalg s_{3/7}$$
for a homeomorphism $\varphi $ which maps the non-round circle 
$\gamma _{2}$ to the round circle $\gamma _{1}$ and maps 
$Y_{0}(s_{3/7})$ to $Y_{0}(s_{1/7})$.  This construction can be
generalised.  If we replace $1/7$ by any odd denominator rational $p$ 
with
$\mu _{p}\geq \mu _{1/7}$, then, similarly to the above, we obtain an
second non-round invariant circle for $s_{3/7}\Amalg s_{p}$which
separates the two periodic orbits, and hence gives rise to an
equivalence
$$s_{3/7}\Amalg s_{p}\simeq s_{1/7}\Amalg s_{q}$$
for some $q$.  This construction occurs in Wittner's thesis.

There is similar construction for captures. Take any type III
capture $\sigma _{p}\circ s_{3/7}$ for $p\in (5/7,6/7)$.  We are 
going to
show that this is Thurston equivalent to a capture $\sigma _{q}\circ 
s_{1/7}$ for some $q$.  The idea is to find
a tree which connects up the points of $Z_{n}(s_{3/7})$, if $e^{2\pi 
ip}$ is in the boundary of a gap of $L_{3/7}$ of preperiod $n$, and 
such that this tree is
naturally isomorphic to the $n$'th preimage of the Hubbard tree of
$s_{1/7}$.  Once we have this, the Thurston equivalence to a capture
$\sigma _{q}\circ s_{1/7}$ is automatic. Let $T_{0}$ be the union of 
$ET$
and the points of $Y_{0}(s_{3/7})$, joined to the triangle by arcs
which do not intersect $\gamma _{2}$.  Then $T_{0}\subset (\sigma
_{\beta }\circ s_{3/7})^{-1}(T_{0})=T_{1}$ up to isotopy preserving
$Y_{1}(s_{3/7})$.

Let $\Delta _{0}$ be a small tubular neighbourhood in the unit disc of
the triangular region $\Delta $ which features in the symbolic
dynamics for $s_{1/7}$ in \ref{2.9}, so that $Z_{0}(z_{1/7})\cap
\Delta _{0}=\emptyset $.  Write $\Delta _{n}=s_{1/7}^{-n}(\Delta 
_{0})$. Then
there is a homeomorphism $\varphi _{1}$ which maps $T_{1}$ to a subset
of the unit disc, consisting of $\Delta _{1}$ and arcs joining the
two components of $\Delta _{1}$ to the five points of 
$Y_{1}(s_{1/7})$,
and mapping $Y_{1}(s_{3/7})$ to $Y_{1}(s_{1/7})$. Similarly, for each 
$1<m\leq n$,
$$T_{0}\subset T_{m-1}\subset (\sigma _{p}\circ
s_{3/7})^{-1}(T_{m-1})=T_{m}{\rm{\ rel\ }}Y_{m}(s_{3/7}),$$
the set 
$T_{m}$ is connected, and there is a homeomorphism $\varphi _{m}$ 
mapping $T_{m}$ to the union of $\Delta _{m}$ and arcs joining up
components of $\Delta _{m}$ and points of $Y_{m}(s_{1/7})$.  If we
consider $n+1$, we still have
$$T_{n}\subset T_{n+1}{\rm{\ rel\ }}Y_{n+1}(s_{3/7}),$$
and there is a homeomorphism 
$\varphi _{n+1}$ such that 
$$\varphi _{n+1}(T_{n})=\varphi _{n}(T_{n}){\rm{\ rel\ 
}}Y_{n}(s_{1/7}),$$ and 
$$\varphi 
_{n+1}(T_{n})\subset \varphi _{n+1}(T_{n+1}){\rm{\ rel\ 
}} Y_{n}(s_{1/7}),$$
but, this time, $\infty \in T_{n+1}$ and $\infty \in \varphi 
_{n+1}(T_{n+1})$.  
Then the component $C$ of $\varphi _{n+1}(T_{n+1}\setminus T_{n})$ 
which 
contains $\infty $ also contains exactly two points of   
$Y_{n+1}(s_{1/7})\setminus Y_{n}(s_{1/7})$. If we join one of these 
points (it 
does not matter which) to $\infty $ by  arc $\zeta '$ in $C$, 
then $\beta '=s_{1/7}(\zeta ')$ joins $\infty $ to a point of 
$Z_{n}(s_{1/7})$, has only one $S^{1}$-crossing up 
to isotopy preserving $Y_{n}(s_{1/7})$, at some point $e^{2\pi iq}$. So $\beta '=\beta _{q}$ is a capture path, 
and 
$$\sigma _{p}\circ s_{3/7}\cong _{\varphi _{n}}\sigma _{q}\circ 
s_{1/7}.$$ 

Note that it is not stated that $\beta _{q}=\varphi _{n}(\beta _{p})$ 
up to 
isotopy preserving $Y_{n}(s_{1/7})$, and this is 
certainly not true in general. However, the proof 
above implies an algorithm for computing $q$ 
from $p$. The tree $T_{n}$ is defined using the map $\sigma _{p 
}\circ s_{3/7}$. We therefore sometimes use the more precise notation 
$T_{n}(\sigma _{\beta }\circ s_{3/7})$ or $T_{n}(p)$. We also use 
$T_{n}'(p)$ to denote the branches of $T_{n}$ which connect $ET$ to 
the forward orbit of $e^{2\pi ip}$, that is, the forward orbit of the 
gap with $e^{2\pi ip}$ in its boundary, the gap containing the second 
critical value. The map 
$\varphi _{n}$ also depends on $p$, and defines a correspondence 
between components of $(\sigma _{p }\circ s_{3/7})^{-n}(ET)$ and 
components of $s_{1/7}^{-n}(\Delta _{0})$. It also gives a 
correspondence between the gaps of $L_{3/7}$ of preperiod $n$ and the 
gaps of $L_{1/7}$ of preperiod $n$, in the full orbit of the critical 
gap. So $\varphi _{n}$ gives a relabelling of the gaps of $L_{3/7}$, 
in terms of the symbolic dynamics of $s_{1/7}$. 

We start with an example of the smallest possible preperiod,
$\sigma _{23/28}\circ s_{3/7}$, of preperiod two.  Note that
$\frac{23}{28}\in (\frac{5}{7},\frac{6}{7})$, and that the forward
orbit under $x\mapsto 2x{\rm{\ mod\ }}1$ is given by 
$$\frac{23}{28}\mapsto \frac{9}{14}\mapsto \frac{2}{7}.$$
With respect to the symbolic dynamics of $s_{3/7}$, the word of the 
gap
with $e^{2\pi i(23/28)}$ in its boundary is $R_{3}L_{2}C$.  Now we 
need to compute the word of this gap with respect to the symbolic
dynamics of $s_{2/7}$, using the tree $T_{2}$ and homeomorphism
$\varphi _{2}$ constructed above. In the following array, each column 
indicates a path of gaps in the tree $T_{n}$, starting from the gap 
$C$ at the bottom and working upwards. The word of each gap in a 
column is obtained by reading the word starting from that letter in 
the column, and including all letters to the right of it.  We are
using the convention introduced in \ref{2.9}, which allows us to use
finite words to represent gaps of $L_{3/7}$. 
The first row is the word of the gap
containing $e^{2\pi ip}$ in its boundary.  The $\uparrow $ symbol
indicates a connection across a component of $(\sigma _{p}\circ
s_{3/7})^{-n}(ET)$.  The gaps in each column are the inverse images
under $s_{3/7}$ (and also under $\sigma _{p}\circ s_{3/7})$ of the
gaps to the right.  The last two rows indicate symbolic
dynamics for $s_{1/7}$.  The last is a refined lettering of the
last-but-one.  The last-but-one row of letters $L$, $R$, $BC$
and $UC$ is obtained from the last connection in the row above.  The
letter $L$ is chosen when the letters above are $C$ with $R_{2}$ or
$L_{1}$ above this, and $R$ is chosen if the letter $C$ above has
$L_{2}$ or $R_{1}$ immediately above it.  The letters $BC$ and $UC$ do
not occur in this first example, but when they do, they are given by
$BC$ and $UC$ respectively in the row above.  This can be seen to be
correct by examining the initial picture of the two invariant circles
$\gamma _{1}$ and $\gamma _{2}$ for $s_{3/7}\Amalg s_{1/7}$.  
$$\begin{array}{lll}R_{3}&\ L_{2}&\ C\\
\uparrow &\ \uparrow &\ \\
R_{2}&\ C&\ \\
\uparrow &\ &\ \\
C&\ &| \\
L&\ R&\ C\\
L_{2}&\ R_{2}&\ C\\
 \end{array}$$
This gives
$$\sigma _{23/28}\circ s_{3/7}\simeq s_{\beta '}\circ s_{1/7}$$
where the gap containing the endpoint of $\beta '$ is encoded by $L_{2}R_{2}C$, using the symbolic dynamics of \ref{2.9} for $s_{1/7}$.

A capture $\sigma _{p}\circ s_{3/7}$ with $p\in
(\frac{5}{7},\frac{6}{7})$ must end in a gap whose word in the symbolic
dynamics for $s_{3/7}$ starts with the letter $R_{3}$ or $BC$. Now we 
consider
a capture path with endpoint in a gap whose word in the symbolic
dynamics of $s_{3/7}$ starts with $BC$.  The shortest possible word is
$BCL_{1}R_{1}R_{2}C$.  The gap has $e^{2\pi i(41/56)}$ in its
boundary.  The other end of the leaf in the gap boundary is at 
$e^{2\pi i(43/56)}$.  The gap has preperiod $3$. The computation is 
as follows.

$$\begin{array}{lllll}BC&\ L_{1}&\ R_{1}&\ R_{2}&\ C\\
\uparrow&\ \uparrow &\ \uparrow &\ \\
R_{3}&\ L_{2}&\ C&\ \\
\uparrow &\ \uparrow &\ &\ &\ \\
R_{2}&\ C&\ &\ &\ \\
\uparrow &\ &\ &\ &\ \\
C &\ &\ &\ &\ \\
L&\ R&\ R&\ L& C\\
L_{2}&\ R_{2}&\ R_{1}&\ L_{2}&\  C\\

\end{array}$$
So
$$\sigma _{41/56}\circ s_{3/7}\simeq \sigma _{\beta '}\circ s_{1/7}$$
where the gap of $L_{1/7}$ containing the second endpoint of $\beta '$ is encoded by $L_{2}R_{2}R_{1}L_{2}C$.

We next consider an example chosen at random. 
Take $\sigma _{p}\circ s_{3/7}$ with $p\in ({5\over 7},{6\over 7})$
with endpoint coded, with respect to the symbolic dynamics for
$s_{3/7}$, by $w(\beta ,{3\over 7})=R_{3}L_{3}L_{2}R_{3}L_{2}C$. 
The computation then uses the following array.

$$\begin{array}{llllll}\ R_{3}&\ L_{3}&\ L_{2}&\ R_{3}&\
L_{2}&\ C\cr \ \uparrow{\rm{top}}&\ \uparrow{\rm{top}}&\
\uparrow{\rm{bot}}&\ \uparrow{\rm{top}}&\ \uparrow{\rm{top}}&\ \cr \
R_{2}&\ UC&\ L_{1}&\ R_{2}&\ C&\ \cr \ \uparrow &\ \uparrow &\
\uparrow &\ \uparrow &\ &\ \cr \ C&\ C&\ C&\ C&\ &\ \cr \ &\ &\ &\ &\
&\ \cr \ L&\ UC&\ L&\ L&\ R&\ C\cr \ L_{2}&\ UC&\ L_{1}&\ L_{2}&\
R_{2}&\ C\cr  \end{array}$$
  So 
$$\sigma _{127/224}\circ s_{1/7}\simeq \sigma _{\beta '}\circ
s_{3/7},$$
where the gap of $L_{1/7}$ containing the second endpoint of $\beta '$ is encoded by $L_{2}UCL_{1}L_{2}R_{2}C$.

These calculations are fairly routine, and can easily be automated. 
They are perhaps not, however, quite as simple as they seem.  In the
examples above, $T_{n}'(p)$ coincides with path of 
$T_{n}(s_{3/7}\Amalg
s_{1/7})$, the tree used to compute the periodic-orbit-exchanging 
self-Thurston
equivalence of $s_{3/7}\Amalg s_{1/7}$ up to isotopy preserving the
$n$'th preimages of periodic orbits.  But $T_{n}'(p)$ does not always
so coincide.  A relatively small example of preperiod $8$ is given by 
the
gap with word $R_{3}L_{3}^{3}L_{2}R_{3}L_{3}L_{2}C$.  There is at least one point $e^{2\pi
ip}$  in the boundary of this gap with
$p\in (\frac{5}{7},\frac{6}{7})$.  This gap is between the gaps with words
$R_{3}L_{3}L_{2}C$ and $R_{3}L_{2}C$.  It follows that the component
$E'$ of $(\sigma _{p}\circ s_{3/7})^{-3}(ET)$ which joins the gaps
with words $R_{3}L_{3}L_{2}C$ and $R_{3}L_{2}C$ separates $e_{2\pi
ip}$ from $\infty $ in $\{ z:\vert z\vert \geq 1\} $.  Then $(\sigma
_{p}\circ s_{3/7})^{-1}(E')$ consists of two long thin triangles,
which are not simply tubular neighbourhoods of triangles of
$L_{1/7}^{-1}$, as was the case with all previous examples.  These
triangles are intersected by $T_{4}'(p)$.  It follows that $T_{4}'(p)$
is not contained in $T_{4}(s_{3/7}\Amalg s_{1/7})$ up to isotopy, and
a fortiore $T_{8}'(p)$ is not contained in $T_{8}(s_{3/7}\Amalg 
s_{1/7})$.

In summary, every type III capture $\sigma _{p}\circ s_{3/7}$ 
 with $p\in ({5\over 7},{6\over 7})$ is 
Thurston 
equivalent to a type III capture $\sigma_{\beta '}\circ s_{1/7}$, 
which can be algorithmically computed. Similarly, every type III 
capture $\sigma _{p}\circ s_{3/7}$ with 
$p\in ({1\over 7},{2\over 7})$ is Thurston equivalent to a type III 
capture $\sigma_{\beta '}\circ s_{5/7}$. So we have the following.

\begin{theorem}\label{3.2} Every type III capture  in $V_{3}$ is 
Thurston equivalent to one of the form $\sigma _{\beta }\circ s_{q}$, 
where $\beta $ has second endpoint  $x\in s_{q}^{-n}(\{ 
s_{q}^{j}(0):j=0,1,2\} )$, for some least $n>0$, and the gap containing 
$x$ has longest leaf with endpoint at $e^{2\pi ip}$, where $e^{2\pi 
ip}$ is the $S^{1}$-crossing of $\beta $ and
one of the following holds: 
\begin{itemize}
\item $q={1\over 7},\ \ p\in
({\textstyle{2\over 7}},{\textstyle{8\over
7}})$;
\item $q={6\over 7}$, $p\in ({\textstyle{-1\over 7}},\textstyle{5\over 
7})$;
\item $q={3\over 7}$, $p\in ({\textstyle{2\over 7}},{\textstyle{1\over 
3}})\cup 
({\textstyle{2\over 3}},{\textstyle{5\over 7}})\cup 
(-{\textstyle{1\over 
7}},{\textstyle{1\over 7}}) $.
\end{itemize}
In the cases $q={1\over 7}$ and ${6\over 7}$, and $q=\frac{3}{7}$ with $p\in (-\frac{1}{7},\frac{1}{7})$, we only need to include 
one path $\beta $ for each $x$. In the case of $q=\frac{3}{7}$ with $p\in (\frac{2}{7},\frac{1}{3})\cup (\frac{2}{3},\frac{5}{7})$ we also need to include at most one path, except in the case when the boundary of the gap containing $x$ intersects both the upper and lower halves of the unit circle, when we only need to include one path crossing the upper half circle and one crossing the lower half circle into the gap containing $x$.  \end{theorem}

In the cases $q={1\over 7}$ and  
${6\over 7}$, this theorem is 
subsumed in the statement of the main theorem in \ref{2.10}.

\section {Equivalences between captures $\sigma _{p}\circ 
s_{3/7}$, for $p\in ({2\over 7},{1\over 3})\cup ({2\over 3},{5\over 
7})$}\label{3.3}

 First we 
simplify 
the notation. We write  $L$ for  
invariant lamination $L_{3/7}$ and $s=s_{3/7}$.  Note that the only
gaps of $L$ lie in the full orbit of the critical gap. We also write 
$\sigma
_{q}$ for $\sigma _{\beta _{q}}$ (as in \ref{2.7}) or $\sigma _{\zeta
_{q}}^{-1}\circ \sigma _{\beta _{q}}$ if $\beta _{q}$ is a path which
crosses the unit circle at the point $e^{2\pi iq}$ into a gap $G$ of
$L$, with all of the path being in $G\cup \{ z:\vert z\vert >1\} $
apart from the point $e^{2\pi iq}$.  Here, we use $\sigma _{\zeta
_{q}}$ if and only if $G$ is in the periodic orbit of the critical 
gap, in
which case $\zeta _{q}$ is the path in $s^{-1}(\beta _{q})$ mapping
homeomorphically under $s$ and the endpoint of $\zeta $ is in the
periodic component of $s^{-1}(G)$.  We are interested only in 
captures for $q\in
[2/7,1/3)\cup (2/3,5/7]$.  Then we have
$$ 
\sigma _{5/7}\circ s\simeq \sigma _{2/7}\circ s,$$
and, for any $k\geq 1$,
$$\sigma _{1-q_{k}}\circ s\simeq \sigma _{q_{k}}\circ s,$$
where 
$$q_{k}={1\over 3}-2^{-2k}{1\over 21}.$$

These equivalences are obtained as follows.  Let $D_{q}$ be the
disc whose boundary intersects the
unit circle only at the points $e^{\pm 2\pi iq}$ and running along the
lamination leaf in the unit disc between these points, which contains the smaller segment of the unit disc cut off by this leaf, but which  
is also taken to include $v_{2}=c_{2}=s(c_{2})=\infty $ in its interior.  Let $\alpha =\alpha 
_{q}$ be the
clockwise loop in $D=D_{q}$ based at $v_{2}$ crossing $S^{1}$ only at
the points $e^{\pm 2\pi iq}$, running along the lamination leaf
between them.  Thus, $\alpha _{q}$ is homotopic to an arbitrarily
small perturbation of $\beta _{q}*\overline{\beta _{1-q}}$.  First we
claim that for any $q$,
\begin{equation}\label{3.3.1}(\sigma _{q_p}\circ 
s,Y_{0}(s))\simeq _{{\rm{identity}}}(\sigma _{1-q_p}\circ s,Y_{0}(s)).\end{equation}
We see this as follows.  The disc $D'=\overline{\mathbb C}\setminus
{\rm{int}}(D)$ contains no critical values.  So
$$s^{-1}(D')=D_{+}\cup D_{-}$$
has two components, homeomorphic preimages under $s$ of $D'$, where $D_{+}=D_{+,q}$ contains 
$v_{1}$
and $D_{-}=D_{-,q}$ contains $s(v_{1})$.  These discs are 
peripheral in
$Y_{0}(s)$.  So Dehn twists round their boundaries are trivial up to
isotopy preserving $Y_{0}(s)$. A clockwise Dehn twist round a simple closed loop $\gamma $ on an orientable surface $S$ means a homeomorphism up to isotopy which can be chosen to be the identity outside an arbitrarily small annulus neighbourhood of $\gamma $, and twists one boundary component $\gamma _{1}$ of the annulus clockwise relative to the other one, $\gamma _{2}$. Since the annulus inherits an orientation from $S$, a clockwise twist of $\gamma _{1}$ relative to $\gamma _{2}$ is the same as a clockwise twist of $\gamma _{2}$ relative to $\gamma _{1}$.  Note that $\partial D=\partial D'$. So (anti)clockwise Dehn twist round $\partial D$ is the same  as (anti)clockwise Dehn
twist round  $\partial D'$.  Now $\sigma _{\alpha }$ is isotopic to a composition
of clockwise Dehn twist round  $\partial D$ and anticlockwise Dehn twist round the boundary of a disc
strictly inside $\alpha $.  Since this inner disc is peripheral, we
obtain
\begin{equation}\label{3.3.2}(\sigma _{\alpha }\circ 
s,Y_{0}(s))\simeq _{{\rm{identity}}}(s,Y_{0}(s)).\end{equation}
In the case when $q=q_{0}=\frac{2}{7}$, since $\alpha $ is a 
perturbation of $\beta _{2/7}*\overline{\beta _{5/7}}$, and since 
$\zeta _{2/7}$ and $\zeta _{5/7}$ are homotopic via a homotopy 
preserving $Y_{0}(s)$, we have, composing on the left with $\sigma 
_{\zeta _{2/7}}^{-1}\circ \sigma _{\beta _{5/7}}$,
$$(\sigma _{2/7}\circ s, Y_{0}(s))=(\sigma _{\zeta _{2/7}}^{-1}\circ 
\sigma _{\beta _{2/7}}\circ s,Y_{0}(s))\simeq 
_{{\rm{identity}}}(\sigma _{\zeta _{2/7}}^{-1}\circ \sigma _{\beta 
_{5/7}}\circ s,Y_{0}(s))$$
$$\simeq _{{\rm{identity}}}(\sigma _{\zeta _{5/7}}^{-1}\circ \sigma 
_{\beta _{5/7}}\circ s,Y_{0}(s))=(\sigma_{5/7}\circ s,Y_{0}(s)).$$
So this gives (\ref{3.3.1}) in the case $q=q_{0}$. For any $k>0$, 
and $q=q_{k}$, the disc inside $\beta _{q}*\overline{\beta _{1-q}}$ 
and $D_{-,q}$, are both disjoint from the forward orbit of the second 
endpoint of $\beta _{q_{k}}$. So in these cases also, composing 
(\ref{3.3.2}) on the left with $\sigma _{\beta _{1-q}}$ gives 
(\ref{3.3.1})

We can express (\ref{3.3.2}) more accurately.  We write $\psi 
_{0,q}$ for the anticlockwise twist round the boundary of a disc  which is strictly inside 
$\alpha _{q}$, but with boundary arbitrarily close to $\alpha _{q}$. 
Then $\psi _{0,q}$ commutes with $\sigma _{\alpha _{q}}$. Also, 
$\sigma _{\alpha _{q}}^{-1}\circ \psi _{0,q}$ is isotopic, relative to 
$Y_{n}(s)$ for any $n$, to the anticlockwise Dehn twist round  $\partial D=\partial D_{q}$, 
which as already
mentioned, is isotopic to the anticlockwise Dehn twist round  $\partial D'$.
Let $\xi _{0,q}$ denote 
the anticlockwise
Dehn twist round  $\partial D_{-}$.  
Note that $D_{-,q}$ is peripheral in $Y_{0}(s)$, but not in
$Y_{1}(s)$, and $D_{+,q}$ and $D$ are isotopic in
$\overline{\mathbb C}\setminus Y_{1}(s)$, in fact in 
$\overline{\mathbb C}\setminus Y_{2k+1}(s)$ if $q=q_{k}$ -- but not in 
$\overline{\mathbb C}\setminus Y_{2k+2}(s)$.  So
$\xi _{0,q}$ is isotopic to the identity relative to $Y_{0}(s)$, and
$$\psi _{0,q}\circ s=\sigma _{\alpha _{q}}\circ s\circ \psi _{0,q}\circ
\xi _{0,q}{\rm{\ rel\ }}Y_{1}(s).$$
For any $q=q_{k}$, write
$$\psi _{1,q}=\xi _{0,q}\circ \psi _{0,q}=\psi _{0,q}\circ \xi _{0,q}.$$
Then we have
$$(s,Y_{1}(s))\simeq _{\psi _{1,q}}(\sigma _{\alpha }\circ 
s,Y_{1}(s)).$$
  For $m<2k+1$, we then define $\xi _{m,q}$ inductively by 
\begin{equation}\label{3.3.3}s\circ
\xi _{m,q}=\xi _{m-1,q}\circ s,\end{equation}
choosing the lift $\xi
_{m,q}$ of $\xi _{m-1,q}$ which is isotopic to the identity relative
to $Y_{m}(s)$.  This is possible inductively, because $\xi _{0,q}$ is
isotopic to the identity relative to $Y_{0}(s)$: we just lift this
isotopy.  Then we define $\psi _{m+1,q}$ inductively by
$$\psi _{m+1,q}=\xi _{m,q}\circ \psi _{m,q}.$$
It is then true, for $m\leq 2k+1$, that
\begin{equation}\label{3.3.4}\psi _{m,q}\circ s=\sigma _{\alpha _{q}}
\circ s\circ \psi _{m+1,q}{\rm{\ rel\ }}Y_{m+1}(s).\end{equation}
For $m=2k$, since the support of $\psi _{2k}$ is disjoint from $\beta
_{1-q_{k}}$,  we obtain, by composing on the left with $\sigma
_{1-q_{k}}^{-1}$, for $q=q_{k}$,
\begin{equation}\label{3.3.5}\sigma
_{1-q}\circ s\simeq _{\psi _{2k,q}}\sigma _{q}\circ s.\end{equation}

For $m=2k+1$ and $q=q_{k}$, we define $\xi _{2k+1,q}$  slightly 
differently. Write
$$A_{q}=A_{2k+1,q}=D_{q}\setminus D_{+,q}.$$
The index $2k+1$ is chosen because $A_{2k+1,q}$ is disjoint from 
$Y_{2k+1}(s)$. Let $\xi _{2k+1,1,q}$ denote the clockwise twist 
of 
the annulus $A_{2k+1,q_k}$ relative to $\overline{\mathbb C}\setminus 
A_{2k+1,q}$. Then we define $\xi _{2k+1,q}$ by
$$\xi _{2k,q}\circ s\circ \xi _{2k+1,1,q}=s\circ \xi _{2k+1,q}$$ 
and $\xi _{2k+1,q}$ is isotopic to the identity relative to 
$Y_{2k+1,q_k}$. Hence, (\ref{3.3.4}) holds for $m=2k+1$, and this yields
$$\psi _{2k+2,q}\circ s=\sigma _{\alpha _{q}}\circ s\circ \psi 
_{2k+2,q}{\rm{\ rel\ }}Y_{2k+2}(s).$$
Composing on the right with $\psi _{2k+2,q}^{-1}$, this gives
\begin{equation}\label{3.3.6}(s,Y_{2k+2}(s))\simeq _{\psi 
_{2k+2,q}}(\sigma _{\alpha }\circ s,Y_{2k+2}(s)).\end{equation}

The support of $\xi _{m,q}$ is as follows, up to isotopy preserving 
$Y_{m}(s)$. Let $C_{q}= C_{0,q}$ be the annulus which is the 
complement, in $D_{-}(q)$, of a small neighbourhood of $s(v_{1})$. 
For $m\leq 2k+1$, define 
$$C_{m,q}=s^{-m}(C_{q}).$$
Note that $C_{m,q}$ is disjoint from $Y_{m}(s)$ for all $m$, and, 
for $m\leq 2k+1$, $C_{m,q}$ is disjoint from $\alpha _{q}$. 
For $m<2k+1$, the support of $\xi _{m,q}$ is $C_{m,q}$, up to 
isotopy. The support of $\xi _{2k+1,q}$ is $C_{2k+1,q}\cup 
A_{2k+1,q}$, again up to isotopy.

It is also interesting to compare $\psi _{2k+2,q_{k}}$ and $\psi _{2k+2,q_{k+1}}$. 
Note that
$$D_{q_{k}}=A_{q_{k}}\cup D_{q_{k+1}}.$$
The overlap between $A_{q_{k}}$ and $D_{q_{k+1}}$ is trivial in $Y_{2k+2}$.
Therefore
\begin{equation}\label{3.3.7}\psi _{0,q_{k+1}}=\xi _{2k+2,2,q_{k}}\circ \psi _{0,q_{k}},\end{equation}
where the support of $\xi _{2k+2,2,q_{k}}$ is strictly smaller than the support $A_{q_{k}}$ of $\xi _{2k+1,1,q_{k}}$, and is trivial in $Y_{2k+2}$. So 
$$[\psi _{0,q_{2k+2}}]=[\psi _{0,q_{k}}]{\rm{\ in\ MG}}(\overline{\mathbb C},Y_{2k+2}(s)).$$
Also, for $m\leq 2k$,
$$\xi _{m,q_{k}}=\xi _{m,q_{k+1}}{\rm{\ rel\ }}Y_{2k+2}(s).$$
and
$$\xi _{2k+1,q_{k}}\circ \xi _{2k+1,1,q_{k}}^{-1}=\xi _{2k+1,q_{k+1}}{\rm{\
rel\ }}Y_{2k+2}(s).$$
So, for $m\leq 2k+2$, 
\begin{equation}\label{3.3.8}[\psi _{m,q_{k+1}}]=[\psi _{m,q_{k}}]{\rm{\ in\ MG}}(\overline{\mathbb C}, Y_{2k+2}(s)).\end{equation}
In particular, using this for $m=2k+2$, $q=q_{k}$ and $r=q_{k+1}$, we have, from
(\ref{3.3.5}) for $r$ replacing $q$,
\begin{equation}\label{3.3.9}\sigma _{1-r}\circ s\simeq _{\psi
_{2k+2,q}}\sigma _{r}\circ s.\end{equation}

We can generalise  to some other $q$, as follows. Let $G$ be the gap 
 containing the endpoint of $\beta _{q}$. Suppose that the word $w$ 
 for $G$ (see \ref{2.9}) 
starts $L_{3}^{2k+1}L_{2}$ for some integer $k\geq 0$, and has no 
other occurrence of $L_{3}^{2k+1}$, and no occurrence of $BC$ or $UC$ 
before the final $C$. Then (\ref{3.3.3}), (\ref{3.3.4}) hold for 
all $m$ less than the preperiod of the endpoint 
of $\beta _{q}$ (equivalently of $e^{2\pi iq}$), and (\ref{3.3.4}) 
also holds when $m$ is this preperiod, and (\ref{3.3.6}) holds with 
$2k+2$ replaced by this preperiod.

\section{Numbers of type III components}\label{3.4}

In this section, for each $m$, we calculate the number of type III 
components of preperiod $m$ and, and give an upper bound on the 
number of captures of preperiod $m$ in $V_{3}=\{ h_{a}:a\in {\mathbb 
C}\} $. One basic point of interest is that the second number is less 
than the first, for all $m\geq 3$, and asymptotically the proportion 
of preperiod $m$ hyperbolic components which are captures is strictly 
less than $1$. The exact proportion is an open question, as it seems 
to be a hard problem to identify all Thurston equivalences between 
captures. Anyway, these  calculations give a useful check on the 
validity of Theorem 
\ref{2.10}.
Each type III component in $V_{3}$ has a unique centre $h_{a}\in 
P_{3,m}\setminus P_{3,m-1}$ for which 
the second critical value 
$$v_{2}(h_{a})={-(a-1)^{2}\over 4a}\in 
Z_{m}(h_{a})\setminus 
Z_{m-1}(h_{a})=h_{a}^{-(m-1)}(h_{a}^{-1}(Z_{0}(a))\setminus 
Z_{0}(a)),$$
where $Z_{0}(h_{a})=\{ 0,1,\infty \} $, for all $a\in \mathbb 
C\setminus  \{ 0\} $, for a 
unique integer $m>0$. We then have $h_{a}^{m}(v_{2})=0$ or $1$, and 
we write $P_{m}'(0)$ or $P_{m}'(1)$ for the corresponding sets of 
punctures, so that 
$$P_{3,m}\setminus P_{3,m-1}=P_{m}'(0)\cup P_{m}'(1).$$
Now
$$h_{a}^{-1}(\{ 0\} )=\{ 1,a\} ,\ \ h_{a}^{-1}(\{ 1\} =\left \{ 
\infty 
,{a\over a+1}\right \} .$$
So
$$P_{m}'(0)=\left \{ a:h_{a}^{m-1}\left( {-(a-1)^{2}\over 
4a}\right)-a=0,\ 
a\neq 0,1\right \} ,$$
$$P_{m}'(1)=\left \{ a:h_{a}^{m-1}\left( {-(a-1)^{2}\over 
4a}\right) -{a\over a+1}=0,\ a\neq 0,-1\right \} .$$
So each of $P_{m}'(0)$ and $P_{m}'(1)$ is a zero set of a polynomial, 
and $\# (P_{m}'(0))$, 
$\# (P_{m}'(1))$ are the respective degrees of these polynomials, if 
all zeros of the polynomials are simple. In fact, the zeros are 
necessarily simple, as pointed out to me by Jan Kiwi. Intersections 
between $V_{3}$ and 
$$W=\{ f_{c,d}:f_{c,d}^{m}(v_{2})=c_{1}\} $$
are transversal and similarly for $c_{1}$ replaced by 
$f_{c,d}(v_{1})$. Given an intersection point between $W$ and 
$V_{3}$, a transversal $F:\{ \lambda :\vert \lambda \vert <1\} \to W$ 
can be constructed so that $F(0)$ is the intersection point, and the 
multiplier of a period $3$ point of $F(\lambda )$ is $\lambda $. We 
omit the details, but this is a standard technique which uses the 
Measurable Riemann Mapping Theorem, and used by Douady and Hubbard, 
for example. 
The polynomials are divisors of, respectively,
$$p_{m}(a)-aq_{m}(a),\ \ (a+1)p_{m}(a)-aq_{m}(a),$$
where
$$h_{a}^{m}\left({2a\over a+1}\right)={p_{m}(a)\over 
q_{m}(a)}.$$
and $p_{m}$ and $q_{m}$ have no common factors. 
Note that if $a=0$ or $1$, then $\displaystyle{a={2a\over a+1}}$, so 
that 
$\displaystyle{h_{a}^{m}\left( {2a\over a+1}\right) =a}$ whenever $m$ 
is 
divisible by $3$. Similarly if $a=0$ or $-1$ then 
$\displaystyle{{a\over 
a+1}={2a\over 
a+1}}$, so that $\displaystyle{h_{a}^{m-1}\left({-(a-1)^{2}\over 
4a}\right) 
={a\over a+1}}$ whenever $m$ is divisble by $3$. So $\# 
(P_{m}'(0))$, 
$\# (P_{m}'(1))$ are the respective degrees of these polynomials when 
$m$ is not divisible by $3$, and the respective degrees minus $2$ if 
$m$ is divisible by $3$.
Now if $a$ is not a factor of $p_{m}$, we have
$$p_{m+1}(a)=(p_{m}(a))^{2}-(a+1)p_{m}(a)q_{m}(a)+a(q_{m}(a))^{2},$$
$$q_{m+1}(a)=(p_{m}(a))^{2}.$$
If $a$ is a factor of $p_{m}$, then the only possible common factor 
of 
the two polynomial expressions above is $a$, and so these expressions 
give $ap_{m}(a)$ and $aq_{m}(a)$. Note that $a\vert p_{0}$, or 
$p_{0}=0{\rm{\ mod\ }}a$. A simple induction shows that if 
$p_{m}=0{\rm{\ mod\ }}a$ then $q_{m+1}=0{\rm{\ mod\ }}a$, 
$p_{m+2}=q_{m+2}{\rm{\ mod\ }}a$ and $p_{m+3}=0{\rm{\ mod\ a}}$. So 
$a$ is a factor of $p_{m}$ if and only if $3\vert m$.
Let $c_{m}$, $d_{m}$ be the coefficients of the highest order terms 
of $p_{m}$ and $q_{m}$ respectively. 
For all $m$, we trivially have
$$d_{m+1}=c_{m}^{2}>0.$$
Then a simple induction shows that
$$\begin{array}{llll}{\rm{deg}}(p_{m})={\rm{deg}}(q_{m}),\ 
&0<2d_{m}\leq c_{m},\ &c_{m+1}=d_{m}^{2}-c_{m}d_{m}&{\rm{\ if\ 
}}m{\rm{\ 
is\ 
even,}}\cr 
{\rm{deg}}(p_{m})={\rm{deg}}(q_{m})+1,\ 
& -d_{m}<c_{m}<0,\ &c_{m+1}=c_{m}^{2}-c_{m}d_{m}&{\rm{\ if\ }}m{\rm{\ 
is\ 
odd.}}\cr 
\end{array}$$
Let $a_{m}={\rm{deg}}(p_{m})$. Then we have
$$a_{m+1}=\begin{cases}2a_{m}{\rm{\ if\ }}m=0,\ 1,\ 5{\rm{\ mod\ 
}}6,\\ 
2a_{m}-1{\rm{\ if\ }}m=3{\rm{\ mod\ }}6,\\
2a_{m}+1{\rm{\ if\ }}m=2,\ 4{\rm{\ mod\ }}6. \end{cases}$$
Then
$$a_{6k+6}=2a_{6k+5}=4a_{6k+4}+2=8a_{6k+3}-2=16a_{6k+2}+6=64a_{6k}+6$$
$$=2^{6k+6}\left(  1+{2\over 21}\right) -{2\over 21}.$$
So for all $m$,
$$a_{m}=2^{m}\left(  1+{2\over 21}\right) +O(1).$$
More precisely,
$$a_{m}=2^{m}\left(  1+{2\over 21}\right) \begin{cases}-{2\over 
21}& {\rm{\ if\ }}m=0{\rm{\ mod\ }}6,\\
-{4\over 21}{\rm{\ if\ }}m=1{\rm{\ mod\ }}6,\\
-{8\over 21} {\rm{\ if\ }}m=2{\rm{\ mod\ }}6,\\
+{5\over 21} {\rm{\ if\ }}m=3{\rm{\ mod\ }}6,\\
-{11\over 21} {\rm{\ if\ }}m=4{\rm{\ mod\ }}6,\\
-{1\over 21} {\rm{\ if\ }}m=5{\rm{\ mod\ }}6. \end{cases}$$
For example:
$$a_{1}=2,\ a_{2}=4,\ a_{3}=9,\ a_{4}=17,\ a_{5}=35,\ a_{6}=70\cdots 
$$
We also have 
$${\rm{deg}}(p_{m}(a)-aq_{m}(a))=\begin{cases}a_{m}+1{\rm{\ if\ 
}}m{\rm{\ is\ even,}}\\ a_{m}{\rm{\ if\ }}m{\rm{\ is\ odd.}}
\end{cases}$$

So 
\begin{equation}\label{3.4.1}\# (P_{m}'(1))=2^{m}\left(  1+{2\over 
21}\right) 
    \begin{cases}-{23\over 
21} {\rm{\ if\ }}m=0{\rm{\ mod\ }}6,\\
+{17\over 21} {\rm{\ if\ }}m=1{\rm{\ mod\ }}6,\\
+{13\over 21} {\rm{\ if\ }}m=2{\rm{\ mod\ }}6,\\
-{16\over 21} {\rm{\ if\ }}m=3{\rm{\ mod\ }}6,\\
+{10\over 21}{\rm{\ if\ }}m=4{\rm{\ mod\ }}6,\\
+{20\over 21} {\rm{\ if\ }}m=5{\rm{\ mod\ }}6. 
\end{cases}\end{equation}

Similarly, 
\begin{equation}\label{3.4.2}\# (P_{m}'(0))=2^{m}\left(  1+{2\over 
21}\right) 
    \begin{cases}-{23\over 
21}{\rm{\ if\ }}m=0{\rm{\ mod\ }}6,\\
-{4\over 21} {\rm{\ if\ }}m=1{\rm{\ mod\ }}6,\\
+{13\over 21} {\rm{\ if\ }}m=2{\rm{\ mod\ }}6,\\
-{37\over 21} {\rm{\ if\ }}m=3{\rm{\ mod\ }}6,\\
+{10\over 21} {\rm{\ if\ }}m=4{\rm{\ mod\ }}6,\\
-{1\over 21} {\rm{\ if\ }}m=5{\rm{\ mod\ }}6. 
\end{cases}\end{equation}
For example, starting from $m=1$, the values of 
$\# (P_{m}'(1))$ are 
$3$, $5$, $8$, $18$, $36$, $69$ \ldots 
while the values of $\# 
(P_{m}'(0))$ are $2$, $5$, $7$, $18$, 
$35$, $69$ \ldots 

So, adding the values of $\# (P_{m}'(1))$ and $\# 
(P_{m}'(0))$, the number of type III hyperbolic components 
of 
preperiod $m$ is $2^{m+1}(1+{2\over 21})+O(1)$. 
We note in passing that centres of type IV components in $V_{3}$ of 
periods dividing $m$ are zeros of $p_{m}$. So the  number of 
type IV hyperbolic components of period dividing $m$, and even of 
period exactly $m$, is $2^{m}(1+{2\over 21})(1+o(1))$.

We now consider the number of type III capture hyperbolic components. 
In fact we shall not do exactly this, because as already stated in 
\ref{3.2}, we do not have an asymptotic formula for the number of 
distinct Thurston equivalence classes of preperiod $m$ among captures 
of the form $\sigma _{p}\circ s_{3/7}$ for $p\in ({2\over 7},{1\over 
3})\cup ({2\over 3},{5\over 7})$. We shall actually give a 
recursive 
formula for the number 
of points in $Z_{m}(s)\setminus Z_{m-1}(s)$ in specified subsets of 
the unit disc, and also give asymptotes of the number as $m\to \infty 
$. This will at least give a bound on the number of capture 
components, by Theorem \ref{3.2}.
We only have three lamination maps to consider: 
$s_{1/7}$, $s_{3/7}$ and $s_{6/7}$. The descriptions of captures up 
to equivalence in sections \ref{3.1} to \ref{3.3} suggest that we 
should compute the numbers 
$$a_{m,x,q}=\# (\{ z\in D(q):s_{q}^{m}(z)=x,\ s_{q}^{m-1}(z)\neq 
s_{q}^{2}(x)\} ),$$
where $x=c_1$ or $s_q(v_1)$, and  $D(q)$ is the following subset of the unit disc, for each of 
$q=1/7$, $3/7$, $6/7$. If $q=1/7$ or $q=6/7$, we take 
$$D(q)=D(L_{2})\cup D(C)\cup \Delta \cup D(R_{1})\cup 
D(R_{2}),$$ 
for $C$, 
$\Delta $, $R_{1}$, $R_{2}$ as in \ref{2.9} for $s_{1/7}$ (or 
$s_{6/7}$). For $q={3\over 7}$, we take $D(3/7)=D(3/7,1)\cup 
D(3/7,-1)$ where $D(3/7,-1)=D(R_{1})\cup D(R_{2})$, using the 
Markov 
partition of \ref{2.9} for $s_{3/7}$, and $D(3/7,1)\subset D(L_{3})$ is 
the 
region bounded by the leaves of $L_{3/7}$ with endpoints at $e^{\pm 
2\pi i(1/3)}$ and  $e^{\pm 2\pi i(2/7)}$. In terms of the symbolic 
dynamics,
$$D(3/7,1)=\cup _{k=0}^{\infty }D(L_{3}^{2k+1}L_{2}).$$
We then define, for $j=\pm 1$, $x=s_{3/7}(v_{1})$ or $c_{1}$,
$$a_{m,x,3/7,j}=\# (\{ z\in D(q,j):s_{q}^{m}(z)=x,\ 
s_{q}^{m-1}(z)\neq 
s_{q}^{2}(x)\} ).$$

First we consider the calculation of $a_{m,x,q}$ for $q=\frac{1}{7}$, $\frac{6}{7}$, and $x=c_1$ or $s_q(v_1)$. It turns out that the number $a_{m,x,q}$ is the same for all four of these choices. We 
consider the Markov partition into the sets
$$L_{1},\ L_{2},\ C\cup \Delta \cup R_{1}\cup R_{2}=X.$$
The map $s_{q}$, for $q={1\over 7}$ or ${6\over 7}$,  maps these 
partition elements as follows:
$$L_{1}\to L_{2},\ L_{2}\to X,\ X\to X\cup L_{2}\cup 2L_{1},$$
where $X\to 2L_{1}$ means that each point in $L_{1}$ has two 
preimages 
under $s$ in $X$, while $X\to X$ (for example) means that each point 
in $X$ has just one preimage in $X$. The corresponding matrix is
$$A_{q}=\begin{pmatrix} 0& 1& 0\cr 0& 0& 1\cr 2& 1&1\cr \end{pmatrix}{\rm{\ for\ which\ 
}}A_{q}^{2}=\begin{pmatrix}0&0&1\cr 2& 1& 1\cr 2& 3& 2\cr \end{pmatrix},$$
where the first row of $A_{q}$ indicates the degree of the map of 
$L_{1}$ to each 
of $L_{1}$, $L_{2}$ and $X$, and so on. The only nonzero entry of 
$A_{q}$ in the 
first row, for example, is the second one, because $L_{1}$ maps onto 
$L_{2}$, and the $1$ indicates that the map onto $L_{2}$ is degree 
$1$. Note that all columns sum to $2$, because $s_{q}$ is degree two 
overall. Therefore, $2$ is an eigenvalue, with eigenvector
$$v_{q}=\begin{pmatrix}1\cr 2\cr 4\cr \end{pmatrix},$$
where this is easily deduced from either straight computation or from 
the relative lengths of $L_{1}\cap S^{1}$, $L_{2}\cap S^{1}$ and 
$X\cap S^{1}$. The characteristic polynomial of $A_{q}$ is 
$$(\lambda -2)(\lambda ^{2}+\lambda +1).$$
The nullspace of $A_{q}^{2}+A_{q}+I$, which, of course, is also the 
nullspace of $A_{q}^{3}-I$, is
$$\left\{ \underline{x}=\begin{pmatrix}x_{1}\cr x_{2}\cr 
x_{3}\cr\end{pmatrix}:x_{1}+x_{2}+x_{3}=0.\right\} .$$
Now $a_{m,x,q}$ is the sum of the second and third entries of
$$A_{q}^{m-1}\begin{pmatrix} 0\cr 0\cr 1\cr\end{pmatrix}=2^{m-1}\begin{pmatrix}{1\over7}\cr 
{2\over 7}\cr {4\over 7}\cr \end{pmatrix}+A_{q}^{m-1}\begin{pmatrix}-{1\over7}\cr 
-{2\over 7}\cr {3\over 7}\cr \end{pmatrix}.$$
Now
$$A_{q}\begin{pmatrix}-{1\over7}\cr 
-{2\over 7}\cr {3\over 7}\cr \end{pmatrix}=\begin{pmatrix}-{2\over 7}\cr {3\over 7}\cr 
-{1\over 7}\cr\end{pmatrix},\ \ A_{q}^{2}\begin{pmatrix}-{1\over7}\cr 
-{2\over 7}\cr {3\over 7}\cr \end{pmatrix}=\begin{pmatrix}{3\over 7}\cr -{1\over 7}\cr 
-{2\over 7}\cr \end{pmatrix}.$$
So
$$a_{m,x,q}=2^{m}.{3\over 7}\begin{cases}-{3\over 7}{\rm{\ if\ 
}}m=0{\rm{\ mod\ }}3,\\ +{1\over 7}{\rm{\ if\ 
}}m=1{\rm{\ mod\ }}3,\\ +{2\over 7}{\rm{\ if\ 
}}m=2{\rm{\ mod\ }}3, \end{cases}$$

Now we consider the calculation of $a_{m,x,3/7,1}$ and $a_{m,x,3/7,-1}$ for $x=c_1$ or $s_{3/7}(v_1)$. We use the partition of $D(3/7)$ as the union of the sets  $D(3/7,1)$ and $D(3/7,-1)$, and, 
correspondingly $a_{m,x,3/7}=a_{m,x,3/7,1}+a_{m,x,3/7,-1}$. 
It would be possible to do the computation using a different Markov 
partition from the one in \ref{2.9}, and, given the nature of the set 
$D(3/7,1)$, it might seem more natural to do so, But in fact, 
because of a calculation we shall do later, in \ref{3.5}, it makes 
sense to use the partition we already have. In fact, we shall use the 
coarser partition into the four sets
$$L_{1},\ L_{2}\cup L_{3}=X_{1},\ C\cup R_{3}=X_{2},\ R_{1}\cup 
R_{2}=X_{3}.$$
This time, $s_{3/7}$  maps these 
partition elements as follows:
$$L_{1}\to X_{3},\ X_{1}\to X_{1}\cup X_{2},\ X_{2}\to 2L_{1}\cup 
X_{1}, X_{3}\to X_{3}\cup X_{2}.$$
The corresponding matrix is
$$A_{3/7}=\begin{pmatrix}0& 0& 0& 1\cr 0&1&1&0\cr 2&1&0&0\cr 0& 0& 1& 
1\cr\end{pmatrix},{\rm{\ with\ }}A_{3/7}^{2}=\begin{pmatrix}0&0&1&1\cr 2&2&1&0\cr 
0&1&1&2\cr 2&1&1&1\cr \end{pmatrix}.$$
Once again, $2$ is an eigenvalue, this time with eigenvector
$$v_{3/7}=\begin{pmatrix}1\cr 2\cr 2\cr 2\cr\end{pmatrix}.$$
The characteristic polynomial is
$$(A_{3/7}-2)(A_{3/7}^{3}-1).$$
Once again, the nullspace of $A_{3/7}^{3}-I$ is the set of vectors 
with coefficients summing to $0$. Now $a_{m,3/7,s(v_{1}),-1}$ and 
$a_{m,c_{1},3/7,-1}$ are the fourth entries of, respectively,
$$A_{3/7}^{m-1}\begin{pmatrix}0\cr 0\cr 0\cr 1\cr\end{pmatrix},\ \ 
A_{3/7}^{m-1}\begin{pmatrix}0\cr 1\cr 0\cr 0\cr\end{pmatrix}.$$ 
These column vectors are chosen because the nonperiodic preimages of $s(v_1)$ and $c_1$ are in $X_3$ and $X_1$ respectively, that is, the fourth and second elements of the Markov partition currently being used. Now 
$$A_{3/7}^{m-1}\begin{pmatrix}0\cr 0\cr 0\cr 1\cr\end{pmatrix}=2^{m-1}\begin{pmatrix}{1\over 
7}\cr {2\over 7}\cr {2\over 7}\cr {2\over 7}\cr 
\end{pmatrix}+A_{3/7}^{m-1}\begin{pmatrix}-{1\over 
7}\cr -{2\over 7}\cr -{2\over 7}\cr {5\over 7}\cr \end{pmatrix},$$
$$A_{3/7}^{m-1}\begin{pmatrix}0\cr 1\cr 0\cr 0\cr\end{pmatrix}=2^{m-1}\begin{pmatrix}{1\over 
7}\cr {2\over 7}\cr {2\over 7}\cr {2\over 7}\cr 
\end{pmatrix}+A_{3/7}^{m-1}\begin{pmatrix}-{1\over 
7}\cr {5\over 7}\cr -{2\over 7}\cr -{2\over 7}\cr\end{pmatrix}.$$
Now 
$$A_{3/7}\begin{pmatrix}-{1\over 
7}\cr -{2\over 7}\cr -{2\over 7}\cr {5\over 7}\cr \end{pmatrix}=\begin{pmatrix}{5\over 
7}\cr -{4\over 7}\cr -{4\over 7}\cr {3\over 7}\cr\end{pmatrix},\ \ 
A_{3/7}^{2}\begin{pmatrix}-{1\over 
7}\cr -{2\over 7}\cr -{2\over 7}\cr {5\over 7}\cr \end{pmatrix}=\begin{pmatrix}{3\over 
7}\cr {-8\over 7}\cr {6\over 7}\cr -{1\over 7}\cr \end{pmatrix},$$
$$A_{3/7}\begin{pmatrix}-{1\over 
7}\cr {5\over 7}\cr -{2\over 7}\cr -{2\over 7}\cr \end{pmatrix}=\begin{pmatrix}-{2\over 
7}\cr {3\over 7}\cr {3\over 7}\cr -{4\over 7}\cr \end{pmatrix},\ \  
A_{3/7}^{2}\begin{pmatrix}-{1\over 
7}\cr {5\over 7}\cr -{2\over 7}\cr -{2\over 7}\cr \end{pmatrix}=\begin{pmatrix}-{4\over 
7}\cr {6\over 7}\cr -{1\over 7}\cr -{1\over 7}\cr \end{pmatrix}.$$
So
$$a_{m,s(v_{1}),3/7,-1}={1\over 7}.2^{m}\begin{cases}-{1\over 
7}{\rm{\ if\ }}m=0{\rm{\ mod\ }}3,\\ +{5\over 
7}{\rm{\ if\ }}m=1{\rm{\ mod\ }}3,\\ +{3\over 
7}{\rm{\ if\ }}m=0{\rm{\ mod\ }}3,\end{cases}$$
$$a_{m,c_{1},3/7,-1}={1\over 7}.2^{m}\begin{cases}-{1\over 
7}{\rm{\ if\ }}m=0{\rm{\ mod\ }}3,\\ -{2\over 
7}{\rm{\ if\ }}m=1{\rm{\ mod\ }}3,\\ -{4\over 
7}{\rm{\ if\ }}m=2{\rm{\ mod\ }}3.\end{cases}$$

For future reference, we record
\begin{equation}\label{3.4.3}
    a_{m,1/7,s(v_{1})}+a_{m,6/7,s(v_{1})}+a_{m,3/7,s(v_{1}),-1}
    =2^{m}
    \begin{cases}-1{\rm{\ if\ }}m=0{\rm{\ mod\ }}3,\\
	+1{\rm{\ if\ }}m=1{\rm{\ mod\ }}3,\\
	+1{\rm{\ if\ }}m=2{\rm{\ mod\ }}3, \end{cases}\end{equation}
	
\begin{equation}\label{3.4.4}
    a_{m,1/7,c_{1}}+a_{m,6/7,c_{1}}+a_{m,3/7,c_{1},-1}
    =2^{m}
    \begin{cases}-1{\rm{\ if\ }}m=0{\rm{\ mod\ }}3,\\
	+0{\rm{\ if\ }}m=1{\rm{\ mod\ }}3,\\
	+0{\rm{\ if\ }}m=2{\rm{\ mod\ }}3. \end{cases}\end{equation}
	
The calculation of $a_{m,x,3/7,1}$ is a little less direct. First 
we note that, for $m\geq 2$,
$$b_{m,s(v_{1})}=\# (\{ z\in D(L_{3}L_{2}): 
s_{3/7}^{m}(z)=s_{3/7}(v_{1}),s_{3/7}^{m-1}(z)\neq v_{1}\} )$$
$$=
\# (\{ z\in X_{2}: 
s_{3/7}^{m-2}(z)=s_{3/7}(v_{1}),s_{3/7}^{m-3}(z)\neq v_{1}\} ),$$
where the last condition $s_{3/7}^{m-3}(z)\neq v_{1}$ is dropped if 
$m=2$. Then $b_{1,s(v_{1})}=0=b_{2,s(v_{1})}$ and, for $m\geq 
3$, $b_{m,s(v_{1})}$ is the third entry of 
$$A_{3/7}^{m-3}\begin{pmatrix}0\cr 0\cr 0\cr 1\cr\end{pmatrix} .$$
So if $m\geq 3$,
$$b_{m,s(v_{1})}={1\over 28}.2^{m}\begin{cases}-{2\over 
7}{\rm{\ 
if\ }}m=0{\rm{\ mod\ }}3,\\ -{4\over 7}{\rm{\ 
if\ }}m=1{\rm{\ mod\ }}3,\\ + {6\over 7}{\rm{\ 
if\ }}m=2{\rm{\ mod\ }}3.\\ \end{cases}$$
If we define $b_{m,c_{1}}$ similarly then $b_{1,c_{1}}=0$, 
$b_{2,c_{1}}=1$, and for $m\geq 3$ we have
$$b_{m,c_{1}}={1\over 28}.2^{m}+\begin{cases}-{2\over 7}{\rm{\ 
if\ }}m=0{\rm{\ mod\ }}3,\\ +{3\over 7}{\rm{\ 
if\ }}m=1{\rm{\ mod\ }}3,\\-{1\over 7}{\rm{\ 
if\ }}m=2{\rm{\ mod\ }}3.\end{cases}$$
Then we can obtain $a_{m,3/7,x,1}$ from
$$a_{m,x,3/7,1}=\sum _{m-2-2k\geq 0}b_{m-2k,x}.$$

So we obtain
$$a_{m,s(v_{1}),3/7,1}={1\over 21}.2^{m}\begin{cases}-{22\over 
21}{\rm{\ if\ }}m=0{\rm{\ mod\ }} 6,\\ -{2\over 21}{\rm{\ if\ 
}}m=1{\rm{\ mod\ }} 6,\\ -{4\over 21}{\rm{\ if\ 
}}m=2{\rm{\ mod\ }}6,\\ -{8\over 21}{\rm{\ if\ 
}}m=3{\rm{\ mod\ }}6,\\ -{16\over 21}{\rm{\ if\ 
}}m=4{\rm{\ mod\ }}6,\\ +{10\over 21}{\rm{\ if\ 
}}m=5{\rm{\ mod\ }}6, \end{cases}$$ 
$$a_{m,c_{1},3/7,1}={1\over 21}.2^{m}\begin{cases}+{22\over 
21}{\rm{\ if\ }}m=0{\rm{\ mod\ }}6,\\ -{2\over 21}{\rm{\ if\ 
}}m=1{\rm{\ mod\ }}6,\\ +{17\over 21}{\rm{\ if\ 
}}m=2{\rm{\ mod\ }}6,\\ -{8\over 21}{\rm{\ if\ 
}}m=3{\rm{\ mod\ }}6,\\ +{26\over 21}{\rm{\ if\ 
}}m=4{\rm{\ mod\ }}6,\\ -{11\over 21}{\rm{\ if\ 
}}m=5{\rm{\ mod\ }}6.\end{cases}$$ 

So we obtain, if $m\geq 2$, if
$$t(s(v_{1}))= a_{m,1/7,s(v_{1})}+a_{m,6/7,s(v_{1})} 
+a_{m,s(v_{1}),3/7,-1}
+a_{m,s(v_{1}),3/7,1},$$
then:
\begin{equation}\label{3.4.5}
  t(s(v_{1}))=\left( 1+\displaystyle{{1\over 21}}\right) 
.2^{m}\begin{cases}-{43\over 21}{\rm{\ if\ }}m=0{\rm{\ mod\ }}6,\\ 
 +{19\over21}{\rm{\ if\ }}m=1{\rm{\ mod\ }} 6,\\
+{17\over 21}{\rm{\ if\ }}m=2{\rm{\ mod\ }}6,\\
 -{29\over 21}{\rm{\ if\ }}m=3{\rm{\ mod\ }}6,\\
  +{5\over 21}{\rm{\ if\ }}m=4{\rm{\ mod\ }}6,\\
   +{31\over 21}{\rm{\ if\ }}m=5{\rm{\ mod\ }}6.\cr \end{cases}\end{equation}

The first $6$ numbers are
$$3,\ \ 5,\ \ 7,\ \ 17,\ \ 35\ \ 65,$$
compared with the values $3$, $5$, $8$, $18$, $36$, $69$
of $\# (P_{m}'(1))$

Similarly, if
$$t(c_{1})=a_{m,1/7,c_{1}}+a_{m,6/7,c_{1}}+a_{m,c_{1},3/7,-1}+
a_{m,c_{1},3/7,1},$$
then:
\begin{equation}\label{3.4.6}
    t(c_{1})
=\left( 1+{1\over 21}\right) .2^{m}\begin{cases}
-{1\over 
21}{\rm{\ if\ }}m=0{\rm{\ mod\ }}6,\\
 -{2\over 21}{\rm{\ if\ }}m=1{\rm{\ mod\ }}6,\\  
 +{17\over 21}{\rm{\ if\ }}m=2{\rm{\ mod\ }}6,\\
  -{29\over 21}{\rm{\ if\ }}m=3{\rm{\ mod\ }}6,\\
 +{26\over 21}{\rm{\ if\ }}m=4{\rm{\ mod\ }}6,\\
 -{11\over 21}{\rm{\ if\ }}m=5{\rm{\ mod\ }}6\end{cases}\end{equation} 

The first $6$ numbers are 
$$2,\ \ 5,\ \ 7,\ \ 18,\ \ 33,\ \ 67,$$
compared with the values $2$, $5$, $7$, $18$, $35$, $69$
of $\# (P_{m}'(0))$. The deficit of the captures from the total
in this case is slightly less than 
it appears to be, because in the cases of preperiods $5$, the two
captures $\sigma _{r}\circ s_{3/7}$ with $e^{2\pi ir}$ in the 
boundary 
of the gap coded by $L_{3}L_{2}R_{3}L_{3}L_{2}C$ are not Thurston 
equivalent, and are not Thurston equivalent to any other captures. 
Similar properties hold for the two preperiod $6$ captures 
corresponding to the gap coded by $L_{3}L_{2}R_{3}L_{3}^{2}L_{2}C$.

Overall, the deficit of $t(c_{1})$ and $t(s(v_{1}))$ compared to the 
values of $\# (P_{m}'(0))$ and $\# (P_{m}'(1))$ is equal to 
$${1\over 21}.2^{m}+O(1).$$

This calculation shows 
that, 
asymptotically, at most half of the points of $P_{3,m}(a_{1},+)$ are 
Thurston equivalent to captures, that is, at most ${1\over 
21}2^{m}(1+o(1))$ out of ${2\over 21}2^{m}(1+o(1))$. In fact, it is 
probably considerably less. I would guess that the number is 
$c2^{m}(1+o(1))$ for some $c>0$ which can be determined. It is 
interesting to note that the average image size of a map from a set 
of $n_{1}$ elements to a set of $n_{2}$ elements is 
$$n_{2}\left( 1-\left( 1-{1\over n_{2}}\right) \right) ^{n_{1}}.$$
I thank my colleague Jonathan Woolf for producing and deriving this 
formula. If 
$n_{2}=2n_{1}(1+o(1))$ then this gives 
$$n_{2}(1-e^{-1/2})(1+o(1)).$$

\section{Check on the numbers in Theorem \ref{2.10}}\label{3.5}

Theorem \ref{2.10} includes a precise formula for the number of type 
III 
hyperbolic components of preperiod $m$ of each of the two possible 
types, that is, with $\displaystyle{{2a\over a+1}}$ in the backward 
orbit of 
$\displaystyle{{a\over a+1}}$ or $a$. Let  $U^{p}$ be as in 
\ref{2.10}. Throughout this section, let $s=s_{3/7}$.
Similarly to the numbers $a_{m,x,q}$ as in  \ref{3.4}, we define
$$d_{m,x}=\# (\{ z\in D(L_{2}):s_{3/7}^{m}(z)=x,\ s_{3/7}^{m-1}(z)\neq 
s_{3/7}^{2}(x)\} ),$$
$$e_{m,x}=\# (\{ z\in D(BC):s_{3/7}^{m}(z)=x,\ s_{3/7}^{m-1}(z)\neq 
s_{3/7}^{2}(x)\} ),$$
for each of $x=s_{3/7}(v_{1})$ and $x=c_{1}$, and let
$$f_{m,p,x}=\# (\{ z\in U^{p}:s_{3/7}^{m}(z)=x,\ s_{3/7}^{m-1}(z)\neq 
s_{3/7}^{2}(x)\} )$$
and
$$c_{m,x}=\sum _{p\geq 0}f_{m,p,x}.$$
In the terminology of \ref{3.5}, Theorem \ref{2.10} implies that the 
number of type III hyperbolic 
components with ${2a\over a+1}$ in the backward orbit 
of 
${a\over a+1}$ or $a$ is
\begin{equation}\label{3.5.1}
    a_{m,x,1/7}+a_{m,x,6/7}+a_{m,x,3/7,-1}+c_{m,x}\end{equation}
where $x$ is, respectively, $s(v_{1})$ or $c_{1}$. The numbers 
of hyperbolic components are bounded by, respectively, $\# 
(P_{m}'(1))$ 
and $\# (P_{m}'(0))$. So a useful check on the validity of the 
theorem 
is to show that the numbers in (\ref{3.5.1}) are, for 
$x=s_{3/7}(v_{1})$ and $x=c_{1}$ respectively, $\# (P_{m}'(1))$ and 
$\# (P_{m}'(0))$.  The numbers
$a_{m,x,1/7}+a_{m,x,6/7}+a_{m,x,3/7,-1}$ are given in (\ref{3.4.3}) 
and (\ref{3.4.4}). So we only need to compute the numbers $c_{m,x}$, which are  computed in terms of the numbers $f_{m,p,x}$. Now we claim that, if $p\geq 0$,
$$f_{m,p+1,x}=f_{m-2,p,x}=f_{m-2p-2,0,x}.$$

There is a lot of cancellation in $U^{p}$ for $p\geq 1$:
\begin{itemize}
\item[.] Points in $S_{1,p,k+1}D(L_{2})$ cancel with points in $S_{1,p,k}D(v_{p,0})$ for $k\geq 0$. The length of the word ending in $L_{2}$ in both cases is $4p+4+k(2p+3)$;
\item[.] Points in $S_{2,p,k}D(v_{t,n})$ cancel with points in $S_{1,p,k}D(v_{t-1,n})$ for all $n\geq 0$ and $k\geq 0$ and $2\leq t\leq p$;
\item[.] Points in $S_{2,p,k+1}D(v_{1,0})$ cancel with points in $S_{1,p,k}D(v_{p,1})$. In both cases the lengths of words are $4p+8+k(2p+3)$;
\item[.]  Points in $S_{2,p,k+1}D(v_{1,n-1})$  cancel with points in $S_{1,p,k}D(v_{p,n})$ for $n\geq 2$ and $k\geq 0$. Both have length $4p+5+3n+k(2p+3)$;
\item[.] Points in $S_{2,p,k+1}D(u_{n})$ cancel with points in $S_{2,p,k}D((L_{2}R_{3})^{2}L_{3}^{2p-1}u_{n})$ for $k\geq 0$; 
\item[.] Points in $S_{1,p,k}D((L_{2}R_{3})^{2}L_{3}^{2p-1}u_{n})$ cancel with points in $S_{2,p,k+1}D(u_{n})$ for $k\geq 0$.
\end{itemize}

So we are left with points in $S_{2,p,0}D(v_{1,0})$ and points in $S_{2,p,0}D(u_{n})$ minus points in $S_{1,p,0}D(u_{n})$ for $n\geq 0$ and points in $S_{2,p,0}D(v_{1,k})$ for $k\geq 0$. Note that $S_{2,p,0}D(u_{k})$ and $S_{1,p,0}D(u_{k})$ map under $s^{2p+3k}$ and $s^{2p+2+3k}$ to $D(BC)$ and $S_{2,0}D(v_{1,k})$ maps under $s^{2p+3k}$ to $D(L_{3}L_{2})$, while $D(L_{3}^{3k-1}L_{2}BC)$ and $D(L_{3}^{3k+1}L_{2})$ map under $s^{3k}$ to $D(BC)$ and $D(L_{3}L_{2})$ respectively.

So we deduce that, for $p\geq 0$,
$$f_{m,p,x}=f_{m-2,p-1,x}=f_{m-2p,0,x}.$$
So
$$c_{m,x}=\sum _{2p\leq m}f_{m-2p,0,x}.$$
Note that
$$f_{m,0,x}=\sum _{k\geq 0,3k\leq m}(e_{m-3k,x}+d_{m-3k-1,x}-e_{m-3k-2,x}).$$
Here we define $e_{r,x}=0$ if $r\leq 0$.

Then in the 
terminology we have previously used, in particular in \ref{3.4}, 
$D(L_{2})$ maps homeomorphically under $s$ to $D(C)\cup D(R_{3})$, 
 and $D(BC)$ maps homeomorphically under $s^{2}$ to $D(R_{1})\cup 
D(R_{2})$, and $D(L_{3}L_{2}BC)=S_{1,0,0}D(u_{0})$ maps homeomorphically under $s^{4}$ to $D(R_{1})\cup 
D(R_{2})$. 
Let
$$\underline{v}=\begin{pmatrix}0\cr 0\cr 0\cr 1\cr\end{pmatrix}{\rm{\ or\ }}
\begin{pmatrix}0\cr 1\cr 0\cr 0\cr\end{pmatrix},$$
depending on whether $x=s(v_{1})$ or $c_{1}$.
Then similarly to the calculation of $b_{m,x}$ in \ref{3.4}, with 
matrix $A=A_{3/7}$ as in \ref{3.4},  
$d_{m-1,x}$ and $e_{m,x}$ are the third and fourth entries respectively  of 
$A^{m-3}(\underline{v})$for $m\geq 3$.
 Using the 
decomposition of $\underline{v}$ as a sum of an eigenvector of $A$ 
with 
eigenvalue $2$ and a vector in the sum of the other eigenspaces of 
$A$,  we obtain, for $m\geq 4$,
$$e_{m,s(v_{1})}+d_{m-1,s(v_{1})}-e_{m-2,s(v_{1})}=$$
$${4\over 7}.2^{m-3}-{2\over 7}.2^{m-5}+\begin{cases}-\frac{2}{7}+\frac{5}{7}-(\frac{3}{7}){\rm{\ (if\ }}m=0{\rm{\ mod\ }}3{\rm{)}}\\ -\frac{4}{7}+\frac{3}{7}-(-\frac{1}{7}){\rm{\ (if\ }}m=1{\rm{\ mod\ }}3{\rm{)}}\\ \frac{6}{7}-\frac{1}{7}-(\frac{5}{7}){\rm{\ (if\ }}m=2{\rm{\ mod\ }}3{\rm{)}}\end{cases} =2^{m-4}+0,$$
and
$$e_{m,c_{1}}+d_{m-1,c_{1}}-e_{m-2,c_{1}}=$$
$${4\over 7}.2^{m-3}-{2\over 7}.2^{m-5}+\begin{cases}-\frac{2}{7}-\frac{2}{7}-(-\frac{4}{7}){\rm{\ (if\ }}m=0{\rm{\ mod\ }}3{\rm{)}}\\ \frac{3}{7}-\frac{4}{7}-(-\frac{1}{7}){\rm{\ (if\ }}m=1{\rm{\ mod\ }}3{\rm{)}}\\ -\frac{1}{7}-\frac{1}{7}-(-\frac{2}{7}){\rm{\ (if\ }}m=2{\rm{\ mod\ }}3{\rm{)}}\end{cases} =2^{m-4}+0.$$
We can also check that
$$d_{1,s(v_{1})}=d_{2,s(v_{1})}=0,$$
$$d_{1,c_{1}}=0,\ d_{2,c_{1}}=1,$$
$$e_{r,x}=0{\rm{\ for\ }}r<3{\rm{\ and\ }}x=s(v_{1}){\rm{\ or\ }}c_{1},$$
$$e_{3,s(v_{1})}=1,\ \ e_{3,c_{1}}=0.$$
So
$$e_{3,s(v_{1})}+d_{2,s(v_{1})}-e_{1,s(v_{1})}=1,$$
$$e_{3,c_{1}}+d_{2,c_{1}}-e_{1,c_{1}}=0.$$

We can then check that
\begin{equation}\label{3.5.2}
    f_{m,0,s(v_{1})}={1\over 14}2^{m}\begin{cases}+{3\over 
7}{\rm{\ if\ }}m=0{\rm{\ mod\ }}3,\\ 
-{1\over 7}{\rm{\ if\ }}m=1{\rm{\ mod\ }}3,\\ 
-{2\over 7}{\rm{\ if\ }}m=2{\rm{\ mod\ }}3.
\end{cases}\end{equation}
So then we obtain
\begin{equation}\label{3.5.3}
    c_{m,s(v_{1})}={2\over 21}2^{m}\begin{cases}
    -{2\over 21}{\rm{\ if\ }}m=0{\rm{\ mod\ }}6,\\ 
    -{4\over 21}{\rm{\ if\ }}m=1{\rm{\ mod\ }}6,\\ 
    -{8\over 21}{\rm{\ if\ }}m=2{\rm{\ mod\ }}6,\\ 
    +{5\over 21}{\rm{\ if\ }}m=3{\rm{\ mod\ }}6,\\ 
    -{11\over 21}{\rm{\ if\ }}m=4{\rm{\ mod\ }}6,\\
    -{1\over 21}{\rm{\ if\ }}m=5{\rm{\ mod\ }}6. 
    \end{cases}\end{equation} 
From this we obtain that
$$a_{m,s(v_{1}),1/7}+a_{m,s(v_{1}),6/7}+a_{m,s(v_{1}),3/7,-1}
+c_{m,s(v_{1})}=\# (P_{m}'(1)),$$
as we wanted to check. 
The calculations for $f§_{m,0,c_{1}}$ and $c_{m,c_{1}}$ are similar, 
taking onto account the different values of $d_{m,c_{1}}$.
We have
\begin{equation}\label{3.5.4}
    f_{m,0,c_{1}}={1\over 14}2^{m}\begin{cases}-{4\over 
7}{\rm{\ if\ }}m=0{\rm{\ mod\ }}3,\\ 
-{1\over 7}{\rm{\ if\ }}m=1{\rm{\ mod\ }}3,\\ 
+{5\over 7}{\rm{\ if\ }}m=2{\rm{\ mod\ }}3.\\ 
\end{cases}\end{equation}
So
\begin{equation}\label{3.5.5}
    c_{m,c_{1}}={2\over 21}2^{m}\begin{cases}
    -{2\over 21}{\rm{\ if\ }}m=0{\rm{\ mod\ }}6,\\ 
    -{4\over 21}{\rm{\ if\ }}m=1{\rm{\ mod\ }}6,\\ 
    +{13\over 21}{\rm{\ if\ }}m=2{\rm{\ mod\ }}6,\\
    -{16\over 21}{\rm{\ if\ }}m=3{\rm{\ mod\ }}6,\\ 
    +{10\over 21}{\rm{\ if\ }}m=4{\rm{\ mod\ }}6,\\
    -{1\over 21}{\rm{\ if\ }}m=5{\rm{\ mod\ }}6.
    \end{cases}\end{equation} 
From this we obtain that
$$a_{m,c_{1},1/7}+a_{m,c_{1},6/7}+a_{m,c_{1},3/7,-1}
+c_{m,c_{1}}=\# (P_{m}'(0)),$$
 as we 
wanted to check. 

\chapter{The Resident's View}\label{4}
\section{}\label{4.1}

In this section we recall the theory of the Resident's View from 
\cite{R3}
as we need it. The basic idea is that we shall describe $V_{3}$ in 
terms of the dynamical planes of the three maps $h_{a_{0}}$, 
$h_{\overline{a_{0}}}$ and $h_{a_{1}}$ of $V_{3}$ which are M\"obius 
conjugate to polynomials, and which are Thurston equivalent to the 
maps $s_{1/7}$, $s_{6/7}$ and $s_{3/7}$. Recall from the diagram of 
\ref{2.2}  that $V_{3}(a_{0})$ and $V_{3}(\overline{a_{0}})$ are the 
regions in the upper and lower half-planes respectively bounded by 
the hyperbolic coponents $H_{\pm 1}$, that $V(a_{1},-)$ is bounded by 
$H_{a_{1}}$ and $H_{-1}$, and $V(a_{1},+)$ by $H_{a_{1}}$ and 
$H_{1}$.  We shall use the dynamical 
planes of 
$$s_{1/7},\ \ s_{6/7},\ \ s_{3/7},$$
to describe the regions 
$$V_{3}(a_{0}),\ \ V_{3}(\overline{a_{0}}),\ \ V_{3}(a_{1},\pm ).$$
This means that we identify type III hyperbolic components of 
preperiod $\leq m$ in $V_{3}$ with points in the set $Z_{m}(s)$,
where $s=s_{1/7}$ or $s_{6/7}$, or $s_{3/7}$. Since the Resident's View 
is actually an identification between universal covers, that is, by 
paths to points rather than simply the endpoints of the paths, there 
is no 
general reason why we should 
have a one-to-one correspondence between the set of  type III 
hyperbolic components of preperiod $\leq m$ in $V_{3}$ with points in 
the set and subsets of the sets $Z_{m}(s_{q})$, for $q={1\over 7}$ or 
${6\over 7}$ or ${3\over 7}$. But in fact we do find such a 
correspondence in $V_{3}(a_{0})$ and  $V_{3}(\overline{a_{0}})$ and  
$V_{3}(a_{1},-)$. The correspondence is  
many-to-one for $V_{3}(a_{1},+)$. 
The precise statement will have to 
wait until the main theorem in \ref{7.8}. 
All the theory described in this section works for
any $V_{n,m}$ or $B_{n,m}$, but we shall be applying it only in the 
case $n=3$ in the present paper. So, to simplify notation a bit, we 
restrict to the case $n=3$.

\section{The maps $\Phi _{1}$, $\Phi _{2}$, $\rho $}\label{4.2}

We fix $m$ and $g_{0}\in B_{3,m}$ for which $g_{0}(v_{2})=v_{2}$.
If $g_{0}\in V_{3,m}$ then $g_{0}=h_{a}$ for  $a=a_{0}$ or 
$\overline{a_{0}}$ or $a_{1}$. We shall always take $g_{0}$ to be one 
of 
these, or $g_{0}=s_{q}$ for $q={1\over 7}$ or  ${6\over 7}$ or ${3\over 
7}$, where these are as in \ref{2.3}. Although these maps are not in 
$V_{3,m}$,  the map $s_{1/7}$ is Thurston equivalent to $h_{a_{0}}$, and  $s_{6/7}$ 
to $h_{\overline{a_{0}}}$  and $s_{3/7}$  to $h_{a_{1}}$.
We write  $v_{1}$ and $v_{2}$ 
for the critical 
values of $g_{0}$, with $v_{1}$ of period $3$. Because Thurston 
equivalence 
classes of these maps are simply connected (using our conventions, as 
noted in \ref{2.4}) there is a unique path,  up to homotopy in 
$B_{3,m}$, 
joining $s_{1/7}$ and $h_{a_{0}}$ in their common Thurston 
equivalence class, and similarly for the other pairs. Inclusion of 
$V_{3,m}$ in $B_{3,m}$ therefore gives rise to a natural homomorphism 
from $\pi _{1}(V_{3,m},h_{a})$ to $\pi _{1}(B_{3,m},s_{q})$, for each 
of $(a,q)=(a_{0},{1\over 7})$ or $(\overline{a_{0}},{6\over 7})$ or  
$(a_{1},{3\over 7})$. One of the main results of \cite{R3} is that 
this homomorphism is injective. By abuse of notation we therefore 
write $\pi _{1}(V_{3,m},s_{q})$ for the resulting subgroup of $\pi 
_{1}(B_{3,m},s_{q})$. For any pair of spaces $(X,A)$ with $A\subset 
X$, and $x\in X$ we write $\pi _{1}(X\setminus A,A,x)$ for the set of homotopy 
classes of paths from $x$ to $A$, intersecting $A$ only in the second endpoint, using homotopies preserving $A$ and 
$x$. Then, similarly to the above, we write $\pi 
_{1}(V_{3,m},P_{3,m},s_{q})$ for the subset of $\pi 
_{1}(B_{3,m},P_{3,m},s_{q})$ of homotopy classes which can be 
represented as a path in the Thurston equivalence class of 
$s_{q}$ and $h_{a_{0}}$ from $s_{q}$ to $h_{a_{0}}$, followed by an 
element of $\pi _{1}(V_{3,m},P_{3,m},h_{a_{0}})$.

So now we let $g_{0}$ be any of these six maps $h_{a}$ or $s_{q}$. We 
also write $Z_{m}=Z_{m}(g_{0})$ and $Y_{m}=Z_{m}\cup \{ v_{2}\} $. The 
universal 
cover of 
$\overline{\mathbb C}\setminus Z_{m}$ is conformally the unit disc 
$D$, 
for  all $m\geq 0$. The {\em{Resident's View}} 
of 
\cite{R3} identifies the universal cover of $V_{3,m}$ with a subset 
of the universal cover of $\overline{\mathbb C}\setminus Z_{m}$ which 
is the disc.
This is done using set-theoretic injections
$$\rho =\rho (.,g_{0}):\pi _{1}(B_{3,m},g_{0})\to \pi 
_{1}(\overline{\mathbb 
C}\setminus 
Z_{m},v_{2}),$$
$$\rho_{2}=\rho _{2}(.,g_{0}):\pi _{1}(B_{3,m},N_{3,m},g_{0})\to \pi 
_{1}(\overline{\mathbb C}\setminus 
Z_{m},Z_{m},v_{2})$$
which are defined in 1.12 of \cite{R3}.
These combine with the injective-on-$\pi _{1}$ inclusions of 
$V_{3,m}$ 
in $B_{3,m}$ and of $(V_{3,m},P_{3,m})$ in $(B_{3,m},N_{3,m})$ (of 
which the proof takes up perhaps half of \cite{R3}) to give
set-theoretic injections
$$\rho:\pi _{1}(V_{3,m},g_{0})\to \pi _{1}(\overline{\mathbb 
C}\setminus 
Z_{m},v_{2}),$$
$$\rho_{2}:\pi _{1}(V_{3,m},P_{3,m},g_{0})\to \pi 
_{1}(\overline{\mathbb C}\setminus 
Z_{m},Z_{m},v_{2}).$$ 
 From now on in this work, we shall simply write $\rho $ for the maps 
which were called $\rho $ and $\rho _{2}$ in \cite{R3}. The 
definitions are very similar. Only the domains are different.  We now 
recall 
the definitions in the present context, considering just the 
definitions on 
$\pi _{1}(V_{3,m},g_{0})$ and $\pi _{1}(V_{3,m},P_{3,m},g_{0})$. Let 
$t\mapsto g_{t}:[0,1]\to V_{3,m}$ be either a closed path starting 
and ending at $g_{0}$ (for the 
definition of $\rho $) or a path from $g_{0}$ to a small 
neighbourhood of a point 
of $P_{3,m}\setminus P_{3,0}$ (for the definition of $\rho _{2}$). 
Then the path $g_{t}$ defines a path 
$\varphi 
_{t}$ of homeomorphisms of $\overline{\mathbb C}$ where $\varphi 
_{0}$ 
is the identity and $\varphi _{t}$ maps $Z_{m}=Z_{m}(g_{0})$ to 
$Z_{m}(g_{t})$ and $v_{2}(g_{0})$ to $v_{2}(g_{t})$. There is also a 
lifted 
path of homeomorphisms $\psi _{t}$ defined by $g_{t}\circ \psi 
_{t}=\varphi _{t}\circ g_{0}$, such that $\psi _0$ is the identity with $\psi _t(c_j(g_0))=c_j(g_t)$ for $j=1$, $2$. Standard covering space theory gives the lift of $\varphi _t$ to a homeomorphism $\psi _t$of $\overline {\mathbb C}\setminus \{ c_1(g_0),c_2(g_0)\} $ to $\overline{\mathbb C}\setminus \{ c_1(g_t),c_2(g_t)\} $, but this extends to a homeomorphism mapping $c_j(g_0)$ to $c_j(g_t)$ for $j=1$, $2$. By continuity in $t$, since $\psi _t(Z_m)\subset g_t^{-1}(Z_m(g_t))$, it follows that $\psi _t=\varphi _t$ on $Z_m$, and $\varphi _{t}$ and 
$\psi _{t}$ are isotopic via an isotopy which is constant on 
$Z_{m}$. It follows that there is a path $\alpha _{t}$ starting 
from $v_{2}(g_{0})$ in $\overline{\mathbb C}\setminus 
Z_{m}$ such that $\varphi _{t}$ and $\psi _{t}\circ \sigma _{\alpha 
_{t}}$ are isotopic via an isotopy constant on $Z_{m}\cup 
\{ v_{2}(g_{0})\} $. Here, for any path $\gamma $, the homeomorphism $\sigma _{\gamma }$ 
is as in \ref{2.5}.
In the case when $t\mapsto g_{t}$ is a closed path, $\alpha =\alpha 
_{1}$ is also a closed path, whose homotopy class in $\pi 
_{1}(\overline{\mathbb C}\setminus 
Z_{m},v_{2})$ depends only on the homotopy class of $t\mapsto 
g_{t}$ 
in $\pi _{1}(V_{3,m},g_{0})$. So in this case we take 
$$\rho ([t\mapsto g_{t}])=[\alpha ].$$
We also define
$$\Phi _{1}([t\mapsto g_{t}])=[\varphi _{1}]$$
and 
$$\Phi _{2}([t\mapsto g_{t}])=[\psi _{1}].$$
Then $\Phi_{1}$ and $\Phi _{2}$ are anti-homomorphisms into the 
modular groups ${\rm{MG}}(\overline{\mathbb C},Y_{m} )$, 
${\rm{MG}}(\overline{\mathbb C},Y_{m+1})$ respectively. It is shown 
in 1.11 of \cite{R3} that they 
are injective on $\pi _{1}(B_{3,m},g_{0})$. It is pointed out in 1.13 
of 
\cite{R3} that
\begin{equation}\label{4.2.1}\sigma _{\alpha }\circ g_{0}
    \simeq _{\psi _{1}}g_{0},
    \end{equation}
for $\alpha =\rho ([t\mapsto g_{t}])$ and 
$\psi _{1}=\Phi _{2}([t\mapsto g_{t}])$ as above, where, as earlier, 
$\simeq $ 
denotes Thurston equivalence of critically finite branched coverings 
and $\simeq _{\psi _{1}}$ is as in \ref{2.4}. In fact (\ref{4.2.1}) 
is 
the defining equation of 
$$G_{2}=\Phi _{2}(\pi _{1}(B_{3,m},g_{0})),$$
or at least one form of the defining equation. The group 
$$G_{1}=\Phi _{1}(\pi _{1}(B_{3,m},g_{0}))$$
is equivalently characterised by 
\begin{equation}\label{4.2.2}G_{1}=\left \lbrace 
\begin{array}{l}[\varphi ]
    \in {\rm{MG}}(\overline{\mathbb 
C},Y_{m} ):[\varphi ]=[\sigma _{\alpha }\circ  \psi ] ,\cr 
{\rm{\ for\ at\ least\ one\ }}\psi 
{\rm{\ with\ }}\varphi \circ 
g_{0}=g_{0}\circ \psi ,\cr  \alpha \in \pi _{1}(\overline{\mathbb 
C}\setminus Z_{m},v_{2})\cr \end{array}\right \rbrace .\end{equation}
Here, part of the definition given by (\ref{4.2.2}) is that at least 
one  lift 
$\psi $ of $\varphi $ by the equation $g_{0}\circ \psi =\varphi \circ 
g_{0}$ satisfies $\psi =\varphi $ on $Z_{m}$. For further details, 
see 
1.11 to 1.13 of \cite{R3}.

We now consider the definition on $\pi _{1}(V_{3,m},P_{3,m},g_{0})$. 
We shall again want to obtain (\ref{4.2.1}), with $\alpha =\rho 
({g_{t}})$. Recall that the Polynomial-and-Path Theorem of \cite{R1} 
says that every critically finite type III quadratic rational map is 
Thurston-equivalent to a map of the form $\sigma _{\alpha }\circ 
g_{0}$ for some $g_{0}$ which is Thurston equivalent to a polynomial, 
and for some path $\alpha $ from $v_{2}(g_{0})$ to $Z_{m}(g_{0})$. 
First, we consider a path $g_{t}$ with $g_{0}$ as before, and with 
$g_{1}$  close to an 
element of $P_{3,m}\setminus P_{3,0}$, rather than equal to such an 
element. Then we can take $\varphi 
_{t}^{-1}$ bounded except near $v_{2}(g_{t})$, 
for $t$ near $1$. 
This means that $\psi _{t}^{-1}$ is bounded near 
$c_{2}(g_{t})$, and the second endpoint of $\alpha 
_{t}$ is near the corresponding point of 
$Z_{m}$ for $t$ near $1$. We then extend $\alpha _{1}$ by a small 
arc to a path $\alpha $ from $v_{2}(g_{0})$ to $Z_{m}$. Once again, the 
class of $\alpha $ in $\pi _{1}(\overline{\mathbb C}\setminus 
Z_{m},Z_{m},v_{2})$ depends only on the isotopy class of the path 
in $\pi _{1}(V_{3,m},P_{3,m},g_{0})$ which is a small extension of 
$t\mapsto g_{t}$. So this defines the map $\rho $ on $\pi 
_{1}(V_{3,m},P_{3,m},g_{0})$. We then have (\ref{4.2.1}) as before.

If $(a,q)=(a_{0},{1\over 7})$ or $(\overline{a_{0}},{6\over 7})$ or 
$(a_{1},{3\over 7})$, the images 
under $\rho ( .,h_{a})$ and $\rho (., s_{q})$ in homotopy 
classes in the dynamical planes of $h_{a}$ or $s_{q}$ correspond under 
any 
suitably close homeomorphism approximation to the semiconjugacy 
between $s_{q}$ and $h_{a}$, since the quotient map of $s_{q}$ on 
the quotient of the invariant lamination $L_{q}$ is conjugate to 
$h_{a}$ on 
the  Julia set $J(h_{a})$. If we define $\beta =\rho (t\mapsto 
g_{t},s_{q})$, then
$$\sigma _{\beta }\circ s\simeq g_{1}.$$
Moreover, we can recover $g_{1}$ up to conjugacy from $\beta $ and 
$s_{q}$. This is described in the Lamination Map Conjugacy and 
Lamination Map Equivalence Theorems of \cite{R1}, and works as 
follows. Let $L_{\beta }$ be the closure of the  set of geodesics 
homotopic to the paths of 
$$\cup _{n=0}^{\infty }(\sigma _{\beta }\circ s_{q})^{-n}(L_{q}).$$
Then an easy homotopy gives a Thurston equivalence between $\sigma 
_{\beta }\circ s_{q}$ and a map (unfortunately, in view of the use of 
$\rho $ in \cite{R3} and in the current paper, called $\rho 
_{L_{\beta }}$ in \cite{R1}) which preserves $L_{\beta }$, and is 
semiconjugate to $g_{1}$, under a semiconjugacy which collapses each 
leaf of $L_{\beta }$ to a point, but little more: the preimages of 
points are equivalence classes of the smallest closed equivalence 
relation generated by: each leaf of $L_{\beta }$ is in a single 
equivalence class. 

The point of the current work, and the way to prove Theorem 
\ref{2.10}, is to choose a fundamental domain for 
$V_{3,m}$. This means choosing an open topological disc $U$ in 
$V_{3,m}$, containing $g_{0}\in V_{3,m}$, whose complement is a union 
of edges joining up the points of $P_{3,m}$ to form a tree. Up to 
homotopy preserving $g_{0}$ and $P_{3,m}$, there is then a one-to-one 
correspondence between points in $V_{3,m}$, and the set of paths in 
$U$  from $g_{0}$ to $P_{3,m}$. Using $\rho $, this will then give a 
one-to-one correspondence between $P_{3,m}$ and a certain set of paths 
in $\pi _{1}(\overline{\mathbb C}\setminus Z_{m},Z_{m},v_{2})$. Using 
(\ref{4.2.1}), this will give a single chosen representation of the 
centre of each type III hyperbolic component of preperiod $\leq m$ in 
the form $\sigma _{\beta }\circ g_{0}$. 

\section{The Resident's View of Rational Maps Space}\label{4.3}

Now $\pi _{1}(V_{3,m},g_{0})$ and $\pi _{1}(\overline{\mathbb 
C}\setminus 
Z_{m},v_{2})$ are naturally embedded in the universal covers of 
$V_{3,m}$ and $\overline{\mathbb C}\setminus 
Z_{m}$ respectively, both of which identify conformally with the 
unit disc $D$. To make the natural embeddings, of course we have to 
fix preimages $\tilde{g_{0}}$ and $\tilde{v_{2}}$ of $g_{0}$ and $v_{2}$ 
in the universal covers. In the same way,
$\pi _{1}(V_{3,m},P_{3,m},g_{0})$ and $\pi 
_{1}(\overline{\mathbb C}\setminus 
Z_{m},Z_{m},v_{2})$ identify with subsets of $\partial D$. So we can 
regard $\rho $ as a map from a subset of $\partial D$ to a subset 
of $\partial D$. The 
content of the 
{\em{Resident's View of Rational Maps Space}} of \cite{R3} was 
essentially that $\rho $ extends monotonically to $\partial D$, 
with just countably many discontinuities, with continuous inverse on 
the set
$$\partial _{1}=\overline{\rho (\pi 
_{1}(V_{3,m},P_{3,m},g_{0}))}.$$

Since $\pi _{1}(V_{3,m},g_{0})$ acts naturally on the left of $\pi 
_{1}(V_{3,m},P_{3,m},g_{0})$, this action can be transferred using 
$\rho $ to one on  $\pi _{1}(\overline{\mathbb C}\setminus 
Z_{m},Z_{m},v_{2})$. It can be shown quite easily (1.13 of 
\cite{R3}) that the action is by homeomorphisms of $\partial D$ 
which, of course, restrict to homeomorphisms on $\partial 
_{1}$. The 
homeomorphisms are lifts of elements of ${\rm{MG}}(\overline{\mathbb 
C},Z_{m}\cup \{ v_{2}\} )$. As a 
consequence of the Resident's View of Rational Maps Space, the action 
can be extended (non-canonically) to an action on the 
Poincar\'e-metric-convex-hull
$$D'={\rm{Conv}}(\partial _{1})\subset D$$ 
of
$\partial _{1}$, which is 
conjugate to the action of $\pi _{1}(V_{3,m},g_{0})$ on the universal 
cover of $V_{3,m}$. Thus, $D$ 
is the universal cover of $\overline{\mathbb C}\setminus Z_{m}$, and 
the universal cover of $V_{3,m}$ has been identified with a 
Poincar\'e-metric-convex subset $D'=D'(g_{0})$ of $D$, with the $\pi 
_{1}(V_{n,m},g_{0})$ action transferring to a natural action 
on the boundary $\partial D$
of the universal cover of $\overline{\mathbb C}\setminus Z_{m}$, 
and preserving the subset 
$\partial D'\cap \partial D=\partial _{1}$. 

The action of $\pi _{1}(V_{3,m},g_{0})$ on 
$\pi _{1}(\overline{\mathbb C}\setminus 
Z_{m},Z_{m},v_{2})\subset \partial D
$ has an interpretation in terms of Thurston 
equivalence. If $\omega _{2}\in \pi _{1}(
V_{3,m},P_{3,m},g_{0})$ and $\omega _{1}\in \pi _{1}(V_{3,m},g_{0})$, 
and 
\begin{equation}\label{4.3.1}
    \beta _{1}=\rho (\omega _{1}),\ \beta _{2}=\rho(\omega 
_{2}),\ \beta _{3}=\rho (\omega _{1}*\omega _{2}),\end{equation}
\begin{equation}\label{4.3.2} [\psi _{1}]=[\Phi _{2}(\beta 
_{1})],\end{equation}
then (see also 1.13 of \cite{R3})
\begin{equation}\label{4.3.3}\beta _{3}=
    \beta _{1}*\psi _{1}^{-1}(\beta _{2}),\end{equation}
and hence, using (\ref{4.2.1}),
\begin{equation}\label{4.3.4}\sigma _{\beta _{3}}\circ g_{0}
    \simeq _{\psi _{1}}\sigma _{\beta  
_{2}}\circ g_{0}.\end{equation}
For later use, it is worth pointing out that (\ref{4.3.3}) and 
(\ref{4.3.4}) also hold 
if $\omega _{2}\in \pi _{1}(V_{3,m},g_{0})$, and we modify the 
definitions of $\beta _{2}$ and $\beta _{3}$ and add to them:
$$ [\psi _{2}]=[\Phi _{2}(\beta _{2})].$$
In this case we also have:
\begin{equation}\label{4.3.5} \sigma _{\beta _{3}}\circ 
g_{0}\simeq_{\psi _{2}\circ \psi _{1}}g_{0},\ \ g_{0}\simeq_{\psi 
_{1}^{-1}\circ \psi _{2}^{-1}}\sigma _{\beta _{3}}\circ 
g_{0}.\end{equation}

 The action of $\pi 
_{1}(V_{3,m},g_{0})$ on $D$, regarded as the universal covering space 
of $\overline{\mathbb C}\setminus Z_{m}$, extends to an action of 
$\pi _{1}(B_{3,m},g_{0})$ on $D$, which can equally be regarded as an 
(anti)action 
of $\Phi _{1}(\pi _{1}(B_{3,m},g_{0}))$ or $\Phi _{2}(\pi 
_{1}(B_{3,m},g_{0}))$. 
The groups   $G_{j}=\Phi _{j}(\pi _{1}(B_{3,m},g_{0}))$ are 
identified in 
\ref{4.2} using (\ref{4.2.1}) and (\ref{4.2.2}). There is 
also a natural action of  
  $\Phi _{i}(\pi _{1}(B_{3,m},g_{0}))$ on $\partial D$, which extends 
the 
action of 
 $\Phi _{i}(\pi _{1}(V_{3,m},g_{0}))$, $i=1$, $2$, and is given by 
the same 
 formulae as in (\ref{4.3.1}) to (\ref{4.3.4}).
 
 \section{Identifying $D'$}\label{4.4}

The best way to identify  $D'$ of \ref{4.3}, and  
$\Phi _{1}(\pi _{1}(V_{3,m},g_{0}))$, is by identifying $\partial 
D'\cap 
D$. The boundary is always a union 
of (Poincar\'e) geodesics  in the {\em{Levy convex 
hulls}}  $C(g_{0},\Gamma )$ (3.13 of \cite{R3}) of pairs 
$(g_{0},\Gamma )$ 
with $g_{0}\in B_{3,m}$, and for $\Gamma $ a set of simple loops in 
$\overline{\mathbb C}\setminus Z_{m}$, where $(g_{0},\Gamma )$ 
satisfies 
the {\em{Invariance and Levy Conditions}} (2.2 of \cite{R3}), and is 
{\em{minimal isometric satisfying the Edge Condition}} (Resident's 
View of Rational Maps Space in 5.10 of \cite{R3}). 
These conditions mean the following.
\begin{itemize} 
\item $\Gamma $ is a set of simple, non-trivial, non-peripheral, isotopically disjoint and isotopically distinct closed loops in $\overline{\mathbb C}\setminus Y_m$, such that $\Gamma \subset g_0^{-1}(\Gamma )$ up to isotopy preserving $Z_m$, and every loop in $g_0^{-1}(\Gamma )\setminus \Gamma $ is either trivial or peripheral in $\overline{\mathbb C}\setminus Z_m$. This is the {\em{Invariance Condition}}.
\item Every loop $\gamma \in \Gamma $ is in the backward orbit of a {\em{Levy cycle}} $\gamma _1,\cdots \gamma _r$ for $g_0$ and $\Gamma $, that is, $\gamma =\gamma _t$ for some $t\geq 1$, and $\gamma _j\in \Gamma $ for $1\leq j\leq t$ with $\gamma _{j+1}\subset g_0^{-1}(\gamma _j)$ up to $Z_m$-preserving isotopy  for $1\leq j\leq t-1$ and also $\gamma _1\subset g_0^{-1}(\gamma _r)$ up to $Z_m$-preserving isotopy, and if $\gamma _{j+1}'$ is the component of $g_0^{-1}(\gamma _j)$ which isotopic to $\gamma _{j+1}$, and $\gamma _1'\subset g_0^{-1}(\gamma _r)$ is similarly defined, then $g_0\vert \gamma _k'$ is a homeomorphism for $1\leq j\leq k$. This is the {\em{Levy Condition}}.
\item 

  Let $E_{2}$ denote the component 
of 
$\overline{\mathbb C}\setminus (\cup \Gamma )$ containing $v_{2}$. Then the {\em{fixed set}} $P=P(g_0,\Gamma )$ of $(g_0,\Gamma )$ (see 2.8 of \cite{R3}) is either a boundary component of $E_2$ --- and hence in $\Gamma $ --- or is a component of $\overline{\mathbb C}\setminus (\cup \Gamma )$ which is adjacent to $E_2$. This is the {\em{minimal condition}}. The fixed set $P$ satisfies a number of conditions:
\begin{itemize}
\item $P$ is not a disc, 
\item there is a component $P_1$ of $g_0^{-1}(P)$ such that $P_1=P$ up to $Z_m$-preserving isotopy, and $g_0\vert P_1$ is an homeomorphism, reversing orientation if $P$ is a loop and cyclically permuting boundary components, up to $Z_m$-preserving isotopy, if $P$ is surface (with boundary). \item The {\em{isometric condition}} means that $g\vert P_1$ is an isometry, up to isotopy.\end{itemize}\end{itemize}

In the case of $B_{3,m}$, and any pair $(g_0,\Gamma )$ satisfying the above conditions, the fixed set $P$ is either a loop or a pair of pants. So the isometric condition is automatic.  Now we define the set $C(g_0,\Gamma )$.  We use the notation of 
3.13 of \cite{R3}.
 Let $\Delta _{2}$ be the closed 
disc 
containing $E_{2}$, with interior disjoint from $P$. Then 
$\tilde{\Delta} _{2}$ is the 
preimage of $\Delta _{2}$ in $D$ with $\tilde{v_{2}}\in 
\tilde{E_{2}}$, where $\tilde{v_{2}}$ is a chosen preimage of 
$v_{2}$ in $D$, and 
$C(g_{0},\Gamma )$ is the union of geodesics with the same endpoints 
as 
the components of $\partial \tilde{\Delta} _{2}$. Thus, 
$C(g_{0},\Gamma )$ is a straightening of components of preimages in 
$D$ of $P$ or $\partial P$, depending on whether $P$ is a single loop 
or a pair of pants which are adjacent to, or 
in, preimages of $P$.  There are 
only two 
orbits of sets $C(g_{0},\Gamma )$ for $B_{3,m}$ under the action of $\pi _1(V_{3,m},g_0)$ on $D'$, because there are 
only 
two equivalence classes of loop sets $(g_{0},\Gamma )$ satisfying the 
Invariance, Levy and Edge conditions, and being minimal and 
isometric.

 So computing $\partial D'$ means computing sets 
$C(g_{0},\Gamma )$, which means computing the fixed  sets $P=P(g_0,\Gamma )$. Now this is 
easier if we can ensure that $\partial \tilde{\Delta} _{2}$ is not 
too far 
from geodesic. This means ensuring that  $P$  or $\partial P$ is not 
too far from geodesic in $\overline{\mathbb C}\setminus Z_{m}$. We 
specify to the case $g_{0}=s_{q}=s$ for some $q$, in our case with 
$q={1\over 7}$, ${6\over 7}$ or ${3\over 7}$. A simple closed loop in 
$\overline{\mathbb C}\setminus Z_{m}$ which has only essential 
intersections with $S^{1}$ can be isotoped, via an isotopy preserving 
$S^{1}$, to a geodesic. Our 
set $P$ might have nonessential intersections because it is defined 
up 
to isotopy preserving $Y_{m}=Z_{m}\cup \{ v_{2}\} $. But then we can 
choose $P$ up to isotopy preserving $Y_{m}$, 
and 
a closed loop $\alpha \in \pi _{1}(\overline{\mathbb C}\setminus 
Z_{m},v_{2})$ so 
that  $Q=\sigma _{\alpha  }(P)$ has only essential intersections with 
$S^{1}$. Our set $C(s,\Gamma )$ is then a geodesic with the same 
endpoints as a suitable lift of $Q$ or $\partial Q$, depending on 
whether $Q$ is a single loop or a pair of pants. We also have 
$Q_{1}=Q$ up to isotopy preserving $Z_{m}$, for a suitable 
component $Q_{1}$ of $s^{-1}(Q)$. It is a fact that 
$\alpha $ is bound to have at least one essential intersection with 
$\partial Q$

Now this gives a criterion for checking whether  $\beta \in 
\pi _{1}(\overline{\mathbb C}\setminus Z_{m},Z_{m},v_{2})$ has lift 
$\tilde{\beta }$ in $D'$. A necessary and sufficient condition for 
$\tilde{\beta }$ not to be in $D'$ is the following 
{\em{Inadmissibility criterion}}.
\begin{itemize}
\item There is a pair $(Q,\alpha )$ as 
above and such that the following holds. Assume without loss of 
generality that $\alpha $ and $\beta $ have only essential 
intersections with $\partial Q$. Let $\alpha _{1}$ be the portion of 
$\alpha $ from the start point at $v_{2}$ up until the last 
intersection with $\partial Q$. Then there is an initial portion 
$\beta _{1}$ of $\beta $ such that $\alpha _{1}$ and $\beta _{1}$ are 
homotopic via a homotopy fixing the start points at $v_{2}$ and 
keeping the second endpoint in $\partial Q$.
\end{itemize}
 This  criterion is essentially
considered in \cite{R2}, and we can use the Nonrational Lamination 
Map Theorem there to check it, using the invariant lamination 
$L_{\beta }$ defined in \ref{4.2}. In our case, this
criterion involves analysing the leaves of $L_{\beta }$ which have 
period two or three, and determining whether the closure of their 
union can be homotopic to a set such as $Q$. There are very few such 
periodic leaves, and therefore not much to check. This inadmissibility
criterion is a generalisation of the criterion for a mating to be 
not Thurston equivalent to a rational map \cite{TL}, and the results 
about Thurston equivalence of capture to rational maps, quoted in 
\ref{2.6}, can be easily derived from it.  

\section{The arc set $\Omega _{0}$ and the image under $\rho
$}\label{4.5}

We consider the picture of $V_{3}$ with an arc set $\Omega _{0}$:
\begin{center}
\includegraphics[width=8cm]{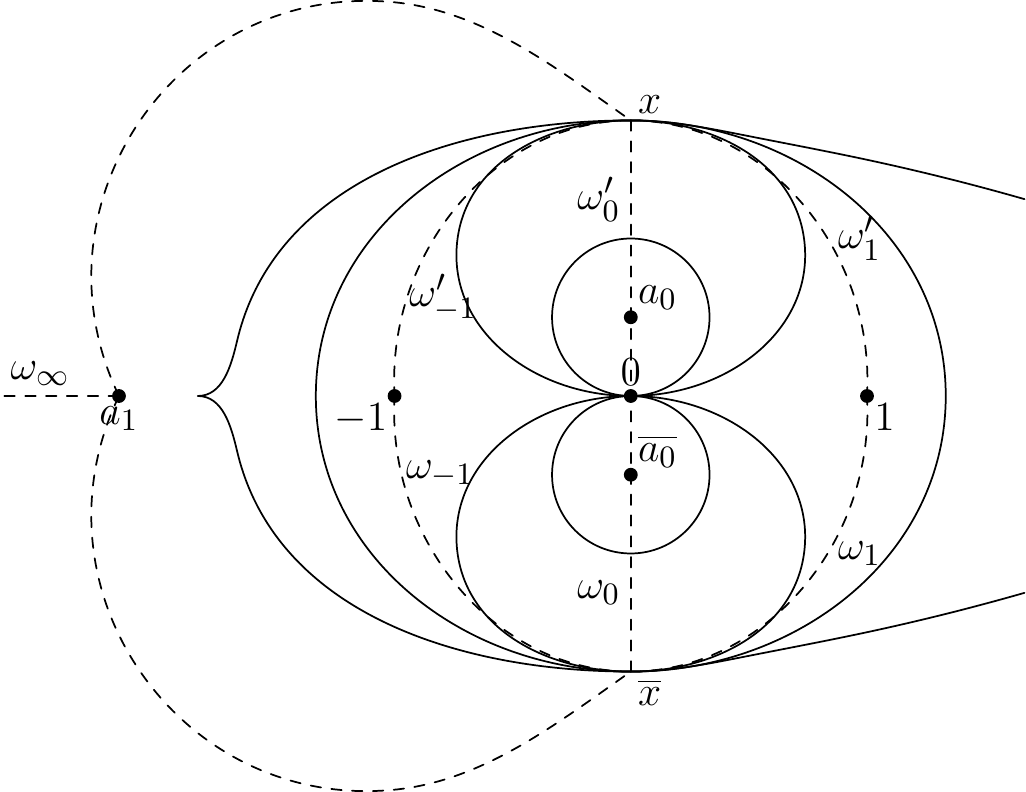}

$V_3$ and $\Omega $
\end{center}

The paths $\omega _{\infty }$, $\omega _{0}$, $\omega _{1}$ and $\omega _{-1}$ --- which are, in fact, arcs --- all start 
from 
$a_{1}$ and end at the points $\infty $ (the straight line path  $\omega _{\infty }$
heading 
to the left), $\pm 1$ and  $0$, and are in the closed lower half-plane.  
The paths $\omega _{0}$, $\omega _{1}$ and $\omega _{-1}$ all run close to the same path from $a_{1}$ 
to 
$x$, 
which is a path in the hyperbolic component of $h_{a_{1}}$. The paths 
from  near $\overline{x}$ to $\pm 1$ are then within the hyperbolic 
components 
of 
$h_{\pm 
1}$ respectively. The path from $\overline{x}$ to $0$ passes through 
$\overline{a_{0}}$ and 
is always in a copy of the Mandelbrot set, corresponding to matings  (\ref{2.6})
of 
quadratic polynomials  in the Mandelbrot set, minus a limb, with 
$h_{\overline{a_{0}}}$.  We also define paths $\omega _{0}'$, $\omega _{1}'$ and $\omega _{-1}'$, which are
 simply the complex 
conjugates of the paths $\omega _{0}$, 
$\omega _{1}$,  and $\omega _{-1}$. We write 
$\omega _{i}'$ for the conjugate of $\omega _{i}$, which is 
an arc of $\Omega _{0}'$, since, if $\beta $ is a path 
$\overline{\beta }$ is our notation for the reverse of $\beta $, and 
we shall be using this notation. Note that $\omega _{\infty 
}'=\omega _{\infty }$.

 Now we describe $\rho
(\omega ,s_{3/7})$ for each of $\omega =\omega _{\pm 1}'$, $\omega
_{\pm 1}$, $\omega _{\infty }$. The path $\rho (\omega _{\infty })$ 
is a bit different from the others because its second endpoint is not 
at a point of $Z_{0}(s_{3/7})$, but on a closed loop which is itself 
only defined up to isotopy relative to $Z_{m}(s_{3/7})$. This closed 
loop cuts the unit disc at precisely two points, $e^{\pm 2\pi 
i(3/7)}$. Recall that the path $\omega _{\infty }$ was along the 
negative real axis. It is not possible to represent $\rho (\omega _{\infty })$ by a path along the real axis without 
hitting $Z_{m}(s_{3/7})$ (even if $m=0$), so a choice of 
representative has been made, which meets the closed loop $\gamma $ 
at $e^{-2\pi i(3/7)}=e^{2\pi i(4/7)}$. Other possible representatives 
of $\rho (\omega _{\infty })$ are obtained from the chosen representative by 
homotopy keeping $Z_{m}(s_{3/7})$ fixed, and keeping the second 
endpoint on $\gamma $. So moving along $\gamma $ outside the unit 
disc, another possible representative of $\rho (\omega _{\infty })$ is the 
complex conjugate of the path drawn. 

  The picture is as follows. 
The unit circle and some lamination leaves of $L_{3/7}$ are indicated
by dashed lines
\begin{center}
\includegraphics[width=6cm]{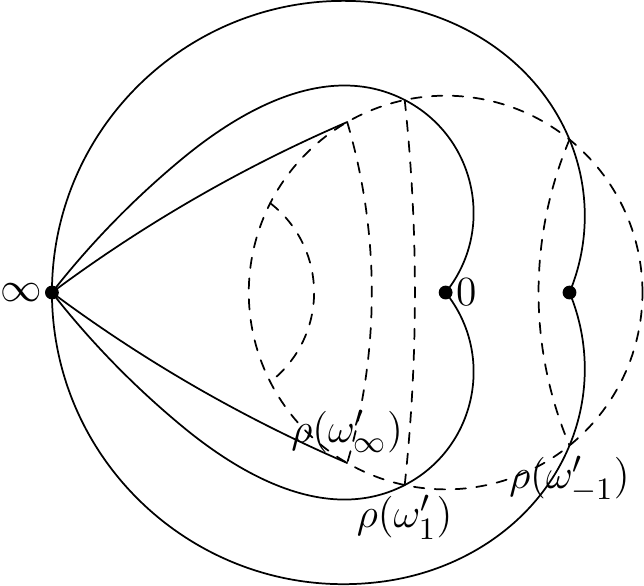}

Image under $\rho (.,s_{3/7})$\end{center}


The images under $\rho $ of $\omega _{\pm 1}'$ and $\omega 
_{\infty }'$ have been drawn with solid lines and labelled, but not 
$\omega _{0}'$, 
because the regions $V_{3}(a_{0})$ and $V_{3}(\overline{a_{0}})$ are 
best described using $\rho (.,s_{1/7})$ or $\rho 
(.,s_{6/7})$. (The images under $\rho $ of $\omega _{\pm 1}$ and $\omega _{\infty }$ have also been drawn, with solid lines, but not labelled.) The pictures are correct up to isotopy preserving 
$Z_{m}(s)$ for any $m\geq 0$. To see that this is the correct 
picture, we 
just consider $\rho (\omega _{1}',s_{3/7})$ since the other cases 
are exactly similar. Apart from one common boundary point between 
hyperbolic components, the  path $\omega _{1}'$ is entirely in 
$H_{a_{1}}\cup H_{1}$. The common boundary 
point represents the parabolic point on $\partial H_{a}$, where the 
centre of $H_{a}$ is Thurston 
equivalent to $s_{3/7}\Amalg s_{1/7}$, this boundary point also being 
in  $\partial H_{a_{1}}\cap \partial H_{1}$. If we 
parametrise $\omega _{1}'$ by $t\in [0,1]$ with $\omega 
_{1}'(0)=a_{1}$ 
and $\omega _{1}'(1)=1$, then we can choose $\omega _{1}'$ in its 
homotopy class so that $\omega _{1}'({1\over 2})$ is the parabolic 
parameter value which is a common boundary point, and for all 
$t<{1\over 2}$ or $t=\frac{1}{2}$ or $t>\frac{1}{2}$ respectively, 
$v_{2}(\omega
_{1}'(t))$ is in the fixed attractive basin of $h=h_{\omega 
_{1}'(t)}$,
or in the parabolic basin of $c_{1}(h_{\omega _{1}'(t)})$, or in the
attractive basin of $c_{1}(h_{\omega _{1}'(t)})$.  We then see that
$\rho (\omega _{1}',s_{3/7})=\beta _{5/7}$, as shown in the figure, and this
realises the Thurston equivalence between $h_{1}$ and the type II capture $\sigma _{5/7}\circ s_{3/7}$ defined in \ref{2.7}.

Now we describe parts of the domain of the maps $\rho 
(.,s_{1/7})$ and $\rho (.,s_{6/7})$.
Recall that the arcs $\omega 
_{0}'$ and  $\omega _{\pm 1}'$ run close to a common path until they reach 
the 
upper-half-plane common boundary component $x$ of the hyperbolic
components containing $\pm 1$.  A similar statement holds for 
$a_{0}$ and  $x$
replaced by $\overline{a_{0}}$ and $\overline{x}$, and with $\omega 
_{\pm 1}'$ and
$\omega _{0}'$ replaced by $\omega _{\pm 1}$ and $\omega _{0}$.  Now we
take adjustments $\omega _{0,a_{0}}'$ and $\omega _{\pm
1,a_{0}}'$ of the arcs $\omega _{0}'$ and $\omega _{\pm 1}'$,
arcs in $V_{3,0}$ which all start at $a_{0}$, because we
shall want to determine other arcs by their images under $\rho (.,s_{1/7})$.  To do this, we simply choose an arc $\gamma _{0}'$
from $a_{0}$ to $a_{1}$, and define $\omega
_{i,a_{0}}$ to be $\gamma _{0}'*\omega _{i}'$, as an element
of $\pi _{1}(V_{3,0},P_{3,0},\overline{a_{0}})$.  So it remains to
define $\gamma _{0}'$, as an element of $\pi
_{1}(V_{3,0},a_{1},a_{0})$.    Note that we can then take $\omega _{0,a_{0}}'$ to lie
in the hyperbolic component of $h_{a_{0}}$, except for the
endpoint at $0$, and we shall do so.

The arc $\gamma _{0}'$ 
can be determined up to isotopy in $V_{3,m}$, for any $m$, by 
describing the 
image of the reverse 
path $\overline{\gamma _{0}'}$ under $\rho 
(.,s_{3/7})$ up to isotopy preserving $Y_{m}(s_{3/7})$ for any $m$. 
Even though $\overline{\gamma _{0}'}$ is not in the domain of $\rho 
(.,s_{3/7})$
as prevously defined, any path in $V_{3,m}$ between $a_{1}$ and either
$a_{0}$ or $\overline{a_{0}}$ does determine an element $[\beta ]$ of 
$\pi
_{1}(\overline{\mathbb C}\setminus Z_{m}(s_{3/7}),\infty )$,
and a Thurston equivalence
$$s_{3/7}\simeq \sigma _{\beta }\circ s_{3/7}.$$
 The figure in \ref{3.1} should help: this figure gives
two invariant circles under $s_{3/7}\Amalg s_{1/7}$, the usual unit
circle being $\gamma _{1}$ and the other being $\gamma _{2}$, which is
drawn up to isotopy preserving $Y_{0}(s_{3/7})$ in the figure.  Up to
isotopy preserving $Y_{0}(s_{3/7})$, the path $\rho (\overline{\gamma
_{0}'},s_{3/7})$ crosses both $\gamma _{1}$ and $\gamma _{2}$ between
$e^{2\pi i(5/7)}$ and $e^{2\pi i(6/7)}$, and has no further crossings
with $\gamma _{2}$ before returning to $v_{2}$.  We simply refine this
description to give $\rho (\overline{\gamma _{0}'},s_{3/7})$ up to
isotopy preserving $Y_{m}(s_{3/7})$: just one essential crossing with
$\gamma _{2}$.  It follows that $\rho (\overline{\gamma _0'}*\omega _{i,a_{0}}',s_{3/7})$
also just has one essential crossing with $\gamma _{2}$, for $i=\pm
1$, and hence $\rho (\omega _{i,a_{0}}',s_{1/7})$ has just one
crossing with the unit circle.  Examining the diagram in \ref{2.10},
we see that the crossing is at $e^{2\pi i(4/7)}$, up to isotopy
preserving $Z_{m}(s_{1/7})$, for any $m$.  Write $\beta _{p,x}$ for
an arc from $v_{2}=\infty $ to $x$, which crosses the unit circle only
at $e^{2\pi i(p)}$, passing from there into the gap of $L_{1/7}$ 
containing $x$,
for $x=c_{1}$ or $s_{1/7}(v_{1})$ and $p=1/7$ or $2/7$ or $4/7$.  Then we have shown that
 \begin{equation}\label{4.5.3}
     \rho (\omega _{1,a_{0}}',s_{1/7})=\beta _{4/7,c_{1}},\ \ 
     \rho (\omega _{-1,a_{0}}',s_{1/7})=\beta 
     _{4/7,s(v_{1})}.\end{equation}
We write $\beta _{4/7,c_{1},s_{1/7}(v_{1})}$ for an arc joining 
$c_{1}(=0)$ and
$s_{1/7}(v_{1})$, passing between the gaps containing $c_{1}$ and
$s_{1/7}(v_{1})$ via the common boundary point $e^{2\pi i(4/7)}$. 
  Then it follows immediately from (\ref{4.5.3})
that \begin{equation}\label{4.5.2}\overline{\rho (\omega
_{1,a_{0}}',s_{1/7})}* \rho (\omega _{-1,a_{0}}',s_{1/7})=\beta
_{4/7,c_{1},s(v_{1})}.\end{equation} 
Similar formulae to (\ref{4.5.3})
and (\ref{4.5.2}) hold for $\overline{a_{0}}$ and $s_{6/7}$ replacing
$a_{0}$ and $s_{1/7}$, for loops $\omega _{i,\overline{a_{0}}}=\gamma
_{0}*\omega _{i}$, for $\gamma _{0}$ defined similarly to $\gamma
_{0}'$.  Even without the precise definition of $\gamma _{0}'$ given
above, it is true that the homotopy class of the arc $\overline{\omega
_{-1,a_{0}}}'*\omega _{1,a_{0}}'$ is uniquely determined in $V_{3,m}$,
up to homotopy fixing endpoints, but one cannot canonically define the
image under $\rho $ of $\overline{\omega _{-1,a_{0}}'}*\omega
_{1,a_{0}}'$.

We also define two more arcs $\omega _{0,\pm 1,a_{0}}'$
from $a_{0}$ to $\pm 1$, by taking a path close to $\omega 
_{0,a_{0}}'$
followed by an arc in $H_{\pm 1}$.  So these arcs can be taken to lie
in the union of the hyperbolic component of $h_{a_{0}}$ and either
$H_{1}$ or $H_{-1}$, apart from arbitrarily small neighbourhoods of
$0$.  So the arcs $\omega _{0,a_{0}}'$ and $\omega _{0,\pm 1,a_{0}}'$
are uniquely determined as elements of $\pi (V_{3,m},P_{3,m},a_{0})$,
for any $m$.  We write $\beta _{1/7,c_{1}}$ for an arc from
$v_{2}=\infty $ to $c_{1}$, which crosses the unit circle only at
$e^{2\pi i(1/7)}$, passing from there into the gap containing $c_{1}$,
and similarly for $\beta _{2/7,s(v_{1})}$, which crosses the unit
circle at $e^{2\pi i(2/7)}$, passing into the gap containing
$s_{1/7}(v_{1})$.  Then the homotopy classes of these two paths in
$\overline{\mathbb C}\setminus Z_{m}(s_{1/7})$ are uniquely
determined.  Now $\omega _{0,1,a_{0}}'$ is an arc in $H_{a_{0}}\cup
H_{1}$, apart from an arbitrarily small neighbourhood of $0$, which is
a deleted common boundary point, while $\beta _{1/7,c_{1}}$ is the
path defining the corresponding type II capture, and similarly with
$1$ replaced by $-1$.  So 
\begin{equation}\label{4.5.1}\rho (\omega _{0, 1,a_{0}}',s_{1/7}) =\beta _{1/7,c_{1}},\ \rho (\omega _{0, -1,a_{0}}',s_{1/7}) =\beta
_{2/7,s_{1/7}(v_{1})}.\end{equation}
 
It therefore makes sense to 
use $\rho (.,s_{3/7})$ to describe 
$V_{3,m}(a_{1},\pm )$, to use $\rho (.,s_{1/7})$ to describe 
$V_{3,m}(a_{0})$, and to use $\rho (.,s_{6/7})$ to describe 
$V_{3,m}(\overline{a_{0}})$.

\chapter{ Fundamental Domains}\label{5}

\section{Fundamental domains: a restricted class}\label{5.1}

The proof of Theorem \ref{2.10} is simply the 
construction of a fundamental domain for the action of $\pi 
_{1}(V_{3,m})$ on $D'$. In this section we consider the general 
problem of constructing a fundamental domain for the action of a 
finitely-generated discrete group $\Gamma $ of M\" obius 
transformations, 
acting freely on the open unit disc $D$. We specialise to $V_{3,m}$ 
in \ref{5.4}. We are interested in fundamental domains only up to 
homeomorphism. Then $F\subset D$ is a {\em{fundamental domain 
for $\Gamma $}} if 
$$D=\cup \{ \gamma .F:\gamma \in \Gamma \} ,$$
and $\gamma .({\rm{int}}(F))\cap {\rm{int}}(F)=\emptyset $ for all 
$\gamma \neq {\rm{identity}}$, $\gamma \in \Gamma $. A fundamental 
domain always exists. If $\Gamma $ is a free group, and has at least 
one 
parabolic 
element, $D/\Gamma $ is known to be a compact surface minus at least 
one puncture, and so $\Gamma $ is a free group containing parabolic 
elements. Then we can, and do, choose  $F$ bounded by finitely 
many smooth arcs 
with both ends at lifts of punctures, with 
$\overline{F}$ intersecting $\partial D$ only in lifts of parabolics. We shall give a vague restatement of the Main Theorem \ref{2.10} in \ref{5.7}, and refine this, in restatements in two separate cases, in \ref{6.1}, \ref{7.8}.

\section{Fundamental Domains: construction from graphs and from 
matching pairs of adjacent path pairs}\label{5.2}

The union of the punctures and the  projection of $\partial F$ to 
$D/\Gamma $ is a graph whose complement is a topological disc. 
Conversely, let $G$ be any graph in $V=D/\Gamma $ which is a union of 
smooth arcs and such that, allowing vertices at punctures, the 
complement is a topological disc. Then any component of the lift to 
the universal 
cover of the 
complement  of $G$ is a fundamental domain for $\Gamma $. 

Let $P$ denote 
the set of punctures of $V=D/\Gamma $, and fix a basepoint $x_{0}\in 
V$. Then another equivalent way of choosing a fundamental domain with 
vertices at lifts of punctures is to choose a set of arcs 
$\Omega $
from $x_{0}$ to $P$, such that the interiors of the arcs are disjoint 
and they
all represent distinct elements of $\pi 
_{1}(V,P,x_{0})$, such that each component $U$ of $V\setminus \cup 
\Omega $ is bounded by either three or four arcs of $\Omega $, and is 
disjoint 
from $P$, but has exactly two points of $P$ in its boundary, provided 
that $\# (P)>2$, as we are assuming. Here is an example. The arcs of 
$\Omega $ are shown as solid lines, and the edges of the 
corresponding graph $G_{\Omega }$ -- which is a tree --- by dashed 
lines. In this rather small example, there is just one component of 
$V\setminus \cup \Omega $ which has four arcs of $\Omega $ in its 
boundary. All the other components have just three arcs.
\begin{center}
\includegraphics[width=4cm]{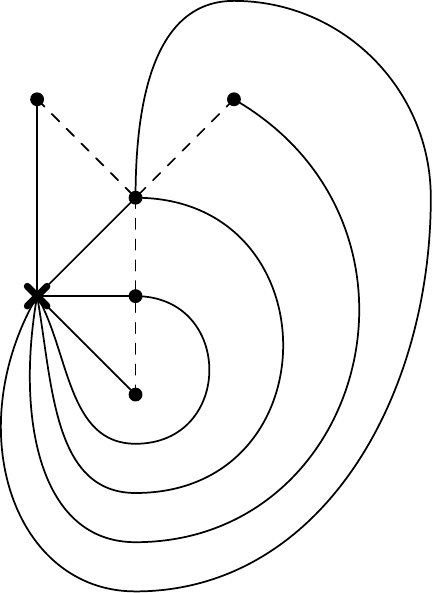}

An example of a graph
\end{center}
There is exactly one edge of $G_{\Omega }$ in each component $U$ of 
$V\setminus \cup \Omega $, joining the 
unique pair of points in $P\cap \partial U$, and these are all the 
edges of $G_{\Omega }$. Two of the arcs of 
$\Omega $ in $\partial U$ are accessible from each side of $e(U)$ in 
$U$. This gives two pairs of arcs $(\omega _{1},\omega _{2})$, 
$(\omega _{1}',\omega _{2}')$, in no specified order, such that one 
pair is accessible from $U$ from one side of $e(U)$ and the other 
pair from the other side. One of the arcs $\omega _{1}$, $\omega 
_{2}$ 
might coincide with one of $\omega _{1}'$, $\omega _{2}'$ but at 
least 
three of the arcs are distinct and these comprise all the arcs of 
$\Omega $ in $\partial U$. We call such a pair of pairs 
$((\omega _{1},\omega _{2}), (\omega _{1}',\omega _{2}'))$ 
a {\em{matching pair of adjacent pairs}}. Note that $\omega _{1}$ and  
$\omega _{2}$ are indeed adjacent in $\Omega $, and similarly for 
$\omega _{1}'$ and $\omega _{2}'$. By this we mean the following. 
Take lifts of all  arcs in $\Omega $ to the 
universal cover, with lifts all starting from the same primage 
$\tilde{x}_{0}$ of $x_{0}$. Let $\tilde{\omega }$ denote the lift of 
$\omega $ and $\tilde{\Omega }=\{ \tilde{\omega }:\omega \in \Omega 
\} $. Then $\omega _{1}$ and $\omega _{2}$ are said to be adjacent 
in $\Omega $ if $\tilde{\omega }_{1}$ and $\tilde{\omega }_{2}$ are 
not separated by any other arc in $\tilde{\Omega }$. We can 
use this concept of adjacency if $\Omega $ is any subset of $\pi 
_{1}(V,P,x_{0})$, not necessarily a set of disjoint arcs.  
Conversely, to construct 
$\Omega =\Omega _{G}$ from $G$, we need $x_{0}\notin G$, and we then 
take $\Omega $ to be the set of arcs, starting from $x_{0}$ and which,
up to homotopy, intersect 
$G$ only at endpoints in $P$, with one path ending at each 
sector 
of the complement of $G$, at each point of $P$. The set of matching 
pairs of adjacent pairs in $\Omega $ is then determined from $G$ as 
above. Any adjacent pair in $\Omega $ does occur as one of a matching 
pair of adjacent pairs for exactly one edge, and is thus matched with 
exactly one other adjacent pair. The correspondence 
$G\mapsto \Omega $ is therefore a bijection.

Given a fundamental domain $F$ for some group $\Gamma $ which has all 
vertices on $\partial D$, and given a vertex $v$ of $F$, we can 
define the {\em{group element of $v$}} to be the 
product $g_{r}\cdots g_{1}$ defined as follows. Put 
an anticlockwise orientation on $\partial F$ round $F$. Let $\ell 
_{1}$ be the geodesic in $\partial F$ starting at $v$. Now 
inductively define a finite
sequence of geodesics $\ell _{i}$ in $\partial F$, and group 
elements $g_{i}\in \Gamma \setminus \{ 1\} $, such that 
$g_{i}\ell _{i}$ is also in $\partial 
F$. Thus $\ell _{i}$ uniquely determines $g_{i}$. Note that each 
$g_{i}:\ell _{i}\to g_{i}\ell _{i}$ reverses orientation.  
The inductive definition of $\ell 
_{i+1}$ is that $\ell _{i+1}$  is the geodesic which starts at the 
vertex where $g_{i}\ell _{i}$ finishes. Thus, $\ell _{i+1}$ starts at 
the vertex $g_{i}\cdots g_{1}v$. We take $r$ to be the least integer 
such that $\ell _{r+1}=\ell _{1}$, equivalently, the least integer 
such that $g_{r}\cdots g_{1}v=v$. Then $g_{r}\cdots g_{1}$ 
 is automatically parabolic.

Let $F$ be a fundamental domain in $D$, $G$ the corresponding graph 
in 
the quotient space $V$ with puncture set $P$, $\Omega \subset \pi 
_{1}(V,x_{0},P)$ the corresponding arc set. Let $A\subset \pi 
_{1}(V,x_{0})$ be the set of elements $g$ such that $g.F$ has a 
common edge with $F$. Then $A=A^{-1}$ generates $\pi _{1}(V,x_{0})$. 
The elements of $A$ are 
in one-to-one correspondence with the edges of $F$, in two-to-one 
correspondence with the edges of $G$, with the two elements 
associated to an edge $e$ being of the form $g^{\pm 1}$, and in 
one-to-one correspondence with ordered matched pairs of adjacent 
pairs of elements of $\Omega $. This is the set of all $(\omega _{1},
\omega _{2},\omega _{1}',\omega _{2}')$, where $\omega _{1}$ and 
$\omega _{2}$ are adjacent in $\Omega $ with $\omega _{2}$ 
anticlockwise from $\omega _{1}$, and  $\omega _{1}'$ and 
$\omega _{2}'$ are adjacent in $\Omega $ with $\omega _{2}'$ 
clockwise from $\omega _{1}'$,  and $\omega _{i}$ and $\omega _{i}'$ end 
at the same point of $P$ for $i=1$, $2$, and the disc with 
anticlockwise boundary made up of  
$\omega _{1}$, $\overline{\omega _{1}'}$, $\omega _{2}'$, 
$\overline{\omega _{2}}$ is disjoint from $P$. The corresponding 
element of $A$ then satisfies $g.\omega _{i}=\omega _{i}'$ for $i=1$, 
$2$, using the usual action of $\pi _{1}(V,x_{0})$ on $\pi 
_{1}(V,x_{0},P)$. There is only one such element $g\in \pi 
_{1}(V,x_{0})$. For suppose there are two such, $g_{1}$ and $g_{2}$. 
Then $g_{2}^{-1}g_{1}\neq 1$ fixes  the geodesic in $D$ which is the 
side of $F$ joining the endpoints of lifts of $\omega _{1}$ and 
$\omega _{2}$. It also fixes the endpoints. But there are also 
parabolic elements of $\pi _{1}(V,x_{0})$ fixing these endpoints. For 
a discrete group of hyperbolic isometries, this is impossible.

The vertex group elements can then be computed from the set of 
matching pairs of adjacent pairs, and the group elements matched 
pairs, by taking a cycle $(\omega _{1},\omega _{2},\omega 
_{1}',\omega 
_{2}')$, $(\omega _{1}',\omega _{3}',\omega _{1}'',\omega _{2}'')$ 
\ldots 
and taking the product of group elements for each successive matching 
pair of adjacent pairs in the cycle.

\section{Examples of Fundamental domains: the graph from a Julia 
set}\label{5.3}

In this paper, we are concerned with the punctured spheres $V_{3,m}$. 
As we have seen, this means we are also interested in the punctured 
sphere $\overline{\mathbb C}\setminus Z_{m}(g_{0})$ for various maps 
$g_{0}\in B_{3,m}$, special consideration being given to the maps 
$g_{0}=s_{p}$ for $p=\frac{1}{7}$, $\frac{6}{7}$, $\frac{3}{7}$. A 
simple way to 
construct a fundamental domain $F$ for 
$\overline{\mathbb C}\setminus Z_{m}(s_{p})$, or, more accurately, 
for 
its fundamental group, is to take as the graph $G$ the tree 
$T_{m}(s_{p})$ which has as vertices the points of $Z_{m}(s_{p})$, and 
is defined as follows. Two points $z_{1}$ and 
$z_{2}\in Z_{m}(s_{p})$ are joined by an edge in $T_{m}(s_{p})$ if 
and only if the gaps $G_{1}$ and $G_{2}$ of $L_{p}$ containing them are adjacent in the 
unit disc, that is, not separated by any other gap containing a point 
of $Z_{m}(s_{p})$, and either one of $G_{1}$ or $G_{2}$ is the  gap containing the critical value, or separates the other one from this gap. Then all edges are taken in the unit disc and 
transverse to $L_{p}$. This completely determines $T_{m}(s_{p})$, up 
to isotopy. Then $T_{m}(s_{p})$ determines an arc set, 
together with a matching of pairs of 
adjacent pairs. This example will be important when we come to restate  
our main theorem \ref{2.10} in \ref{5.7}, \ref{6.1}, \ref{7.8}.

\section{Fundamental domain for the whole group}\label{5.4}
The following gives a sufficient condition for a subset  $F$ of the 
unit disc, with a pairing of edges via elements of $\Gamma $,
to be a fundamental domain for a group $\Gamma $, rather than for a 
subgroup of $\Gamma $. The lemma implies that, in order to construct 
a fundamental domain for $V$, we only need to find  
$\Omega \subset \pi _{1}(V,P,x_{0})$, together with a set of matching 
pairs of adjacent pairs, such that each adjacent pair is matched to 
exactly one other, such that the vertex elements 
are all parabolic, and a rather mild extra condition. We do not need 
to know that the elements of $\Omega $  are 
all represented by arcs in $V$. That is a consequence, but not a 
prior 
condition. The construction of $F$ from $\Omega $ in \ref{4.2} uses 
lifts of $\Omega $ from a chosen lift $\tilde{x_{0}}$ of $x_{0}$. The 
lifts of paths in $\Omega $ emanating from $\tilde{x_{0}}$ are 
homotopic to geodesic rays in $D$, whether or not the original paths 
in $\Omega $ are arcs.

\begin{ulemma} Let $\Gamma $ be a discrete group of M\" 
obius transformations acting on the unit disc $D$. 
Let $F$ be a topological disc in $D$ bounded by $2r$ geodesics $\ell 
_{i}$ ($1\leq i\leq r$)  and $\gamma _{i}.\ell _{i}$, where $\gamma 
_{i}\in \Gamma $, each $\ell _{i}$ has endpoints at parabolic 
points of $\Gamma $, and $\gamma _{i}.F\cap F=\emptyset $. Then $F$ 
is a 
fundamental domain for a finite index subgroup $\Gamma _{1}$ of 
$\Gamma $, 
which is a free group on the $r$ generators $\gamma _{i}$. 

If, in 
addition, for each vertex $v$ of $F$, the group element of the vertex 
is a parabolic representing a puncture of $D/ \Gamma $ , and 
either 
$D/\Gamma $ 
is a punctured sphere, or at least one puncture on $D/\Gamma $ is 
represented by only one $\Gamma _{1}$-conjugacy-class of a group 
element of a vertex of $F$, then $F$ is a 
fundamental domain for $\Gamma $.\end{ulemma}

\noindent {\em{Proof.}} Let $\Gamma _{1}$ be the group generated by 
$\gamma _{i}$ for $1\leq i\leq r$. Let $w$, $w'$ be any two words 
in $\gamma _{i}^{\pm 1}$ for $1\leq i\leq r$ in which $\gamma _{i}$ 
and $\gamma _{i}^{-1}$ are never adjacent. Then by induction, $w.F$ 
and $w'.F$ are disjoint if $w\neq w'$. In fact $w.F$ and $w'.F$ 
are separated by $w''.F$, where $w''$ is the longest common prefix 
(possibly trivial) of $w$ 
and $w'$, unless one of $w$, $w'$ is a prefix of the other.
If $w$ is a prefix of $w'$, then $w.F$ and $w'.F$ are 
adjacent if $w'$ has one more letter than $w$, and separated by 
$w''.F$ 
if $w'$ has at least two more letters than $w$, and $w''$ is the 
prefix 
of $w'$ which has one more letter than $w$. So $\Gamma _{1}$ is a 
free 
group, and discrete, since it is a subgroup of $\Gamma $. Now we 
claim 
that $D=\Gamma _{1}.F$. If not, then there is some sequence $w_{n}$ 
of 
increasing words in $\{ \gamma _{i}^{\pm 1}:1\leq i\leq r\} $ such 
that 
$w_{n}.F$ does not converge to a point. Now the minimum distance of 
$F$ from $w.F$ is bounded from $0$ unless $w=\gamma _{i}^{m}$ for 
some $1\leq i\leq r$ and some $m\in {\mathbb Z}$. So the increasing 
sequence of words $w_{n}$ must be of the form $w_{k}\gamma _{i}^{\pm 
(n-k)}$ for some $w_{k}$ and $\gamma _{i}$. But for such a sequence, 
$w_{n}.F$ converges to a point. So now $D=\Gamma _{1}.F$, and since 
$F$ has finite area, $\Gamma _{1}$ must be of finite index in $\Gamma 
$.

If $\Gamma _{1}\neq \Gamma $, then $D/\Gamma _{1}$ is a finite cover 
of $D/\Gamma $. We call the covering map $\pi $. If the group element 
of each vertex of $F$ 
is a simple parabolic 
element of $\Gamma $, then we can fill in punctures to obtain closed 
surfaces $S_{1}$ and $S$ with marked points, and we can extend $\pi $ to 
a covering map $\pi :S_{1}\to S$, such that marked points map to 
marked points, and the inverse image of every marked point of $S$ is a 
marked 
point of $S_{1}$. If we have an $n$-fold covering for $n>1$, then 
each marked point on $S$ has $n$ preimages in $S_{1}$,
and $\chi (S_{1})=n\chi (S)$, where $\chi $ denotes Euler 
characteristic. So we must have $\chi (S)\leq 0$, and $S$ cannot be a 
sphere, We must have $n=1$ if some puncture is represented by only one $\Gamma _{1}$ conjugacy class of a group element of a vertex of $F$.
\Box

\section{Fundamental domains for increasing numbers of 
punctures}\label{5.5}
 
Let $V$ be a surface with a finite puncture set $P$. Now let 
$Q\subset V$ be 
a finite set with $P\subset Q$. Let $F_{P}$ be a 
fundamental domain 
for 
$V\setminus P$ with all vertices at $P$. Then we can choose a 
fundamental domain for $V\setminus Q$ by the following 
procedure. We can regard $F_{P}$ as a topological disc in $V$, and 
after homotopy on the boundary, we can assume that $\partial 
F_{P}\cap Q=\emptyset$.
 Fix $x_{0}\in V\setminus \partial F_{P}$. We have a graph $G_{P}$ 
 and set  $\Omega 
 _{P}$ of arcs from $x_{0}$ to $P$, derived from $F_{P}$, as explained 
 in \ref{5.2}. We also have matching pairs of adjacent pairs in 
 $\Omega _{P}$, as explained in \ref{5.2}. Let $((\omega _{1},\omega 
 _{2}),(\omega _{1}',\omega _{2}'))$ be any matching pair of adjacent 
 pairs in $\Omega _{P}$. Then up to homotopy, $\omega _{1}\cup \omega 
 _{2}\cup \omega _{1}'\cup \omega _{2}'$ bounds a topological disc 
$U$ 
in $V$, containing a single edge $e(U)$ of the graph $G_{P}$, which is 
also 
the projection of two matched edges of $F_{P}$.
 We then modify $F_{P}$ to a fundamental domain $F_{Q}$ 
for $V\setminus Q$, by changing $e(U)$, 
for each such $U$, to a tree $G(e(U))\subset U\cup 
{\rm{ends}}(e(U))$, satisfying the following properties:
\begin{itemize}
\item $G(e(U))$ is a tree which is contained in $U$ apart from the  two ends of $e(U)$, which are extreme vertices of $G(e(U))$;
\item the other vertices of $G(e(U))$ are all the  points of $Q$ in $U$.\end{itemize}
The graph $G_Q$ is  the union of all the trees $G(e(U))$. Then $G_Q$ determines a fundamental domain for $V\setminus Q$. 
The arc set $\Omega _{Q}$ contains $\Omega _{P}$ up to homotopy. Each 
matching pair 
of adjacent pairs is contained in $U\cup \omega _{1}\cup \omega 
_{2}\cup \omega _{1}'\cup \omega _{2}'$ for some $((\omega _1,\omega _2),(\omega _1',\omega _2'))$ and $U$ as above.

Let $K$ be the kernel of the forgetful homomorphism $\pi 
_{1}(V\setminus Q)\to \pi _{1}(V\setminus P)$. Each edge $e$ of 
$G_{P}$ has two lifts in $\partial F_{P}$, which are identified by 
$g_{e}\in \pi _{1}(V\setminus P)$. The element $g_{e}$ is unique up 
to 
replacing it by its inverse. Let $e'$ be an edge of $G(e)\subset 
G_{Q}$.
If one end of $e'$ is an extreme vertex of $G(e)$, that is, not an 
end of any other edge of $G_{Q}$, then $g_{e'}\in K$. Otherwise, 
$e'$ extends to a path in $G(e)$ between the endpoints of $e$, this 
path being unique up to homotopy in $G(e)$, and $g_{e'}$ projects 
under the forgetful homomorphism to $g_{e}$ or
$g_{e}^{-1}$.

Let $\Sigma _{P}$ be the (finite disjoint) union of all sets $U$ for 
$F_{P}$,
as above. Each $\partial U$ is the 
union of two arcs $\overline{\omega _{1}}*\omega _{2}$ and 
$\overline{\omega _{1}'}*\omega _{2}'$. If we have a sequence of 
fundamental domains $F_{m}$, with each $F_{m+1}$ constructed from 
$F_{m}$ as $F_{Q}$ is constructed from $F_{P}$ above, then we get a 
sequence of sets $\Sigma _{m}$. 

Topological discs of the type $U\subset V\setminus P$ will play an 
important role in the construction of the fundamental domain for  
$V_{3,m}$. 
Let $K$ denote the kernel of the forgetful 
homomorphism $\pi _{1}(V\setminus Q)\to \pi _{1}(V\setminus P)$. 
Then $\partial U$ represents a conjugacy class in $K$, if we perturb 
$U$ isotopically to a set disjoint from $P$, inside the original $U$.
For then $\partial U $ is a simple loop which bounds a loop disjoint 
from $P$, but possibly not disjoint from $Q$. Now $\pi 
_{1}(V\setminus Q)$ acts on 
$D$, and the 
quotient space $D/K$ is homeomorphic to a punctured disc, on which 
$\pi_{1}(V\setminus P)\cong \pi _{1}(V\setminus Q)/K$ acts. 
Filling in the punctures, we obtain the universal cover of 
$V\setminus P$, which is again the unit disc $D$ up to conformal 
isomorphism, and the canonical action of of $\pi _{1}(V\setminus 
P)$ on $D$. Any lift $\tilde{U_{1}}$ of $U$ to $D/K$ 
projects homeomorphically to $U$ and any two distinct lifts 
$\tilde{U_{1}}$ and $\tilde{U_{2}}$ are disjoint, but in the same 
orbit under the action of $\pi _{1}(V\setminus Q)/K$. Any two 
lifts $\tilde{U_{1}}$ and $\tilde{U_{2}}$ of $U$ to the universal 
cover $D$ of $\pi _{1}(V\setminus Q)$ are either equal or 
disjoint, with $\tilde{U_{2}}=g.\tilde{U_{1}}$ for at least one 
$g\in \pi _{1}(V\setminus Q)$. The stabiliser of 
$\tilde{U_{1}}$ in $\pi _{1}(V\setminus Q)$ is a subgroup of 
$K$. 

\section{Specifying to $V_{3,m}$}\label{5.6}

To construct a fundamental domain $F_{m+1}$ for $V_{3,m+1}$ from a 
fundamental domain $F_{m}$ for $V_{3,m}$, we shall use the procedure 
outlined in \ref{5.5}, constructing a fundamental domain for 
$V_{3,m+1}$ from a fundamental domain for $V_{n,m}$. This means that 
we start with a fundamental 
domain for $V_{3,0}$. By \ref{5.2}, this is equivalent to having  a 
set $\Omega 
_{0}\subset \pi _{1}(V_{3,0},a_{1},P_{3,0})$ with matching pairs of 
adjacent pairs. We use the set $\Omega _{0}$ of \ref{4.5}.

Again, by
\ref{5.2}, the construction of a fundamental domain for $V_{3,m}$ is 
equivalent to constructing an arc set $\Omega _{m}$ with similar 
properties, with matching pairs of adjacent pairs. Note that we do 
not 
need to know that $\Omega _{m}$ is a set of disjoint arcs,
if we can deduce it later on, using the lemma in \ref{5.4}. So we only
need a set of paths  $\Omega _{m}\subset \pi 
_{1}(V_{3,m},a_{1},P_{3,m})$ 
such that the 
vertex group elements of the fundamental domain (a priori a 
fundamental domain for a subgroup: see \ref{5.4}) are all simple 
parabolics, 
that is, representing paths which go once round 
punctures of $V_{3,m}$. We shall define
$$\Omega _{m}=\{ \gamma _{0}'*\omega :\omega \in \Omega _{m}(a_{0})\} 
\cup \{ \gamma _{0}*\omega :\omega \in \cup \Omega 
_{m}(\overline{a_{0}})\}  \cup \Omega _{m}(a_{1},-)\cup \Omega _{m}(a_{1},+),$$
where $\gamma _{0}'$ and $\gamma _{0}$ are as in \ref{4.5}. The set $\gamma _{0}'*\Omega 
_{m}(a_{0})$ will consist of the paths of $\Omega _{m}$ with 
endpoints in 
$P_{m}(a_{0})$, and the paths of $\Omega 
_{m}(a_{0})$ will have endpoints in $P_{m}(a_{0})$, and otherwise
lie entirely in 
$V_{3,m}(a_{0})$. Similar properties hold for $\gamma _{0}$ and
$\Omega _{m}(\overline{a_{0}})$. The paths of  $\Omega 
_{m}(a_{1},-)$ will have endpoints in $P_{3,m}(a_{1},-)$, and 
otherwise lie entirely in $H_{1}\cup 
\{ x,\overline{x}\} \cup V_{3,m}(a_{1},-)$ and similarly for $\Omega 
_{m}(a_{1},+)$. We shall construct $\Omega _{m}(a_{0})$ by 
constructing the image under $\rho (.,s_{1/7})$, and $\Omega 
_{m}(\overline{a_{0}})$ by constructing the image under $\rho 
(.,s_{6/7})$ and $\Omega _{m}(a_{1},\pm )$ by constructing the 
images under $\rho (.,s_{3/7})$. The sets $\Omega _{m}(a_{0})$, 
$\Omega _{m}(\overline{a_{0}})$, and $\Omega _{m}(a_{1},-)$ will be 
constructed in Section \ref{6}, and the sets $\Omega _{m}(a_{1},+)$ 
will be constructed in Section \ref{7}.

Now provided we know that $\rho (\Omega 
_{m}(a_{0}),s_{1/7})\subset D'$ for $D'=D'(s_{1/7})$ as in \ref{4.3}, 
and similarly for $\Omega _{m}(\overline{a_{0}})$ and $\Omega 
_{m}(a_{1},\pm )$, the Resident's View implies that adjacency of 
paths 
in $\Omega _{m}$ transfers under $\rho $. We  need to use $\rho 
(.,s_{p})$ for $p=\frac{1}{7}$ or $\frac{6}{7}$ or $\frac{3}{7}$, for 
different subsets of $\Omega _{m}$. So let $\zeta 
_{1}$ and  $\zeta  _{2}\in \Omega 
_{m}(a_{0})$ or $\Omega 
_{m}(\overline{a_{0}})$ or  $\Omega 
_{m}(a_{1},\pm )$, and let $\beta _{i}=\rho (\zeta 
_{i},s_{p})$ for $i=1$ 
and $2$ and $p={1\over 7}$, ${6\over 7}$ or ${3\over 7}$ in the 
respective cases. Then $\beta _{1}$ and $\beta _{2}$ are adjacent in 
$\rho (\Omega _{m},s_{p})$ if and only if $\zeta _{1}$ and $\zeta 
_{2}$ are adjacent in $\Omega _{m}$. So each adjacent pair of paths 
in $\Omega _{m}$ will lie in one of the following sets:
$$\gamma _{0}'*(\Omega _{m}(a_{0})\cup \{ \omega _{0,a_{0}},\omega 
_{\pm 1,a_{0}},\omega _{0,\pm 1,a_{0}}\} ),$$
and similarly for $a_{0}$ replaced by $\overline{a_{0}}$ and $\gamma 
_{0}'$ replaced by the corresponding path $\gamma _{0}$ in the lower 
half-plane, or
$$\Omega _{m}(a_{1},-)\cup \{ \omega _{-1},\omega _{-1}'\} ,$$
$$\Omega _{m}(a_{1},+)\cup \{ \omega _{1},\omega _{1}',\omega _{\infty 
}\} .$$
We recall from \ref{4.5} that $\gamma _{0}'*\omega _{i,a_{0}}=\omega _{i}'$ for $i=0$, 
$\pm 
1$. Also $(\gamma _{0}'*\omega _{0,1,a_{0}},\omega _{1}')$ is an 
adjacent pair 
matched with $(\gamma _{0}*\omega _{0,1,\overline{a_{0}}},\omega 
_{1})$ and similarly with $1$ replaced by $-1$. So to find a complete 
set of  
matching pairs of adjacent pairs in $\Omega _{m}$, it suffices to 
achieve a 
complete set of matching pairs  of adjacent pairs in each of the 
following sets:
\begin{equation}\label{5.6.1}\Omega _{m}(a_{0}),\ \Omega 
_{m}(\overline{a_{0}}),\ \Omega 
_{m}(a_{1},-),\ \Omega _{m}(a_{1},+).\end{equation}
We then define
\begin{equation}\label{5.6.2}\begin{array}{ll}R_{m}(a_{0})=\rho 
(\Omega 
_{m}(a_{0}),s_{1/7}),\ &R_{m}(\overline{a_{0}})=\rho (\Omega 
_{m}(\overline{a_{0}}),s_{6/7}),\cr 
R_{m}(a_{1},-)=\rho (\Omega _{m}(a_{1},-),s_{3/7}),\ &
R_{m}(a_{1},+)=\rho (\Omega _{m}(a_{1},+),s_{3/7}).\cr 
\end{array}\end{equation}

Adjacency is preserved by the mappings $\rho (.,s_{q})$ for 
$q={1\over 7}$, ${6\over 7}$ or ${3\over 7}$.
Matching pairs can also be viewed by considering the images under 
$\rho 
(., )$. Let $(\zeta _{1}',\zeta _{2}')$ be an adjacent pair in $\Omega _{m}$, and we consider how the property of matching with the adjacent pair $(\zeta _{1},\zeta _{2})$ transfers under $\rho $. Let $\gamma \in \pi 
_{1}(V_{3,m},s_{q})$ and let $\rho (\gamma ,s_{q})=\alpha $ and let $\Phi 
_{2}(\gamma )=[\psi ^{-1}]$.  Then, as we recalled in (\ref{4.3.1}) to
(\ref{4.3.3}) (but note that we are now replacing $\psi $ by $\psi 
^{-1}$),
$$\rho (\gamma *\zeta _{i}')=\alpha *\psi (\rho (\zeta 
_{i}')).$$
Hence, $\gamma *\zeta _{i}'=\zeta _{i}$ if and only if 
(\ref{5.6.3}) and (\ref{5.6.4}) hold for $\beta _{i}=\rho (\zeta 
_{i})$ and  $\beta _{i}'=\rho (\zeta _{i}')$:
\begin{equation}\label{5.6.3}\alpha 
*\psi (\beta _{i}')=\beta _{i}{\rm{\ rel\ }} 
Y_{m}(s_{q}),\end{equation}
\begin{equation}\label{5.6.4}(s_{q},Y_{m+1}(s_{q}))\simeq _{\psi 
}(\sigma _{\alpha }\circ s_{q},Y_{m+1}(s_{q})).\end{equation}
In particular, $\psi (\beta _{i}')$ 
and $\beta _{i}$ have the same endpoint, and the paths
$$\overline{\beta _{1}}*\beta _{2},\ \psi (\overline{\beta 
_{1}'}*\beta _{2}')$$
are homotopic via a homotopy preserving $Y_{m}(s_{q})$. For any 
$\beta _{i}$ and $\beta _{i}'\in \pi _{1}(\overline{\mathbb 
C}\setminus Z_{m}(s_{q}),Z_{m}(s_{q}),v_{2})$ for which (\ref{5.6.3}) 
and (\ref{5.6.4}) hold, we have $\alpha =\rho (\gamma )$ and 
$[\psi ]=\Phi _{2}(\gamma )$ for some $\gamma \in \pi _{1}(B_{3,m}, 
s_{p})$. 
But if in addition $\beta _{i}$ and $\beta _{i}'\in D'(s_{q})$ then we 
know that $\gamma \in \pi _{1}(V_{3,m}, s_{q})$, because the 
stabiliser 
of $D'$ in $\pi _{1}(B_{3,m}, s_{q})$ is $\pi _{1}(V_{3,m}, s_{q})$. 

So now we define a matching pair of adjacent pairs in  one of the 
sets of 
(\ref{5.6.2}) to be a pair of adjacent pairs $((\beta _{1},\beta 
_{2}),(\beta _{1}',\beta _{2}'))$ such that each $(\beta _{i},\beta 
_{i}')$ satisfies (\ref{5.6.3}) for some $\alpha \in \pi 
_{1}(\overline{\mathbb C}\setminus Z_{m}(s_{q}),v_{2})$ and  $[\psi 
]\in 
{\rm{MG}}(\overline{\mathbb C},Y_{m+1}(s_{q}))$ satisfying 
(\ref{5.6.4}). Then we can find a complete set of matching pairs of 
adjacent pairs in $\Omega _{m}$, if  we can find a complete set of 
matching pairs of adjacent pairs in each of the sets of (\ref{5.6.1}), 
and each of these sets lies in $D'(s_{q})$ for the appropriate value 
of $q$, for $q={1\over 7}$, ${6\over 7}$, or ${3\over 7}$. 

\section{}\label{5.7}
We are now ready to give a second, rather vague, statement of the Main Theorem.
\begin{maintheorem1.5}
A fundamental domain for $V_{3,m}$ can be 
constructed using a set $\Omega _{m}\subset 
\pi_{1}(V_{3,m},P_{3,m},a_{1})$ with
$$\Omega _{m}=\gamma _{0}'*\Omega _{m}(a_{0})\cup \gamma _{0}*\Omega 
_{m}(\overline{a_{0}})\cup \Omega _{m}(a_{1},-)\cup \Omega 
_{m}(a_{1},+),$$
$$\begin{array}{ll}R_{m}(a_{0})=\rho (\Omega 
_{m}(a_{0}),s_{1/7}),\ &R_{m}(\overline{a_{0}})=\rho (\Omega 
_{m}(\overline{a_{0}}),s_{6/7})\cr 
R_{m}(a_{1},-)=\rho (\Omega _{m}(a_{1},-),s_{3/7}),\ &
R_{m}(a_{1},+)=\rho (\Omega _{m}(a_{1},+),s_{3/7}).\cr 
\end{array}$$
Here: 
\begin{itemize}
\item $\gamma _{0}'$ and  $\gamma _{0}$ are the  paths from $h_{a_{1}}$ 
to $h_{a_{0}}$ $h_{\overline{a_{0}}}$ of \ref{4.5};
\item the 
paths of $\Omega _{m}(a)$ are in $V_{m}(a)$ up to homotopy, apart 
from second endpoints in $P_{m}(a)$, and with first endpoints at 
$h_{a}$, for $a=a_{0}$ and $\overline{a_{0}}$;
\item the paths of $\Omega 
_{m}(a_{1},+)$ are in $V_{m}(a_{1},+)$ up to homotopy, apart from 
second endpoints at $P_{m}(a_{1},+)$, and with first endpoints at 
$h_{a_{1}}$, and similarly for $\Omega _{m}(a_{1},-)$.
\end{itemize}
\end{maintheorem1.5}

\section{Checking vertex group elements}\label{5.8}

To know that the 
set $\Omega _{m}$ of \ref{5.7} gives a fundamental domain for 
$V_{3,m}$, by 
\ref{5.4}, we only need to know that the group vertex elements are 
simple parabolics. We now assume that $\Omega _{n}$ is a sequence of 
path sets, $0\leq n\leq m$, as in \ref{5.5}, with 
$\Omega 
_{n-1}\subset \Omega _{n}$. Then the vertex group elements of the 
fundamental domain associated to $\Omega _{n-1}$ are projections 
under 
the forgetful homomorphism of the vertex group elements of the 
fundamental domain associated to $\Omega _{n}$, by \ref{5.5}. If the 
projection under the forgetful homomorphism of a parabolic element 
$h$ is 
simple, then $h$ is itself simple. So to 
show that all group vertex elements are simple parabolics, we only need to 
show this for all vertices of $\Omega _{n}$ corresponding to points 
of $P_{3,n}\setminus P_{3,n-1}$, for each $n\leq m$. Now we give a 
criterion for this in terms of $\rho (\Omega _{n},s)$.

So let  $\beta _{i}$, $\beta _{i,2}$ and $\beta _{i,3}\in \rho 
(\Omega 
_{n},s)\subset \pi 
_{1}(\overline{\mathbb C}\setminus Z_{n}(s),v_{2},Z_{n}(s))$ end at  
points of $Z_{n}(s)\setminus Z_{n-1}(s)$ for $1\leq i\leq r$ and 
suppose that each $((\beta _{i},\beta _{i,2}),(\beta _{i+1},\beta 
_{i,3}))$ is the image of an adjacent  pair of matching adjacent 
pairs, with $\beta _{r+1}=\beta _{1}$. 
Equivalently for each $1\leq i\leq r$, there is 
$\gamma _{i}\in \pi _{1}(V_{3,n},s)$ such that the following holds 
for $[\psi 
_{i}^{-1}]=\Phi _{2}(\gamma _{i})$:
\begin{equation}\label{5.8.1}(\sigma _{\beta _{i+1}}\circ 
s,Y_{n}(s))\simeq _{\psi _{i}}(\sigma _{\beta 
_{i}}\circ s,Y_{n}(s)),\end{equation}
$$(\sigma _{\beta _{i,3}}\circ s,Y_{n}(s))\simeq _{\psi _{i}}(\sigma 
_{\beta 
_{i,2}}\circ s,Y_{n}(s)).$$
In particular, $\beta _{i}$ and $\psi _{i}(\beta _{i+1})$ share a 
common endpoint 
in $Z_{n}(s)\setminus Z_{n-1}(s)$. Write $\psi _{j,i}=\psi _{j}\circ 
\cdots \circ \psi _{i}$ for $1\leq j\leq i\leq r$ and $\psi 
_{j,j-1}={\rm{identity}}$. Then 
\begin{equation}\label{5.8.2}\Phi 
_{2}(\gamma _{j}*\cdots *\gamma _{i})=[\psi 
_{j,i}^{-1}].\end{equation} 

Let $h_{i}$ be 
the element of the covering group of $V_{3,n}$ corresponding to 
$\gamma _{i}$.  Then $\beta _{1}$ 
and $\psi _{1,i}(\beta _{i})$ share a common endpoint for $1\leq 
i<r$, 
$\psi _{1,r}(\beta _{r+1})=\psi _{1,r}(\beta _{1})=\beta _{1}$, 
and $h_{1}\cdots h_{r}$ is a vertex group element for the fundamental 
domain corresponding to $\Omega _{n}$ under the correspondence 
described in \ref{4.2}.
We claim that the following is sufficient for $h_{1}\cdots h_{r}$ to 
be a simple parabolic element, that is, for $\gamma _{1}*\cdots 
*\gamma 
_{r}$ to be a simple path round a puncture of $V_{3,n}$. 

\medskip

{\em{Simple parabolic criterion}}  
The cyclic order of  paths $\psi _{1,i}(\beta _{i+1})$ 
$0\leq i<r$ round the common endpoint in $Z_{n}(s)\setminus 
Z_{n-1}(s)$ respects the order of the indices.

\medskip

We see this as follows. Let $\zeta _{i}$ be the unique loop which is 
isotopic 
in $\overline{\mathbb C}\setminus Z_{n}(s)$  to an arbitrarily small 
perturbation 
of both 
$\beta 
_{i}*\psi _{i}(\overline{\beta _{i+1}})$ and $\beta 
_{i,2}*\psi _{i}(\overline{\beta _{i,3}})$. Thus $\zeta _{i}=\rho 
(\gamma _{i},s)$ and in fact we can see directly from (\ref{5.8.1}) 
that
\begin{equation}\label{5.8.3}\begin{array}{l}
    s\simeq_{\psi _{i}}\sigma _{\zeta _{i}}\circ s,\cr 
\zeta _{k}=\rho (\gamma _{k},s).\cr \end{array}\end{equation}
    Then for $1\leq j\leq i\leq r$,
    $$s\simeq _{\psi _{j,i}}\sigma _{\zeta _{j,i}}\circ s,$$
where
$$\zeta _{j,i}=\zeta _{j}*\psi _{j,j}(\zeta _{j+1})*\cdots *\psi 
_{j,i-1}(\zeta _{i})=\rho (\gamma _{j}*\cdots *\gamma _{i},s).$$
(See (\ref{4.3.3}), (\ref{4.3.5}), but this also follows from 
(\ref{5.8.2}).) Now $\zeta _{j,i}$ is an arbitrarily small 
perturbation of 
$\beta _{j}*\psi _{j,i}(\overline{\beta _{i+1}})$. So $\zeta _{1,r}$ 
is homotopic to a simple loop once round a point of $Z_{n}\setminus 
Z_{n-1}(s)$. We claim that $\gamma _{1}*\cdots *\gamma _{r}$ must 
also be a simple loop.  For we know that $h_{1}\cdots h_{r}$ 
is parabolic. So $\gamma _{1}*\cdots *\gamma 
_{r}=\delta ^{p}$ for some $p\geq 1$ and simple loop $\delta $ round 
a 
point of $V_{3,n}\setminus V_{3,n-1}$. So $\Phi _{2}(\delta )$ is 
isotopic to the identity via an isotopy preserving $Z_{n}(s)$, and 
by (\ref{4.3.5}) we have 
$$\zeta _{1,r}=\rho (\delta ^{p},s)=(\rho (\delta ,s))^{p}.$$
Since $\zeta _{1,r}$ is simple, we must have $p=1$.

In our cases,  the Simple parabolic
criterion is satisfied rather easily. Let $R$ 
one of the sets of (\ref{5.6.2}). Then in all cases, for a cycle 
of paths $\beta _{i}$ as above with endpoints in $Z_{n}(s)\setminus 
Z_{n-1}(s)$, we shall have $r=1$ or $2$. There is only one cyclic 
order on a set of one or two paths.

\chapter{Easy cases of the Main Theorem}\label{6}

\section{Restatement of the Main Theorem in the easy 
cases}\label{6.1}

We now give detail of the Main Theorem in the easy cases. We 
recall
that for any $q\in (0,1)\cap \mathbb Q$ which is in the boundary of a
gap $G$ of $L_{3/7}$, $\beta _{q}$ is a path from $v_{2}=\infty $ to a
point of $\cup _{m}Z_{m}(s_{3/7})$, which crosses $S^{1}$ at $e^{2\pi
iq}$ into $G$, and ends at the point of $\cup _{m}Z_{m}(s_{3/7})$ in
$G$.  We also use $\beta _{1/3}$ for a path in $\{ z:\vert z\vert \geq
1\} $ from $v_{2}$ to $e^{2\pi i(1/3)}$.  This path is defined up to
homotopy keeping the first endpoint and $Z_{m}(s_{3/7})$ fixed, and
the second endpoint on the loop $\ell _{1/3}\cup \ell _{1/3}^{-1}$,
where $\ell _{1/3}$ is the leaf of $L_{3/7}$ with endpoints $e^{\pm
2\pi i(1/3)}$

\begin{maintheorem2}A fundamental domain for $V_{3,m}$ can be 
constructed using a set $\Omega _{m}\subset 
\pi_{1}(V_{3,m},P_{3,m},a_{1})$ with
$$\Omega _{m}=\gamma _{0}'*\Omega _{m}(a_{0})\cup \gamma _{0}*\cup 
\Omega 
_{m}(\overline{a_{0}})\cup \Omega _{m}(a_{1},-)\cup \Omega 
_{m}(a_{1},+),$$
$$\begin{array}{ll}R_{m}(a_{0})=\rho (\Omega 
_{m}(a_{0}),s_{1/7}),\ &R_{m}(\overline{a_{0}})=\rho (\Omega 
_{m}(\overline{a_{0}}),s_{6/7})\cr 
R_{m}(a_{1},-)=\rho (\Omega _{m}(a_{1},-),s_{3/7}),\ &
R_{m}(a_{1},+)=\rho (\Omega _{m}(a_{1},+),s_{3/7}).\cr 
\end{array}$$
Here, $\gamma _{0}'$ and $\gamma _{0}$ are the  paths from $h_{a_{1}}$ 
to $h_{a_{0}}$ and $h_{\overline{a_{0}}}$ of \ref{5.6}. The 
paths of $\Omega _{m}(a)$ are in $V_{m}(a)$ up to homotopy, apart 
from second endpoints in $P_{m}(a)$, and with first endpoints at 
$h_{a}$, for $a=a_{0}$ and $\overline{a_{0}}$. The paths of $\Omega 
_{m}(a_{1},+)$ are in $V_{m}(a_{1},+)$ up to homotopy, apart from 
second endpoints at $P_{m}(a_{1},+)$, and with first endpoints at 
$h_{a_{1}}$, and similarly for $\Omega _{m}(a_{1},-)$. The structures 
of $\Omega _{m}(a_{0})$ and $\Omega _{m}(\overline{a_{0}})$ and $\Omega 
_{m}(a_{1},-)$ are as follows.
\begin{itemize}
   
\item[1.] Write $s=s_{1/7}$. Let $T_{m}(s)$ be as in \ref{5.3}. 
 Let $T_{m}'(s)$ be the subset of $T_{m}(s)$, obtained by deleting 
all edges and vertices lying entirely in the smaller region of the disc bounded by 
the 
leaf with endpoints at $e^{2\pi i(1/7)}$,  $e^{2\pi i(2/7)}$,  apart from the vertex at $v_{1}$. Then 
$R_{m}(a_{0})$ is the set of arcs in the complement of 
$T_{m}(s_{1/7})$, apart from the second endpoints, 
from $v_{2}$ to the vertices of $T_{m}'(s)$, with one path of 
$R_{m}(a_{0})$ approaching each vertex between each pair of adjacent 
edges ending at that vertex. The matching pairs of adjacent pairs 
$((\beta _{1},\beta _{2}),(\beta _{1}',\beta _{2}'))$ of 
$R_{m}(a_{0})$
are defined by the property that $\overline{\beta _{1}}*\beta _{2}$ 
and $\overline{\beta _{1}'}*\beta _{2}'$ bound a topological disc 
containing an edge of $T_{m}'(s)$, and $\beta _{1}$, $\beta 
_{1}'$ have a common endpoint, as do $\beta _{2}$ and $\beta _{2}'$. 
Suppose that $\beta _{i}=\rho (\zeta _{i})$ and $\beta _{i}'=\rho 
(\zeta _{i}')$, and let $\gamma \in \pi _{1}(V_{3,m},s)$ be the 
element with 
$$\gamma *\zeta _{i}=\zeta _{i}',\ i=1,\ 2.$$
Let $\alpha =\rho (\gamma )\in \pi _{1}(\overline{\mathbb C}\setminus 
Z_{m}(s),v_{2})$ and $\Phi _{2}(\gamma )=[\psi ^{-1}]\in  
{\rm{MG}}(\overline{\mathbb C},Y_{m+1}(s))$, that is, similarly to 
\ref{5.8}:
\begin{equation}\label{6.1.1}
    (s,Y_{m+1}(s))\simeq _{\psi }(\sigma _{\alpha }\circ 
    s,Y_{m+1}(s)),\end{equation}
    \begin{equation}\label{6.1.2}
\alpha *\psi (\beta _{i}')=\beta _{i}{\rm{\ in\ }}\pi 
_{1}(\overline{\mathbb C}\setminus 
Z_{m}(s),Z_{m}(s),v_{2}).\end{equation}
Then $\psi $ fixes the common 
endpoint of $\beta _{i}$ and $\beta _{i}'$ for $i=1$, $2$, and
\begin{equation}\label{6.1.3}
   \overline{\beta _{1}'}*\beta _{2}'=
    \psi (\overline{\beta _{1}'}*\beta _{2}'){\rm{\ rel\ }}Z_{m+1}(s).
\end{equation}
Exactly similar statements hold for $R_{m}(\overline{a_{0}})$, with 
${1\over 7}$, ${2\over 7}$ replaced by 
${6\over 7}$, ${5\over 7}$.

\item[2.] Exactly similar statements hold for $R
_{m}(a_{1},-)$. Here $T_{m}'(s_{3/7})$ is obtained from 
$T_{m}(s_{3/7})$ by deleting all edges and vertices 
lying entirely in the larger region of the disc bounded by the leaf 
with endpoints at $e^{2\pi i(1/7)}$,  $e^{2\pi i(6/7)}$,  apart from the vertex at $s_{3/7}(v_{1})$.
\end{itemize}
\end{maintheorem2}

\section{Derivation of the first version of the Main Theorem from the 
second}\label{6.2}

One claim of \ref{6.1} is that two paths in $V_{3,m}(a_{0})$ or 
$V_{3,m}(\overline{a_{0}})$ or $V_{3,m}(a_{1},-)$ end at the same 
point of $P_{3,m}(a_{0})$ or $P_{3,m}(\overline{a_{0}})$ or 
$P_{3,m}(a_{1},-)$ if and only if the images under $\rho 
_{2}(.,s_{p})$, for $p={1\over 7}$ or ${6\over 7}$ or ${3\over 7}$, have 
the same endpoint in $Z_{m}(s_{p})$. Another claim is that, for any 
path $\omega $ in $V_{3,m}$ with such an endpoint, $\rho 
_{2}(\omega )=\beta $ for $\beta $ with a single $S^{1}$-crossing and 
thus, as we have seen in \ref{2.8},
if $\omega $ ends at $a$, $h_{a}$ is Thurston equivalent to the 
capture $\sigma _{\beta }\circ s_{p}$. So all maps $h_{a}$  with
$$a\in P_{3,m}(a_{0})\cup P_{3,m}(\overline{a_{0}})\cup 
P_{3,m}(a_{1},-)$$
are Thurston equivalent to captures of the types  
claimed in \ref{2.10}, and the paths in $\Omega _{m}(a_{0})$, $\Omega 
_{m}(\overline{a_{0}})$ and $\Omega _{m}(a_{1},-)$ are completely 
determined by their images under $\rho _{2}$ in $R_{m}(a_{0})$, 
$R_{m}(\overline{a_{0}})$ and $R_{m}(a_{1},-)$. Note that by Tan 
Lei's result, summarised in \ref{2.6}, all these captures are indeed 
Thurston equivalent to 
rational maps, so that the sets  $R_{m}(a_{0})$, 
$R_{m}(\overline{a_{0}})$ and $R_{m}(a_{1},-)$ do lift to $D'(s_{p})$ 
for $p={1\over 7}$, ${6\over 7}$, ${3\over 7}$. 
Furthermore, the restatement of \ref{6.1} 
implies that we have a fundamental domain for $V_{3,m}$ such that the 
associated tree intersects $V_{3,m}(a_{0})$ in a tree which is 
naturally homeomorphic to $T_{m}'(s_{1/7})$, and similarly for 
$V_{3,m}(\overline{a_{0}})$ and $V_{3,m}(a_{1},-)$. So we have a 
one-to-one correspondence between $P_{m}(a_{0})$ and 
$T_{m}'(a_{0})\cap 
Z_{m}(s_{1/7})$, and similarly for $P_{m}(\overline{a_{0}})$ and 
$T_{m}'(\overline{a_{0}})\cap Z_{m}(s_{6/7})$, and for 
$P_{m}(a_{1},-)$ 
and $T_{m}'(a_{1})\cap Z_{m}(s_{3/7})$. So captures with the same 
endpoint in any of these regions are Thurston equivalent --- which we 
already saw directly in \ref{2.8} --- and there are no further 
Thurston equivalences. The map $a\mapsto \omega (a)$ in these 
regions is not completely canonically defined, but paths wth the same 
endpoint in $P_{3,m}(a_{0})\cup P_{3,m}(\overline{a_{0}})\cup 
P_{3,m}(a_{1},-)$ 
are mapped under $\rho (.,s_{p})$ to paths with 
the same endpoint,  we can choose $\omega (a)$ to be any path in 
$\Omega $ with endpoint at $a$, for all such $a$.

\section{Proofs in the easy cases}\label{6.3}

The completion of the proofs of Theorem \ref{6.1} in the case of 
$R_{m}(a_{0})$, 
$R_{m}(\overline{a_{0}})$ and $R_{m}(a_{1},-)$ are all very similar, so 
we concentrate on the case of  $R_{m}(a_{0})$. Write $s=s_{1/7}$. We 
only need to show 
that pairs of pairs of  adjacent arcs $((\beta _{1},\beta 
_{2}),(\beta 
_{1}',\beta _{2}'))$ of $R_{m}(a_{0})$ which bound a common edge $e$ 
of 
$T_{m}'(s)$ are indeed matched as claimed, that is, there is 
$\alpha \in \pi _{1}(\overline{\mathbb C}\setminus Z_{m}(s),v_{2})$ 
and 
$[\psi ]\in {\rm{MG}}(\overline{\mathbb C},Y_{m}(s))$ such that
\begin{equation} \label{6.3.1}\alpha *\psi (\beta _{i}')\simeq \beta 
_{i}{\rm{\ rel\ }}Y_{m}(s),\end{equation}
$$( s_{1/7},Y_{m}(s)),\simeq _{\psi 
}(\sigma _{\alpha }\circ s,Y_{m}(s)).$$
 This will imply that
 \begin{equation}\label{6.3.2} \sigma _{\beta _{i}'}\circ s\simeq 
_{\psi }\sigma _{\beta _{i}}\circ s.\end{equation}
We also need to show that $\psi $ fixes the second (common) endpoint 
of $\beta _{i}$ and $\beta _{i}'$. The Simple parabolic criterion of 
\ref{5.8} then follows from the fact that every vertex of $T_{m}'(s)$ 
in $Z_{m}(s)\setminus Z_{m-1}(s)$ is a meeting of at most two edges 
of $T_{m}'(s)$. In the first two cases, such vertices are always 
extreme, and therefore attached to only one edge. 

Assume without loss of generality that $\beta _{1}*\overline{\beta 
_{1}'}$ bounds an open  disc $D(\beta _{1})$ containing 
$\beta _{2}*\overline{\beta 
_{2}'}$, and let $r$ be the greatest integer such that $D(\beta _{1})$ 
does not intersect $Z_{r-1}(s)$. Then $r\geq 1$. In the case of 
$R_{m}(a_{0})$, $r$ is greater than the preperiod of the endpoint 
of $\beta _{1}$. In the case of $R_{m}(a_{1},-)$, this may not be 
true, but this does not matter.
Let $\alpha _{r}$ be the closed loop which is an arbitrarily 
small perturbation of $\partial D(\beta _{1})$, and in $D(\beta _{1})$, 
apart from having endpoint at $v_{2}$. Then $\alpha _{r}$ bounds a a 
disc 
$D(\alpha _{r})$,
whose intersection with $T_{m}(s)$ is a subtree 
of $T_{m}(s)$, with basepoint on $e$. Then we define 
$$\psi _{r}=\sigma _{\alpha _{r}},$$ 
and  
define $\alpha _{p}$, $\psi _{p}$ inductively  for $p\geq r$ by 
\begin{equation}\label{6.3.3}\psi _{p}\circ s= 
    \sigma _{\alpha _{p}}\circ s\circ \psi _{p+1},\end{equation}
and $\alpha _{p+1}$ is an arbitrarily small perturbation of $\beta 
_{1}*\psi _{p}(\overline{\beta _{1}'})$. Write 
\begin{equation}\label{6.3.4}\psi _{p+1}=\xi _{p}\circ \psi 
_{p}.\end{equation}
Then
$$\xi _{p}\circ \psi _{p}\circ s=\sigma _{\alpha _{p+1}}\circ s\circ 
\xi 
_{p+1}\circ \psi _{p+1},$$
giving 
$$\xi _{p}\circ \sigma _{\alpha _{p}}\circ s\circ \psi _{p+1}=
\sigma _{\alpha _{p+1}}\circ s\circ \xi 
_{p+1}\circ \psi _{p+1},$$
and
since  $\xi _{p}(\alpha _{p})$ is a perturbation of $\xi _{p}(\beta 
_{1})*\psi _{p+1}(\overline{\beta _{1}'})$, we obtain
\begin{equation}\label{6.3.5}\xi _{p}\circ \sigma _{\beta _{1}}\circ 
s=\sigma _{\beta _{1}}\circ 
s\circ \xi _{p+1}.\end{equation}
Then the support of $\xi 
_{p}$ is contained in $(\sigma _{\beta _{1}}\circ 
s)^{r-p-1}(D(\zeta ))$, and  does not intersect the edge $e$ for 
$m\geq p>r$, nor the adjacent vertices. So $\psi _{p}$ fixes the common 
endpoint of $\beta _{i}$ and $\beta _{i}'$ for both $i=1$ and $2$ and for $p\geq r$. 
The support 
is 
allowed to intersect $\beta _{1}'$ and $\beta _{2}'$ elsewhere, and 
almost certainly will.  Then for 
$p=m$, $\alpha _{m}=\alpha $, $\psi _{m}=\psi $, we 
have (\ref{6.3.2}) for $i=1$, $2$. So we have all 
the 
required properties. \Box

\chapter{The hard case: final statement and examples}\label{7}

\section{}\label{7.1}
 Throughout this section, we write
$$s=s_{3/7},$$ 
$$Y_{n}=Y_{n}(s)=Y_{n}(s_{3/7}),$$
$$Z_{n}=Z_{n}(s)=Z_{n}(s_{3/7}).$$
Also, $U^{p}$ is the subset of the unit disc defined in \ref{2.10}, and, as in \ref{3.3},
$$q_{p}=\frac{1}{3}-2^{-2p}\frac{1}{21}.$$

\section{Some conjugacy tracks}\label{7.2}

 We recall the notation $\psi _{m,q}$ of \ref{3.3} for the conjugacy, up to $Y_{m}$-preserving isotopy, between $\sigma _{q}\circ s$ and $\sigma _{1-q}\circ s$, remembering that $\psi 
_{m+1,q}=\xi _{m,q}\circ \psi _{m,q}$. In \ref{3.3}, we defined $\psi _{m,q}$ for $q=q_{k}$ and $0\leq m\leq 2k+2$.  We can extend the definition to all $m\geq 0$ by defining $\alpha _{m,q_{k}}$ and $\psi _{m+1,q_{k}}$ inductively for $m\geq 2k+2$ as follows. We define $\alpha _{m,q_{k}}$ to be an arbitrarily small perturbation of $\beta _{q_{k}}*\psi _{m,q_{k}}(\overline{\beta _{1-q_{k}}})$ which bounds a disc containing $v_{1}$ and disjoint from the endpoint of $\beta _{q_{k}}$. Note that $\psi _{2k+2,q_{k}}$ is the identity on $\beta _{1-q_{k}}$, and hence $\alpha _{2k+2,q_{k}}$ is an arbitrarily small perturbation of $\beta _{q_{k}}*\overline{\beta _{1-q_{k}}}$. Writing $q=q_{_k}$, we then define $\psi _{m+1,q}$ for $m\geq 2k+2$ by 
\begin{equation}\label{7.2.3}\sigma _{\alpha _{m,q}}\circ s\circ \psi _{m+1,q}=\psi _{m,q}\circ s,\end{equation}
and 
\begin{equation}\label{7.2.4}[\psi _{m+1,q}]=[\psi _{m,q}]{\rm{\ in\ }}{\rm{MG}}(\overline{\mathbb C},Y_{m}(s)).\end{equation}
As in \ref{3.3}, we then define $\xi _{m,q}$, for all $m$, and for $q=q_{k}$, by
$$\psi _{m+1,q}=\xi _{m,q}\circ \psi _{m,q}.$$
As in \ref{6.3} we then obtain that, for $m\geq 2k+1$,
\begin{equation}\label{7.2.5}\sigma _{\beta _{q}}\circ s\circ \xi _{m+1,q}=\xi _{m,q}\circ \sigma _{\beta _{q}}\circ s.\end{equation}
The support of $\xi _{m,q}$ 
is a union of annuli
\begin{equation}\label{7.2.1}A_{m,q}=(\sigma _{\beta_{q}}\circ s)^{2k+1-m}(A_{2k+1,q})\end{equation} 
for $m\geq 2k+1$, and 
\begin{equation}\label{7.2.2}C_{m,q}=(\sigma _{\beta _{q}}\circ s)^{-m}(C_{0,q})\end{equation}
for $m\geq 0$. In the case $q=\frac{2}{7}$ we still have this, but we also have
$$A_{m,2/7}=(\sigma _{2/7}\circ s)^{1-m}(A_{1,2/7}),$$
and 
$$C_{m,2/7}=(\sigma _{2/7}\circ s)^{1-m}(C_{1,2/7})$$
for $m\geq 1$, remembering that $s_{1/7}=\sigma _{\zeta _{2/7}}^{-1}\circ \sigma _{\beta _{2/7}}\circ s$, and that $\zeta _{2/7}\subset s^{-1}(\beta _{2/7})$, which is disjoint from $A_{m,2/7}$ for all $m\geq 1$.

The definition of $\xi _{m,q}$ from the sets $A_{m,q}$ and $C_{m,q}$ is given, up to isotopy preserving $Y_{n}$ for any $n$, in terms of {\em{beads}} on the sets $A_{m,q}$ and $C_{m,q}$. A bead on $A_{m,q}$ is a component of $A_{m,q}\cap \{ z:\vert z\vert \leq 1\} $. The definition of a bead on $C_{m,q}$ is slightly different, and in addition, there are two different types of bead on $C_{m,q}$ for $q=q_{k}$ and $k\geq 1$, but only one type of bead for $k=0$, that is, for $q_{0}=\frac{2}{7}$. The easiest way to define them seems to be to define the beads on $C_{0,q}$, which we shall do shortly. The beads on $C_{m,q}$ are then the preimages under $s^{m}$ (and also under $(\sigma _{\beta _{q}}\circ s)^{m}$) of the beads on $C_{0,q}$. The homeomorphism $\xi _{m,q}$ sends each bead on $A_{m,q}$ to the next bead on the same annulus component of $A_{m,q}$, in the anticlockwise direction, and sends each bead on $C_{m,q}$ to the next same type bead on the same annulus component, in the clockwise direction. We now describe the sets in some detail.   

First, we draw $A_{m,2/7}$ for $1\leq m\leq 5$. We draw $C_{0,2/7}$ on the same diagram as $A_{1,2/7}$. The intersection of $C_{0,2/7}$ with $\cup _{n}Y_{n}$ is contained in the (non-closed)  disc which is the smaller component of $\{ z:\vert z\vert \leq 1\} \setminus \ell $, where $\ell $ is the leaf of $L_{3/7}$ with endpoints at $e^{\pm 2\pi i(1/7)}$. We can choose a smaller disc within this disc, which is closed, is contained in $C_{0,2/7}$, contains all the points of $\cup _{n}Y_{n}$ within $C_{0,2/7}$, and does not intersect $\ell $, although it does contain the endpoints $e^{\pm 2\pi i(1/7)}$. This is the unique bead on $C_{0,2/7}$.  In the case $m=5$, we only 
draw two of the components of $A_{5,2/7}$, and these on separate diagrams. 
We also draw one component of $A_{6,2/7}$. 
\begin{center}
\includegraphics[width=7cm]{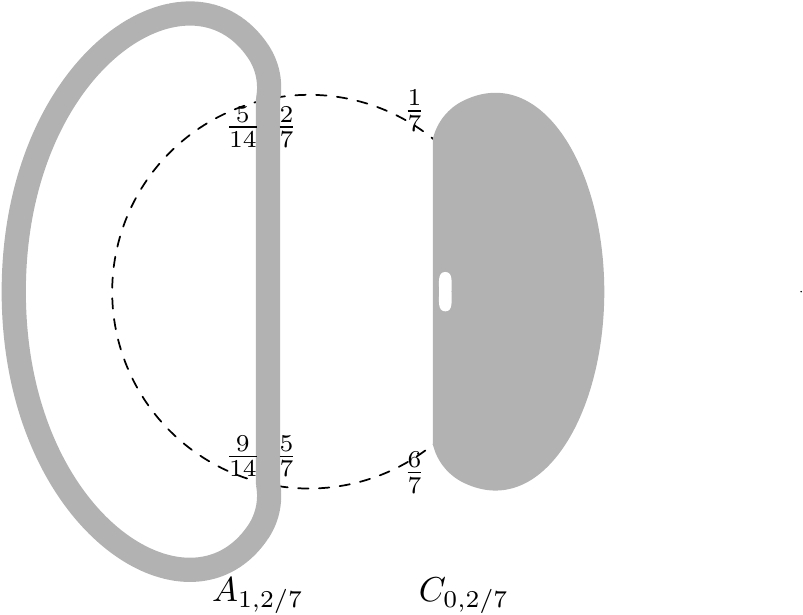}\includegraphics[width=3cm]{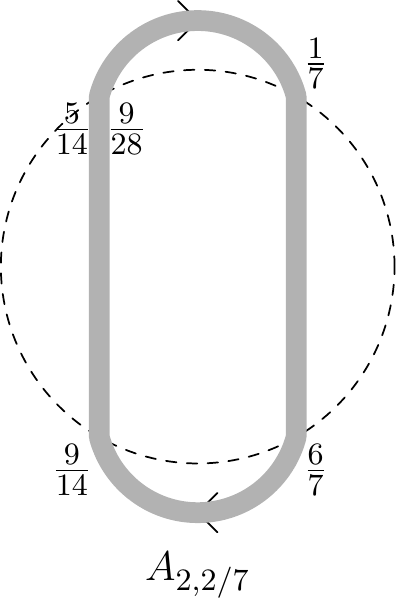}


\includegraphics[width=4cm]{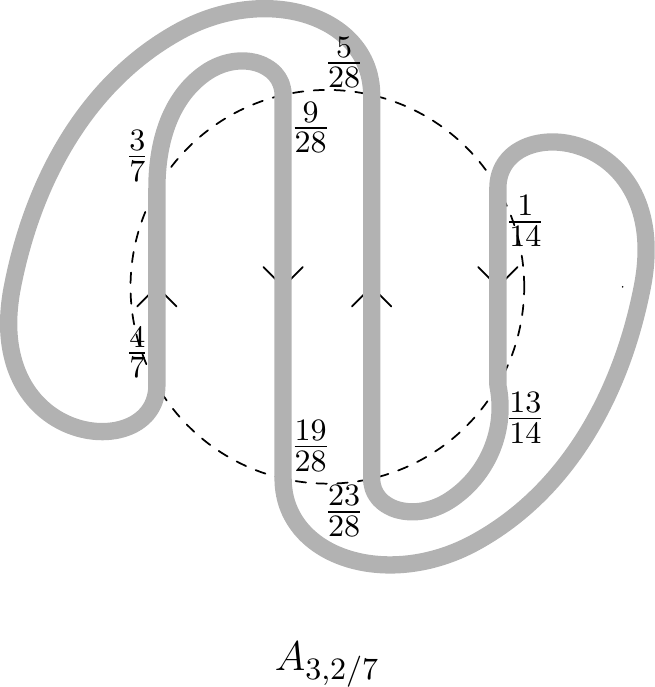}\ \ \ \ \ \includegraphics[width=4cm]{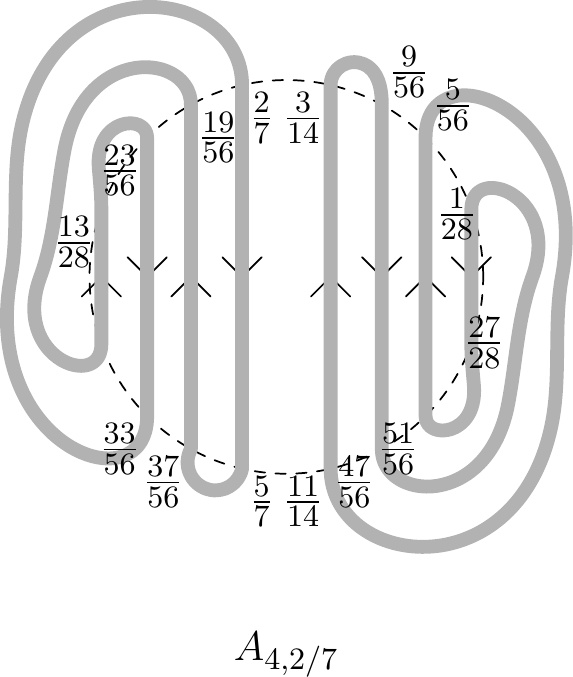}


\includegraphics[width=4cm]{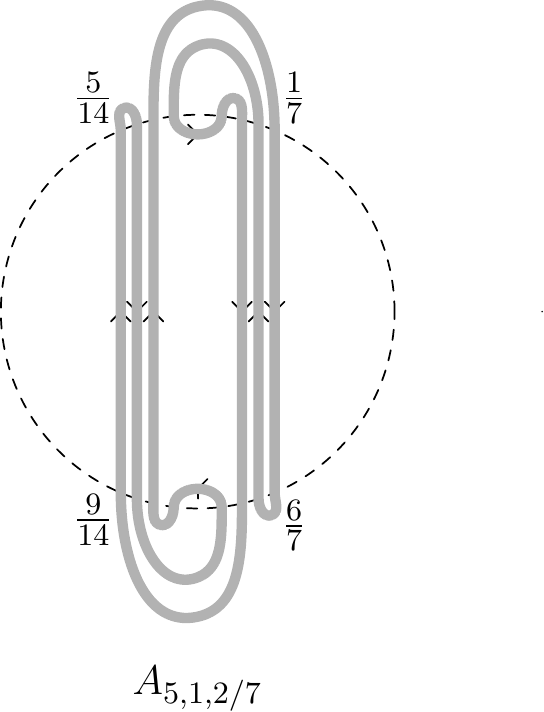}\ \ \ \includegraphics[width=4cm]{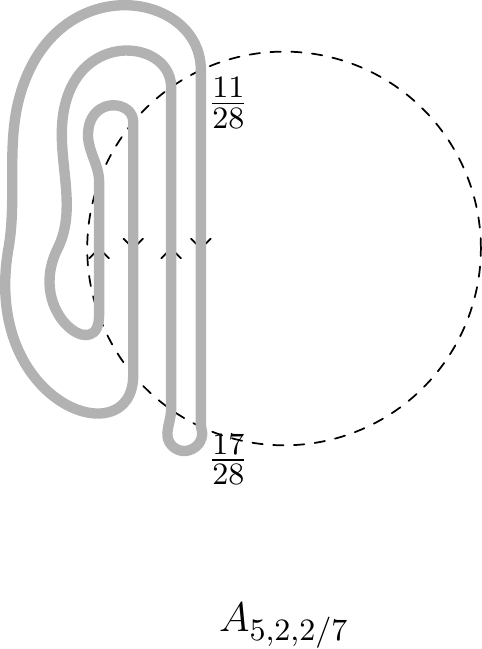}


\includegraphics[width=3cm]{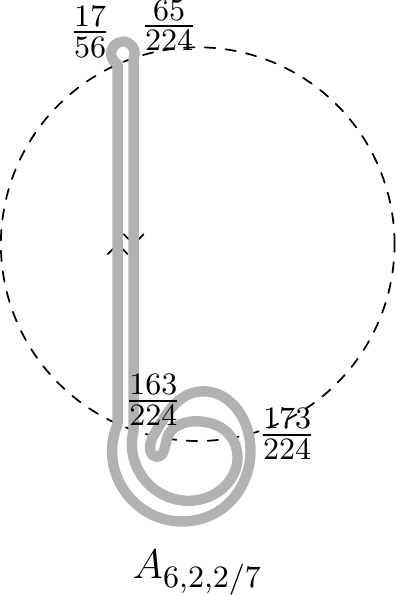}

 \end{center}
   
We see that $A_{m,2/7}$ has just one component for $m=2$ or $3$, 
and 
two for $m=4$. We call the lefthand component $A_{4,1,2/7}$, and the 
righthand component $A_{4,2,2/7}$. Then the map 
$\sigma _{2/7}\circ s:A_{5,1,2/7}\to A_{4,1,2/7}$ is degree two and  
$A_{5,1,2/7}$ is homotopic
to $A_{2,2/7}$ in $\overline{\mathbb C}\setminus Z_{0}(s)$.  We let $A_{5,2,2/7}$ and $A_{5,3,4/7}$ be the two 
components of the preimage of $A_{4,2,2/7}$, with $A_{5,2,2/7}$ to the 
left as drawn. On most of the diagrams, we have labelled 
intersection points of only one boundary of the annulus with the 
unit circle $S^{1}$, using the  convention of writing $x$ for 
$e^{2\pi ix}$.
 
For  $A_{5,1,2/7}$, $A_{5,2,2/7}$ and $A_{6,2,2/7}$, there is 
not space to include all the labels, even on one side of the annulus. 
For $A_{5,1,2/7}$, starting from the outer edge of the top leftmost 
intersection point with $S^{1}$, the points, proceeding in a 
clockwise direction round this edge, are: ${5\over 14}$ (as 
labelled), and then: 
$${37\over 112},\ {75\over 112},\ {43\over 56},\ {41\over 56},\ 
{79\over 112},
 \ {33\over 112},\ {1\over 7},{6\over 7},\ {93\over 112},\ {19\over 
112},\ {15\over 56},\ 
{13\over 56},
\ {23\over 112},\ {89\over 112},\ {9\over 14}.$$
The width of this annulus is ${1\over 224}$. For $A_{5,2,2/7}$, 
starting 
at the top rightmost intersection point 
with $S^{1}$ and proceeding in a clockwise direction the intersection 
points are ${11\over 28}$ (as labelled) and then:
$${17\over 28},\ {65\over 112},\ {47\over 112},\ {29\over 56},\ 
{27\over 56},\ {51\over 112},\ {61\over 112}.$$
The width of this annulus is also ${1\over 224}$. For $A_{6,2,2/7}$, 
starting at the top leftmost intersection with $S^{1}$ and proceeding 
in a clockwise direction on this edge, the points are ${17\over 56}$ 
(as labelled) and then:
$${65\over 224},\ {159\over 224},\ {85\over 112},\ {83\over 112},\ 
{163\over 224},\ {173\over 224}, {39\over 56}.$$
The width of this annulus is ${1\over 448}$.

The beads on $A_{m,2/7}$ can be described by their words: each bead is $D(w)$ for a word $w$, where $D(w)$, as in \ref{2.9}, is the set of points (topologically a closed disc) labelled by $w$.
For $A_{1,2/7}$ we have just one word:
$$L_{3}.$$
The homeomorphism $\xi _{1,2/7}$ maps this word to itself, rotating the annulus $A_{1,2/7}$ in an anticlockwise direction. For $A_{2,2/7}$ we have: 
$$L_{3}^{2}\to {\rm{ (top)}}R_{3}L_{3}\to {\rm{(bot)}}L_{3}^{2}.$$
Here, and subsequently, the arrows denote the direction of movement of beads in $A_{m,2/7}$ under $\xi _{m,2/7}$. For $A_{3,2/7}$, we have:
$$L_{2}R_{3}L_{3}\to  L_{3}^{3}\to R_{2}R_{3}L_{3}\to  
R_{3}L_{3}^{2}\to L_{2}R_{3}L_{3}.$$
For $A_{4,2/7}$ we have two cycles, in $A_{4,1,2/7}$ and $A_{4,2,2/7}$:
$$L_{1}R_{2}R_{3}L_{3}\to  L_{2}R_{3}L_{3}^{2}\to  
L_{3}L_{2}R_{3}L_{3}\to L_{3}^{4}\to L_{1}R_{2}R_{3}L_{3},$$ 
$$R_{1}R_{2}R_{3}L_{3}\to  R_{2}R_{3}L_{3}^{2}\to  
R_{3}L_{2}R_{3}L_{3}\to R_{3}L_{3}^{3}\to R_{1}R_{2}R_{3}L_{3}.$$
For $A_{5,2/7}$ there are three cycles, in $A_{5,i,2/7}$ for $i=1$, $2$ and $3$, with 8 beads in $A_{5,1,2/7}$, as was 
shown. 
The cycle for $A_{5,1,2/7}$ is:
$${\rm{(bot)}}BCL_{1}R_{2}R_{3}L_{3}\to  
{\rm{(bot)}}L_{3}L_{2}R_{3}L_{3}^{2}\to  
{\rm{(top)}}
R_{3}L_{3}L_{2}R_{3}L_{3}\to {\rm{(bot)}}R_{3}L_{3}^{4}\to $$
$${\rm{(top)}}UCL_{1}R_{2}R_{3}L_{3}\to 
{\rm{(top)}}R_{3}L_{2}R_{3}L_{3}^{2}\to  
{\rm{(bot)}}
L_{3}^{2}L_{2}R_{3}L_{3}\to {\rm{(top)}}L_{3}^{5}\to $$
$${\rm{(bot)}}BCL_{1}R_{2}R_{3}L_{3}.$$

Inductively, for $m\geq 6$, we define a component $A_{m+1,1,2/7}$ to be a component of $A_{m+1,2/7}$ which is a preimage under $\sigma _{2/7}\circ s$ of $A_{m,1,2/7}$, and homotopic to $A_{m-2,1,2/7}$ in $\overline{\mathbb C}\setminus Z_{m-4}(s)$. The connection $L_{3}^{2}\to {\rm{(top)}}R_{3}L_{3}$ in $A_{2,2/7}$
is essentially periodic of period $3$. It reappears in $A_{5,1,2/7}$ with: 
$$L_{3}L_{2}R_{3}L_{3}^{2}\to  
{\rm{(top)}}R_{3}L_{3}L_{2}R_{3}L_{3},$$
and similarly in $A_{3n+2,1,2/7}$ for all $n\geq 1$. All other 
preimages remain close under backward iterates. Even if the first 
letter is not the same, as in $L_{3}w_{1}\to BCw_{2}$, there is only 
way letters can remain different under taking preimages, for this 
one, for example, by taking $(L_{2}R_{3}L_{3})^{m}w_{1}\to 
(L_{1}R_{2}BC)^{m}w_{2}$. We call $A_{m,1,2/7}$ the {\em{periodic component}} of $A_{m,2/7}$

Now we draw some of the sets $C_{m,2/7}$.
\begin{center}
\includegraphics[width=5cm]{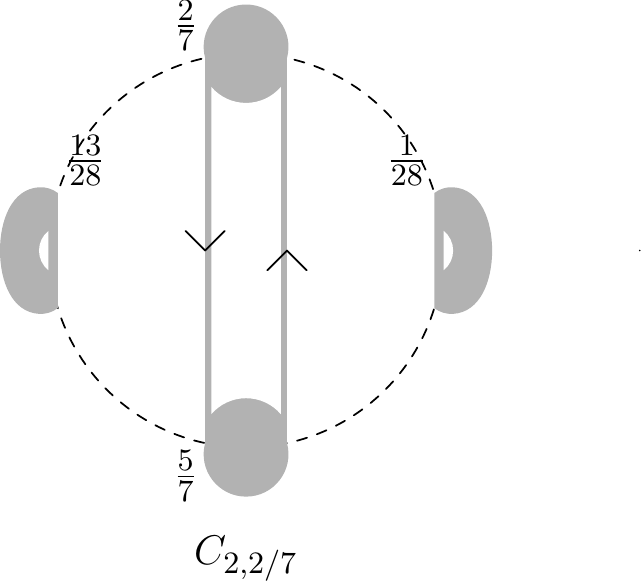}\ \ \ \includegraphics[width=4cm]{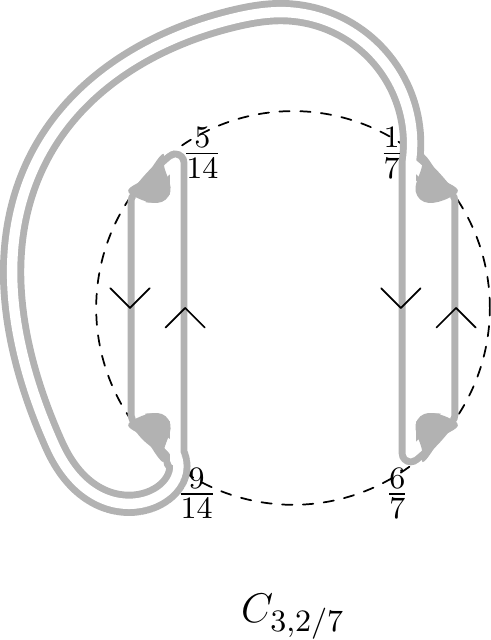}


\includegraphics[width=4cm]{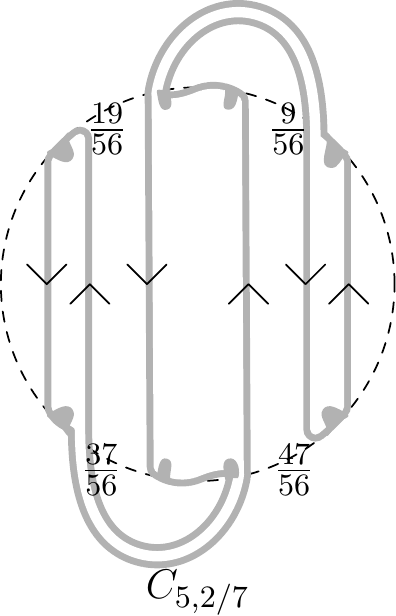}

\end{center}
Now, similarly to what we did for the sequence $A_{m,2/7}$, we consider the words for the beads on  
$C_{m,2/7}$ for $m\geq 2$. For all $m\geq 1$, the set $C_{m,2/7}$ has more 
than 
one component, but one bounds a disc containing a component of $Z_{0}(s)$. We call this the 
{\em{periodic}} component. The beads are again labelled by words. In fact, each bead is the union of the closure of the interior of $D(w)$ for a word $w$, for $D(w)$ as in \ref{2.9}. The word path for the periodic component 
of  $C_{2,2/7}$ is
$$UC\to BC\to UC.$$
This has one preimage in $C_{3,2/7}$, represented by:
$$R_{2}BC\to R_{2}UC\to L_{2}BC\to L_{2}UC\to R_{2}BC.$$
This has two preimages in $C_{4,2/7}$:
$${\rm{(top)}}L_{1}R_{2}BC\to {\rm{(bot)}}L_{1}R_{2}UC\to 
{\rm{(top)}}L_{3}L_{2}BC\to 
{\rm{(bot)}}L_{3}L_{2}UC\to $$
$${\rm{(top)}}L_{1}R_{2}BC,$$
$${\rm{(bot)}}R_{1}R_{2}BC\to {\rm{(top)}}R_{1}R_{2}UC\to 
{\rm{(bot)}}R_{3}L_{2}BC\to 
{\rm{(top)}}R_{3}L_{2}UC\to $$
$${\rm{(bot)}}R_{1}R_{2}BC.$$
The first of these is the periodic one, and its preimage in 
$C_{5,2/7}$ is represented by: 
$$BCL_{1}R_{2}BC\to BCL_{1}R_{2}UC\to 
{\rm{(top)}}L_{3}^{2}L_{2}BC\to {\rm{(bot)}}L_{3}^{2}L_{2}UC\to $$
$$UCL_{1}R_{2}BC\to UCL_{1}R_{2}UC\to 
{\rm{(bot)}}R_{3}L_{3}L_{2}BC\to 
{\rm{(top)}}R_{3}L_{3}L_{2}UC\to $$
$$BCL_{1}R_{2}BC.$$
As with the $A_{m,2/7}$ sequence, there is only one segment of track 
which is periodic, again of period $3$. This is the piece represented 
by
$${\rm{(top)}}R_{3}L_{3}L_{2}UC\to BCL_{1}R_{2}BC,$$
which has third preimage
$$BCL_{1}R_{2}R_{3}L_{3}L_{2}UC\to 
{\rm{(top)}}L_{3}^{2}L_{2}BCL_{1}R_{2}BC.$$

Now we make some remarks about the sets $A_{m,q_{k}}$ and $C_{m,q_{k}}$ for $k\geq 1$. Write $\sigma _{q_{k}}$ for $\sigma _{\beta }$, where $\beta =\beta _{q_{k}}$ The general shape of the sets $A_{m,q_{k}}$, for $m\leq 4$, is the same as for $A_{m,2/7}$, but the beads are thinner. So we have two components $A_{4,1,q_{k}}$ and $A_{4,2,q_{k}}$ of $A_{4,q_{k}}$, with $A_{4,1,q_{k}}$ on the left and $A_{4,2,q_{k}}$ on the right. The  annulus $A_{4,2,q_{k}}$ is homotopically trivial relative to the critical forward orbit of $\sigma _{q_{k}}\circ s$, and consequently the preimages under $(\sigma _{q_{k}}\circ s)^{m-4}$ of $A_{4,2,q_{k}}$ in $A_{m,q_{k}}$ are homotopic to corresponding components of $A_{m,2/7}$ in $Z_{m-3}(s)$. But the preimage under $\sigma _{q_{k}}\circ s$ in $A_{5,q_{k}}$ of $A_{4,1,q_{k}}$ has two components  for $k\geq 1$, and, therefore, subsequent preimages are different. Inductively, for $m\geq 5$, we define a component $A_{m+1,1,q_{k}}$ of $A_{m,q_{k}}$ which is a component of $(\sigma _{q_{k}}\circ s)^{-1}(A_{m,1,q_{k}})$ and homotopically nontrivial in $\overline{\mathbb C}\setminus Z_{0}(s)$. These properties determine $A_{m+1,q_{k}}$ uniquely. Then $A_{m+1,1,q_{k}}$ is a homeomorphic preimage of $A_{m,q_{k}}$ for $5\leq m<2k+4$, and a degree two preimage for $m=2k+4$. Moreover, $A_{2k+5,1,q_{k}}$ is homotopic to $A_{2,q_{k}}$ relative to the critical forward orbit of $\sigma _{q_{k}}\circ s$, and $A_{m+2k+3,1,q_{k}}$ and $A_{m,1,q_{k}}$ are similarly homotopic, for all $m\geq 3$, defining $A_{3,q_{k}}=A_{3,1,q_{k}}$. We call $A_{m,1,q_{k}}$ the {\em{periodic component}}  of $A_{m,q_{k}}$, just as we did in the case $k=0$.

For $C_{0,q_{k}}$ for $k\geq 1$, the annulus has a rectangle next to the leaf with endpoints $e^{\pm 2\pi i (1/7)}$. The other vertical side of the rectangle is the leaf with endpoints $e^{\pm 2\pi i(q_{k}/2)}$. This gives the second type of bead on $C_{0,q_{k}}$, referred to above. The preimage in $C_{3,q_{k}}$ of the periodic component of $C_{2,q_{k}}$ has two components for $k\geq 1$ and so again, of course, subsequent preimages are different. There is a component of $C_{2k+3,q_{k}}$ which is in the backward orbit of the periodic component of $C_{2,q_{k}}$ and which is a degree two preimage of a component of $C_{2k+2,q_{k}}$. This is a natural analogue of the periodic component of $C_{3,2/7}$.

\section{Another sequence of homeomorphisms}\label{7.3}

Let $\beta _{q_{p}}$, $\alpha _{m,q_{p}}$ and $\psi _{m,q_{p}}$ be as in \ref{3.3} and \ref{7.2}. 
We are now going to define new sequences $\alpha _{m,q_{p}}'$ and  $\psi _{m,q_{p}}'$,
for $p\geq 0$ and $m\geq 0$.  For $m=0$,
$$\alpha _{0,q_{p}}'=\alpha _{2/7},\ \ \psi _{0,q_{p}}'=\psi _{0,2/7}.$$
For $m\leq 2p+2$,
$$\psi _{m,q_{p+1}}'=\psi _{m,q_{p}}',$$
and $\beta _{q_{p}}'$ is defined by:
$$\alpha _{2p,q_{p}}*\psi _{2p,q_{p}}'(\beta _{q_{p}}')=\beta _{q_{p}}.$$
and, for $m\leq 2p+1$:
$$\alpha _{m,q_{p+1}}'=\alpha _{m,q_{p}}'.$$
But $\alpha _{2p+2,q_{p+1}}$ is such that 
$$\alpha
_{2p+2,q_{p+1}}'*\overline{\alpha _{2p+2,q_{p}}'}=\alpha
_{q_{p+1}}*\overline{\alpha _{q_{p}}}.$$

Then inductively for $m\geq 2p$, we define $\psi _{m+1,q_{p}}'$ in terms of $\psi _{m,q_{p}}'$, for $q=q_{p}$, by:
\begin{equation}\label{7.3.1}\sigma _{\alpha _{m,q}}\circ s\circ \psi _{m+1,q}'=\psi _{m,q}'\circ s,\end{equation}
\begin{equation}\label{7.3.2}\psi _{m+1,q}'=\psi _{m,q}'{\rm{\ rel\ }}Y_{m},\end{equation}
and for $m\geq 2p$, we also define $\alpha _{m+1,q_{p}}'$  in terms of $\beta _{q_{p}}'$ and $\psi _{m+1,q_{p}}'$ by: $\alpha _{m+1,q_{p}}'$ is an arbitrarily small perturbation of 
$$\beta _{q_{p}}*\psi _{m+1,q_{p}}'(\overline{\beta _{q_{p}}'}).$$
Finally, we define $\xi _{m,q_{p}}'$ by
$$\psi _{m+1,q_{p}}'=\xi _{m,q_{p}}'\circ \psi _{m,q_{p}}'.$$

As an example which we shall consider later, let $p=1$. Then:
$$\psi _{1,9/28}'=\psi _{1,2/7},\ \ \xi _{1,9/28}'=\xi _{1,2/7},$$
$$\psi _{2,9/28}'=\psi _{2,2/7}=\xi _{1,2/7}\circ \psi _{1,2/7}.$$
Recall from \ref{3.3} that $\psi _{1,2/7}$ is a composition of anticlockwise Dehn twists round the boundaries of two discs, which we call $D_{2/7}$, and $D_{-,2/7}$. Meanwhile, $\xi _{1,2/7}$ is a composition of anticlockwise Dehn twists round the boundaries of the two disc components of  $s^{-1}(D_{-,2/7})$ and a clockwise twist in the annulus $A_{1,2/7}$, that is, a clockwise Dehn twist round the outer boundary composed with an anticlockwise Dehn twist round the inner boundary. Then we see that $\beta _{9/28}'=\beta _{19/28}$,
and 
$$\psi _{2,9/28}'=\psi _{2,9/28}{\rm{\ rel\ }}Y_{3}.$$
In fact 
$$\psi _{2,9/28}'=\xi _{2,9/28}''\circ \psi _{2,9/28},$$
where $\xi _{2,9/28}''$ is a composition of clockwise twists in annuli --- in each annulus, the composition of clockwise twist round the outer boundary composed with with anticlockwise twist round the inner boundary --- where each annulus is trivial in $Y_{3}$. In fact there are four annuli involved. One annulus is  the subannulus of $A_{1,2/7}$ bounded by the component of $\partial A_{1,2/7}$ crossing $S^{1}$ at $e^{\pm 2\pi i(5/14)}$ and by an interior loop of $A_{1,2/7}$ whose intersection with the unit disc is the leaf between $e^{\pm 2\pi i(9/28)}$. The other annuli are the one which intersects the unit disc only between $e^{\pm 2\pi i(9/56)}$ and $e^{\pm 2\pi i(1/7)}$, and which bounds a  finite disc, and the two preimages of this one under $s$.
It follows that
$$\psi _{3,9/28}'=\psi _{3,9/28}{\rm{\ rel\ }}Y_{4},$$
and indeed for all $m\geq 2$,
$$\psi _{m,9/28}'=\psi _{m,9/28}{\rm{\ rel\ }}Y_{m+1}.$$
Moreover if we write
$$\psi _{m,9/28}'=\xi _{m,9/28}''\circ \psi _{m,9/28}$$
then the support of $\xi _{2,9/28}''$ is disjoint from $\overline{\beta _{19/28}}*\beta _{5/7}$ and the support of $\xi _{4,9/28}''$ is disjoint from $\overline{\beta _{75/112}}*\beta _{19/28}$. Extending this gives the following lemma.

\begin{ulemma}
\begin{equation}\label{7.3.3}\alpha _{m,q_{p}}'=\alpha _{m,q_{p}}=\alpha _{q_{p}}={\rm{\ in\ }}\pi _{1}(\overline{\mathbb C}\setminus Z_{m+1},v_{2}){\rm{\ for\ }}m\geq 2p\end{equation}
\begin{equation}\label{7.3.4}[\psi _{m,q_{p}}']=[\psi _{m,q_{p}}]{\rm{\ in\ MG}}(\overline{\mathbb C},Y_{m+1}){\rm{\ for\ }}m\geq 2p,\end{equation}
\begin{equation}\label{7.3.5}{\rm{supp}}(\psi _{2p+2,q_{p}}')\cap (\overline{\beta _{1-q_{p}}}*\beta _{1-q_{p+1}}\cup s^{-1}(\overline{\beta _{1-q_{p-1}}}*\beta _{1-q_{p}}))=\emptyset.
\end{equation}
\begin{equation}\label{7.3.6}\beta _{q_{p}}'=\beta _{1-q_{p}}{\rm{\ pointwise. }}\end{equation}

\end{ulemma}
\noindent {\em Proof.}
The proof is by  induction. These  statements are trivially true for $p=0$. Now suppose that they are true for $p$ and we prove them with $p+1$ replacing $p$.
First we prove (\ref{7.3.6}) for $p+1$. Since $\psi _{2p+2,q_{p+1}}'=\psi _{2p+2,q_{p}}'$ by definition, we have, by the definition of $\beta _{q_{p+1}}'$ and $\alpha _{2p+2,q_{p+1}}'$ and $\alpha _{2p+2,q_{p}}'$,
$$\alpha _{2p+2,q_{p+1}}'*\overline{\alpha _{2p+2,q_{p}}'}=\alpha _{2p+2,q_{p+1}}*\overline{\alpha _{2p+2,q_{p}}}.$$
$$\alpha _{2p+2,q_{p+1}}'={\rm{\ perturbation\ of\ }}\beta _{q_{p+1}}*\psi _{2p+2,q_{p}}'(\beta _{q_{p+1}}')$$
$$\alpha _{2p+2,q_{p}}'={\rm{\ perturbation\ of\ }}\beta _{q_{p}}*\psi _{2p+2,q_{p}}'(\beta _{q_{p}}').$$
So, from the definition of $\alpha _{2p+2,q_{p}}$ and $\alpha _{2p+2,q_{p+1}}$,
$$\beta _{q_{p+1}}*\psi _{2p+2,q_{p}}'(\overline{\beta _{q_{p+1}}'}*\beta _{q_{p}}')*\overline{\beta _{q_{p}}}=\beta _{q_{p+1}}*\psi _{2p+2,q_{p}}(\overline{\beta _{1-q_{p+1}}}*\beta _{1-q_{p}})*\overline{\beta _{q_{p}}}.$$

So this gives 
$$\psi _{2p+2,q_{p}}'(\overline{\beta _{q_{p+1}}'}*\beta _{q_{p}}')=\overline{\beta _{1-q_{p+1}}}*\beta _{1-q_{p}}.$$
Then (\ref{7.3.5}) gives $\beta _{q_{p+1}}'=\beta _{1-q_{p+1}}$, as required.
Now (\ref{7.3.3}) for $m=2p+2$ and $q_{p+1}$ follows from the definition of $\alpha _{2p+2,q_{p+1}}'$ and (\ref{7.3.3}) for $m=2p+2$ and $q_{p}$. Also, (\ref{7.3.4}) for $m=2p+2$ and $q_{p+1}$ follows from (\ref{3.3.8}): 

$$[\psi _{2p+2,q_{p+1}}]=[\psi _{2p+2,q_{p}}]{\rm{\ in\ MG}}(\overline{\mathbb C},Y_{2p+3}).$$

Now we prove (\ref{7.3.4}) and (\ref{7.3.3}) for $q_{p+1}$ and $m=2p+3$. Write $r=q_{p+1}$ and $q=q_{p}$. We have, since $\psi _{2p+2,r}'=\psi _{2p+2,q}'$:
$$\sigma _{\alpha _{2p+2,r}'}\circ s\circ \psi _{2p+3,r}'=\psi _{2p+2,q}'\circ s.$$ 
But from what we have  for $m=2p+2$, this gives 
$$\sigma _{\alpha _{2p+2,r}}\circ s\circ \psi _{2p+3,r}'=\psi _{2p+2,q}\circ s{\rm{\ rel\ }}Y_{2p+4}.$$
Since $\psi _{2p+3,r}'=\psi _{2p+2,r}'$ relative to $Y_{2p+3}$ and the same equation is solved by $\psi _{2p+3,r}$ relative to $Y_{2p+3}$, we obtain
$$\psi _{2p+3,q_{p+1}}'=\psi _{2p+3,q_{p+1}}{\rm{\ in\ }}{\rm{MG}}(\overline{\mathbb C},Y_{2p+4}).$$
Then we obtain, from the definition, 
$$\alpha _{2p+3,q_{p+1}}'=\alpha _{2p+3,q_{p+1}}{\rm{\ in\ }}\pi _{1}(\overline{\mathbb C}\setminus Z_{2p+4},v_{2}).$$
In exactly the same way, we obtain (\ref{7.3.3}) and (\ref{7.3.4}) for $q_{p+1}$ and $m\geq 2p+4$. 

Finally, we consider (\ref{7.3.5}) for $q_{p+1}$.  Write 
$$\psi _{m,q_{k}}'=\xi _{m,q_{k}}''\circ \psi _{m,q_{k}},$$
for $k=p$ and $p+1$. Write
$$X_{1}=\overline{\beta _{q_{p+1}}}*\beta _{q_{p}}\cup s^{-1}(\overline{\beta _{q_{p}}}*\beta _{q_{p-1}}),$$
$$X_{2}=\overline{\beta _{q_{p+2}}}*\beta _{q_{p+1}}\cup s^{-1}(\overline{\beta _{q_{p+1}}}*\beta _{q_{p}}).$$
Since the support of $\psi _{2p+4,q_{p+1}}$ is disjoint from $X_{2}$, it suffices to show that $\xi _{2p+4,q_{p+1}}''$ is disjoint from $X_{2}$. Define $\psi _{m,q_{p}}''$ for $2p+2\leq m\leq 2p+4$ by
$$\psi _{2p+2,q_{p}}''=\psi _{2p+2,q_{p}},$$
and for $m=2p+3$ and $m=2p+2$, $\psi _{m,q_{p}}''$ is isotopic to $\psi _{2p+2,q_{p}}$ relative to $Y_{2p+2}$, and for $r=q_{p+1}$,
$$\sigma _{\alpha _{m-1,r}}\circ s\circ \psi _{m,q_{p}}''=\psi _{m-1,q_{p}}''\circ s.$$
But for $k=2p+3$ or $2p+4$, we can write 
$$\xi _{k,q_{p+1}}''=\xi _{k,q_{p+1},1}''\circ \xi _{k,q_{p+1},2}'',$$
where 
$$\psi _{k,q_{p+1}}'=\xi _{k,q_{p+1},1}''\circ \psi _{k,q_{p}}''$$
and 
$$\psi _{k,q_{p}}''=\xi _{k,q_{p+1},2}''\circ \psi _{k,q_{p+1}}.$$
Now the support of $\xi _{k,q_{p+1},1}''$ is contained in the preimage under $(\sigma _{\beta _{r}}\circ s)^{k-2p-2}$ of the support of $\xi _{2p+2,q_{p}}''$. Since the support of $\xi _{2p+2,q_{p}}''$  is disjoint from $X_{1}$,
the support of $\xi _{2p+3,q_{p+1},1}''$ is disjoint from $\overline{\beta _{q_{p+1}}}*\beta _{q_{p}}$. The support of $\psi _{2p+3,q_{p+1}}$ is the preimage under $\sigma _{\beta _{r}}\circ s$ of an annulus which intersects the unit disc between the leaves with endpoints  $e^{\pm 2\pi iq_{p+1}}$ and $e^{\pm 4\pi iq_{p+1}}$. The preimage is disjoint from  $X_{2}$
So now the supports of $\xi _{2p+4,q_{p+1},1}''$ and $\xi _{2p+4,2}''$ are obtained by taking preimages under $\sigma _{\beta _{r}}\circ s$ again, and must be disjoint from $X_{2}$.

 This completes the proof of (\ref{7.3.5}) for $q_{p+1}$.
\Box

\section{Hard part of the fundamental domain: the first few cases}\label{7.5}
 
In this subsection, we describe the part of the fundamental domain corresponding to $V_{3,m}(a_{1},+)$ for $m\leq 5$, with some partial information in the cases $m=6$, $7$. To do this, we shall describe the set $\Omega _{m}(a_{1},+)$ in terms of its image $R_{m}(a_{1},+)$ under $\rho $. Three paths in $\Omega _{m}(a_{1},+)$ (for any $m$) have already been chosen: $\omega _{1}$, $\omega _{1}'$ and $\omega _{\infty }$. The images under $\rho (.,s)$ are, respectively, $\beta _{2/7}$,  $\beta _{5/7}$ and $\beta _{1/3}$, or equivalently, $\beta _{2/3}$, since $\beta _{1/3}$ and $\beta _{2/3}$ are homotopic under a homotopy moving the second endpoint along $\gamma _{1/3}$.  Here, $\gamma _{1/3}$ denotes the closed loop $\ell _{1/3}\cup \ell _{1/3}^{-1}$, where $\ell _{1/3}$ is the leaf of $L_{3/7}$ with endpoints $e^{\pm 2\pi i(1/3)}$. As in the proof of the easy cases in \ref{6.3}, we need to describe matched pairs of adjacent pairs in $R_{m}(a_{1},+)$. As in (\ref{6.3.1}), the adjacent pairs $(\beta _{1},\beta _{2})$ and $(\beta _{1}',\beta _{2}')$ in $R_{m}(a_{1},+)$ are matched if there is $[\psi ]\in {\rm{MG}}(\overline{\mathbb C},Y_{m}(s))$
 and $\alpha \in \pi _{1}(\overline{\mathbb C}\setminus Z_{m}(s),v_{2})$ such that
\begin{equation}\label{7.5.1}\begin{array}{l}(s,Y_{m}(s))\simeq _{\psi }(\sigma _{\alpha }\circ s,Y_{m}(s))\\
\beta _{i}=\alpha *\psi (\beta _{i}'){\rm{\ rel\ }}Y_{m}(s),\ i=1,2\end{array}\end{equation}

In particular, $\overline{\beta _{1}}*\beta _{2}$ and $\psi (\overline{\beta _{1}'}*\beta _{2}')$ are homotopic via a homotopy preserving  $Z_{m}$. We are going to try to make choices so that  $\overline{\beta _{1}}*\beta _{2}$ and $\psi (\overline{\beta _{1}'}*\beta _{2}')$ are disjoint arcs, after arbitrarily small perturbation near $v_{2}$, bounding an open disc disjoint from $m$.

If $(\beta _{2/7},\beta _{2})$ and $(\beta _{5/7},\beta _{2}')$ are adjacent pairs in $R_{m}(a_{1},+)$, then these are always matched by $[\psi ]=[\psi _{m,2/7}]$ and $\alpha _{m,2/7}$, where $\alpha _{m,2/7}$ and $\psi _{m,2/7}$ are as in \ref{7.2}. For $m\leq 2$, $\psi _{m,2/7}$ is the identity on $\beta _{5/7}$ and $\alpha _{m,2/7}$ is an arbitrarily small perturbation of $\beta _{2/7}*\psi _{m,2/7}(\beta _{5/7})$. 

We are now ready to start an inductive construction of $R_{m}(a_{1},+)$ with a matching of pairs of adjacent pairs. After homotopy preserving endpoints if necessary, $\beta _{2/7}$ and $\beta _{5/7}$ are disjoint from $\gamma _{1/3}$ and also from $\beta _{1/3}$, apart from the common endpoint at $v_{2}$. Also, after homotopy preserving endpoints if necessary, $\beta _{1/3}$ is disjoint from $\gamma _{1/3}$, apart from the second endpoint being in $\gamma _{1/3}$. Then $\beta _{2/7}\cup \beta _{5/7}\cup \beta _{1/3}\cup \gamma _{1/3}$ bounds an open topological disc  containing just one point of $Z_{2}(s)$, namely the common endpoint of $\beta _{9/28}$ and $\beta _{19/28}$. (The boundary of this disc is not an embedded circle.) Also, $\psi _{2,2/7}$ fixes both $\beta _{5/7}$ and $\beta _{2/7}$, and $\psi _{2,2/7}(\beta _{2/3})$ is homotopic to $\beta _{2/3}$ via a homotopy fixing endpoints and $Z_{2}$. Moreover, (\ref{7.5.1}) holds for $m=2$, and 
$$(\beta _{1},\beta _{2})=(\beta _{2/7},\beta _{9/28}),\ \ (\beta _{1}',\beta _{2}')=(\beta _{5/7},\beta _{19/28}),\ \ \alpha =\alpha _{2,2/7},\  \ \psi =\psi _{2,2/7},$$
and also for 
$$(\beta _{1},\beta _{2})=(\beta _{9/28},\beta _{1/3}),\ \ (\beta _{1}',\beta _{2}')=(\beta _{19/28},\beta _{1/3}),\ \alpha =\alpha _{2,9/28},\ \ \psi =\psi _{2,2/7}.$$
So then we can define 
$$R_{2}(a_{1},+)=\{ \beta _{2/7},\ \beta _{5/7},\ \beta _{9/28},\ \beta _{19/28},\ \beta _{1/3},\} ,$$
where the adjacent pair $(\beta _{2/7},\beta _{9/28})$ is matched with the adjacent pair $(\beta _{5/7},\beta _{19/28})$, and the adjacent pair $(\beta _{9/28},\beta _{1/3}))$ is matched with the adjacent pair $(\beta _{19/28},\beta _{1/3})$.  The set $U(\beta _{2/7},\beta _{9/28})$, bounded by $\overline{\beta _{2/7}}*\beta _{9/28}$ and 
$$\overline{\beta _{5/7}}*\beta _{19/28}=\psi _{2,2/7}(\overline{\beta _{5/7}}*\beta _{19/28}),$$
is shown below, up to homotopy preserving $Z_{2}$, with 
$\overline{\beta _{2/7}}*\beta _{9/28}$ (up to homotopy) indicated by solid line and $\overline{\beta _{5/7}}*\beta _{19/28}$ (up to homotopy) indicated  by dashed 
line. The unit circle and some lamination leaves are indicated by 
dotted lines. 
 \begin{center}
\includegraphics[width=4cm]{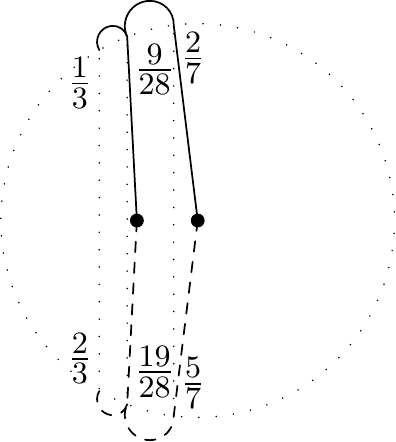}
\end{center}   
\medskip

Now we want to choose $R_{3}(a_{1},+)$, containing $R_{2}(a_{1},+)$ up to homotopy preserving $Y_{3}(s)$. We take the loop $\alpha $ as in (\ref{7.5.1}) to be an arbitrarily small neighbourhood of $\beta _{2/7}*\psi _{2,2/7}(\overline{\beta _{5/7}})$. Then the homeomorphism $\psi $ used as in (\ref{7.5.1}) for matching adjacent pairs between $\beta _{5/7}$ and $\beta _{19/28}$ with adjacent pairs between $\beta _{2/7}$ and $\beta _{9/28}$, is $\psi _{3,2/7}$, up to isotopy preserving $Y_{3}(s)$. Because of this, we consider the region bounded by $\overline{\beta _{9/28}}*\beta _{2/7}$ and $\psi _{3,2/7}(\overline{\beta _{19/28}}*\beta _{5/7})$, up to homotopy relative to $Z_{3}(s)$. So, since $\psi _{3,2/7}=\xi _{2,2/7}\circ \psi _{2,2/7}$, we need to consider the image of $\psi _{2,2/7}(\overline{\beta _{19/28}}*\beta _{5/7})$ under $\xi _{2,2/7}$. The support of $\xi _{2,2/7}$ is $A_{2,2/7}\cup C_{2,2/7}$. The annulus $A_{2,2/7}$ does not intersect  $\psi _{2,2/7}(\overline{\beta _{19/28}}*\beta _{5/7})$ up to homotopy preserving $Z_{3}(s)$, but $C_{2,2/7}$ does. So $\psi _{3,2/7}(\overline{\beta _{19/28}}*\beta _{5/7})$ has a bulge into the lower half of the unit disc between the two leaves of $L_{3/7}$ with endpoints at $e^{\pm 2\pi i(2/7)}$ and $e^{\pm 2\pi i(3/14))}$. This bulge contains a point of $Z_{3}(s)$ with symbolic code $w_{1}=BCL_{1}R_{1}R_{2}C$ (using the symbolic dynamics and conventions of \ref{2.9}). We therefore need to divide the region bounded by $\overline{\beta _{9/28}}*\beta _{2/7}$ and $\psi _{3,2/7}(\overline{\beta _{19/28}}*\beta _{5/7})$ into two, by defining two more paths in $R_{3}(a_{1},+)$, which we call $\beta (w_{1})$ and $\beta '(w_{1}')$, for $w_{1}'=UCL_{1}R_{1}R_{2}C$. It will turn out later that this notation is valid: no other path in $\cup _{p}R_{m,p}$ will have endpoint at the point of $Z_{3}$ labelled by $w_{1}$. Now $R_{3}(a_{1},+)$ is simply the union of $R_{2}(a_{1},+)$ and $\beta (w_{1})$ and $\beta '(w_{1}')$. Below is the sketch of the two regions $U(\beta _{2/7},\beta (w_{1}))$ and $U(\beta (w_{1}),\beta _{9/28})$, one bounded by  $\overline{\beta (w_{1})}*\beta _{2/7}$ and $\psi _{3,2/7}(\overline{\beta '(w_{1}')}*\beta _{5/7})$ and the other bounded by $\overline{\beta _{9/28}}*\beta (w_{1})$ and $\psi _{3,2/7}(\overline{\beta _{19/28}}*\beta '(w_{1}'))$.
\begin{center}
\includegraphics[width=5cm]{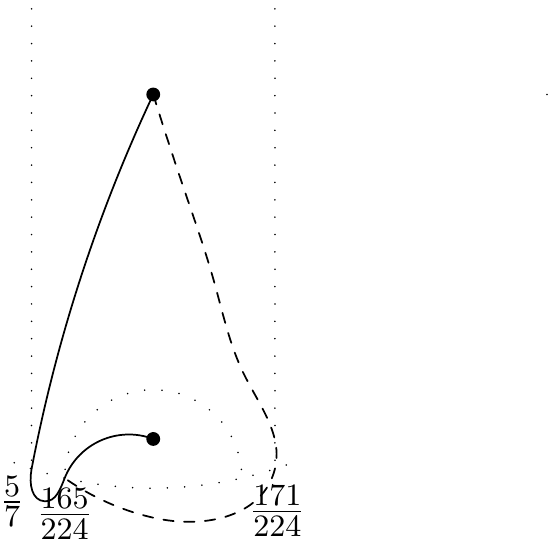}\ \ \ \includegraphics[width=4cm]{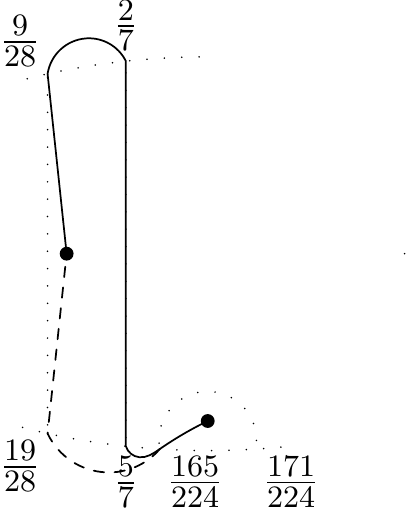}

\medskip

$U(\beta _{2/7},\beta (w_1))$ and $U(\beta (w_1),\beta _{9/28})$\end{center}
 
 \medskip
 
The critically finite map represented 
by $\beta (w_{1})$  is not Thurston-equivalent to any 
capture, 
as is demonstrated by the comparison between the capture numbers and 
expected value of $\# (P_{3}'(1))$, the deficit in the 
third number in the sequence, noted after (\ref{3.4.5}).
 The point $e^{2\pi i(165/224)}$ 
is not a point of lowest preperiod  under $z\mapsto z^{2}$ in the 
boundary of the gap containing the endpoint of $\beta 
(w_{1})$. The points of lowest preperiod are actually 
$e^{2\pi i(41/56)}$ and $e^{2\pi i(43/56)}$. The rule will be 
specified later, but the last $S^{1}$ crossing points of paths of $R_{m}(a_{1},+)$, for any $m$, are always in the backward orbit of $e^{2\pi i(2/7)}$. 
For the matching pair $(\beta _{9/28},\beta _{1/3})$ and $(\beta _{19/28},\beta _{1/3})$, we use the homeomorphism $\psi _{3,9/28}'$ of \ref{7.3} to effect the matching. Thus, 
$$\sigma _{\alpha _{2,9/28}}\circ s\circ \psi _{3,9/28}'=\psi _{2,2/7}\circ s.$$
Then $\psi _{3,9/28}'(\overline{\beta _{1/3}}*\beta _{19/28})$ and  $\psi _{3,9/28}(\overline{\beta _{1/3}}*\beta _{19/28})$ are homotopic up to homotopy preserving $Z_{3}(s)$, and the region bounded by $\overline{\beta _{1/3}}*\beta _{9/28}$ and $\psi _{3,9/28}'(\overline{\beta _{1/3}}*\beta _{19/28})$ contains no points of$Z_{3}(s)$ in its interior. So
we do not have to subdivide the region bounded by $\overline{\beta _{1/3}}*\beta _{9/28}$ and  $\psi _{3,9/28}'(\overline{\beta _{1/3}}*\beta _{19/28})$. The paths $\beta (w_{1})$, and  $\beta '(w_{1}')$ are indeed the only ones we need to add to $R_{2}(a_{1},+)$ to produce $R_{3}(a_{1},+)$. 

So far, $R_{m}(a_{1},+)\setminus \{ \beta _{1/3}\} $ is a union of two sets, such that each path in the first set is matched with a path in the second, and the matching of adjacent pairs is effected by this matching of individual paths. For example, in $R_{3}(a_{1},+)$, the path $\beta _{2/7}$ is matched with $\beta _{5/7}$, and $\beta (w_{1})$ with $\beta '(w_{1}')$, and $\beta _{9/28}$ with $\beta _{19/28}$. This matching will continue. We continue trying to subdivide the regions bounded by $\overline{\beta _{2}}*\beta _{1}$ and $\psi _{m+1}(\overline{\beta _{1}'}*\beta _{2}')$ in $R_{m}(a_{1},+)$ for matching adjacent pairs $(\beta _{1},\beta _{2})$ and $(\beta _{1}',\beta _{2}')$, where $\psi _{m}$ and $\alpha _{m}$ satisfy (\ref{7.5.1}), with $\psi _{m}$ and $\alpha _{m}$ replacing $\psi $ and $\alpha $, and $\psi _{m+1}$ is defined by $[\psi _{m+1}]=[\psi _{m}]$ in ${\rm{MG}}(\overline{\mathbb C},Y_{m})$ and 

\begin{equation}\label{7.5.2}\sigma _{\alpha _{m}}\circ s\circ \psi _{m+1}=\psi _{m}\circ s.\end{equation}
At each stage there is a choice to be made for the loop $\alpha _{m}$ for an adjacent pair $(\beta _{1},\beta _{2})$, where $\beta _{1}$ is nearer to $\beta _{2/7}$ than $\beta _{2}$. We shall always take $\alpha _{m}$ to be an arbitrarily small perturbation of $\beta _{1}*\psi _{m}(\overline{\beta _{1}'})$, so that $\alpha _{m}$ is also  homotopic, up to homotopy preserving $Y_{m}$, to an arbitrarily small perturbation of $\beta _{2}*\psi _{m}(\overline{\beta _{2}'})$. We write $U(\beta _{1},\beta _{2})$ for the region bounded by arbitrarily small pertubations of $\beta _{1}*\overline{\beta _{2}}$ and $\psi _{m}(\beta _{1}'*\overline{\beta _{2}'})$ up to $Z_{m}$-preserving homotopy, where the perturbations are chosen so that $U(\beta _{1},\beta _{2})$ does not contain $v_{2}$ and can be homotoped into $\{ z:\vert z\vert \leq 1\} $ in the complement of $v_{2}$.
We denote by $\partial 'U(\beta _{1},\beta _{2})$ the part of the boundary of $\partial U(\beta _{1},\beta _{2})$ which is homotopic to an arbitrarily small perturbation of $\psi _{m}(\overline{\beta _{1}'}*\beta _{2}')$. 

So, following this procedure, we next insert paths of $R_{4}(a_{1},+)\setminus R_{3}(a_{1},+)$ between adjacent paths in $R_{3}(a_{1},+)$, where necessary. The sets $U(\beta _{1},\beta _{2})$ for adjacent pairs $(\beta _{1},\beta _{2})$ between $\beta _{2/7}$ and $\beta _{9/28}$ are shown in the pictures below. 
\begin{center}
\includegraphics[width=5cm]{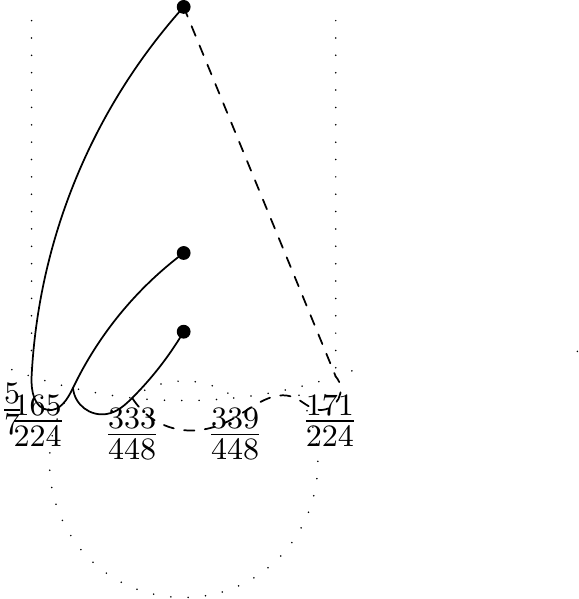}\ \ \ \includegraphics[width=3.5cm]{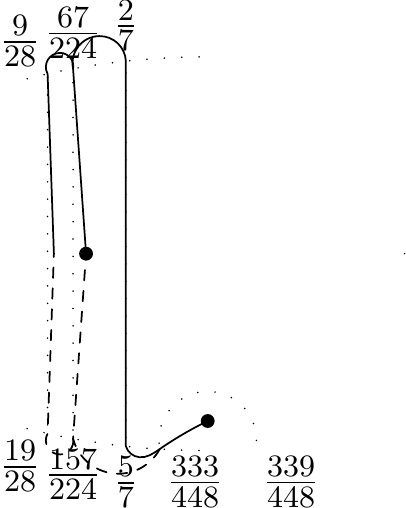}
\end{center}
  \medskip
  
We label each path by the word representing its second
endpoint, using the convention of \ref{2.9}. There are two new paths $\beta (w_{2})$ and $\beta (w_{3})$ between $\beta (w_{1})$ and $\beta _{9/28}$, with $w_{2}=BCL_{1}R_{1}^{2}R_{2}C$ and $w_{3}=L_{3}L_{2}R_{3}L_{2}C$.   The path  $\beta (w_{3})$ is actually $\beta _{67/224}$. The point of lowest preperiod on the top  boundary of the gap is $e^{2\pi i(33/112)}$, but we want to keep to the rule that $S^{1}$-crossing points are in the backward orbit of $e^{2\pi i(2/7)}$.
There are corresponding paths $\beta '(w_{2}')$ and $\beta '(w_{3}')$ between $\beta '(w_{1}')$ and $\beta _{19/28}$ with $w_{2}'=UCL_{1}R_{1}^{2}R_{2}C$ and $w_{3}'=w_{3}$, such that $(\beta (w_{1}),\beta (w_{2}))$ is matched with $(\beta '(w_{1}'),\beta '(w_{2}'))$, and $(\beta (w_{2}),\beta (w_{3}))$ is matched with $(\beta '(w_{2}'),\beta (w_{3}'))$.  The  path $\beta '(w_{3}')$ is $\beta _{157/224}$. We have $w_{3}'=w_{3}$ because $\psi _{4,2/7}$ fixes $D(w_{3})$. Although $\beta (w_{2})$ is between $\beta (w_{1})$ and $\beta _{9/28}$, we have 
$$U(\beta (w_{1}),\beta (w_{2}))\subset U(\beta _{2/7},\beta (w_{1})).$$

There are two remaining paths in $R_{4}(a_{1},+)\setminus R_{3}(a_{1},+)$, which are $\beta _{37/112}$ and $\beta _{75/112}$. The adjacent pair $(\beta _{9/28},\beta _{37/112})$ is matched with the adjacent pair $(\beta _{19/28},\beta _{75/112})$, and the adjacent pair $(\beta _{37/112},\beta _{1/3})$ is matched with the adjacent pair $(\beta _{75/112},\beta _{1/3})$. The dashed boundary $\partial 'U(\beta _{9/28},\beta _{37/112})$ is simply an arbitrarily small perturbation of $\overline{\beta _{19/28}}*\beta _{75/112}$, up to homotopy preserving $Z_{4}$.

So far, the matching of adjacent pairs in $R_{m}(a_{1},+)$ has been induced by a matching of each path in $R_{m}(a_{1},+)\setminus \{ \beta _{1/3}\} $ with exactly one other path in $R_{m}(a_{1},+)\setminus \{ \beta _{1/3}\} $. This means that, so far, the tree in $V_{3,m}(a_{1},+)\cup P_{3,m}(a_{1},+)$ which is dual to $\cup \Omega _{m}(a_{1},+)$ is simply an interval made up of edges between vertices in $P_{m}(a_{1},+)$. It is reasonable to attempt to continue this pattern, and we shall do so. So we aim to define $R_{m}(a_{1},+)$ so that
$$R_{m}(a_{1},+)=\cup _{0\leq p,\ 2p\leq m}R_{m,p}\cup R_{m,p}'\cup \{ \beta _{1/3}\} ,$$
where $R_{m,p}$ is the set of paths between $\beta _{q_{p}}$ and $\beta _{q_{p+1}}$, and 
\begin{equation}\label{7.5.3}R_{m,p}\cap R_{m,q}'=\emptyset \end{equation}
for all $p$ and $q$, and if $2p+2\leq m$, then
\begin{equation}\label{7.5.4}R_{m,p}\cap R_{m,p+1}=\{ \beta _{q_{p+1}}\} ,\end{equation}
\begin{equation}\label{7.5.5}R_{m,p}'\cap R_{m,p+1}'=\{ \beta _{1-q_{p+1}}\} {\rm{\ rel\ }}Y_{2p+2},\end{equation}
\begin{equation}\label{7.5.6}\beta _{2/7}\in R_{m,0},\ \ \beta _{5/7}\in R_{m,0}'.\end{equation}
We shall also aim to have
\begin{equation}\label{7.5.7}R_{m}(a_{1},+)\subset R_{m+1}(a_{1},+){\rm{\ rel\ }}Y_{n}{\rm{\ for\ all\ }}n\geq 0.\end{equation}

 Now we consider $R_{5,0}$. For any pair $(\beta _{1},\beta _{2})$ in $R_{4,0}$, it can be shown that the corresponding homeomorphisms $\psi _{4}$ and $\psi _{5}$ coincide with $\psi _{4,2/7}$ and $\psi _{5,2/7}$ on $\overline{\beta _{1}'}*\beta _{2}'$, for the matching pair $(\beta _{1}',\beta _{2}')$. (We shall prove such results later.) The definition of $\psi _{m+1}$ in terms of $\alpha _{m}$ and $\psi _{m}$ is given in (\ref{7.5.2}) with $\psi _{0}=\psi _{0,2/7}$. The definition of $\alpha _{m}$ in terms of $\psi _{m}$, $\beta _{1}$ and $\beta _{1}'$, for an adjacent pair $(\beta _{1},\beta _{2})$ in $\cup _{p}R_{m,,p}$, matched with an adjacent pair $(\beta _{1}',\beta _{2}')$ in $\cup _{p}R_{m,p}'$, is given immediately after (\ref{7.5.2}). Then, if we write $\psi _{5}=\xi _{4}\circ \psi _{4}$, we have 
 $$\xi _{4}=\xi _{4,2/7}{\rm{\ on\ }}\partial 'U(\beta _{1},\beta _{2}).$$
 The support of $\xi _{4,2/7}$ is the union of $A_{4,2/7}$ and $C_{4,2/7}$, each of these being a union of disjoint annuli. One annulus from $A_{4,2/7}$ intersects $\partial 'U(\beta (w_{2}),\beta (w_{3}))$, and one annulus from $C_{4,2/7}$ intersects $\partial 'U(\beta _(w_{3})\beta _{9/28})$. Otherwise, the support has no intersection with the sets $\partial 'U(\beta _{1},\beta _{2})$ for adjacent pairs $(\beta _{1},\beta _{2})$ in $R_{4,0}$. The part of  $A_{4,2/7}$ which intersects $\partial 'U(\beta (w_{2}),\beta (w_{3}))$ is
$$L_{3}L_{2}R_{3}L_{3}\to L_{3}^{4},$$
and the part of $C_{4,2/7}$ which intersects $\partial 'U(\beta (w_{3}),\beta _{9/28})$ is
$$L_{3}L_{2}UC\to L_{3}L_{2}BC.$$
We now draw the image under $\xi _{4,2/7}$ of the dashed boundary of $U=U(\beta 
(w_{2}),\beta (w_{3}))$, which uses the intersection with $A_{4,2/7}$.
\begin{center}
\includegraphics[width=8cm]{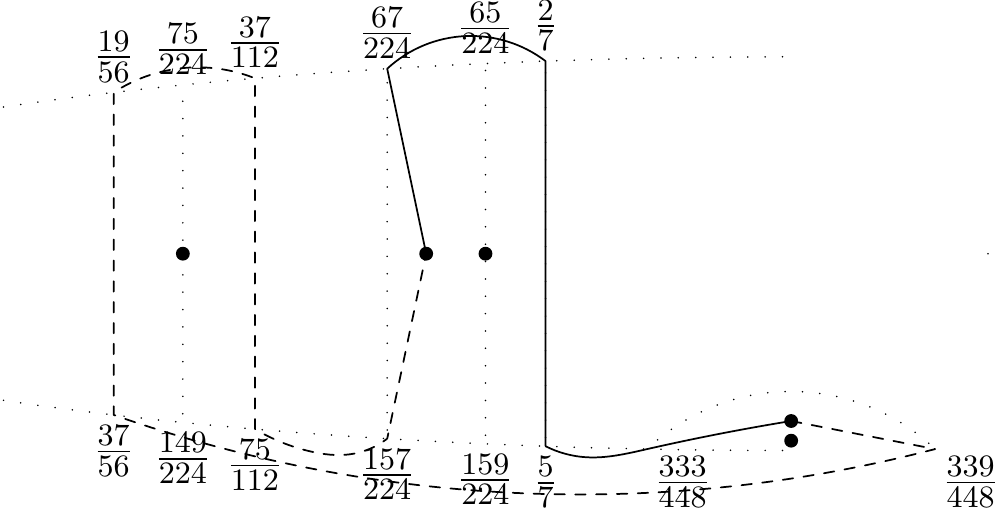}

\medskip

   Image under $\xi _{4,2/7}$
\end{center} 
\medskip
  
The region bounded by $\overline{\beta (w_{2})}*\beta (w_{3}))$ and $\xi _{4,2/7}(\partial 'U)$ contains three points of $Z_{5}\setminus Z_{4}$, with words 
$$w_{5}=BCL_{1}R_{1}^{3}R_{2}C,\ \ w_{6}=L_{3}^{4}L_{2}C,\ \ w_{7}=L_{3}L_{2}R_{3}L_{3}L_{2}C.$$
 The region bounded by $\overline{\beta _{2/7}}*\beta (w_{2})$ and $\xi _{4,2/7}(\partial 'U(\beta _{2/7},\beta (w_{2})))$ also contains one point of $Z_{5}\setminus Z_{4}$, with word $w_{4}=BCL_{1}R_{2}R_{3}L_{2}C$. 
The region bounded by $\overline{\beta (w_{3})}*\beta _{9/28}$ and $\xi _{4,2/7}(\partial 'U(\beta (w_{3}),\beta _{9/28}))$ contains one point of $Z_{5}\setminus Z_{4}$ with the word $w_{8}=L_{3}L_{2}UCL_{1}R_{1}R_{2}C$.
Then
$$R_{5,0}\setminus R_{4,0}=\{ \beta (w_{4}),\beta (w_{5}),\beta (w_{6}),\beta (w_{7}),\beta (w_{8})\} .$$
In all these cases, $\beta (w_{i})$ is the only path in $R_{5,0}$ (and also, it turns out, the only such path in $\cup _{p}R_{5,p}$) with endpoint coded by $w_{i}$. Correspondingly,
$$R_{5,0}'\setminus R_{4,0}'=\{ \beta '(w_{4}'),\beta '(w_{5}'),\beta '(w_{6}'),\beta '(w_{7}'),\beta (w_{8}')\} ,$$
where:
$$\begin{array}{l}w_{4}'=UCL_{1}R_{2}R_{3}L_{2}C,\ \ w_{5}'=UCL_{1}R_{1}^{3}R_{2}C,\ \ w_{6}'=w_{7}=L_{3}L_{2}R_{3}L_{3}L_{2}C,\\ w_{7}'=L_{2}R_{3}L_{3}^{2}L_{2}C,\ \  w_{8}'=L_{3}L_{2}BCL_{1}R_{1}R_{2}C.\end{array}$$

Now we draw $U(\beta _{1},\beta _{2})$ 
for adjacent pairs $(\beta _{1},\beta _{2})$ from $R_{5,0}$.
\begin{center}
\includegraphics[width=9cm]{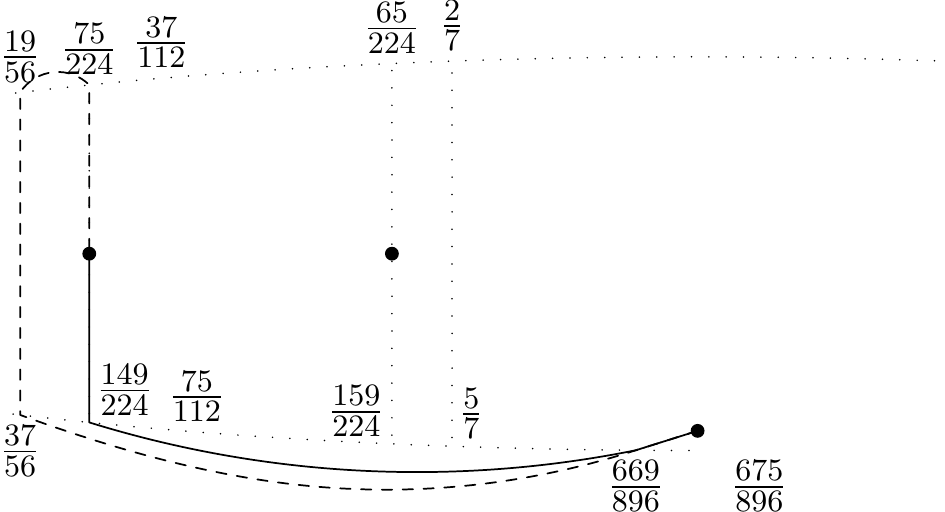}

\medskip

$U(\beta (w_5),\beta (w_6))$

\medskip
\includegraphics[width=7cm]{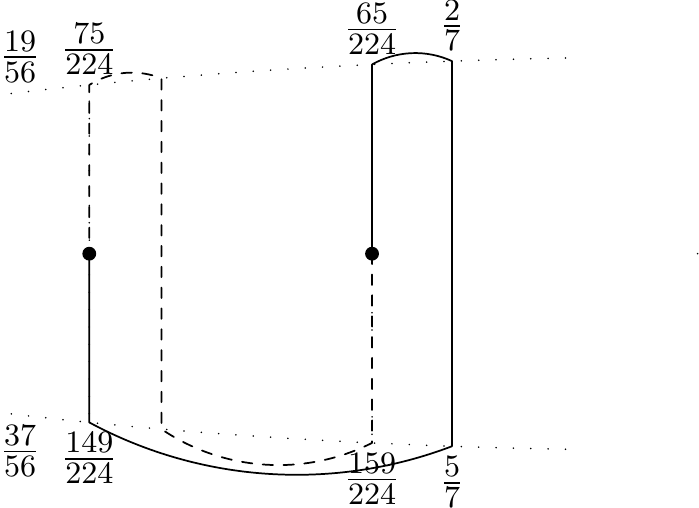}\ \ \ \ \ \includegraphics[width=2.5cm]{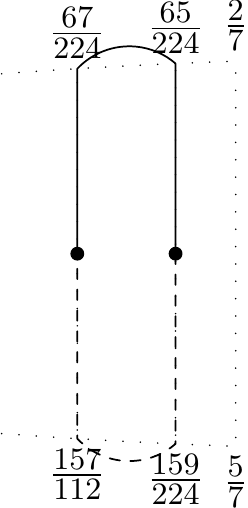}

$U(\beta (w_6),\beta (w_7))$ and $U(\beta (w_7),\beta (w_3))$.
\medskip

\includegraphics[width=4cm]{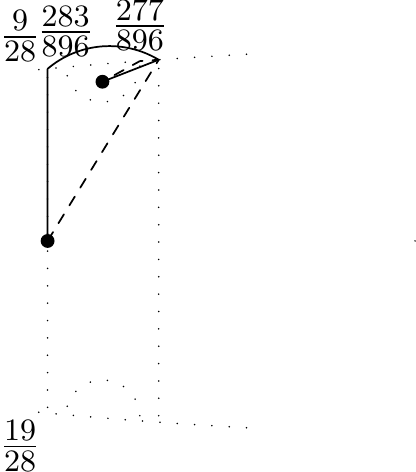}\ \ \ \ \ \ \includegraphics[width=3cm]{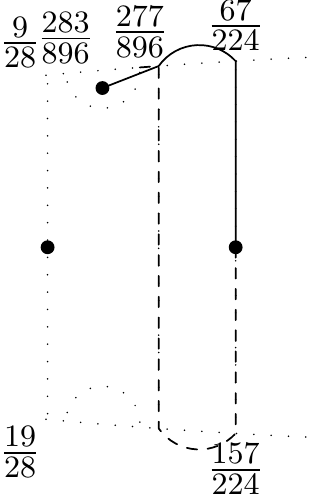}

\medskip

$U(\beta (w_8),\beta _{9/28})$ and $U(\beta (w_7),\beta (w_8))$  
\medskip

\includegraphics[width=5cm]{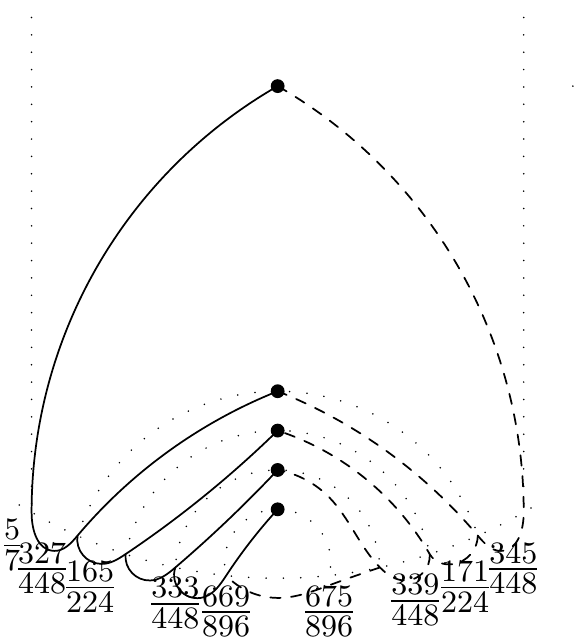}

\medskip 

$U(\beta _{2/7},\beta (w_{4}))$,  $U(\beta (w_{4}),\beta (w_{1}))$,
$U(\beta (w_{1}),\beta (w_{2}))$ and $U(\beta (w_{2}),\beta (w_{5}))$
\end{center}
 Note that, similarly to what happened with $\beta (w_{2})$, although $\beta (w_{5})$ is between $\beta (w_{2})$ and $\beta (w_{3})$, for the adjacent pair $(\beta (w_{2}),\beta (w_{5}))$, we have
$$U(\beta (w_{2}),\beta (w_{5}))\subset U(\beta (w_{1}),\beta (w_{2})).$$ 
Note, also, that $U(\beta (w_{6}),\beta (w_{7}))$ has two essential components of intersection with $\{ z:\vert z\vert <1\} $, which has not happened before. The reason for this becomes partially apparent when considering the image under $\xi _{5}$, which, as with $\xi _{4}$, coincides with $\xi _{5,2/7}$. The region bounded by $\overline{\beta (w_{6})}*\beta (w_{7})$ and $\xi _{5,2/7}(\partial 'U(\beta (w_{6}),\beta (w_{7})))$ again has two components of intersection, but overall the picture is simpler than it would have been if $U(\beta (w_{6}),\beta (w_{7}))$ had been chosen to have just one essential component of intersection with $\{ z:\vert z\vert <1\} $, adjacent to the endpoint of $\beta (w_{5})$, which might seem to be the more natural choice.

Note, also, that $\beta (w_{7})=\beta _{65/224}\in 
R_{5,0}$ is matched with$\beta '(w_{7}')$, where $w_{7}'=L_{2}R_{3}L_{3}^{2}L_{2}C$, and $\beta 
(w_{6})$ is matched with $\beta 
'(w_{6}')=\beta '(w_{7})$, which is, in fact, equal to $\beta 
_{159/224}$. For each matched 
$\beta $ and $\beta '$, the maps $\sigma _{\beta }\circ s$ and 
$\sigma _{\beta '}\circ s$ are Thurston equivalent. In particular, $\sigma _{159/224}\circ s$ is Thurston equivalent to $\sigma _{\beta (w_{6})}\circ s$. Since $R_{5}(a_{1},+)$, together with the sets $R_{5}(a_{1},-)$ and  $R_{5}(a_{0})$ and $R_{5}(\overline{a_{0}})$, determines a fundamental domain for $V_{3,5}$, and since $\beta _{65/224}$ and $\beta _{159/224}$ determine non-matched vertices of the fundamental domain, the two captures $\sigma _{65/224}\circ s$ and $\sigma _{159/225}\circ s$ are not Thurston equivalent.
This is the first stage at which captures $\sigma _{q}\circ s$ and 
$\sigma _{1-q}\circ s$ are not Thurston equivalent even though $q$ 
and $1-q$ are endpoints of the same lamination leaf: for $q={65\over 
224}$. The path $\beta (w_{6})=\beta (L_{3}^{4}L_{2}C)$, which is not a capture path, is drawn below.
\begin{center}
\includegraphics[width=6cm]{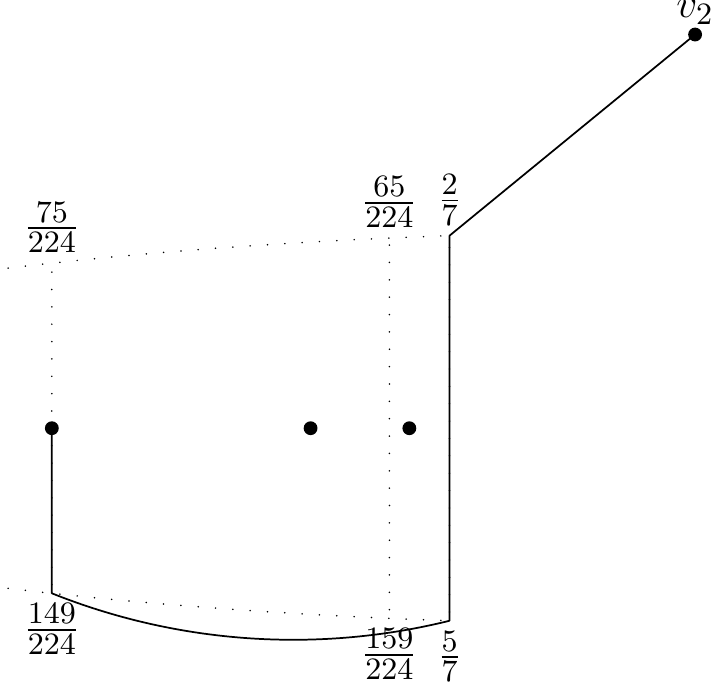}

\medskip

The path $\beta (w_6)$
\end{center}

\medskip

There is one path in $R_{5,1}$ between $\beta _{9/28}$ and $\beta _{37/112}$, which we call $\beta (w_{9})$, with $w_{9}=L_{3}L_{2}BCL_{1}R_{2}C$. This path is matched with $\beta '(w_{9}')$ in $R_{5,1}'$ between $\beta _{19/28}$ and $\beta _{75/112}$, where $w_{9}'=w_{8}=L_{3}L_{2}UCL_{1}R_{2}C$.

We now consider the case $m=6$.
We list all words corresponding to points in $U^{0}$ 
of preperiod $6$, remembering that the word corresponding to 
a point of preperiod $m$ has length $m+1$ or $m+2$ and ends in $C$. For preperiod $6$, it is still true that each word corresponds to exactly one path.
There are nine words, split between those 
ending in $R_{1}R_{2}C$ and those ending in $L_{2}C$:
\begin{equation}\label{7.5.9}\begin{array}{l}w_{10}=BCL_{1}R_{1}^{4}R_{2}C,\ 
w_{11}=BCL_{1}R_{2}BCL_{1}R_{1}R_{2}C,\cr 
w_{12}=BCL_{1}R_{2}UCL_{1}R_{1}R_{2}C\ \ 
w_{13}=L_{3}L_{2}UCL_{1}R_{1}^{2}R_{2}C,\cr  
w_{14}=L_{3}^{2}L_{2}BCL_{1}R_{1}R_{2}C, \end{array}\end{equation}

\begin{equation}\label{7.5.10}\begin{array}{l}w_{15}=BCL_{1}R_{2}R_{3}L_{3}L_{2}C,\ \ 
w_{16}=L_{3}(L_{2}R_{3})^{2}L_{2}C,\cr  
w_{17}=L_{3}L_{2}R_{3}L_{3}^{2}L_{2}C,\ \   
w_{18}=BCL_{1}R_{1}R_{2}R_{3}L_{3}L_{2}C. \end{array}\end{equation}
We note that there are also $2$ preperiod $6$ words corresponding to 
$U^{1}$:
$$w_{19}=L_{3}^{3}L_{2}R_{3}L_{2}C,\ w_{20}=L_{3}L_{2}BCL_{1}R_{1}^{2}R_{2}C,$$
and one corresponding to $U^{2}$, 
$$w_{21}=L_{3}^{5}L_{2}C.$$

As with the previous cases, the construction of $R_{6}(a_{1},+)$ and $R_{6}'(a_{1},+)$ is effected by the action of the maps $\xi _{5}$ on sets $\partial 'U(\beta _{1},\beta _{2})$ for adjacent pairs $(\beta _{1},\beta _{2})$ in $R_{5}(a_{1},+)$. We cannot expect, at this stage, that $\xi _{5}$ is always equal to $\xi _{5,2/7}$ on $\partial 'U(\beta _{1},\beta _{2})$, especially if $\beta _{1}$ and $\beta _{2}$ are relatively far from $\beta _{2/7}$.

 In some cases, we may want to redraw the paths of $R_{m}'(a_{1},+)$ up to finer isotopy preserving $Y_{m+1}$. This is permissible, but will, of course, affect the definitions of homeomorphisms $\xi _{n}$ for $n\geq m+1$.  We have not, of course, drawn the paths of $R_{m}'(a_{1},+)$ explicitly, but their definitions up to isotopy preserving $Y_{m}$ are given by the sets $\partial 'U(\beta _{1},\beta _{2})$, for adjacent pairs $(\beta _{1},\beta _{2})$ in $R_{m}(a_{1},+)$. An example is given by the adjacent pairs $(\beta _{2/7},\beta (w_{4}))$ and $(\beta (w_{5}),\beta (w_{6}))$ in $R_{5,0}'$. Let  $\psi _{5,j}$ and $\psi _{6,j}$, for  $j=1$ and $2$, be the homeomorphisms for these respective adjacent pairs, which are determined up to isotopy preserving $Y_{5}$ and $Y_{6}$. We have $\psi _{5,1}=\psi _{5,2/7}$ up to isotopy preserving $Y_{5}$ and $\psi _{6,1}=\psi _{6,2/7}$ on $\partial 'U(\beta _{2/7},\beta (w_{4}))$ up to isotopy preserving $Y_{6}$, and similarly for $\psi _{5,2}$ and $\psi _{6,2}$ on $\partial 'U(\beta (w_{5}),\beta (w_{6}))$. So we can define $\psi _{5,j}=\psi _{5,2/7}$ pointwise and $\psi _{6,j}=\psi _{6,2/7}$, for both $j=1$ and $2$. So then we have $\psi _{6,j}=\xi _{5,2/7}\circ \psi _{5,j}$. But if we look of the image of the two sets $\partial 'U(\beta _{2/7},\beta (w_{4}))$ and $\partial 'U(\beta (w_{5}),\beta (w_{6}))$, as originally drawn under $\xi _{5,2/7}$, they are  complicated, because of the nature of the intersection of $A_{5,2/7}\cup C_{5,2/7}$ with $\partial 'U(\beta _{2/7},\beta (w_{4}))$ and $\partial 'U(\beta (w_{5}),\beta (w_{6}))$. We concentrate first on the intersection with $C_{5,2/7}$. To simplify the images, we simply change the sets $\partial 'U(\beta _{2/7},\beta (w_{4}))$ and $\partial 'U(\beta (w_{5}),\beta (w_{6}))$, keeping the same sets up to homotopy preserving $Y_{5}$, but not up to homotopy preserving $Y_{6}$. We create an extra ``hook'' in the set $\partial 'U(\beta (w_{5}),\beta (w_{6}))$, so that the set $U(\beta (w_{5}),\beta (w_{6}))$ is increased at the expense of $U(\beta _{2/7},\beta (w_{4})$. This means that we have to redraw the paths $\beta '(w_{5}')$ and $\beta '(w_{4}')$, and, in consequence, the paths in $R_{5}(a_{1},+)$ between them also, that is, $\beta '(w_{1}')$ and $\beta '(w_{2}')$. Below, we draw the modified set $\partial 'U(\beta (w_{5}),\beta (w_{6}))$ shown by the dashed lines, and also the modified path $\beta '(w_{5}')$.
\begin{center}
\includegraphics[width=8cm]{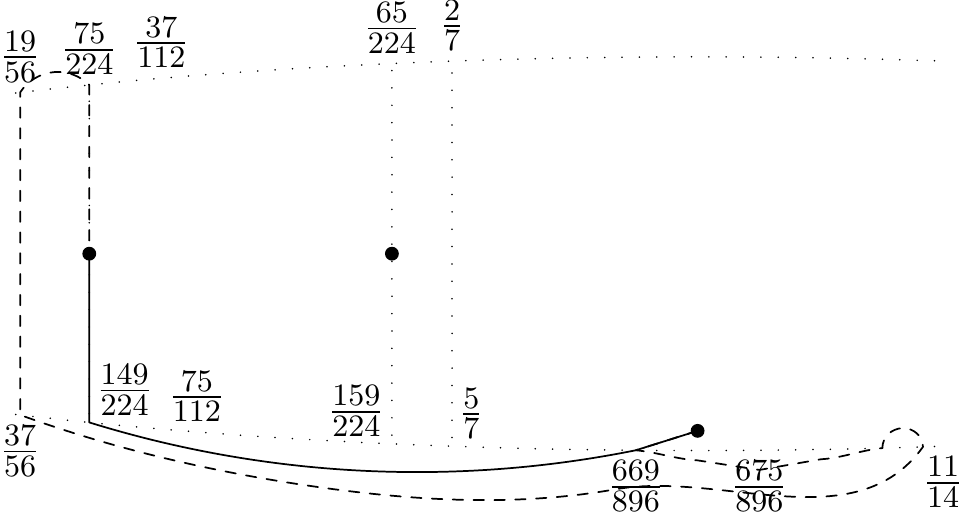}

\medskip

The set $\partial 'U(\beta (w_{5}),\beta (w_{6}))$, after modification

\medskip \medskip
\includegraphics[width=6cm]{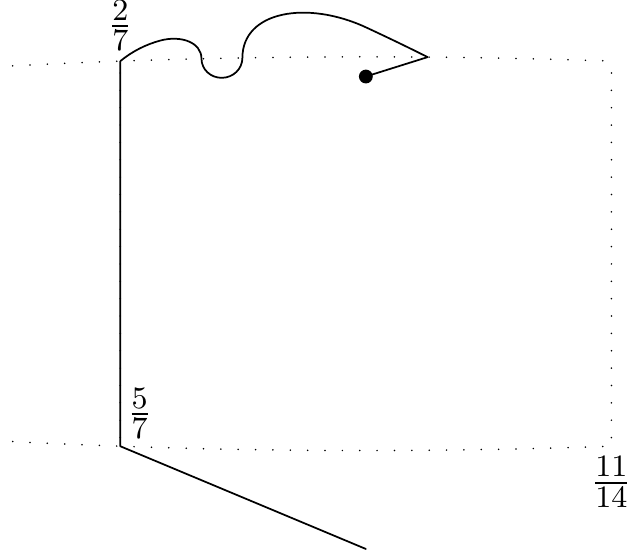}

\medskip

The path $\beta '(w_5')$, after modification.
\end{center}

\medskip

We also consider the intersection of $A_{5,2/7}$ with $U(\beta _{2/7},\beta (w_{4}))$, $U(\beta (w_{5}),\beta (w_{6}))$ and $U(\beta (w_{6}),\beta (w_{7}))$. The image of $\partial 'U(\beta (w_{6}),\beta (w_{7}))$ under $\xi _{5,2/7}$ is then as shown. 
\begin{center}
\includegraphics[width=7cm]{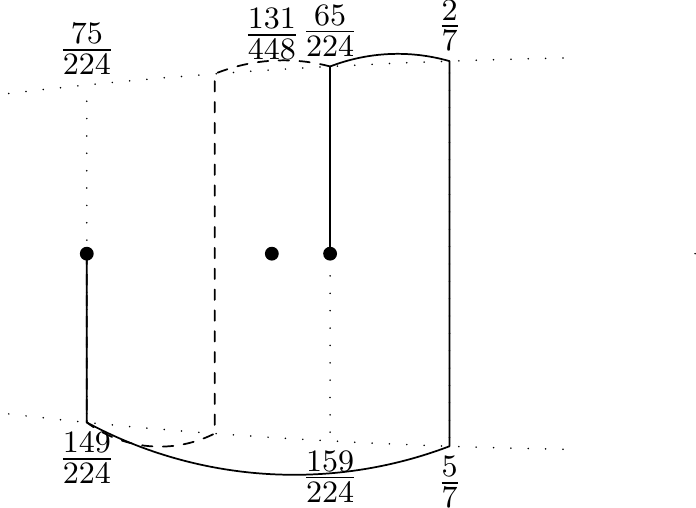}

\end{center}

\medskip
 
  Taking into consideration all these intersections, the images under $\xi _{5,2/7}$ of (modified) $\partial'U(\beta _{2/7},\beta (w_{4})$, (modified) $\partial 'U(\beta (w_{5}),\beta (w_{6}))$, and $\partial 'U(\beta (w_{6}),\beta (w_{7}))$, the sets bounded by these and $\overline{\beta _{2/7}}*\beta (w_{4})$ and $\overline{\beta (w_{5})}*\beta (w_{6})$ and $\overline{\beta (w_{6})}*\beta (w_{7})$ subdivide into sets $U(\beta _{1},\beta _{2})$ for adjacent pairs $(\beta _{1},\beta _{2})$ from $R_{6,0}$. The order of paths $\beta (w)\in R_{6,0}$ is with $w$ as follows:

\begin{equation}\label{7.5.11}\begin{array}{l}C,\ w_{11}=BCL_{1}R_{2}BCL_{1}R_{1}R_{2}C,\ w_{12}=BCL_{1}R_{2}UCL_{1}R_{1}R_{2}C,\cr
w_{4}=BCL_{1}R_{2}R_{3}L_{2}C,\ w_{1}=BCL_{1}R_{1}R_{2}C,\cr w_{2}=BCL_{1}R_{1}^{2}R_{2}C,\ w_{5}=BCL_{1}R_{1}^{3}R_{2}C,\cr w_{10}=BCL_{1}R_{1}^{4}R_{2}C,\ w_{18}=BCL_{1}R_{2}R_{3}L_{3}L_{2}C,\cr
w_{14}=L_{3}^{2}L_{2}BCL_{1}R_{1}R_{2}C,\ w_{6}=L_{3}^{4}L_{2}C,\cr w_{17}=L_{3}L_{2}R_{3}L_{3}^{2}L_{2}C,\ w_{7}=L_{3}L_{2}R_{3}L_{3}L_{2}C,\cr w_{3}=L_{3}L_{2}R_{3}L_{2}C,\ w_{16}=L_{3}(L_{2}R_{3})^{2}L_{2}C,\cr w_{13}=L_{3}L_{2}UCL_{1}R_{1}^{2}R_{2}C,\ w_{8}=L_{3}L_{2}UCL_{1}R_{1}R_{2}C,\cr L_{3}L_{2}C,\  w_{9}=L_{3}L_{2}BCL_{1}R_{1}R_{2}C,\  w_{20}=L_{3}L_{2}BCL_{1}R_{1}^{2}R_{2}C.\end{array}
\end{equation}

For example, we have the following.
\begin{center}
\includegraphics[width=7cm]{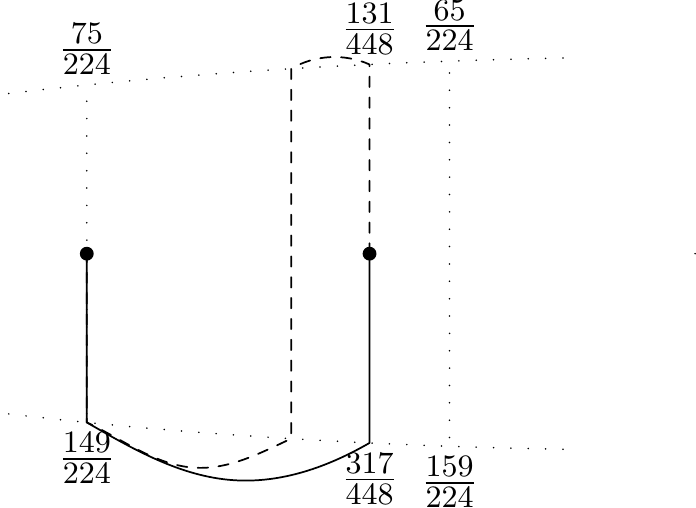}\ \ \ \includegraphics[width=3.5cm]{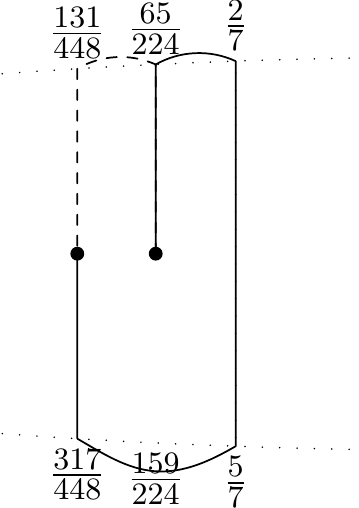}

\medskip

$U(\beta (w_6),\beta (w_17))$ and $U(\beta (w_{17}),\beta (w_7))$.
\end{center}

\medskip

Now, without considering $R_{6}(a_{1},+)$ any further, we consider $R_{7}(a_{1},+)$.  For preperiod $7$ there are $24$ words in all corresponding to points in $\cup _{p}U^{p}$, with $12$ ending in 
$L_{2}C$ and $12$ ending in $R_{1}R_{2}C$. However note that for $w_{22}=L_{3}(L_{2}R_{3}^{2}L_{3}L_{2}C$, $D(w_{22})$ is contained in both $U^{0}$ and $U^{1}$. We therefore expect a path in both $R_{7,0}$ and $R_{7,1}$ with endpoint in $D(w_{22})$. We shall write $\beta (w_{22},0)$ and $\beta (w_{22},1)$ for these paths.

First, for $\beta (w_{22},0)$, we consider the adjacent pair $(\beta (w_{17}),\beta (w_{6}))$ in $R_{6,0}$. We consider $\partial 'U(\beta (w_{17}),\beta (w_{6}))$. On this set, if $\psi _{7}$ is the homeomorphism effecting the matching, we have $\psi _{7}=\psi _{7,2/7}=\xi _{6,2/7}\circ \psi _{6}$. The only part of the support of $\xi _{6,2/7}$ which intersects $\partial 'U(\beta (w_{17}),\beta (w_{6}))$ is $A_{6,2/7}$, with intersection occurring on the piece of track
$$BCL_{1}R_{2}R_{3}L_{3}^{2}\to L_{3}(L_{2}R_{3})^{2}L_{3},$$
and the image of $\partial 'U(\beta (w_{17}),\beta (w_{6}))$ is as shown.
\begin{center}
\includegraphics[width=7cm]{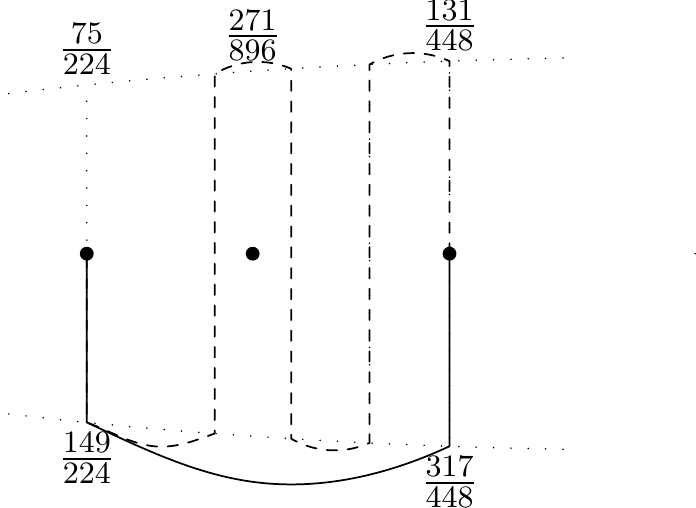}

\medskip 

$\xi _{6,2/7}(\partial 'U(\beta (w_{17}),\beta (w_{6})))$
\end{center}

\medskip

Therefore, for $w_{22}=L_{3}(L_{2}R_{3})^{2}L_{3}L_{2}C$, the path $\beta 
(w_{22},0)$ is between $\beta (w_{6})$ and $\beta (w_{17})$.  Since the region bounded by $\overline{\beta (w_{1})}*\beta (w_{6})$ and $\xi _{6,27}(\partial 'U(\beta (w_{17}),\beta (w_{6})))$ has two essential components of intersection with $\{ z:\vert z\vert <1\} $, up to isotopy preserving $Z_{n}$ for any $n$, it is relatively easy to divide this region up into $U(\beta (w_{17}),\beta (w_{22},0))$ and $U(\beta (w_{22},0),\beta (w_{6}))$.

The path $\beta (w_{22},1)$ is involved in the first example of a set $\partial 'U(\beta _{1},\beta _{2})$ such that  
$$\xi _{n}(\partial 'U(\beta _{1},\beta _{2}))\neq \xi _{n,2/7}(\partial 'U(\beta _{1},\beta _{2})){\rm{\ in\ }}\pi _{1}(\overline{\mathbb C}\setminus Z_{n+1},Z_{n+1},Z_{n+1}),$$
with $n=6$, and  $\beta _{1}=\beta (w_{9})$, $\beta _{2}=\beta (w_{19})$.  These two paths are adjacent in $R_{6,1}$. Here, 
$$\xi _{6}(\partial 'U(\beta _{1},\beta _{2}))=\xi _{6,9/28}(\partial 'U(\beta _{1},\beta _{2}))=\xi _{6,9/28}'(\partial 'U(\beta _{1},\beta _{2})){\rm{\ in\ }}\pi _{1}(\overline{\mathbb C}\setminus Z_{7},Z_{7},Z_{7}).$$
The region bounded by $\overline{\beta _{1}}*\beta _{2}$ and $\xi _{6}(\partial 'U(\beta _{1},\beta _{2}))$ is as shown, with the points of $Z_{7}$ inside the region marked. 
\begin{center}
\includegraphics[width=4cm]{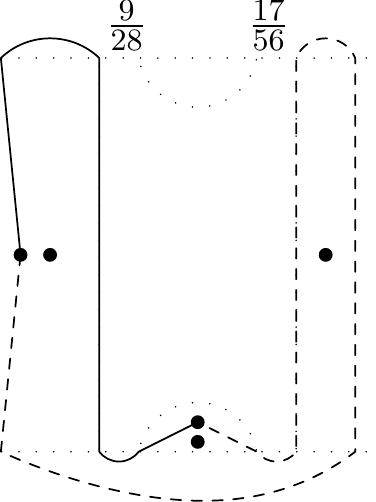}

\medskip

   Image under $\xi _{6,9/28}'$
\end{center} 

\medskip
  
This means that in $R_{7,1}$ we want to insert another three paths between $\beta (w_{19})$ and $\beta (w_{9})$: a path $\beta (w_{24})$ with endpoint in $D(w_{24})$, where
$w_{24}=L_{3}^{3}L_{2}R_{3}L_{3}L_{2}C$,
a path $\beta (w_{23})$ with $w_{23}=L_{3}L_{2}BCL_{1}R_{1}^{3}R_{2}C$,
and a path $\beta (w_{22},1)$ with three $S^{1}$-crossings,a disc crossing from  at $e^{2\pi i(9/28)}$ to $e^{2\pi i(19/28)}$, and then a crossing into $D(w_{22})$, at the lower right-most boundary point. The other path with the same endpoint, $\beta (w_{22},0)$, is between $\beta (w_{6})$ and $\beta (w_{17})$. Note that this is different from the picture in $R_{4,0}$, of the paths between $\beta (w_{2})$ and $\beta (w_{3})$, although there is some resemblance. The paths $\beta (w_{23})$ and $\beta (w_{24})$ are analogues of $\beta (w_{5})$ and $\beta (w_{7})$ and indeed are images of these paths under a local inverse of $s^{2}$, up to homotopy preserving $Y_{6}$. But $\beta (w_{6})$ and $\beta (w_{22},1)$ are not related in this way, even though $\beta (w_{22},1)$ is between $\beta (w_{23})$ and $\beta (w_{24})$ in the ordering, while $\beta (w_{6})$ is between $\beta (w_{5})$ and $\beta (w_{7})$ in the ordering. 

Another example occurs for $n=7$. For this one it is convenient to modify some  paths in $R_{7,1}'$ as elements of $\pi _{1}(\overline{\mathbb C}\setminus Z_{8},Z_{8},v_{2})$, so that, for $\beta _{1}=\beta (w_{9})$ and $\beta _{2}=\beta (w_{22},1)$ as elements of $\pi _{1}(\overline{\mathbb C}\setminus Z_{8},Z_{8},v_{2})$, there is a kink in $\partial 'U(\beta _{1},\beta _{2})$ much like that in $\partial 'U(\beta (w_{5}),\beta (w_{6}))$, effected by changing the definition of paths in $R_{5,0}'$ as elements of $\pi _{1}(\overline{\mathbb C}\setminus Z_{6},Z_{6},v_{2})$. Then
$$\xi _{7}(\partial 'U(\beta _{1},\beta _{2}))\neq \xi _{7,2/7}(\partial 'U(\beta _{1},\beta _{2})){\rm{\ in\ }}\pi _{1}(\overline{\mathbb C}\setminus Z_{8},Z_{8},Z_{8})$$
As before we have
$$\xi _{7}(\partial 'U(\beta _{1},\beta _{2}))=\xi _{7,9/28}(\partial 'U(\beta _{1},\beta _{2})) =\xi _{7,9/28}'(\partial 'U(\beta _{1},\beta _{2})){\rm{\ in\ }}\pi _{1}(\overline{\mathbb C}\setminus Z_{7},Z_{7},Z_{7}).$$
The region bounded by $\overline{\beta _{1}}*\beta _{2}$ and $\xi _{7}(\partial 'U(\beta _{1},\beta _{2}))$ is as shown, with the points of $Z_{8}$ inside the region marked.
\begin{center}
\includegraphics[width=2.5cm]{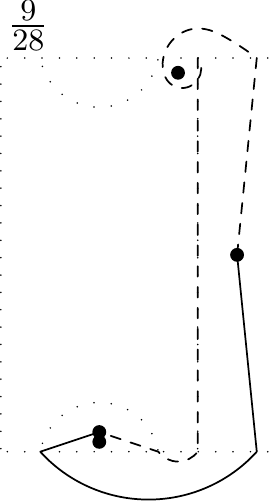}

\medskip

   Image under $\xi _{7,9/28}'$
\end{center}  

\medskip
 
We therefore define $\beta (w_{25})$, for $w_{25}=L_{3}L_{2}UCL_{1}R_{2}BCL_{1}R_{1}R_{2}C$, to have two complete disc crossings, both in the boundary of $D(L_{3}L_{2}C)$, the first from $e^{2\pi i(2/7)}$ to $e^{2\pi i(5/7)}$ and the second from $e^{2\pi i(11/14)}$ to $e^{2\pi i(3/14)}$, and endpoint in $D(w_{25})$. Note that $D(w_{25})\subset U_{1}\setminus  U^{0}$, and in fact $\beta (w_{25})$ is the unique path with its endpoint. It is between $\beta (w_{23})$ and $\beta (w_{22},1)$. It corresponds in some sense to the path $\beta (w_{14})=\beta (L_{3}^{2}L_{2}BCL_{1}R_{1}R_{2}C)$. But clearly it is not an inverse image of this path under $s^{2}$.

The set of paths in $R_{m}(a_{1},+)$, and the ordering given on $\cup _{p}R_{m,p}$ in (\ref{7.5.11}), are obviously rather involved. Nevertheless, a possible strategy is emerging. The most obvious point is that the sets of paths $R_{m,p}$ and $R_{m,p}'$ should be obtained by induction on $m$. In fact,the paths of $R_{m+1,p}$ between an adjacent pair $(\beta _{1},\beta _{2})$ of $R_{m,p}$, and the matching paths in $R_{m+1,p}'$ have been obtained from the set $U(\beta _{1},\beta _{2})$. This means that the definition of the paths in $R_{m+1,p}'$ has been less direct than the definition for paths in $R_{m+1,p}$. Also, it will be noticed that the shape of the sets $U(\beta _{1},\beta _{2})$, and the image of the boundary subset $\partial 'U(\beta _{1},\beta _{2})$ under the homeomorphism $\xi _{m}=\psi _{m+1}\circ \psi _{m}^{-1}$, is strongly influenced by the homeorphism $\xi _{m,2/7}$ and its support $A_{m,2/7}\cup C_{m,2/7}$, in the examples so far considered. 

The later examples show that it is not always true that $\xi _{m}=\xi _{m,2/7}$ on $\partial 'U(\beta _{1},\beta _{2})$,  for adjacent pairs $(\beta _{1},\beta _{2})$ in $\cup _{p}R_{m,0}$. In fact it is not even true if we restrict to $R_{m,0}$, although we have not yet seen any examples. But we shall obtain some control on the variation of the homeomorphisms $\xi _{m}$ on sets $\partial 'U(\beta _{1},\beta _{2})$.
Also, in the examples so far considered, it has been true that the paths in $R_{m,p}$ all have distinct endpoints in  $Z_{m}\cap U^{p}$, and each point of $Z_{m}\cap U^{p}$ has been the endpoint of exactly one path in $R_{m,p}$. This does not remain true for all $m$, although it is always true that the number of points in $Z_{m}\cap U^{p}$ is the same as the number of paths in $R_{m,p}$. There is a clue in the inclusion of $D(w_{22})$ in both $ U^{0}$ and $U^{1}$, and of $D(w_{25})$ in both $U^{0}$ and $U^{1}$. In a somewhat analogous way, there are two paths in $R_{m,0}$ with endpoints in
$D(L_{3}L_{2}R_{3}L_{3}(L_{2}R_{3})^{2}L_{3}L_{2}R_{3}L_{3}L_{2}C)$. But there are no paths with endpoint in $D(L_{3}L_{2}R_{3}L_{3}^{6}L_{2}R_{3}L_{3}L_{2}C)$.

 \section{Definition of  $w_{i}(w,0)$ and  $w_{i}'(w,0)$}\label{7.6}
 
In \ref{2.10}, we defined a set $U^{p}$ for $p\geq 0$. Recall also from \ref{2.10} that there is the same number of paths in $R_{m,p}$ as of points in $U^{p}\cap Z_{m}$, although it is not the case that every point in $U^{p}\cap Z_{m}$ is the endpoint of exactly one path.  A point of $U^{p}\cap Z_{m}$ is determined uniquely by a set $D(w)$ containing it, where $w$ is a word of length $m+1$ ending in $L_{2}C$, or of length $m+2$ ending in $R_{1}R_{2}C$. It turns out that a path in $R_{m,p}$ is determined by its second endpoint and its first $S^{1}$ -crossing. 

From now on, we consider only the case $p=0$. In order to identify the path corresponding to $w$, we shall define a sequence $w_{i}(w,0)$ of prefixes of $w$, and a corresponding sequence $w_{i}'(w,0)$, and finally, a sequence $w_{i}'(w,x,0)$ where $x=w_{1}'(w,0)$. The sequences $w_{i}'(w,0)$ and $w_{i}'(w,x,0)$ are not, in general, sequences of prefixes of $w$. For some $n$, $w_{n}(w,0)=w$, and then $w_{i+1}'(w,0)=w_{i}'(w,0)$, and $w_{i+1}'(w,x,0)=w_{i}'(w,x,0)$, for $i\geq n+1$. 

We recall the definition of $U^{0}$ from \ref{2.10}:
$$U^{0}=\begin{array}{l}D(L_{3}L_{2})\cup D(BC)\\
\cup \cup _{k=1}^{\infty }(D(L_3^{3k+1}L_{2})\cup D(L_{3}^{3k-1}L_{2}BC))\setminus (\cup _{k\geq 1}D(L_{3}L_{2}(UCL_{1}R_{2})^{k}BC)).\end{array}$$

We define $w_{1}(w,0)$ for $w$ with $D(w)\subset U^{0}$ and with $w$ ending in $C$. The definition also works for infinite words.  The first letter of $w$ is 
either $BC$ or $L_{3}$. 
If $BC$ is the first letter, then $w_{1}(w,0)=BC$. 
If $L_{3}$ is the first letter of $w$, and $w$ starts with $L_{3}^{4}$ or $L_{3}^{2}L_{2}BC$, then $w_{1}(w)=L_{3}^{4}$ or $L_{3}^{2}$ respectively. Otherwise, $w$ starts $L_{3}L_{2}$, and 
we look for the first occurrence of one of the following. One of 
these must occur if $w$ is preperiodic, that is, ends in $C$. 

\begin{itemize}
\item[0.] An occurrence of $C$.
\item[1.] An occurrence 
of $xBC$ for any $n\geq 0$, where $x$ is a maximal word in the letters 
$L_{3}$, $L_{2}$, $R_{3}$ with an even number of $L$ letters. 
\item[2.] An occurrence 
of $xUC$, where $x$ is a maximal word in the letters 
$L_{3}$, $L_{2}$, $R_{3}$, with an odd number of $L$ letters. 
\item[3.] An occurrence of $L_{1}R_{2}UC$. 
\item[4.] An occurrence of  $R_{1}R_{2}BC$.  
\item[5.] An occurrence of  a string $v_{1}\cdots v_{n}L_{3}$ which is the end of $w$, or is followed in $w$ by $L_{2}$, where:
\begin{itemize} 
\item $v_{i}=L_{3}(L_{2}R_{3})^{m_{i}}L_{3}^{r_{i}-1}$ for $i\geq 2$, where $m_{i}+r_{i}=3$; 
\item either $v_{1}=L_{3}(L_{2}R_{3})^{m_{1}}L_{3}^{r_{1}-1}$ for $m_{1}\geq 1$ and  an odd $m_{1}+r_{1}\geq 3$, or $v_{1}=XL_{1}R_{1}^{p_{1}}R_{2}R_{3}(L_{2}R_{3})^{m_{1}}L_{3}^{r_{1}-1}$ for $3\leq r_{1}\leq 5$, and $m_{1}\geq 0$, and $p_{1}\geq 0$, and $X=BC$ or $UC$;
\item $v_{n}L_{3}$ is either followed by $L_{2}$, or is at the end of $w$; 
\item the string $v_{1}\cdots v_{n}L_{3}$ is maximal with these properties; 
\item $n$ is odd.
\end{itemize}

\end{itemize}

We define $w_{1}(w,0)$ to be the prefix 
of $w$ ending in the occurrence listed  -- except in case 1 when $n>0$, and in case 4 when $w$ ends in $L_{1}R_{1}R_{2}BC$, and in case 5.  In case 4, we define $w_{1}(w,0)$ to be the word up to, and not including, this occurrence of $L_{1}R_{1}R_{2}BC$, and $w_{2}(w,0)$ to be the prefix of $w$ ending in this occurrence of $L_{1}R_{1}R_{2}BC$. In case 5, if the maximal string of words of this type is $v_{1}\cdots v_{n}$, then we define $w_{1}(w,0)$ to be the prefix ending with $v_{n}L_{3}$.

Now suppose that the first letter of $w$ is $BC$. Then $w$ starts with $(BCL_{1}R_{2})^{k}$ for a maximal $k\geq 0$. 
Then we define 
$w_{2}(w,0)=w$ if there is no occurrence of cases 1 to 5 above. 
We 
define $w_{2}(w,0)$ to be the prefix of $w$ defined exactly  as 
$w_{1}(w,0)$ is 
defined when $w$ starts with $L_{3}$, if there is an occurrence of one of 2 to 5.

It may not be immediately apparent what determines these choices.  Recall that {\em{beads}} on $C_{m,2/7}$ and $A_{m,2/7}$ were defined in \ref{7.2}. In cases 1 to 4, if $w_{1}=w_{1}(w,0)$ has length $m$, then $D(w_{1})$ is a bead on 
$C_{m,2/7}$,  and in case 5 it is a bead on
$A_{m,2/7}\cap \{ z:\vert z\vert \leq 1\} $. Roughly half the number of beads of $A_{m,2/7}$ and $C_{m,2/7}$ in $U^{0}$ are listed. If a bead is listed, then the adjacent bead on $A_{m,2/7}$ or $C_{m,2/7}$ is not listed. If a bead in $U^{0}$ is not listed, then any adjacent beads which are also in $U^{0}$ are listed.

The reason for the rather convoluted additional conditions in 5 is basically that, without them, the order of paths will become extremely complicated. The examples in \ref{7.3} did not get up to high enough preperiod to illustrate the problem, which first appears at preperiod 11. But it is similar to modifications which were made in \ref{7.3} to place $\beta (w_{14})$ between $\beta (w_{5})$ and $\beta (w_{6})$, rather than between $\beta (C)$ and $\beta (w_{5})$. 

Now we define $w_{1}'(w,0)=w_{1}'(w_{1},0)$. If $w=C$ or if $w_{1}$ starts with $BC$ or $L_{3}^{2}$, then we define $w_{1}'(w,0)=L_{3}$. Now suppose that $w$ starts with $L_{3}L_{2}$, and we assume that one of cases 
0 to 5 holds, which, as pointed out, is always true if $w$ ends in 
$C$. 

\begin{itemize}
\item[0.] $w_{1}'$ is the same length as $w_{1}$, ending in 
    $L_{3}$, if preceded by a maximal subword in the letters $L_{3}$, $L_{2}$ or $R_{3}$ and with an even number of $L$ letters, or ending in $R_{3}$ if preceded by a maximal subword in the letters $L_{3}$, $L_{2}$ or $R_{3}$ and with an even number of $L$ letters, such that $D(w_{1}')$ and $D(w_{1})$ are adjacent.
    \item[1.]  Let $w_{1}''$ be obtained from $w_{1}$ by deleting the last $(UCL_{1}R_{2})^{n}BC$, and adding $C$. Then $w_{1}'$ is the same length as $w_{1}''$, ending in $L_{3}$, such that $D(w_{1}')$ and $D(w_{1}'')$ are adjacent. If $n>0$ then, $w_{2}'$ is the same length as $w_{1}''$, ending in $R_{3}$, such that $D(w_{2}')$ and $D(w_{1}'')$ are adjacent.
    \item[2.] Let $w_{1}''$ be obtained from $w_{1}$ by deleting the last $UC$ and adding $C$. Then  $w_{1}'$ is defined similarly to 1, but ending in 
    $R_{3}$.
    \item[3.] $w_{1}'$ is obtained 
from $w_{1}$ by replacing the last letter by $BC$.
\item[4.] If $w_{1}$ ends in $R_{1}^{2}R_{2}BC$, then $w_{1}'$ is obtained 
from $w_{1}$ by replacing the last letter by $R_{3}L_{2}R_{3}L_{3}$.  If $w_{1}$ is extended in $w$ by $L_{1}R_{1}R_{2}BC$, then $w_{1}'(w,0)$ is the word of the same length as $w_{1}$, with $D(w_{1}')$ to the left of $D(v')$, where $v'$ is obtained from $w_{1}$ by replacing the last letter by $C$. We then also define $w_{2}'(w,0)=v''$, where $v''$ is of the same length as $v'$, with $D(v')$ between $D(w_{1}')$ and $D(v'')=D(w_{2}')$.
\item[5.] The word $w_{1}'$ is obtained 
from $w_{1}$ by deleting all but the first letter of $v_{n}L_{3}$ if $n>1$, or if $n=1$, $m_{1}+r_{1}=3$ and the first letter of $v_{1}$ is $L_{3}$. In the other cases with $n=1$, the word  $w_{1}'$ is obtained from $w_{1}$ by replacing the last $L_{3}^{3}$ by $L_{2}R_{3}L_{3}$ if $r_{1}\geq 3$, by deleting the last $L_{2}R_{3}L_{3}^{2}$ if $r_{1}=2$, and by deleting the last $(L_{2}R_{3})^{2}L_{3}$ if $r_{1}=1$.  

\end{itemize}

In many of these cases, there is a connection between $w_{1}'$ and the bead $D(w_{1})$ on $A_{n,2/7}$ or $C_{n,2/7}$. In parts of cases 1, 2, and 4, $D(w_{1}')$ is adjacent to the string on $C_{n,2/7}$ between the bead preceding $D(w_{1})$, which is mapped to $D(w_{1})$ by $\xi _{n,2/7}$, and the bead $D(w_{1})$ itself. In the other part of case 4, $D(w_{1}')$ is adjacent to an earlier piece of string on the same component of $C_{n,2/7}$. In case 3, $D(w_{1}')$ is the bead preceding $D(w_{1})$ on $C_{n,2/7}$ and is mapped to it by $\xi _{n,2/7}$. In case 5, $D(w_{1}')$ contains the bead on $A_{n,2/7}$ which is mapped by $\xi _{n,2/7}$ to $D(w_{1})$.

Now we 
define $w_{i+1}(w,0)$ if $w_{i}=w_{i}(w,0)$ has been defined.  If 
$w_{i}(w,0)=w$ then $w_{i+1}(w,0)=w$. Now suppose that
$w_{i}\neq w$. Write 
$w=w_{i}v_{i}$. If $w_{i}$ 
ends with $L_{3}$, then
$$w_{i+1}=w_{i}v_{i}',{\rm{\ where\ }}L_{3}v_{i}'=w_{1}(L_{3}v_{i},0).$$
If $w_{i}$ ends with $BC$ or $UC$ then
$$w_{i+1}=w_{i}v_{i}',{\rm{\ where\ }}BCv_{i}'=w_{2}(BCv_{i},0).$$

We refer to {\em{cases 1 to 5 of $w_{i+1}$}}, if $L_{3}v_{i}$ or $BCv_{i}$ 
is as in case 1 to 5 for the definition of $w_{1}(L_{3}v_{i})$ or 
$w_{2}(BCv_{i})$. 

Once again, if $w_{i}$ has length $n$, then $D(w_{i})$ is a bead on $A_{n,2/7}$ or $C_{n,2/7}$. Note that the length of $w_{i}(w,0)$ is strictly increasing with $i$ until $w_{i}(w,0)=w$.

We define $w_{i}'$ from $w_{i}$ in exactly the same way as $w_{1}'$ is defined from $w_{1}$, in each of cases 0 to 5. Once again, there is a connection between $w_{i}'$ and the bead and piece of string preceding the bead $D(w_{i})$ on $A_{n,2/7}$ or $C_{n,2/7}$. 

\section{Definition of $U^{x}$ and $w_{i}(w,x,0)$ and $w_{i}'(w,x,0)$}\label{7.9}

Let $W_{1}(0)$ be on the set of words of the form $w_{1}'(w,0)$, and all words of the form $w_{2}'(w,0)$ for $w$ starting with $BC$. So  any word in $W_{1}(0)$ ends in either $L_{3}$ or $R_{3}$.

 It seems best not to go into too much detail at this stage, but the definitions in this section are related to a sequence of homeomorphisms $\xi _{n}$ related to $x$, in much the same way as the definitions of $w_{i}(w,0)$ and $w_{i}'(w,0)$ in \ref{7.6} were related to the family of homeomorphisms $\xi _{n,2/7}$. Many of the definitions are the same, simply because the homeomorphism $\xi _{n}$ is reasonably close to $\xi _{n,2/7}$. 

We put an ordering on the words of $W_{1}(0)$ which start with $L_{3}$, using the anticlockwise ordering on the upper half of the unit circle.  If $w$ is any finite word ending in $L_{3}$ or $R_{3}$, we define $\overline{w}=w(L_{2}R_{3}L_{3})^{\infty }$.  Then $\overline{w}$ labels a leaf $\ell (w)$ of $L_{3/7}$, which is in the boundary of a unique gap $G(w)$. If $w$ contains no letter $BC$ or $UC$, then $\ell (w)$ is a vertical leaf. If $w\in W_{1}(0)$ starts with $L_{3}$, then $\ell (w)$ has at least one endpoint on the upper half unit circle. For such $w$, we define $P(w)$ to be the right-most point of $\ell (w)$ on the upper half unit circle. Then we define $w<w'$ if $P(w)$ is to the right of $P(w')$. This means that either all of $\ell (w)$ is to the right of $\ell (w')$, or the disc $D(w)$, which is bounded by $\ell (w)$, contains $D(w')$.

Now we define $U^{x}$ for $x\in W_{1}(0)$. If $x=L_{3}$, then $U^{x}$ is the union of all $D(w)$ with $w_{1}'(w,0)=L_{3}$, where the union is taken over all words not ending in $L_{3}$. Now suppose that the word $x$ starts with $L_{3}$  and that there is more than one word $w$ with $w_{1}(w,0)=x$. Then there are just two possibilities.
\begin{itemize}
\item If $x$ ends in $L_{3}(L_{2}R_{3})^{2j}L_{3}^{2}$ for some $j\geq 0$, then we define $x_{1}$ by replacing this by  $L_{3}(L_{2}R_{3})^{2j+1}$
\item If $x$ ends in $L_{3}(L_{2}R_{3})^{2j}$ for some $j\geq 1$, then we define $x_{1}$ by replacing this by $L_{3}(L_{2}R_{3})^{2j-1}L_{3}^{2}$.

\end{itemize}

Then $U^{x}$ is the union of all $D(w)$  for which $w$ satisfies  $w_{1}'(w,0)=x$, except for those  of the form $xL_{3}^{3}L_{2}v$ or $x(L_{2}R_{3})^{2}L_{3}L_{2}v$  or $xL_{2}R_{3}L_{3}^{2}L_{2}R_{3}L_{3}L_{2}v$  with $L_{3}L_{2}v\leq x$ . In these cases:  
\begin{itemize}
\item $D(xL_{3}^{3}L_{2}v)$ is replaced in $U^{x}$ by $D(x_{1}L_{2}R_{3}L_{3}L_{2}v)$;

\item $D(x(L_{2}R_{3})^{2}L_{3}L_{2}v$ is replaced in $U^{x}$ by $D(xL_{3}^{5}L_{2}v)$.

\item $D(x(L_{2}R_{3}L_{3}^{2}L_{2}R_{3}L_{3}L_{2}v)$ is replaced in $U^{x}$ by $D(x(L_{2}R_{3})^{2}L_{3}^{3}L_{2}v)$.
\end{itemize}

Note that the number of words extending  $xL_{3}^{3}L_{2}$ and ending in $L_{2}C$ or $R_{1}R_{2}C$, of any fixed length, is the same as the number of words of the same type extending $x_{1}L_{2}R_{3}L_{3}L_{2}$. A similar statement holds for $x(L_{2}R_{3})^{2}L_{3}L_{2}$ and $xL_{3}^{5}L_{2}$, and for $xL_{2}R_{3}L_{3}^{2}L_{2}R_{3}L_{3}L_{2}$ and $x(L_{2}R_{3})^{2}L_{3}^{3}L_{2}$. Therefore, the number of words $w$ of preperiod $m$ with $D(w)\subset U^{x}$ is the same as the number of words of preperiod $m$ with $w_{1}'(w,0)=x$.

Note also that for $x$ and $x_{1}$ as above, 
$$w_{1}(xL_{3}^{3}L_{2},0)=xL_{3}^{3}{\rm{\ and\ }}w_{1}(x_{1}L_{2}R_{3}L_{3}L_{2})=x_{1}L_{2}R_{3}L_{3},$$
$$w_{1}(x(L_{2}R_{3})^{2}L_{3}L_{2})=x(L_{2}R_{3})^{2}L_{3}{\rm{\ and\ }}w_{1}(xL_{3}^{5}L_{2})=xL_{3}^{5},$$
$$w_{1}(xL_{2}R_{3}L_{3}^{2}L_{2}R_{3}L_{3}L_{2})=xL_{2}R_{3}L_{3}^{2},\ \ w_{1}(x(L_{2}R_{3})^{2}L_{3}^{3}L_{2})=x(L_{2}R_{3})^{2}L_{3}^{3}.$$

Now suppose that $x=w_{2}'(w,0)$, for at least one $w$ starting with $BC$. Then $U^{x}$ is the union of all $D(w)$ with $w_{2}'(w,0)=x$. In this case, $w_{i}(w,x,0)=w_{i}(w,0)$ and $w_{i}'(w,x,0)=w_{i}'(w,0)$ for all $w$ with $D(w)\subset U^{x}$.

Now suppose that $x$ starts with $L_{3}$. We define $w_{1}'(w,x,0)=x$ for all $w$ with $D(w)\subset U^{x}$.  So, in most, but not all, cases,  $w_{1}'(w,x,0)=w_{1}'(w,0)$. We always define $w_{i}(w,x,0)=w_{i}(w,0)$. Now we define  $w_{i}'(w,x,0)$ for all $D(w)\subset U^{x}$ and all $i\geq 2$. In most cases, $w_{i}'(w,x,0)=w_{i}'(w,0)$. 

{\em{The exceptions}} are when $w$ has a prefix $yL_{2}$ or $y_{1}L_{2}$ with $w_{i}(yL_{2})=y$ or $w_{i}(y_{1}L_{2})=y_{1}$ and  
$$(y,y_{1})=(uxL_{3}^{3},ux_{1}L_{2}R_{3}L_{3}L_{2}){\rm{\ \ or\ \ }}(uxL_{3}^{5}),ux(L_{2}R_{3})^{2}L_{3}){\rm{\ \ or\ \ }}w_{1}(ux(L_{2}R_{3})^{2}L_{3}^{3}L_{2}),ux(L_{2}R_{3})^{2}L_{3}^{3}).$$
In all these cases, $w_{i-1}(yL_{2})=w_{i-1}(y_{1}L_{2})$ and --- as is easily checked --- $w_{i}(yL_{2})=y$ if and only if $w_{i}(y_{1}L_{2})=y_{1}$. We then define
$$w_{i}'(w,x,0)=y_{1}{\rm{\ or\ }}y,$$
depending on whether $yL_{2}$ or $y_{1}L_{2}$ is a prefix of $w$.

\section{Properties with prefixes}\label{7.10}

The definition of the sequence $w_{j}(w,0)$, and of the other sequences, is obviously rather involved. But a general principle which seems to hold is that definitions depend only on certain suffixes of $w$, although the length of these suffices might be rather long. The most immediate properties are the following.

If 
$i>1$, then $w_{i-1}(w,0)$ is a prefix of $w_{i}'(w,0)$. Also, 
$$w_{j}'(w_{i}'(w,0))=w_{j}'(w,0){\rm{\ for\ }}j\leq i,$$ 
$$w_{j}(w_{i}'(w,0))=w_{j}(w,0){\rm{\ for\ }}j<i.$$ 

Suppose that $w=uv$ where the following hold:
\begin{itemize}
\item $D(w)\subset U^{0}$ and $D(v)\subset U^{0}$; 
\item $w_{i}(w,0)$ is a prefix of $u$, possibly equal to $u$, and $w_{i+1}(w,0)$ is not a prefix of $u$ (and not equal to $u$);
\item The word $v$ does not start in the middle of a word $v_{1}\cdots v_{n}L_{3}$ as in 5 of the definition of $w_{1}(.,0)$ of \ref{7.6}.
\end{itemize}

Then the following hold:
\begin{equation}\label{7.6.1}w_{j}(w,0)=uw_{j-i}(v,0){\rm{\ for\ }}j>i;\end{equation}

\begin{equation}\label{7.6.2}w_{j}'(w,0)=uw_{j-i}'(v,0){\rm{\ for\   
}}j>i.
\end{equation}
 If $w_{1}'(w,0)=w_{1}'(v,0)=x$, then we have
 \begin{equation}\label{7.6.4}w_{j}'(w,x,0)=uw_{j-i}'(v,x,0){\rm{\ for\   
}}j>i.
\end{equation}

Now suppose that $w_{1}'(w,0)=x=x_{0}L_{3}$, still with $w=uv$ as above. Suppose in addition that  $v$ contains no subword $y$ such that 
$$D(y)\subset ((U^{x}\setminus x_{0}U^{0})\cup (x_{0}U^{0}\setminus U^{x})).$$
Then
\begin{equation}\label{7.6.6}w_{j}'(w,x,0)=uw_{j-i}'(v,0){\rm{\ for\   
}}j>i.
\end{equation}

 If $v$ starts with $UC$, then we have somewhat similar statements to these. Let $v'$ be obtained from $v$ by replacing the first letter by $BC$, and simply define $w_{j-i}(v,0)$ and $w_{j-i}'(v,0)$ to be obtained from $w_{j-i}(v',0)$ and $w_{j-i}'(v',0)$ by replacing the first letter by $UC$. Then the statements as above hold, except that we should replace $j-i$ by $j-i+1$.
 
 In a slightly different vein, we have the following properties as regards inserting or deleting strings $R_{1}^{n}$. Suppose that $w=u_{1}R_{1}^{n}v_{1}$ for $n\geq 1$, and $w_{j}(w,0)=u_{1}R_{1}^{n}v_{2})$. Then, for any $m\geq 1$,
 $$w_{j}(u_{1}R_{1}^{m}v_{1},0)=u_{1}R_{1}^{m}v_{2}.$$
Similar properties hold for $w_{j}'(w,0)$ and $w_{j}'(u_{1}R_{1}^{m}v_{1},0)$, and also for $w_{j}'(w,x,0)$ and $w_{j}'(u_{1}R_{1}^{m}v_{1},0)$, if $x=w_{1}'(w,0)$ is a prefix of $u_{1}$. 

\section{Definition of the paths $\beta (w,x,0)$}\label{7.7}

Now we define a path $\beta (w,x,0)$ for each word $w$ ending in $C$ 
with $w_{1}'(w,0)=x$ --- that is, with 
$D(w)\subset U^{x}$ --- and each $x\in W_{1}(0)$. The set $R_{m,0}$ is then the set 
of all such paths for $x\in W_{1}(0)$ and $w$ ending in $L_{2}C$ of length $m+1$, or in 
$R_{2}C$ of length $m+2$.  
So fix $x$, and $w$ with $D(w)\subset U^{x}$. We use the sequence 
$w_{i}'(w,x,0)$ of \ref{7.9}.  
Let $n$ be the least integer with $w_{n}'(w,x,0)=w_{n+1}'(w,x,0)$. For $i\leq n$,  the word $w_{i}'(w,x,0)$ determines the $i$'th 
crossing by $\beta $ of $\{ z:\vert z\vert \leq 1\} $. The word $w_{n}'(w,x,0)$ determines the last 
$S^{1}$-crossing of $\beta (w,x,0)$. 

Recall that the last letter of $w_{i}'$ is either $L_{3}$ or $R_{3}$ or $BC$, where the last possibility only occurs in case 3 of \ref{7.6}. First suppose that the last letter of $w_{i}$ is $L_{3}$ or $R_{3}$. As in \ref{7.9}, let $\ell (w_{i}')$ be the unique leaf labelled by the word $w_{i}'(L_{2}R_{3}L_{3})^{\infty }$. This leaf is in the boundary of a unique gap of $L_{3/7}$, which we call $G(w_{i}')$. 
The $i$'th 
 crossing by $\beta (w,x,0)$ of the unit disc is along $\ell (w_{i}')$, directed in an anticlockwise direction along the boundary of $G(w_{i}')$.   The definition of $\ell (w_{i}')$ is made so that, if $\ell (w_{i}')$ is a vertical leaf, then it is on the left of $G(w_{i}')$ if $i$ is odd, and on the right if $i$ is even, so that an anticlockwise direction around $G(w_{i}')$ is always possible. If $i=n$, then we do not have a complete 
crossing of the disc, but just travel halfway along $\ell (w_{n}')$ before entering $G(w_{n}')$, and ending at $G(w_{n}')\cap Z_{m}(s)$.

Now suppose that the last letter of $w_{i}'$ is $BC$.  Let $u'$ be the word obtained from $w_{i}'$ by replacing the last letter by $UC$. Then $i>1$, and  $u'$ and $w_{i}'$ both have $w_{i-1}'$ as a prefix, while $u'$ is a prefix of $w_{j}'$ for $i<j\leq n$.
The words $u'$ and $w_{i}'$ then label disjoint 
subsets $\partial D(w_{i})$ and $\partial D(w_{i}')$ of $\{ z:\vert 
z\vert \leq 1\} $. 
The $i$'th crossing by $\beta 
(w,x,0)$ of $\{ z:\vert z\vert \leq 1\} $ is then from the first 
point of 
$\partial D(w_{i}')$ to the first point of  $\partial D(u')$, 
using anticlockwise order on these intervals.

If $w_{i}=w_{i}(w,0)$ ends in $L_{3}$  
and there is no $BC$ or $UC$ strictly between $w_{i}$ and $w_{i+1}$, 
then
there is also no $BC$ or $UC$  strictly between $w_{i}'$ and 
$w_{i+1}'$, and 
$w_{i}'$ is a prefix of $w$ and of $w_{i+1}'$, all letters strictly 
between 
$w_{i}'$ and $w_{i+1}'$ are $L_{2}$, $L_{3}$ or $R_{3}$, and the 
number of $L_{i}$ is even or odd, depending on whether $w_{i+1}'$ 
ends 
in $L_{3}$ or $R_{3}$. Then the distance along $S^{1}$ 
between the end of the $i$'th unit-disc-crossing of $\beta (w)$ 
(determined by $w_{i}'$) and the start of $i+1$'st crossing 
(determined by 
$w_{i+1}'$) is shorter than the distance between the ends of the 
$i$'th and $i+1$'st disc crossings. This is not always true if $w_{i}$ ends 
in $BC$ or $UC$. If $w$ starts with 
$w_{i}L_{1}R_{1}^{2}R_{2}R_{3}L_{3}L_{2}$, 
then the distance along $S^{1}$ between the end of the $i$'th 
disc-crossing of 
$\beta (w,x,0)$ 
 and the start of $i+1$'st crossing is longer than the distance 
 between the ends of the two crossings.

The first $S^{1}$ crossing of 
any path $\beta \in R_{m,0}$, for any $m$, is always 
between $e^{2\pi i(2/7)}$ and $e^{2\pi i(9/28)}$, because if 
$w_{1}(w)$ 
ends in $BC$ or $UC$, the set 
$D(w_{1}'(w))$ is in 
the part of the disc bounded by leaves with these endpoints, and if 
$w_{1}(w)$ has no occurrence of $BC$ or $UC$ then $w_{1}'(w)$ has an 
even or odd number of $L_{i}$, depending on whether it ends in $L_{3}$ 
or $R_{3}$.

We saw in \ref{7.9} that the number of $w$ of preperiod $m$ with $w_{1}'(w,0)=x$ coincides with the number of $w$ of preperiod $m$ with $D(w)\subset U^{x}$. The function $w_{1}'(w,0)$ is defined for all finite words $w$ with $D(w)\subset U^{0}$. Therefore, the number of paths in $R_{m,0}$ coincides with the number of points in $U^{0}$. We claim, but do not prove, that the number of paths in $R_{m,p}$ coincides with the number of points in $U^{p}\cap Z_{m}$. In fact, we have not even defined the set of paths $R_{m,p}$, since $w_{i}(w,p)$, $w_{i}'(w,p)$ and $w_{i}'(w,x,p)$ have not been defined for $D(x)\subset U^{p}$. But assuming that this is true, this is consistent with the counts made in \ref{3.5}.

For any path $\beta \in \cup _{m}R_{m,0}$, we write $w(\beta )$ for the word ending in $C$ encoding the second endpoint of $\beta $, and $w_{1}'(\beta )$ for the word ending in $L_{3}$ or $BC$ encoding the first $S^{1}$-crossing of $\beta $. 
We also write $w_{i}'(\beta )$ for the word encoding the $i$'th $S^{1}$-crossing. Thus, for $\beta \in \cup _{m}R_{m,0}$,
$$w_{i}'(\beta )=w_{i}'(w(\beta ),w_{1}'(\beta ),0).$$
We also define
$$w_{i}(\beta )=w_{i}(w(\beta ),0).$$

\section{Path-crossings and adjacency: preservation under inverse images}\label{7.11}

We obtain some information about preservation of path segments and of adjacency of paths under inverse images of $s^{n}$ from \ref{7.10}. A word $u$ describes a local inverse by 
$$SD(v)=D(uv)$$
whenever $uv$ is admissible. We fix a domain of $u$ which is contained in either $D(BC)$ or $D(L_{3})$, and such that the domain does not include any $D(v)$ such that  the end of $u$ starts earlier than the last two letters of a word as in 5 of \ref{7.6}. 

For a local inverse $S$ determined by $u$ satisfying the conditions of (\ref{7.6.4}), or with domain restricted to satisfy the conditions of (\ref{7.6.6}), the identities of \ref{7.10}, and the definitions of paths in this section, imply the following. If $\beta _{1}$ and $\beta _{2}$ are adjacent paths in $R_{m,0}$ with $\overline{\beta _{1}}*\beta _{2}$ in the domain of $S$, then

\begin{equation}\label{7.11.1}S(\overline{\beta _{1}}*\beta _{2})=\overline{\beta _{3}}*\beta _{4}\end{equation}
up to homotopy, where $(\beta _{3},\beta _{4})$ is an adjacent pair in $R_{m,0}$. Conversely, if $(\beta _{3},\beta _{4})$ is an adjacent pair in $R_{m,0}$ in the image of such a local inverse $S$, then there is an adjacent pair $(\beta _{1},\beta _{2})$ in $R_{m,0}$ such that (\ref{7.11.1}) holds.

From the definitions in this section, and the identities in \ref{7.10},  (\ref{7.11.1}) holds for some other local  homeomorphisms $S$ of the form 
$$S_{1}S_{2}s^{n}$$

It holds if $S_{1}$ is of length $n-1$ determined by the letter $UC$, and  $S_{2}$ is the local inverse of $s$ defined by the letter $UC$, and the domain of $s^{n}$ is the image of $S_{1}S_{3}$, where $S_{3}$ is the local inverse of $s$ defined by the letter $BC$. 

It also holds if $S_{1}$ is a local inverse of $s^{n}$ and $S_{2}$ is a local inverse of $s^{t}$  determined by $R_{1}^{t}$ for some $t\geq 2$, and the domain of $s^{n}$ is taken to be ${\rm{Im}}(S_{1})$.

\section{Notes on the definition of paths in $R_{m,p}$ for $p>0$}\label{7.12}

The definition of paths in $R_{m,p}$ for $p>0$ follows the same lines as the definition of paths in $R_{m,0}$. For words $w$ with $D(w)\subset U^{p}$, we define finite sequences $w_{i}(w,p)$ and $w_{i}'(w,p)$, and then $w_{i}'(w,x,p)$, for $x=w_{1}'(w,p)$. The definitions are suggested by the structure of the set $U^{p}$ and the sets $A_{m,q_{p}}$ and $C_{m,q_{p}}$, just as the definitions in the case $p=0$ are suggested by the structure of $U^{0}$ and $A_{m,2/7}$ and $C_{m,2/7}$. 
If $w$ does not start with $L_{3}^{2p+1}L_{2}$, then we define $w_{1}(w,p)=L_{3}^{2p+1}$. If $w$ does start with $L_{3}^{2p+1}L_{2}$, then we look for an occurrence of prefixes of a certain form. These occurrences can be of 1 to 4 as in \ref{7.6}, or an analogue of 5 or \ref{7.6}, or an extra alternative, which we call 6. In the analogue of 5 of \ref{7.6}, the form of the words $v_{i}$ is $L_{3}^{2p+1}(L_{2}R_{3})^{m_{i}}L_{3}^{2p+r_{i}-1}$, and otherwise the conditions are as before. These occurrences arise for parallel reasons to the case $p=0$: strings on sets $A_{m,q_{p}}$ of the form
$$L_{3}^{4p+5}\to L_{3}^{2p-1}L_{2}BCL_{1}R_{2}R_{3}L_{3}^{2p+1},$$
$$L_{3}^{2p-1}L_{2}BCL_{1}R_{2}R_{3}L_{3}^{2p+1}\to L_{3}^{2p+1}L_{2}R_{3}L_{3}^{2p+2},$$
$$L_{3}^{2p-1}L_{2}BCL_{1}R_{2}R_{3}L_{3}^{2p+2}\to L_{3}^{2p+1}(L_{2}R_{3})^{2}L_{3}^{2p+1}.$$
The extra alternative 6 arises from the extra thickness of $C_{2,q_{p}}$ in the case $p>0$. We look for a first occurrence of $xR_{3}L_{2}R_{3}L_{3}^{2t-1}L_{2}$ or $xL_{3}L_{2}R_{3}L_{3}^{2t-1}$  for $t\leq p$, depending on whether the longest suffix of $x$ containing only letters $L_{3}$, $L_{2}$ and $R_{3}$ has an even or odd number of $L$ letters. There is a connection, here, with the definition of the set $U^{p}$, and the inclusion of the sets $S_{2,p,k}D(w_{t})$ in \ref{2.10}.

\section{Final statement  in the hard and 
interesting  case}\label{7.8}
We can now be more precise about the structure of the path sets $R_{m,0}$ and $R_{m,0}'$ --- and hence also of $\Omega _{m,0}$ and $\Omega _{m,0}'$. So the final statement of the Main Theorem in the hard and interesting case is as follows. We have given only sketchy detail of the definition of $R_{m,p}$, and hence also of $\Omega _{m,p}$, for $p\geq 1$. So the detailed statement of the theorem in the case $p\geq 1$ is not contained in this chapter, and we shall prove the detail of the theorem only in the case $p=0$, that is, we shall prove that a piece of the fundamental domain can be chosen as claimed, bounded by leaves whose images under $\rho (.,s_{3/7})$ are $\beta _{2/7}$, $\beta _{5/7}$, $\beta _{9/28}$ and $\beta _{19/28}$.
\begin{maintheorem3}

Recall that $s=s_{3/7}$. We have
$$\Omega _{m}(a_{1},+)=\cup _{p=0}^{\infty }\Omega _{m,p}\cup \Omega 
_{m,p}'\cup \{ \omega _{\infty }\} ,$$
where this is a disjoint union, and
$$R_{m}(a_{1},+)=\cup _{p=0}^{\infty 
}R_{m,p}\cup R_{m,p}'\cup \{ \beta _{1/3}\} ,$$
where
$$R_{m,p}=\rho (\Omega _{m,p},s),\ \ R_{m,p}'=\rho 
(\Omega _{m,p}',s).$$
 
Any adjacent pair in $\Omega _{m}(a_{1},+)$ is an adjacent pair in 
exactly one of the sets $\Omega _{m,p}$ or $\Omega _{m,p}'$. Every 
adjacent pair in $\Omega _{m,p}$ is matched with one adjacent pair in 
$\Omega _{m,p}'$, and with no other. Write $Z_{\infty }=\cup _{n}Z_{n}$. The paths in $R_{m,p}$ and $R_{m,p}'$ are defined as elements of $\pi _{1}(\overline{\mathbb C}\setminus Z_{\infty },v_{2},Z_{\infty })$, with 
$$R_{m,p}\subset R_{m+1,p},{\rm{\ in\ }}\pi _{1}(\overline{\mathbb C}\setminus Z_{\infty },v_{2},Z_{\infty }),$$
$$R_{m,p}'\subset R_{m+1,p}'{\rm{\ in\ }}\pi _{1}(\overline{\mathbb C}\setminus Z_{m},v_{2},Z_{m}).$$

The sets of paths $R_{m,p}$, $R_{n,q}'$ are all 
disjoint, except that, if $m>2p+2$,
 $$R_{m,p}\cap R_{m,p+1}=\{ \beta _{q_{p+1}}\} ,$$
and 
$$R_{m,p}'\cap R_{m,p+1}'=\{ \beta _{1-q_{p+1}}\} .$$ 
Also, 
$\beta _{1/3}\in R_{m,p}\cap R_{m,p}'$, for the largest $p$ such that 
these are nonempty. 

A path $\beta $ in $R_{m,p}$ has first $S^{1}$ crossing at $e^{2\pi ix}$ for some $x\in [q_{p},q_{p+1})$, where $x=w_{1}'(\beta )=w_{1}'(w(\beta ),p)$
and $e^{2\pi ix}$ is an endpoint of a leaf of $L_{3/7}$,  and the first unit-disc crossing by $\beta $ is along this disc. The later crossings are encoded by the words $w_{i}'(w,x,p)$. Thus, up to homotopy, $\beta $ is a union of segments along the upper or lower half of the unit disc, leaves of $L_{3}$ in $U^{p}$, some segments in gaps of $L_{3/7}$ and a final segment in a gap of $L_{3/7}$ ending at a point of $U^{p}\cap Z_{m}$. Also, $\beta $ is determined uniquely by its first $S^{1}$-crossing and its final endpoint, but not always by its final endpoint alone. 

Nevertheless, the paths in $R_{m,p}\setminus \{ \beta _{q_{p}}\} $ are in one-to-one correspondence with the points in $Z_{m}\cap U^{p}$.
 
\end{maintheorem3}

Because of the matching between adjacent pairs in $\Omega _{m,p}$ and $\Omega _{m,p}'$, the part of the tree in $V_{3,m}(a_{1},+)$ associated to 
the fundamental region which has been constructed is an interval. One 
should not read too much into this. There is a lot of choice in 
construction of the fundamental domain. But this case is 
very different from the easy cases. It is certainly not possible to 
choose a fundamental domain which  inherits the structure of part of 
$T_{m}(s_{3/7})$ --- except for $m\leq 2$, as can be seen from the computations in \ref{7.5}.

It will have  noticed that $R_{m,0}$ has been described in exhausting detail, but nothing has been said about $R_{m,0}'$. The set of paths $R_{m,0}'$ will be completely described in the next chapter.

\chapter{Proof of the hard and interesting case}\label{8}

\section{The inductive construction}\label{8.1}

We continue with the notation established at the start of Section \ref{7}. We prove the most detailed part of the theorem only in the case $p=0$, that is, for $R_{m,0}$ and $R_{m,0}'$.

The path sets $\Omega _{m,0}$ and $\Omega _{m,0}'$ are described by the  images $R_{m,0}$ and $R_{m,0}'$,  under $\rho $ in the statement of the Main Theorem \ref{7.8}. An adjacent pair in $\Omega _{m,0}$ or $\Omega _{m,0}'$ corresponds to an adjacent pair in $R_{m,0}$ or $R_{m,0}'$. Matching of an adjacent pair in $\Omega _{m,0}$ with an adjacent pair in $\Omega _{m,0}'$ is viewed by considering the image under $\rho $ and $\Phi _{2}$ of the element of $\pi _{1}(V_{3,m},a_{1})$ effecting the matching, as described in \ref{5.6}, in particular in (\ref{5.6.3}) and (\ref{5.6.4}). So for each adjacent pair $(\beta _{1},\beta _{2})$ in $R_{m,0}$, we shall choose $\alpha _{m}\in \pi _{1}(\overline{\mathbb C}\setminus Z_{m},v_{2})$ and $[\psi _{m}]\in {\rm{MG}}(\overline{\mathbb C},Y_{m})$ so that
\begin{equation}\label{8.1.3}(s,Y_{m})\simeq _{\psi _{m}}(\sigma _{\alpha _{m}}\circ s,Y_{m}),\end{equation}
and, for $i=1$, $2$, 
\begin{equation}\label{8.1.4}\beta _{m,i}=\alpha _{m}*\psi _{m}(\beta _{m,i}'){\rm{\ rel\ }}Y_{m}.\end{equation}
 
 The paths in $R_{m,0}$ have been defined as elements of $\pi _{1}(\overline{\mathbb C}\setminus Z_{\infty },Z_{\infty },v_{2})$, and, if $m\leq n$ we have $R_{m,0}\subset R_{n,0}$ in $\pi _{1}(\overline{\mathbb C}\setminus Z_{\infty },Z_{\infty },v_{2})$. As indicated in \ref{7.8}, the paths in  $R_{m,0}'$ will be defined as elements of $\pi _{1}(\overline{\mathbb C}\setminus Z_{m},Z_{m},v_{2})$, so that, if $m\leq n$,
$$R_{m,0}'\subset R_{n,0}'{\rm{\ in\ }}\pi _{1}(\overline{\mathbb C}\setminus Z_{m},Z_{m},v_{2}).$$
 The total ordering on $R_{m,0}$ is easily seen, by fixing lifts to the universal cover, all with the same first 
endpoint, and using the ordering on the second endpoints. The 
universal cover is the unit disc up to holomorphic equivalence, and 
under this equivalence the lifts of the paths in $R_{m,0}$ 
have endpoints in the unit circle. We take the order on the unit 
circle, with the endpoint of $\beta _{2/7}$ as minimal, and $\beta _{9/28}$ is maximal.   Since the elements of $R_{m,0}'$ are matched with the elements of $R_{m,0}$, they, too,
are totally ordered, by the condition that matching preserves order. 
The minimal and maximal elements of $R_{m,0}'$ are $\beta _{5/7}$ and $\beta _{19/28}$. Similarly, the sets $R_{m,p}$ and $R_{m,p}'$ are totally ordered. The minimal and maximal elements of $R_{m,p}$ are $\beta _{q_{p}}$ and $\beta _{q_{p+1}}$, and the minimal and maximal elements of $R_{m,p}'$ are $\beta _{1-q_{p}}$ and $\beta _{1-q_{p+1}}$.  Since the sets $\{ \beta _{q_{p}}:p\geq 0\} $ and $\{ \beta _{1-q_{p}}:p\geq 0\} $ are also totally ordered, the sets $\cup _{p\geq 0}R_{m,p}$ and $\cup _{p\geq 0}R_{m,p}'$ are also totally ordered

Let $(\beta _{1},\beta _{2})$ be an adjacent pair in $R_{m,0}$, matched with an adjacent pair $(\beta _{1}',\beta _{2}')$ in $R_{m,0}'$. In general, we shall choose $\alpha _{m}\in \pi _{1}(\overline{\mathbb C}\setminus Z_{\infty},v_{2})$ and $\psi _{m}$ representing an element of ${\rm{MG}}(\overline{\mathbb C},Y_{\infty })$ inductively for $m\geq 2$. For each $2\leq k\leq m$, there is an adjacent pair $(\beta _{k,1},\beta _{k,2})$ in $R_{k,0}$ such that $\beta _{1}$ and $\beta _{2}$ are between $\beta _{k,1}$ and $\beta _{k,2}$. The induction starts with $\beta _{2,1}=\beta _{2/7}$ and $\beta _{2,2}=\beta _{9/28}$, and ends with $\beta _{m,1}=\beta _{1}$ and $\beta _{m,2}=\beta _{2}$. Similarly, there is an adjacent pair $(\beta _{k,1}',\beta _{k,2}')$ in $R_{k,0}'$ which is matched with $(\beta _{k,1},\beta _{k,2})$, and such that $\beta _{1}'$ and $\beta _{2}'$ are between $\beta _{k,1}'$ and $\beta _{k,2}'$. For intermediate $k$, $\beta _{k,1}=\beta _{m,1}$ and $\beta _{k,2}=\beta _{m,2}$ are possible. If $\beta _{1}=\beta _{k,1}$ in $\pi _{1}(\overline{\mathbb C}\setminus Z_{k},Z_{k},v_{2})$ then $\beta _{1}=\beta _{k,1}$ in $\pi _{1}(\overline{\mathbb C}\setminus Z_{\infty },Z_{\infty }v_{2})$. The corresponding statement is not quite true for $\beta _{k,1}'$ and $\beta _{1}'$. It can happen that $\beta _{1}'=\beta _{k,1}'$ in $\pi _{1}(\overline{\mathbb C}\setminus Z_{k},Z_{k},v_{2})$ and $\beta _{1}'=\beta _{k+1,1}'$ in $\pi _{1}(\overline{\mathbb C}\setminus Z_{k+1},Z_{k+1},v_{2})$ but $\beta _{k,1}'\neq \beta _{k+1,1}'$ in $\pi _{1}(\overline{\mathbb C}\setminus Z_{k+1},Z_{k+1},v_{2})$. In fact this happens very often, but is not very visible, because of the indirect definition of $R_{m,0}'$.

 For any choice of adjacent pair $(\beta _{1},\beta _{2})$ in $R_{m,0}$, we have
$$\beta _{2,1}=\beta _{2/7},\ \ \beta _{2,2}=\beta _{9/28},$$
and, for $0\leq i\leq 2$, $\alpha _{i}$ is an arbitrarily small perturbation of $\beta _{2/7}*\overline{\beta _{5/7}}$, and
$$\psi _{i}=\psi _{i,2/7},$$
for $\psi _{i,2/7}$ as in \ref{7.3}, which is the identity on $\beta _{5/7}$ for $i\leq 2$. So the loop $\alpha _{0}=\alpha _{1}=\alpha _{2}$ bounds a disc not containing $c_{1}$. The defining equations of $\psi _{m}$ and $\alpha _{m}$ are
\begin{equation}\label{8.2.1}\psi _{m}\circ s=\sigma _{\alpha _{m}}\circ s\circ \psi _{m+1},\end{equation}
and $\alpha _{m}$ is homotopic to an arbitrarily small perturbation of $\beta _{m}*\psi _{m}(\overline{\beta _{m,1}})$ in $\pi _{1}(\overline{\mathbb C}\setminus Z_{m},v_{2})$. We then have
\begin{equation}\alpha _{m}*\psi _{m+1}(\beta _{m,i}')=\beta _{m,i}{\rm{\ in\ MG}}(\overline{\mathbb C}\setminus Z_{m},Z_{m},v_{2}).\end{equation}
and
\begin{equation}(s,Y_{m+1})\simeq _{\psi _{m+1}}(\sigma _{\alpha _{m}}\circ s,Y_{m+1}).\end{equation}

We  also define $[\xi _{m}]\in {\rm{MG}}(\overline{\mathbb C},Y_{\infty })$ by
\begin{equation}\label{8.2.4}\psi _{m+1}=\xi _{m}\circ \psi _{m}.\end{equation}

An example of this inductive construction is given (for general p) in \ref{7.3}, for the adjacent pair $(\beta _{q_{p}},\beta _{q_{p+1}})$ in $R_{p,2p+2}$. We saw in \ref{7.3} that the homeomorphism $\psi $ effecting the matching between $(\beta _{q_{p}},\beta _{q_{p+1}})$ and $(\beta _{1-q_{p}},\beta _{1-q_{p+1}})$ satisfies
\begin{equation}
\psi (\beta _{1-q_{p}})=\beta _{q_{p}}{\rm{\ in\ }}\pi _{1}(\overline{\mathbb C}\setminus Z_{2p+2},Z_{2p+2},v_{2})\end{equation}
An arbitrarily small perturbation of $\beta _{q_{p}}*\psi (\beta _{1-q_{p}})$ is freely homotopic, in $\overline{\mathbb C}\setminus Z_{2p+1}$, to a simple closed loop which intersects the unit disc in the leaf of $L_{3/7}$ joining the points $e^{\pm 2\pi i(1/3)}$. 

Now we introduce another object related to the adjacent pair $(\beta _{m,1},\beta _{m,2})$. We denote the set bounded by $\overline{\beta _{m,1}}*\beta _{m,2}$ and $\psi _{m}(\overline{\beta _{m,1}'}*\beta _{m,2}')$ by $U(\beta _{m,1},\beta _{m,2})$. These two paths with endpoints in $Z_{m}$ are then defined as elements of $\pi _{1}(\overline{\mathbb C}\setminus Z_{\infty },Z_{\infty })$.  The boundary of this set is a union of two arcs: the {\em{solid boundary}} $\overline{\beta _{m,1}}*\beta _{m,2}$ and the {\em{dashed boundary}} $\psi _{m}(\overline{\beta _{m,1}'}*\beta _{m,2}')$. We denote the set bounded by $\overline{\beta _{m,1}}*\beta _{m,2}$ and $\psi _{m+1}(\overline{\beta _{m,1}'}*\beta _{m,2}')$ by $V(\beta _{m,1},\beta _{m,2},m+1)$, and, once again, the boundary is the union of two arcs: the {\em{solid boundary}} $\overline{\beta _{m,1}}*\beta _{m,2}$ and the {\em{dashed boundary}} $\psi _{m+1}(\overline{\beta _{m,1}'}*\beta _{m,2}')$. The definitions then give
$$\partial 'V(\beta _{m,1},\beta _{m,2})=\psi _{m+1}(\overline{\beta _{m,1}'}*\beta _{m,2}')=\xi _{m}(\partial 'U(\beta _{m,1},\beta _{m,2}))$$
We make definitions so that
\begin{equation}\label{8.1.5}\begin{array}{l}U(\beta _{1},\beta _{2}){\rm{\ is\ a\ topological\ disc,}}\\ V(\beta _{1},\beta _{2},m+1){\rm{\ is\ also\ a\ disc}}\end{array}\end{equation}

 In general, for $n\geq m$ for an adjacent pair $(\beta _{1},\beta _{2})$ in $R_{m,p}$, we write $V(\beta _{1},\beta _{2},n)$ for the union of sets $U(\beta _{3},\beta _{4})$ for adjacent pairs $(\beta _{3},\beta _{4})$ in $R_{n,p}$ between $\beta _{1}$ and $\beta _{2}$. We shall see that the two definitions coincide for $n=m+1$, up to isotopy preserving $Z_{m+1}$.

Then it is clear that an inductive definition of $R_{m+1,0}$ and $R_{m+1,0}'$ is given by defining a subdivision of the sets $V(\beta _{1},\beta _{2},m+1)$ into sets $U(\beta _{3},\beta _{4})$, for all adjacent pairs $(\beta _{1},\beta _{2})$ of $R_{m,0}$, thus defining all adjacent pairs $(\beta _{3},\beta _{4})$ in $R_{m+1,0}$. We caution that there will be cases when the set $V(\beta _{1},\beta _{2},m+1)$ is only a union of sets $U(\beta _{3},\beta _{4})$ up to homotopy in $Z_{m}$. In these cases the definition of the elements $\beta _{1}'$, $\beta _{2}'$ matched with $\beta _{1}$ and $\beta _{2}$ will be different, as elements of $\pi _{1}(\overline{\mathbb C}\setminus Z_{\infty },v_{2},Z_{\infty })$, in $R_{m,0}$ and in $R_{m+1,0}$.

In order to prove the Main Theorem \ref{7.8}, 
we have to prove two subsidiary theorems. We state the first theorem only in the case $p=0$, and we shall only prove it only  in this case. However, there is also a version in the case of $p>0$.

\begin{theorem}\label{8.2}
The sets $U(\beta _{1},\beta _{2})$ for adjacent pairs $(\beta _{1},\beta _{2})$ can be defined so that the following holds. 
Fix any $x\in W_{1}(0)$. Let $\beta _{1}\in R_{m,0}$ be minimal with $w_{1}'(\beta _{1})=x$, and $m$ minimal with $\beta _{1}\in R_{m,0}$. Let $\beta _{1}<\beta _{n,2}$ be adjacent in $R_{n,0}$ for any $n\geq m$. Let $\xi _{n,1}$ be the homeomorphism associated to the pair $(\beta _{1},\beta _{n,2})$, like the sequence $\xi _{n}$ in (\ref{8.2.4}). Then for any adjacent pair $(\beta _{3},\beta _{4})$ in $R_{n,0}$ with $w_{1}'(\beta _{n,3})=x$ and associated homeomorphism $\xi _{n}$ as in (\ref{8.2.4}), up to  $Z_{n+1}$-preserving isotopy, the region bounded by $\overline{\beta _{3}}*\beta _{4}$ and $\xi _{n,1}(\partial 'U(\beta _{3},\beta _{4}))$ is, up to isotopy preserving $Z_{n+1}$:
\begin{equation}\label{8.2.5}\cup \{ U(\beta _{5},\beta _{6}):(\beta _{5},\beta _{6}){\rm{\ adjacent\ in\ }}R_{n+1,0},\ \beta _{3}\leq \beta _{5}<\beta _{6}\leq \beta _{4}\} .\end{equation}
Also,
\begin{equation}\label{8.2.6}{\rm{supp}}(\xi _{n}\circ \xi _{n,1}^{-1})\cap \xi _{n,1}(\partial 'U(\beta _{3},\beta _{4}))=\emptyset .\end{equation}
Also, for any $x$ of the form $w_{2}'(w)$ for some $w$ starting with $BC$, a similar result holds for the minimal $\beta _{1}$ with $w_{2}'(\beta _{1})$ =$x$ and any pair $(\beta _{3},\beta _{4})$ in $R_{n,0}$ with $w_{2}(\beta _{3})=x$.
\end{theorem}

\begin{theorem}\label{8.3}
The paths of $R_{m,p}$ and $R_{m,p}'$ are in the set $D'$ of 
\ref{3.3} for all $m\geq 1$ and all $p\geq 0$.\end{theorem}
That is,  the obverse of the
Inadmissibility Criterion of \ref{3.4} holds. This will be proved  in \ref{8.21}. 

\section{Adjacent elements in $R_{m,0}$: how many 
$S^{1}$-crossings are the same?}\label{8.4}

We have seen that a path $\beta $ in $R_{m,0}$ is determined by a sequence 
$$(x,\{ w_{i}'(w,x,0) :2\leq i\leq n\} ).$$
We write $\beta =\beta (w,x)$ if $w$ ends in $C$ and $\beta $ ends at the point of $Z_{m}$ encoded by this word.
Because of the rules on the definitions of $w_{1}'(w,0)$, $w_{i}(w,x,0)$ and 
$w_{i}'(w,x,0)$, the paths adjacent to $\beta (w)\in R_{m,0}$ are 
determined by the end part of the word. The examples calculated in 
\ref{7.5} are therefore highly indicative of the general picture.  The following holds. 

\begin{ulemma}
  Let $\beta _{1}$ and $\beta _{2}$ be 
adjacent in $R_{m,0}$ with $\beta _{1}<\beta _{2}$. 
 If $w_{1}'(\beta _{1})\neq w_{1}'(\beta _{2})$, define $i=0$. Otherwise, let $i$ be such that $w_{j}'(\beta _{1})=w_{j}'(\beta _{2})$ 
for $j\leq i$  but $w_{i+1}'(\beta _{1})\neq 
w_{i+1}'(\beta _{2})$. 
\begin{itemize}
\item[Case 1.] Suppose that it is not the case that one of $w_{i+1}(\beta _{1})$ and  $w_{i+1}(\beta _{2})$ is of the form $vL_{1}R_{2}UC$ and the other of the form $vL_{1}R_{2}(BCL_{1}R_{2})y$. 
\begin{itemize}
\item[a)]If $w_{i}'(\beta _{1})$ ends in $L_{3}$, then  $w_{i+2}(\beta _{1})=w(\beta _{1})$ and  
$w_{i+1}(\beta _{2})=w(\beta _{2})$. 
\item[b)]If $w_{i}'(\beta _{1})$ ends in 
$R_{3}$, and $w_{i+1}(\beta _{1})\neq w(\beta _{1})$, 
then $w_{i+3}(\beta _{1})=w(\beta _{1})$ and $w_{i+1}(\beta _{2})=w(\beta _{2})$.
\end{itemize}
\item[Case 2] Suppose that $w_{i+1}(\beta _{1})$ and  $w_{i+1}(\beta _{2})$ are of the form excluded in Case 1. Then $w_{i+1}(\beta _{1})=w(\beta _{1})$, 
$w_{i+1}(\beta _{2})\neq w(\beta _{2})$ and $w_{i+2}(\beta _{2})=w(\beta _{2})$. 
\end{itemize}

\end{ulemma}

\noindent{\em{Proof.}} 

\par \noindent {\em{Case 1.}}  We are supposing that 
$w_{i+1}(\beta _{2})$ is not of the form $vL_{1}R_{2}UC$. Then 
$w_{i+1}'(\beta _{2})$ ends in $L_{3}$ or $R_{3}$, and there is $z$ of the 
same length as $w_{i+1}'(\beta _{2})$ ending in $C$, with $D(z)$ and 
$D(w_{i+1}'(\beta _{2}))$ adjacent. Let $\beta _{3}$ be the path with $w_{i+1}'(\beta _{3})=w_{i+1}'(\beta _{2})$ and  
  $w_{i+1}(\beta _{3})=z$.  Then since 
$w_{i+1}(\beta _{1})\neq w_{i+1}(\beta _{2})$, we have $\beta 
_{1}\leq \beta _{3}\leq \beta _{2}$, and $\beta _{3}=\beta _{2}$. So  $w_{i+1}(\beta _{2})=w(\beta _{2})$. 

\noindent {\em{Case 1.1: $w_{i}(\beta _{1})\neq w_{i}(\beta _{2})$.}} 
 
 If $i=0$ and $w_{1}(\beta _{2})=w(\beta _{2})$, then $w_{1}'(\beta _{1})\neq w_{1}'(\beta _{2})$, and so $\beta _{1}$ is maximal in $R_{m,0}$ with $w_{1}'(\beta _{1})=x$ (for some $x$). We return to this case shortly.  

 Suppose that $i>0$ and that $u=w_{i}(\beta _{1})\neq 
w_{i}(\beta _{2})=u'$. The possibilities for 
$u=w_{i}(\beta _{1})$ given $w_{i}'(\beta _{1})$ are given in \ref{7.6}: in all cases, $w_{i}'$ is obtained from $w_{i}$ either by deleting an end part of $w_{i}$, or, if $w_{i}$ ends in $C$ or $BC$ or $UC$, then $w_{i}$ and $w_{i}'$ have the same length and the discs $D(w_{i})$ and $D(w_{i}')$ are adjacent. At any rate, the possibilities for $w_{i}$ all end in $L_{3}$, $C$, $BC$ or 
$UC$, and these are ordered in terms of the order 
of the sets $D(u)$, with $u$ being minimal for the unique 
such word ending in $C$. If $u$ ends in $C$, then $u'$ ends in $BC$ or $UC$ --- depending on whether $w_{i}'(\beta _{1})=w_{i}'(\beta _{2})$ ends in $L_{3}$ or $R_{3}$ --- and $\beta _{2}$ must be minimal among $\beta \in R_{m,0}$ with $w_{i}'(\beta )=w_{i}'(\beta _{1})=w_{i}'(\beta _{2})$ and $w(\beta )$ of the form $u'y$. If $u$ ends in $BC$ or $UC$, then $\beta _{1}$ must be maximal among 
$\beta \in R_{m,0}$ with $w(\beta )$ of the form $uy$ and $w_{i}'(\beta )=w_{i}'(\beta _{1})$. We shall return to these cases  
shortly. If $u=w_{i}(\beta _{1})$ ends in 
$L_{3}$, for example, $u=L_{3}^{4}$, then $\beta _{1}$ must be maximal among $\beta \in R_{m,0}$ with $w_{i}(\beta )=u$ and $w_{i}'(\beta )=w_{i}'(\beta _{1})$. We then have 
$w_{i+1}'(\beta _{1})=uL_{3}^{2}$, and $\beta _{1}$ is maximal among $\beta 
\in R_{m,0}$ with $w_{i+1}'(\beta )=uL_{3}^{2}$ and $w_{i}'(\beta )=w_{i}'(\beta _{1})$. Again, we return to 
this case shortly.

\noindent {\em{Case 1.2: $w_{i}(\beta _{1})=w_{i}(\beta _{2})$.}} 
Suppose that
$u'=w_{i+1}'(\beta _{1})\neq w_{i+1}'(\beta _{2})=w(\beta _{2})$. 
Then $\beta _{1}$ is maximal among $\beta \in R_{m,0}$ 
with $w_{i+1}'(\beta )=u'$ and $w_{i}'(\beta )=w_{i}'(\beta _{1})$. This is simply a more general version of the second maximal property of Case 1.1. We return to it below.

\noindent {\em{Case 2.}} Suppose that  
$w_{i+1}(\beta _{2})= vL_{1}R_{2}UC$ for some $v$.  Then since 
$w_{i+1}(\beta _{1})\neq w_{i+1}(\beta _{2})$, we must have $w(\beta _{1})=v'$, where $v'$ is obtained by replacing the last letter of $v$ by $C$, or 
$w(\beta _{1})=vL_{1}R_{2}BCL_{1}R_{2}BCz$ for some $z$, depending on the length of 
$v$ in comparison with $m$, and $\beta _{1}$ must be 
maximal in $R_{m,0}$ among $\beta $ with $vL_{1}R_{2}BCL_{1}R_{2}BC$ as a prefix of $w(\beta )$. This is the same property as 
the first maximal property in Case 1.1.

So suppose that $u$ ends in $BC$ or $UC$. 
 Maximal paths $\beta \in R_{m,0}$  with $w(\beta )$ of the form $uy$  
have  
 $w(\beta )=uL_{1}R_{1}^{p}R_{2}C$ for 
$1\leq p\leq 3$, or $uL_{1}R_{2}R_{3}L_{3}L_{2}C$. Minimal paths with $w(\beta )$ of the form $uy$ have 
$w(\beta )=u(L_{1}R_{2}BC)^{n}L_{1}R_{1}R_{2}C$ for some $n\geq 0$.

Now we return to the other possibilities in Case 1. Words $v$ such that $\beta 
\in R_{m,0}$ is maximal with respect to the property that $w_{i+1}'(\beta )=u$,  
with $u$ ending in $L_{3}$ or $R_{3}$ 
are, respectively, $u''C$ of the same length as $u$
with $D(u''C)$ and $D(u)$ adjacent, then, for increasing $m$, $v=u''XL_{1}R_{1}R_{2}C$ 
where 
$X=BC$ or $UC$, depending on whether $u$ ends in $L_{3}$ 
or 
$R_{3}$, then 
$uL_{3}^{3}L_{2}C$, 
then $uL_{2}R_{3}L_{3}^{2}L_{2}C$, then different possibilities 
depending on whether $u$ ends in $L_{3}$ or $R_{3}$. If $u$ ends in 
$L_{3}$, then the additional possibilities, for increasing $m$, are of the form 
$$u(L_{2}R_{3}L_{3})^{q}L_{2}y,$$
$$y=BCL_{1}R_{1}R_{2}C{\rm{\ or\ }}
R_{3}L_{3}^{4}L_{2}C{\rm{\ or\ }}R_{3}L_{3}L_{2}R_{3}L_{3}^{2}L_{2}C$$
for any $q\geq 0$. If 
$u$ ends in $R_{3}$ then the additional possibiliities  are 
$$uL_{2}R_{3}L_{3}^{2}(L_{2}R_{3}L_{3})^{q}(L_{2}{\rm{\ or\ }}L_{2}R_{3}L_{2})C$$
 for any $q\geq 0$.

In all but the last case, $w_{i+2}(\beta _{1})=\beta _{1}$, and in the last case,  
$w_{i+3}(\beta _{1})=\beta _{1}$. So $w_{i+2}(\beta _{1})=\beta _{1}$ if $w(\beta _{1})$ ends in $L_{3}$ and 
$w_{i+3}(\beta _{1})=\beta _{1}$ if $w$ ends in $R_{3}$.

 \Box
 
 \section{Definition of $\xi_{m,0}$ and $\xi _{m,1}$}\label{8.5}
We define $k_{0}$, and sequences $\psi _{m,0}$, $\alpha _{m,0}$ and $\xi _{m,0}$  associated to an adjacent pair $(\gamma _{1},\gamma _{2})$ in $R_{m,0}$, as follows. If $w(\gamma _{1})$ starts with $BC$, then $\psi _{m,0}=\psi _{m}$, $\alpha _{m,0}=\alpha _{m}$ and $\xi _{m,0}=\xi _{m}$, where these are the usual sequences associated to $(\gamma _{1},\gamma _{2})$. Otherwise, let $k_{0}$ be the largest index such that $w(\beta _{k_{0}})$ starts with $BC$, where $(\beta _{k},\beta _{k,2})$ is the adjacent pair in $R_{k,0}$ with $\gamma _{1}$ and $\gamma _{2}$ between $\beta _{k}$ and $\beta _{k,2}$. Also, let $\beta _{k_{0}}'$ be the element of $R_{k_{0}}'$ which is matched with $\beta _{k_{0}}$. Then $k_{0}=3t$ for some $t\geq 1$ and $w(\beta _{k_{0}})=(BCL_{1}R_{2})^{t-1}BCL_{1}R_{1}R_{2}C$. Then we define 
$$\psi _{k_{0},0}=\psi _{k_{0}},\ \ \alpha _{k_{0},0}=\alpha _{k_{0}},\ \ \xi _{k_{0},0}=\xi _{k_{0}}.$$
Then for $k\geq k_{0}$, we define $\psi _{k+1,0}$ and $\xi _{k,0}$ inductively by
$$\sigma _{\alpha _{k,0}}\circ s\circ \psi _{k+1,0}=\psi _{k,0}\circ s.$$
and 
$$\psi _{k+1,0}=\xi _{k,0}\circ \psi _{k,0},$$
and $\alpha _{k+1,0}$ is the perturbation of $\beta _{k_{0}}*\psi _{k+1,0}(\overline{\beta _{k_{0}}'})$ which is equal to $\alpha _{k,0}$ in $\pi _{1}(\overline{\mathbb C}\setminus Z_{k},v_{2})$.
The support of $\xi _{m,0}$ is relatively easy to compute directly. It can be written as 
$$A_{m,2/7}\cup C_{m},$$
where $C_{m}$ is disjoint from $A_{m,2/7}$, coincides with $C_{m,2/7}$ for $k\leq k_{0}$ and only differs from $C_{m,2/7}$ for $m>k_{0}$ in preimages of the periodic component of $C_{k_{0},2/7}$. The difference happens because components of the backward orbit of the periodic component of $C_{k_{0},2/7}$ which intersect the central gap now pull back homeomorphically, not with degree two, and pieces of string which are close to the central gap, pull back close to the preimage gap. 

Similarly, we define sequences $\psi _{m,1}$, $\xi _{m,1}$ and $\alpha _{m,1}$, using $\beta _{k_{1}}$ for some $k_{1}\geq k_{0}$. If $w_{1}'(\beta _{1})=L_{3}$, then $k_{1}=k_{0}$. Otherwise, $k_{1}$ is the least index such that $w_{1}'(\beta _{k_{1}})=w_{1}'(\beta _{1})$.

We shall use \ref{8.8} to analyse $\xi _{m,1}\circ \xi _{m,0}^{-1}$ and $\xi _{m}\circ \xi _{m,1}^{-1}$.

\section{Comparing $\xi _{n}$ and $\xi _{n,0}$ and $\xi _{n,1}$}\label{8.8}

To prove \ref{8.2}, we need  to compare  a sequence of homeomorphisms $\{ \xi _{n}:n\geq m\} $ (as in \ref{8.2.4})  associated to an element $\beta \in R_{m,0}$ with $\{ \xi _{n,0}:n\geq m\} $.  We shall relate $\xi _{n}\circ \xi _{n,0}^{-1}$ to $\xi _{n-1}\circ \xi _{n-1,0}^{-1}$. So let $(\beta _{n},\beta _{n,2})$ be an adjacent pair in $R_{n,0}$, such that:
\begin{itemize}
\item $\beta _{n}'\in R_{n,0}$ is matched with $\beta _{n}$, 
\item $\beta _{n}$ and $\beta _{n,2}$ are between $\beta _{n-1}$ and $\beta _{n-1,2}$, 
\item $\beta _{0}=\beta _{2/7}$, $\alpha _{0}=\alpha _{0,2/7}$, $\psi _{0}=\psi _{0,2/7}$,
\item $\alpha _{n}$  is homotopic, up to homotopy preserving $n$, to a perturbation of $\beta _{n}*\psi (\overline{\beta _{n}'})$ which separates the second endpoint of $\beta _{n}$ from the second endpoint of $\beta _{n,2}$.
\item As in (\ref{8.2.1}) and (\ref{8.2.5}),
\begin{equation}\label{8.8.1}\sigma _{\alpha _{n}}\circ s\circ \psi _{n+1}=\psi _{n}\circ s,\end{equation}
\begin{equation}\label{8.8.2}\psi _{n+1}=\xi _{n}\circ \psi _{n},\end{equation}
\begin{equation}\label{8.8.3}
\sigma _{\alpha _{n}}\circ s\circ \xi _{n}=\xi _{n-1}\circ \sigma _{\alpha _{n-1}}\circ s.\end{equation}
The {\em{support}} of a homeomorphism $\xi $ is the set where it is not the identity, written ${\rm{supp}}(\xi )$.
\end{itemize}

\begin{ulemma} The following holds.

Let $n$ be such that $\psi _{n}=\psi _{n,0}$ and $\alpha _{n-1}=\alpha _{n-1,0}$. Then
$$s\circ \xi _{n}\circ \xi _{n,0}^{-1}=\sigma _{\alpha _{n}}^{-1}\circ \sigma _{\alpha _{n,0}}\circ s.$$
So in this case, 
\begin{equation}\label{8.8.4}{\rm{supp}}(\xi _{n}\circ \xi _{n,0}^{-1})\subset s^{-1}(E_{n}),\end{equation}
where $E_{n}$ is any disc which contains $\overline{\alpha _{n}}*\alpha _{n,0}$.

For any other $n$,
\begin{equation}\label{8.8.10}\sigma _{\beta _{n}}\circ s\circ \xi _{n}\circ \xi _{n,0}^{-1}=\xi _{n-1}\circ \xi _{n-1,0}^{-1}\circ \sigma _{\zeta _{n}}\circ \sigma _{\beta _{n}}\circ s,\end{equation}
where 
$$\zeta _{n}=\zeta _{n,1}*\zeta _{n,2},$$
and $\zeta _{n,1}$ and $\zeta _{n,2}$ are arbitrarily small perturbations of, respectively,  
$$\overline{\beta _{n}}*\beta _{k_{0}}*\xi _{n-1,0}(\overline{\beta _{k_{0}}}*\beta _{n})$$
and
$$\xi _{n-1,0}(\overline {\beta _{n}}*\beta _{n-1}*\psi _{n-1}(\overline{\beta _{n-1}'}*\beta _{n}'))$$
 where $\zeta _{n,1}$ does not contain $v_{2}$ or the endpoints of either $\beta _{n}$ or $\beta _{k_{0}}$, and $\zeta _{n,2}$ contains $v_{2}$ and, if $\beta _{n}\neq \beta _{n-1}$, also contains the endpoint of $\beta _{n}$, but not of   $\beta _{n-1}$. (Otherwise, it does not contain the common endpoint.) So
$${\rm{supp}}(\xi _{n}\circ \xi _{n,0}^{-1})\subset \cup _{k<n}E_{n,k},$$
where $E_{k,k}$ is any suitably chosen disc which contains $\zeta _{k}$, $E_{k+1,k}=s^{-1}(E_{k,k})$, and for $r>k+1$, $E_{r,k}$ is defined inductively as
$$E_{r+1,k}=(\sigma _{\beta _{r}}\circ s)^{-1}(E_{r,k}).$$

Similar statements hold for $\xi _{n,1}$ replacing $\xi _{n,0}$ for $n\geq k_{1}$, and for $\beta _{k_{1}}$ replacing $\beta _{k_{0}}$, for any fixed $k_{1}$, and $\xi _{n,1}$ replacing $\xi _{n,0}$, where $\xi _{n,1}$ is the sequence obtained by using $\beta _{k_{1}}$ for $n\geq k_{1}$.
\end{ulemma}

\noindent {\em{Proof.}} We need to start the induction and find an expression for $\xi _{n}\circ \xi _{n,0}^{-1}$ when $\psi _{n}=\psi _{n,0}$ and $\alpha _{n-1}=\alpha _{n-1,0}$. We start from the expressions
$$\sigma _{\alpha _{n}}\circ s\circ \xi _{n}\circ \psi _{n}=\psi _{n}\circ \sigma _{\alpha _{n-1}}\circ s$$
and
$$\sigma _{\alpha _{n,0}}\circ s\circ \xi _{n,0}\circ \psi _{n,0}=\psi _{n,0}\circ \sigma _{\alpha _{n-1,0}}\circ s.$$
The two right-hand sides are equal, and $\psi _{n}=\psi _{n,0}$. So  we obtain $$\sigma _{\alpha _{n}}\circ s\circ \xi _{n}=\sigma _{\alpha _{n,0}}\circ s\circ \xi _{n,0}.$$
So 
$$s\circ \xi _{n}\circ \xi _{n,0}^{-1}=\sigma _{\alpha _{n}}^{-1}\circ \sigma _{\alpha _{n,0}}\circ s,$$
which gives (\ref{8.8.4}).

Now we aim to prove (\ref{8.8.10}). We use (\ref{8.8.3}) and the corresponding statement for the sequences $\xi _{n,0}$ and $\alpha _{n,0}$, which is:
\begin{equation}\label{8.8.7}
\sigma _{\alpha _{n,0}}\circ s\circ \xi _{n,0}=\xi _{n-1,0}\circ \sigma _{\alpha _{n-1,0}}\circ s.\end{equation}
This can be rewritten as
\begin{equation}\label{8.8.8}\sigma _{\alpha _{n-1,0}}\circ s\circ \xi _{n,0}^{-1}=\xi _{n-1,0}^{-1}\circ \sigma _{\alpha _{n,0}}\circ s.\end{equation}
So then from (\ref{8.8.3}) and (\ref{8.8.8}), we obtain
$$\sigma _{\alpha _{n}}\circ s\circ \xi _{n}\circ \xi _{n,0}^{-1}=\xi _{n-1}\circ \sigma _{\alpha _{n-1}}\circ s\circ \xi _{n,0}^{-1}$$
$$=\xi _{n-1}\circ \sigma _{\alpha _{n-1}}\circ \sigma _{\alpha _{n-1,0}}^{-1}\circ \xi _{n-1,0}^{-1}\circ \sigma _{\alpha _{n,0}}\circ s.$$
With $\beta =\beta _{k_{0}}$, the right-hand side can be written as
$$\xi _{n-1}\circ \sigma _{\alpha _{n-1}}\circ \sigma _{\beta }^{-1}\circ \xi _{n-1,0}^{-1}\circ \sigma _{\beta }\circ \sigma _{\beta _{n}}^{-1}\circ \sigma _{\beta _{n}}\circ s.$$

This gives
$$\sigma _{\beta _{n}}\circ s\circ \xi _{n}\circ \xi _{n,0}^{-1}=\xi _{n-1}\circ \sigma _{\psi _{n-1}(\beta _{n}')}\circ \sigma _{\alpha _{n-1}}\circ \sigma _{\beta }^{-1}\circ \xi _{n-1,0}^{-1}\circ \sigma _{\beta }\circ \sigma _{\beta _{n}}^{-1}\circ \sigma _{\beta _{n}}\circ s,$$
which gives
$$\sigma _{\beta _{n}}\circ s\circ \xi _{n}\circ \xi _{n,0}^{-1}=\xi _{n-1}\circ \xi _{n-1,0}^{-1}\circ \sigma _{\zeta _{n}}\circ \sigma _{\beta _{n}}\circ s,$$
where $\zeta _{n}$ is as stated. The subsequent bound on ${\rm{supp}}(\xi _{n}\circ \xi _{n,0}^{-1})$ is then clear. 

The similar statements with $\beta _{k_{0}}$ and $\xi _{n,0}$ replaced by $\beta _{k_{1}}$ and $\xi _{n,1}$ are proved in the same way. 

\Box

\section{ Definition of $U(\beta _{1},\beta _{2})$}\label{8.6}

We now define the sets $U(\beta _{1},\beta _{2})$ for adjacent pairs $(\beta _{1},\beta _{2})$ in $R_{m,0}$, for any $m\geq 0$. The structure is suggested by the examples in \ref{7.5}. In all cases, $U(\beta _{1},\beta _{2})$ will have at most two components of intersection with $\{ z:\vert z\vert 
\leq 1\} $ which intersect $Z_{\infty }$ nontrivially.  We start by defining, for $X\subset \overline{\mathbb C}$,
$${\rm{ess}}(X)=X\cap \{ z:\vert z\vert \leq 1\} .$$
We are really only interested in ${\rm{ess}}(X)$ up to isotopy preserving $Z_{\infty }$. 

The basic idea is to choose $U(\beta _{1},\beta _{2})$ to be the union of $D(w(\beta _{2}))$, and everything between $D(w(\beta _{2}))$ and the nearest disc-crossing of $\beta _{1}$, but this might not be quite appropriate, if the nearest disc crossing by $\beta _{1}$ is not close by. So we define $i$  as in \ref{8.4}, that is, $i$ is the least integer $\geq 0$ such that $w_{i+1}'(\beta _{1})\neq w_{i+1}'(\beta _{2})$.  Note that $D(w(\beta _{1}))$ and $D(w(\beta _{2}))$ might not be adjacent in the unit disc. They will be adjacent if: 
$$w_{i+1}(\beta _{1})=w(\beta _{1}),\ \ w_{i+1}(\beta _{2})=w(\beta _{2}).$$

 We consider two conditions:
\begin{itemize}
\item[1.] 
$w_{i}(\beta _{1})=w_{i}(\beta _{2})$;
\item[2.] if $w(\beta _{2})$ has a prefix of length $\geq m-6$ ending in $BC$ or $UC$, then $\beta _{2}$ is not minimal or maximal among paths $\beta $ in $R_{m,0}$ such that $w(\beta )$ has this prefix.
\end{itemize}

\noindent {\em{Case that 1 and 2 hold.}}   Referring to \ref{8.4}, we see that, since $w(\beta _{2})$ does not end in $BC$ or $UC$, we are in case 1 of \ref{8.4}, and that:
\begin{itemize}
\item $w_{i+1}(\beta _{2})=w(\beta _{2})$;  
\item  $w_{i+3}(\beta _{1})=w(\beta _{1})$, that is, $\beta _{1}$ has at most two more crossings than $\beta _{2}$ -- and at least as many crossings as $\beta _{2}$, since the assumption that $w_{i}(\beta _{1})=w_{i}(\beta _{2})$ means that $w(\beta _{1})\neq w_{i}(\beta _{1})$. 

\end{itemize}
   The minimal path $\beta _{3}$ with $w_{i+1}'(\beta _{3})=w_{i+1}'(\beta _{1})$ is such that $w(\beta _{3})$ is of the same length as $w_{i+1}'(\beta _{1})$, ending in $R_{1}R_{2}C$ or $L_{2}C$, with $D(w(\beta _{3}))$ adjacent to $D(w_{i+1}'(\beta _{1}))$.  The closure of the set between $D(w(\beta _{3}))$ and $D(w(\beta _{2}))$ is $D(v)$, where $v$ is of the same length as either $w(\beta _{2})$ or $w(\beta _{3})$ and is obtained from $w(\beta _{2})$ (or $w(\beta _{3})$) by replacing the last letter by $R_{3}$. 

First suppose that $w(\beta _{1})$ does not both have length $m+2$ and end in $R_{1}R_{2}C$.  Then
\begin{equation}\label{8.6.0}{\rm{ess}}(U(\beta _{1},\beta _{2}))=D(w(\beta _{2}))\cup D(v).\end{equation} 
An example is given in \ref{7.5} by $w=L_{3}L_{2}R_{3}L_{3}L_{2}C$ and $w'=L_{3}L_{2}R_{3}L_{2}C$. If $w$ does have length $m+2$ and end in $R_{1}R_{2}C$, then write $w(\beta _{1})=w'R_{2}C$, and define
\begin{equation}\label{8.6.1}{\rm{ess}}(U(\beta _{1},\beta _{2}))=D(w(\beta _{2}))\cup D(v)\cup D(w'R_{1}).\end{equation}
In cases considered so far, ${\rm{ess}}(U(\beta _{1},\beta _{2}))$ is connected.

\noindent {\em{Case that  1 holds and 2 does not.}}   Now suppose  that  $w(\beta _{2})$ has a prefix $x$ of length $\geq m-6$ ending in $BC$ or $UC$, and such that $\beta _{2}$ is minimal or maximal among paths $\beta $ with $x$ as a prefix.  This includes all cases of $w_{i}(\beta _{1})$ and $w_{i}(\beta _{2})$ differing in just the last letter $BC$ or $UC$. If $\beta _{2}$ is maximal among such $\beta $, then let $v$ be as above, and we define 
\begin{equation}\label{8.6.2}{\rm{ess}}(U(\beta _{1},\beta _{2}))=D(w(\beta _{2}))\cup D(v)\setminus  D(xL_{1}R_{2}UC).\end{equation}
Having excluded $D(xL_{1}R_{2}UC)$ from $(\beta _{1},\beta _{2})$, we have to decide where to add it in. Let $x_{0}$ be obtained from $w(\beta _{2})$ by deleting the maximal suffix containing only the letters $BC$, $UC$, $L_{1}$ and $R_{2}$, apart from the last $R_{1}R_{2}C$.
In the case when $\beta _{2}$ is minimal among $\beta $ such that $w(\beta )$ has $x$ as a prefix, but not minimal among $\beta $ such that $w(\beta )$ has $x_{0}(BC{\rm{\ or\ }}UC)$ as a prefix and the end word of  $x$ includes at least one $UC$ and has length at least $4$, we let $x'$ be obtained from $x$ by replacing the end $UC(L_{1}R_{2}BC)^{k}$ by $BC(L_{1}R_{2}UC)^{k}$,  and then
\begin{equation}\label{8.6.7}{\rm{ess}}(U(\beta _{1},\beta _{2}))=D(w(\beta _{2}))\cup D(v)\cup D(x'L_{1}R_{2}UC).\end{equation}
If $\beta _{2}$ is minimal among $\beta $ such that $w(\beta )$ has $x_{0}(BC{\rm{\ or\ }}UC)$ as a prefix, then let $y$ be the word ending in $L_{3}L_{2}R_{3}L_{3}$ with $D(y)$ adjacent to $D(x_{0}C)$, and with $D(y)$ between $D(w(\beta _{1}))$ and $D(w(\beta _{2}))$, with $y$ of length $\geq m$, and, subject to these conditions, of the least possible length. Then we define 
\begin{equation}\label{8.6.8}{\rm{ess}}(U(\beta _{1},\beta _{2}))=D(w(\beta _{2}))\cup D(v)\cup D(y).\end{equation}

If the prefix $x$ does not include any letter $UC$, and is not minimal among $\beta $ such that $w(\beta )$ has $x_{0}(BC{\rm{\ or\ }}UC)$ as a prefix, then usually we use (\ref{8.6.0}). But if $x$ is of the form $L_{3}^{2}L_{2}BC(L_{1}R_{2}BC)^{k}$ and $w(\beta _{2})$ (and $w(\beta _{1})$) start with\\ $(BCL_{1}R_{2})^{k}(BC{\rm{\ or\ }}UC)$, then we define $x'=BC(L_{1}R_{2}UC)^{k+1}$ and again use (\ref{8.6.7}). 

The purpose of these definitions (\ref{8.6.2}) and (\ref{8.6.7}) is to ensure that  the set bounded by $\overline{\beta _{1}}*\beta _{2}$ and $\xi _{m,0}(\partial 'U(\beta _{1},\beta _{2}))$ is contained in $D(x)$. The addition of $D(xL_{1}R_{2}UC)$, or subtraction of $D(x'L_{1}R_{2}UC)$, is made so as to counteract the effect of $\xi _{m,0}$, for $\xi _{m,0}$ as in \ref{8.5}.

\noindent{\em{Case when 1 does not hold and 2 does hold.}}
Now suppose that $i>0$, $w_{i}'(\beta _{1})=w_{i}'(\beta _{2})$ and $w_{i}(\beta _{1})\neq w_{i}(\beta _{2})$, that is, that 1 above does not hold, but assume for the moment that 2 above does hold. If $w_{i}(\beta _{1})=w(\beta _{1})$, then $w_{i+1}(\beta _{2})=w(\beta _{2})$ is of the form $w'R_{2}C$. In this case, we define
\begin{equation}{\rm{ess}}(U(\beta _{1},\beta _{2}))=D(w(\beta _{2}))\cup D(w'R_{1}).\end{equation}

Now suppose that $w_{i}(\beta _{1})\neq w(\beta _{1})$. Then both $w_{i+1}'(\beta _{1})$  and $w_{i+1}'(\beta _{2})$ end in $L_{3}$. Write $w_{i+1}'(\beta _{1})=uL_{3}$. As before, define $v$  to be  the word of the same length as $w(\beta _{2})$, with $D(v)$ adjacent to $D(w(\beta _{2}))$, and ending in $R_{3}$. Then if $i=1$ and $w_{1}'(\beta _{1})=L_{3}$, we have:
\begin{equation}\label{8.6.3}{\rm{ess}}(U(\beta _{1},\beta _{2}))=D(w(\beta _{2}))\cup D(v)\cup D(u),\end{equation}
and in all other cases, we have:
\begin{equation}\label{8.6.4}{\rm{ess}}(U(\beta _{1},\beta _{2}))=D(w(\beta _{2}))\cup D(v)\cup D(u)\setminus D(vL_{3}).\end{equation}
The set $vL_{3}$ is of the form $y$ for some other pair $(\beta _{1},\beta _{2})$ as in (\ref{8.6.8}). So the set which is excluded in (\ref{8.6.4}) is added in in (\ref{8.6.8}).

\noindent{\em{Case when 1 does not hold and 2 does not hold.}}
This can only happen when $i=1$ and $w_{1}'(\beta _{1})=w_{1}'(\beta _{2})=L_{3}$, and $w(\beta _{2})$ starts $L_{3}^{2}L_{2}UC$. Then $\beta _{2}$ is minimal among paths $\beta $ which start with $L_{3}^{2}L_{2}BC$. Then $w(\beta _{2})$ starts with $L_{3}^{2}L_{2}(BCL_{1}R_{2})^{n}BCL_{1}R_{1}$ for some $n\geq 0$. The definition of $U(\beta _{1},\beta _{2})$ in this case is similar to (\ref{8.6.8}). In that case we define $v$ to be the word such that $D(v)$ is sandwiched between $D(L_{3}^{2}L_{2}C)$ and $D(w(\beta _{2}))$. We take $y=BCL_{1}R_{2}(UCL_{1}R_{2})^{n}UCL_{1}$ if this has length $m-2$, and otherwise $y$ is empty, and $D(y)$ is empty. Then we define
\begin{equation}{\rm{ess}}(U(\beta _{1},\beta _{2}))=D(w(\beta _{2}))\cup D(v)\cup D(y).\end{equation}

We define $U(\beta _{1},\beta _{2})$ itself
up to isotopy preserving $Z_{n}$ for any $n$, by adding handles 
round arcs of $\overline{\beta _{1}}*\beta _{2}$ to ${\rm{ess}}(U(\beta _{1},\beta _{2}))$.

\begin{lemma}\label{8.15} The  sets ${\rm{ess}}(U(\beta _{1},\beta _{2}))$, for $(\beta _{1},\beta _{2})$ running over adjacent pairs in $R_{m,0}$ with $w_{1}'(\beta _{2})=x$, have disjoint interiors, and the union contains all sets $D(w)\subset U^{x}$ with $w$ of length $\leq m+1$. 

Also, if $\beta $ is any path in $R_{m,0}$ with $w_{1}'(\beta )=x$ and $(\beta _{1},\beta _{2})$ is adjacent in $R_{n,0}$ for some $n\geq m$ with $w_{1}'(\beta _{1})=x$, then $\beta $ has no transversal intersections with $\partial 'U(\beta _{1},\beta _{2})$.

\end{lemma}
\noindent{\em{Proof.}} The sets $D(w(\beta _{2}))$ are disjoint. The sets $D(v)$ are associated with just one pair $(\beta _{1},\beta _{2})$, as are the sets 
$$D(u),\ \ D(w'R_{1}),\ \ D(y),$$
 where defined. So interiors are disjoint. Now we show that the union is as claimed. The first case of the definition of $U(\beta _1,\beta _2)$ is when it is the region  of the disc between the $i$'th unit-disc-crossings of $\beta _1$ and $\beta _2$, where these are the first crossings which are distinct. If $w$ is a word of length $\leq m+1$, the points in $U^x\cap D(w)$ which are not in such a region are in sets of the form $D(w'R_1)$ or $D(u)$ or $D(y)$, where there are as in \ref{8.6}. The definition of $U(\beta _1,\beta _2)$ is made to include all of these. 
 
For the last part, if $w_{1}'(\beta _{2})\neq x$, then the claim of no transversal intersections follows immediately, because $\beta $ is contained in $U^{x}$ after the first $S^{1}$-crossing, and $U(\beta _{1},\beta _{2})\subset U^y$ for $y\neq x$. There are some intersections between sets $U^x$ and $U^y$, which arise from the exceptional definitions of $w_1'(w,x)$ (and $w_1'(w,y)$) listed in \ref{7.9}: but none of these occurs adjacent to $\beta _1$, if $\beta _1$ is the maximal path with $w_1'(\beta _1)=x$.  So we can assume that $w_{1}'(\beta _{1})=w_{1}'(\beta _{2})=x$. Suppose that  $\beta $ intersects $\partial 'U(\beta _{1},\beta _{2})$. Then for some $j$, $w_{j}'(\beta )$ must be strictly  between $w_{i+1}'(\beta _{1})$ and $w_{i+1}'(\beta _{2})$, and $w_{i-1}(\beta _{1})=w_{i-1}(\beta _{2})$ must be a prefix of $w_{j}'(\beta )$ also. If $w_{i}(\beta _{1})\neq w_{i}(\beta _{2})$, then, since $n\geq m$, and $\beta _{1}$ and $\beta _{2}$ are adjacent, $w_{i}(\beta )$ is not between $w_{i}(\beta _{1})$ and $w_{i}(\beta _{2})$, and hence also not between $w_{i+1}'(\beta _{1})$ and $w_{i+1}'(\beta _{2})$, giving a contradiction. If $w_{i}(\beta _{1})=w_{i}(\beta _{2})$, then $w_{i}(\beta _{1})$ is also a prefix of $w(\beta )$, and hence $w_{i}(\beta )=w_{i}(\beta _{1})$ and $w_{j}'(\beta )=w_{j}'(\beta _{1})$ for $j\leq i$. Then because $\beta _{1}$ and $\beta _{2}$ are adjacent, $w_{i+1}'(\beta )$ is not between $w_{i+1}'(\beta _{1})$ and $w_{i+1}'(\beta _{2})$, and $w_{i+1}(\beta )$ is not between  $w_{i+1}'(\beta _{1})$ and $w_{i+1}'(\beta _{2})$. Then since $j\geq i+2$, we see that $w_{i+1}(\beta _{1})$ is a prefix of $w_{j}'(\beta )$, and  $w_{j}'(\beta )$ is also not between $w_{i+1}'(\beta _{1})$ and $w_{i+1}'(\beta _{2})$, giving the final contradiction.
\Box

\section{Preservation of $U(\beta _{1},\beta _{2})$ under local inverses}\label{8.16}
We are especially interested in the cases in which a set $U(\beta _{1},\beta _{2})$ is mapped by a local inverse of $s^{n}$ to a set $U(\beta _{3},\beta _{4})$. From the definition given of $U(\beta _{1},\beta _{2})$, it is clear that  this happens precisely when  $\overline{\beta _{1}}*\beta _{2}$ is mapped by a local inverse of $s^{n}$ to $\overline{\beta _{3}}*\beta _{4}$. 
So for this, we need to refer to \ref{7.10},  and, more fundamentally, to (\ref{7.6.1}) to (\ref{7.6.6}). The definitions are such that, if $S$ is a local inverse of $s^{n}$ as in \ref{7.11}, or more generally a local homeomorphism as in \ref{7.11}, and $(\beta _{1},\beta _{2})$ are adjacent in $R_{m,0}$, then there is an adjacent pair $(\beta _{3},\beta _{4})$ in $R_{m+n,0}$ such that
\begin{equation}\label{8.16.1}SU(\beta _{1},\beta _{2})=U(\beta _{3},\beta _{4}),\end{equation}
and conversely, if $(\beta _{3},\beta _{4})$ is an adjacent pair in $R_{m+n,0}$, in the image of $S$, then there is an adjacent pair $(\beta _{1},\beta _{2})$ in $R_{m,0}$ such that (\ref{8.16.1}) holds.

\section{Definition of $V(\beta _{1},\beta _{2},n)$}\label{8.17}

Let $(\beta _{1},\beta _{2})$ be an adjacent pair in $R_{m,0}$. Let $\xi _{m,0}$ and $\xi _{m,1}$ be as in \ref{8.5}. The following lemma means that if we define $V(\beta _{1},\beta _{2},n)$ to be the union of sets $U(\beta _{3},\beta _{4})$ with $(\beta _{3},\beta _{4})$ adjacent in $R_{n,0}$ and between $\beta _{1}$ and $\beta _{2}$, then, up to isotopy preserving $Z_{m+1}$, and usually up to isotopy preserving $Z_{\infty }$, the set $V(\beta _{1},\beta _{2},n+1)$ is also the region bounded by $\overline{\beta _{1}}*\beta _{2}$ and $\xi _{m,0}(\partial 'U(\beta _{1},\beta _{2}))$, except in a few cases, when we shall see later that it is the region bounded by $\overline{\beta _{1}}*\beta _{2}$ and $\xi _{m,1}(\partial 'U(\beta _{1},\beta _{2}))$.

\begin{ulemma} 
Let $(\beta _{1},\beta _{2})$ be an adjacent pair in in $R_{m,0}$ and let $i$ be the largest index with $w_{i}'(\beta _{1})=w_{i}'(\beta _{2})$, setting $i=0$ if $w_{1}'(\beta _{1})\neq w_{1}'(\beta _{2})$. If $i\geq 1$, set $x=w_{1}'(\beta _{1})=w_{1}'(\beta _{2})$.  Then the region bounded by $\overline{\beta _{1}}*\beta _{2}$ and $\xi _{m,0}(\partial 'U(\beta _{1},\beta _{2}))$ is the union of arbitrarily small neighbourhoods of sets $U(\beta _{3},\beta _{4})$ for $(\beta _{3},\beta _{4})$ adjacent in $R_{m+1,0}$, and $\beta _{1}\leq \beta _{3}<\beta _{4}\leq \beta _{2}$, up to homotopy preserving $Z_{m+1}$ {\em{except}} in the following cases.  In cases 2 and 3, we compensate by taking the pointwise definitions of some  $\beta _{1}'$ and $\beta _{2}'$ to be different in $R_{m+1,0}'$ from the definitions in $R_{m,0}'$. Apart from these exceptional cases, the region and the union are equal up to homotopy preserving $Z_{\infty }$.

\begin{itemize}
\item[1.] There is $\beta _{3}\in R_{m,0}$ between $\beta _{1}$ and $\beta _{2}$ such that $w=w(\beta _{3})$ and $x=w_{1}'(\beta _{3})$, and  $w_{i+1}'(w,x,0)$ is one of
 the exceptional cases of definition in \ref{7.9}.
\item[2.] The difference is in a set $D(y)$, where $y$ ends in $(L_{3}{\rm{\ or\ }}R_{3})L_{2}R_{3}L_{3}$. This happens   when $w_{i}(\beta _{3})\neq w_{i}(\beta _{4})$ for the largest $i$ with $w_{i}'(\beta _{3})=w_{i}'(\beta _{4})$, except when $i=1$ and $w_{1}'(\beta _{3})=w_{1}'(\beta _{4})=L_{3}$. In all the cases with $i>1$, there  is a $y$ as above with $D(y)$ between $D(w(\beta _{3}))$ and $D(w(\beta _{4}))$ and with $w_{i}(y,w_{1}'(\beta _{3}))=y$. It also happens when such a $D(y)$ is between $D(w(\beta _{3}))$ and $D(w(\beta _{4}))$, adjacent to $D(v)$, with $\beta _{3}<\beta _{4}$, and $\beta _{4}$ is minimal with $v(BC{\rm{\ or\ }}UC{\rm{\ or\ }}C)$ as a prefix for some $v$.
\item[3.] The difference is in a word ending in $D(y)$, where $y$ ends in \\ $(UCL_{1}R_{1})^{n}BCL_{1}R_{1}R_{2}$, and $3n+6>m$. This happens when such a  $D(y)$ is between $D(w_{i}(\beta _{3}))$ and $D(w_{i}(\beta _{4}))$, again for the largest $i$ such that $w_{i}'(\beta _{3})=w_{i}'(\beta _{4})$. It also happens if $y=v(UCL_{1}R_{1})^{n}BCL_{1}R_{1}R_{2}$ and $\beta _{3}$ is maximal  with $w(\beta _{3})$ having $v(BC{\rm{\ or\ }}UC{\rm{\ or\ }}C)$ as a prefix.
\end{itemize}\end{ulemma}

\noindent {\em{Proof.}}

The first exceptional case occurs precisely when, if paths were defined using the sequences $w_{i}'(w,0)$ rather than $w_{i}'(w,w_{1}'(w,0),0)$, the paths between $\beta _{1}$ and $\beta _{2}$ would be different. In all other cases, the paths can equally well be defined using the sequences $w_{i}'(w,0)$. So now we consider all these other cases

Recall that the support of $\xi _{m,0}$ is $A_{m,2/7}\cup C_{m,0}$. The region bounded by $\overline{\beta _{1}}*\beta _{2}$ and $\xi _{m,2/7}(\partial 'U(\beta _{1},\beta _{2}))$ is
$$U(\beta _{1},\beta _{2})\cup D_{1}\setminus D_{2},$$
where $D_{1}$ and $D_{2}$ are as follows. The set $D_{1}$ is the union of any beads on $A_{m,2/7}$ or $C_{m,0}$ which follow on from beads inside $U(\beta _{1},\beta _{2})$ and such that the piece of string between the two beads first crosses $\partial U(\beta _{1},\beta _{2})$ along $\partial 'U(\beta _{1},\beta _{2})$. The set $D_{2}$ is the union of any beads of $A_{m,2/7}$ or $C_{m,0}$ inside $U(\beta _{1},\beta _{2})$ which follow on from pieces of string which cross $\partial U(\beta _{1},\beta _{2})$, and such that the crossing of this piece of string into $U(\beta _{1},\beta _{2})$ is across $\partial 'U(\beta _{1},\beta _{2})$.

 Most pieces of string on $C_{m,0}$ trace round inside gaps of $L_{3/7}$, or segments of the unit circle close to the boundary of a single gap, so that they have limited scope for intersecting sets $\partial 'U(\beta _{1},\beta _{2})$. The definitions in (\ref{8.6.2}) and (\ref{8.6.7}) ensure that alternate strings, those from beads with words ending in $UC$ to beads with words ending in $BC$, usually intersect precisely two sets $\partial 'U(\beta _{1},\beta _{2})$ and $\partial 'U(\beta _{2},\beta _{3})$. The exceptions are words ending  $L_{2}BC(L_{1}R_{2}UC)^{k+1}$ and $L_{3}^{3}L_{2}BC(L_{1}R_{2}BC)^{k}$. But such a piece of string is the image under $S$, for some local inverse of some $s^{n}$ as in  \ref{7.11}, of the piece of string on $C_{m,2/7}$ connecting $D(BC(L_{1}R_{2}BC)^{k})$ and $D(L_{3}^{2}L_{2}BC(L_{1}R_{2}BC)^{k})$. This piece of string also intersects precisely two sets $\partial 'U(\beta _{1},\beta _{2})$ and $\partial 'U(\beta _{3},\beta _{4})$. So the action of $C_{m,0}$ gives rise to the exceptions listed in 3 and the last exception in 2.

Now we consider $A_{m,2/7}$. Note that ${\rm{ess}}(U(\beta _{1},\beta _{2}))$ includes $D(w(\beta _{2}))$ and also  a connected component $D(v)$ which is adjacent to $D(w(\beta _{2}))$, except that in some cases $D(vL_{3})$ may be excluded. Let $i$ be the index with $w_{i}'(\beta _{1})=w_{i}'(\beta _{2})$ but $w_{i+1}'(\beta _{1})\neq w_{i+1}'(\beta _{2})$. There is a path $\beta _{3}\in R_{m+1,0}$ between $\beta _{1}$ and $\beta _{2}$ if and only if $w_{i}'(\beta _{3})=w_{i}'(\beta _{2})=w_{i}'(\beta _{1})$ and $w_{i+1}'(\beta _{3})$ is between $w_{i+1}'(\beta _{1})$ and $w_{i+1}'(\beta _{2})$. Then $D(w(\beta _{3}))$ is in $D(v)$, unless $D(z)$ is in $D(v)$, for some $z$ such that $D(z)$ and $D(z')$ are consecutive beads on some component of $A_{m,2/7}$, with $D(z)$ preceding $D(z')$,  and $w(\beta _{3})=z'L_{2}C$. In that case, $zL_{2}C=w(\beta _{4})$, where $\beta _{4}\in R_{m+1,0}$ is also between $\beta _{1}$ and $\beta _{2}$. The solid boundary of $U(\beta _{1},\beta _{2})$ is $\overline{\beta _{1}}*\beta _{2}$. So the string $z\mapsto z'$ crosses $\partial 'U(\beta _{1},\beta _{2})$, unless either $D(z)$ and $D(z')$ are both excluded from ${\rm{ess}}(U(\beta _{1},\beta _{2}))$ or $D(z)$ and $D(z')$ are both included. If $D(z)$ and $D(z')$ are both excluded, then $z=vL_{3}$ for a nontrivial $v$, and we are in case (\ref{8.6.4}). In that case, the conditions of 5 of \ref{7.6} imply  that $\beta _{4}$ is between $\beta _{1}$ and $\beta _{2}$, and $\beta _{3}$ is not. The definition of (\ref{8.6.4}) ensures that $D(z)$ is in the region bounded by $\overline{\beta _{1}}*\beta _{2}$ and $\xi _{m,0}(\partial 'U(\beta _{1},\beta _{2}))$, and $D(z')$ is not. If $D(z)$ and $D(z')$ are both included in $U(\beta _{1},\beta _{2})$, then we are in case (\ref{8.6.8}), and then $D(z')$ is in the region bounded by $\overline{\beta _{1}}*\beta _{2}$ and $\xi _{m,0}(\partial 'U(\beta _{1},\beta _{2}))$, and $D(z)$ is not. Condition 5 of \ref{7.6} confirms that $\beta _{3}$ is between $\beta _{1}$ and $\beta _{2}$ and $\beta _{4}$ is not. 

The proof of \ref{8.2} is then completed by the following lemma. 

\begin{lemma}\label{8.9}
Let $(\gamma _{1},\gamma _{2})$ be an adjacent pair in $R_{m,0}$. Let $i$ be maximal with $w_{i}'(\beta _{1})=w_{i}'(\beta _{2})$ (possibly with $i=0$) and let $x=w_{1}'(\gamma _{1},0)$. Once again, let $\xi _{m,0}$ and $\xi _{m,1}$ be as in \ref{8.5}. Then, for $m\geq k_{0}$, intersections between  $\xi _{m,0}(\partial 'U(\gamma _{1},\gamma _{2}))$  and the support of $\xi _{m,1}\circ \xi _{m,0}^{-1}$ occur only when $i\geq 1$, and there is $\gamma _{3}$ between $\gamma _{1}$ and $\gamma _{2}$ in $R_{m+1,0}$ such that $w_{i+1}'(w(\gamma _{3}),x,0)$ is one of the exceptional cases of definition in \ref{7.9}. There are no intersections at all with the support of $\xi _{m}\circ \xi _{m,1}^{-1}$.\end{lemma}

\noindent {\em{Proof.}} Let $(\beta _{k},\beta _{k,2})$, $E_{n}$ and $E_{n,k}$ be as in \ref{8.8}, with $(\gamma _{1},\gamma _{2})=(\beta _{m},\beta _{m,2})$. So $E_{n,k}$ is obtained by taking successive preimages under maps $\sigma _{\beta _{j}}\circ s$   a small disc neighbourhood $E_{k,k}$  of  the loop $\zeta _{k}$ as in \ref{8.8}. We need to show that for $r\leq m$ there are only limited intersections  between $E_{m,r}$ and $\xi _{m,0}(\partial 'U(\beta _{1},\beta _{2}))$ for $r\leq k_{1}$, with these intersections corresponding to the exceptional definitions of $w_{i}'(w,x,0)$ in \ref{7.9}, and no intersections at all between $E_{m,r}$ and $\xi _{m,1}(\partial 'U(\gamma _{1},\gamma _{2}))$ for $r> k_{1}$.

The set $E_{k,k}$ is nonempty only when either $\beta _{k}\neq \beta _{k-1}$ or when $\beta _{k}'\neq \beta _{k-1}'$, or when the support of $\xi _{k-1}\circ \xi _{k-1,0}^{-1}$ intersects $\overline{\beta _{k}}*\beta _{2/7}$. Note that the possibility of $\beta _{k-1}'\neq \beta _{k}'$ does occur, even when $\beta _{k}=\beta _{k-1}$. These are the cases listed in \ref{8.17}, when a set bounded by $\overline{\beta _{k-1}}*\beta _{k-1,2}$ and $\xi _{k-1,k_{0}}(\partial 'U(\beta _{k-1}\beta _{k-1,2}))$ (for $k\leq k_{1}$) or 
$\xi _{k-1,k_{1}}(\partial 'U(\beta _{k-1}\beta _{k-1,2}))$ (for $k>k_{1}$) is not the same up to homotopy preserving $Z_{\infty }$ as the union of sets $U(\gamma _{3},\gamma _{4})$ for pairs $(\gamma _{3},\gamma _{4})$ of $R_{k,0}$ between $\beta _{k-1}$ and $\beta _{k-1,2}$.

The loop $\zeta _{m-r,2}$ traces round the union of sets  $\xi _{m,0}(\partial U(\gamma _{3},\gamma _{4}))$  for adjacent pairs $(\gamma _{3},\gamma _{4})$ between $\beta _{m-r-1}$ and $\beta _{m-r}$. It is usually the case that $(\beta _{m-r-1},\beta _{m-r})$ is an adjacent pair, but there can be a path between $\beta _{m-r-1}$ and $\beta _{m-r}$. This means that components of $s^{-r}(\zeta _{m-r,2})$ do not have transversal intersections with any set $\partial 'U(\gamma _{5},\gamma _{6})$. We claim that the same is true for  $\zeta _{m-r,1}$. This is because the preimage under $s^{-t}$ of any piece of string in the support of $\xi _{m-r-1,0}$ which crosses $\overline{\beta _{m-r}}*\beta _{k_{0}}$, and which intersects $U^{0}$, will cross the union of solid boundaries of sets $U(\gamma _{3},\gamma _{4})$, for adjacent pairs $(\gamma _{3},\gamma _{4})$ in $R_{m-r+t,0}$, represented by a component of $s^{-t}(\overline{\beta _{m-r}}*\beta _{k_{0}})$, and will not cross any dashed boundary before the next bead.   So the boundary of the set $s^{-t}(\zeta _{m-r,2})$ has no transversal intersections with any set  $\xi _{m,0}(\partial 'U(\gamma _{5},\gamma _{6}))$.

We now consider a component of $E_{m,m-r}$ such that the corresponding component of $E_{n,m-r}$, for $m-r+2\leq n<m$, does not intersect $\beta _{n+1}$. Then $E_{m,m-r}$ is a component of $s^{1-r}(E_{m-r+1,m-r})$, and intersections with $\xi _{m,0}(\partial 'U(\gamma _{1},\gamma _{2}))$ can only occur as intersections between $\xi _{m,0}(\partial 'U(\gamma _{1},\gamma _{2}))$  and components of 
$s^{-r}(\beta _{m-r})$.
For intersections between $\xi _{m,1}(\partial 'U(\gamma )_{1},\gamma _{2}))$ and components of $s^{1-r}(E_{m-r+1,m-r})$, for $m-r\geq k_{1}$, under the assumption that we are currently making on these components, we only need to consider components of the support of $s^{-r}(\beta _{m-r})$
 which intersect $U^{0}$. In fact, we shall consider, more generally, intersections with $\xi _{m,0}(\partial 'U(\gamma _{3},\gamma _{4}))$, for any adjacent pair $(\gamma _{3},\gamma _{4})$ with $w_{1}'(\gamma _{3})=w_{1}'(\gamma _{1})$.

First we consider intersections between $s^{-r}(\beta _{m-r})$ and $\xi _{m,0}(\partial 'U(\gamma _{3},\gamma _{4}))$.  Any component of  $\xi _{n,1}\circ \xi _{n,0}^{-1}$, for any $n$, which has support close to $s^{-1}(\beta _{m-r})$, can be broadly represented by
$$L_{3}^{2}L_{2}R_{3}{\rm{(bot)}}\leftrightarrow {\rm{(top)}}R_{3}L_{3}L_{2}R_{3}$$
This means that the part of the homeomorphism supported near $s^{-1}(\beta _{m-r})$ interchanges two discs in $D(L_{3}^{2}L_{2}R_{3}L_{3})$ and $D(R_{3}L_{3}L_{2}R_{3})$ along the line of the diagonal $s^{-1}(\beta _{m-r})$. We get corresponding components of $\xi _{n+r,1}\circ \xi _{n+r,0}^{-1}$ by taking preimages under $s^{r}$. The second preimages are 
\begin{equation}\label{8.9.1}L_{3}^{4}L_{2}R_{3}{\rm{(bot)}}\leftrightarrow {\rm{(bot)}}R_{3}L_{2}R_{3}L_{3}L_{2}R_{3},\end{equation}
\begin{equation}\label{8.9.2}R_{3}L_{3}^{3}L_{2}R_{3}{\rm{(top)}}\leftrightarrow {\rm{(top)}}L_{3}L_{2}R_{3}L_{3}L_{2}R_{3},\end{equation}
\begin{equation}\label{8.9.3}L_{2}R_{3}L_{3}^{2}L_{2}R_{3}{\rm{(top)}}\leftrightarrow {\rm{(bot)}}L_{1}R_{2}R_{3}L_{3}L_{2}R_{3},\end{equation}
\begin{equation}\label{8.9.4}R_{2}R_{3}L_{3}^{2}L_{2}R_{3}{\rm{(bot)}}\leftrightarrow {\rm{(top)}}R_{1}R_{2}R_{3}L_{3}L_{2}R_{3}.\end{equation}

 We claim that only components of support of $\xi _{m,1}\circ \xi _{m,0}^{-1}$ represented by further preimages of the first two will ever intersect $\xi _{m,0}(\partial 'U(\gamma _{3},\gamma _{4}))$, for any $m$. We can restrict to looking at any set  $\xi _{m,0}(\partial 'U(\gamma _{3},\gamma _{4}))$ which is not in the inverse orbit of some other $\xi _{n,0}(\partial 'U(\beta _{7},\beta _{8}))$, and a component of support of  $\xi _{m,1}\circ \xi _{m,0}^{-1}$ which intersects $U^{0}$, and which is not in the backward orbit of any other component of support of $\xi _{n,1}\circ \xi _{n,0}^{-1}$ which also intersects $U^{0}$.  For if $\xi _{m,0}(\partial 'U(\gamma _{3},\gamma _{4}))$ is  in the inverse orbit of some other $\xi _{n,0}(\partial 'U(\gamma _{5},\gamma _{6}))$, it can only intersect a component of support of $\xi _{m,1}\circ \xi _{m,0}^{-1}$ if a component of support of $\xi _{n,1}\circ \xi _{n,0}^{-1}$ intersects $\xi _{n,0}(\partial 'U(\gamma _{5},\gamma _{6}))$.
 
So now we consider the backward orbit of (\ref{8.9.1}) and adjacent pairs $(\gamma _{3},\gamma _{4})$ in $R_{m,0}$ with $w_{1}'(\gamma _{3})=w_{1}'(\gamma _{1})$. The two immediate preimages are given by adjoining either $(L_{3},L_{2})$ or $(R_{3},R_{2})$, that is, we have
\begin{equation}\label{8.9.5}L_{3}^{5}L_{2}R_{3}{\rm{(top)}}\leftrightarrow {\rm{(top)}}L_{2}R_{3}L_{2}R_{3}L_{3}L_{2}R_{3}\end{equation}
and
\begin{equation}\label{8.9.6}R_{3}L_{3}^{4}L_{2}R_{3}{\rm{(bot)}}\leftrightarrow {\rm{(bot)}}R_{2}R_{3}L_{2}R_{3}L_{3}L_{2}R_{3}.\end{equation}
Any backward orbit of (\ref{8.9.6}) passes through one of the following:
\begin{equation}\label{8.9.7}L_{3}L_{2}R_{3}L_{3}^{4}L_{2}R_{3}{\rm{(bot)}}\leftrightarrow {\rm{(bot\ left\ )}}BCL_{1}R_{2}R_{3}L_{2}R_{3}L_{3}L_{2}R_{3},\end{equation}
$$R_{3}L_{2}R_{3}L_{3}^{4}L_{2}R_{3}{\rm{(bot)}}\leftrightarrow {\rm{(top\  right)}}UCL_{1}R_{2}R_{3}L_{2}R_{3}L_{3}L_{2}R_{3},$$
\begin{equation}\label{8.9.8}BCL_{1}R_{1}^{n-1}R_{2}R_{3}L_{3}^{4}L_{2}R_{3}{\rm{(bot\ left)}}\leftrightarrow {\rm{(bot\ left)}}BCL_{1}R_{1}^{n}R_{2}R_{3}L_{2}R_{3}L_{3}L_{2}R_{3},\end{equation}
$$UCL_{1}R_{1}^{n-1}R_{2}R_{3}L_{3}^{4}L_{2}R_{3}{\rm{(top\ right)}}\leftrightarrow {\rm{(top\ right)}}UCL_{1}R_{1}^{n}R_{2}R_{3}L_{2}R_{3}L_{3}L_{2}R_{3}.$$

There is no intersection of the component of $\xi _{m,1}\circ \xi _{m,0}^{-1}$ identified by (\ref{8.9.7}) with $\xi _{m,0}(\partial 'U(\gamma _{3},\gamma _{4}))$, for any adjacent pair $(\gamma _{3},\gamma _{4})$ in $R_{m,0}$. We see this as follows. Note that most of the sets $U(\gamma _{3},\gamma _{4})$ have solid boundary on the left side of $D(BC)$ and dashed boundary on the right. The only set $\xi _{m,0}(\partial 'U(\gamma _{3},\gamma _{4}))$ which spans both $D(BC)$ and $D(L_{3})$ is the longest one of all, which passes from the bottom righthand side of $D(BC)$ to $D(L_{3}^{4})$. Similarly, there is no intersection of the component of $\xi _{m,1}\circ \xi _{m,0}^{-1}$ identified by (\ref{8.9.8}) with $\xi _{m,0}(\partial 'U(\gamma _{3},\gamma _{4}))$, for any adjacent pair $(\gamma _{3},\gamma _{4})$ in $R_{m,0}$.  Applying the prefix and substitution  rules of (\ref{8.16}) for the shape of sets $\partial 'U(\gamma _{3},\gamma _{4})$, there are no intersections in the other cases also.

So we consider the backward orbit of (\ref{8.9.5}) obtained by adjoining   $(uL_{3},uL_{2})$. Because of the conditions of 5 of \ref{7.6}, the presence of  intersections, or not, is determined by the longest suffix of $u$ which  consists only of the letters $L_{2}$, $L_{3}$ and $R_{3}$, of the form $u_{1}\cdots u_{n}L_{3}$ described in 5 of \ref{7.6}, that is, $n\geq 1$, $u_{i}=L_{3}(L_{2}R_{3})^{m_{i}}L_{3}^{r_{i}-1}$ with $m_{i}+r_{i}=3$ for $i\geq 2$, $m_{1}+n_{1}$ is odd and $\geq 3$, and $m_{1}\geq 1$. Then intersection occurs with a set $\xi _{m,0}(\partial 'U(\gamma _{3},\gamma _{4}))$  if, and only if, $n$ is odd.  The definition of $w_{i}'(\gamma _{3})$ in these cases was chosen precisely so that the region bounded by $\overline{\gamma _{3}}*\gamma _{4}$ and $\xi _{m,1}(\partial 'U(\gamma _{1},\gamma _{2}))$ is the union of the sets $U(\gamma _{5},\gamma _{6})$ for adjacent pairs $(\gamma _{5},\gamma _{6})$ in $R_{m+1,0}$ between $\gamma _{3}$ and $\gamma _{4}$. These are the exceptions identified in \ref{7.9}.

Now we consider (\ref{8.9.2}). It is probably worth pointing out that, whenever this case arises, $\beta _{m-r}\neq \beta _{k_{0}}$, and therefore the immediate preimage is definitely under $s^{-1}$. We have two immediate preimages that we can adjoin immediately: 
\begin{equation}\label{8.9.12}R_{2}R_{3}L_{3}^{3}{\rm{top}}\leftrightarrow{\rm{top}}R_{3}L_{3}L_{2}R_{3}L_{3}\end{equation}
and
\begin{equation}\label{8.9.13}L_{2}R_{3}L_{3}^{3}{\rm{bot}}\leftrightarrow{\rm{bot}}L_{3}^{2}L_{2}R_{3}L_{3}\end{equation}

Preimages of (\ref{8.9.13}) are similar to (\ref{8.9.5}) but, with top and bottom reversed. This gives rise to some intersection with sets $\xi _{m,0}(\partial 'U(\gamma _{3},\gamma _{4}))$, which, as with (\ref{8.9.5}), arise exactly in the circumstances identified in (\ref{7.9}). 

For (\ref{8.9.12}) we have further preimages defined by similar words to (\ref{8.9.6}):
\begin{equation}\label{8.9.14}BCL_{1}R_{2}R_{3}L_{3}^{3}{\rm{(bot\ right)}}\leftrightarrow{\rm{(bot)}}R_{3}L_{2}R_{3}L_{3}L_{2}R_{3}L_{3},\end{equation}
\begin{equation}\label{8.9.15}UCL_{1}R_{2}R_{3}L_{3}^{3}{\rm{(top\ left)}}\leftrightarrow{\rm{(top)}}L_{3}L_{2}R_{3}L_{3}L_{2}R_{3}L_{3},\end{equation}
\begin{equation}\label{8.9.16}BCL_{1}R_{1}^{n}R_{2}R_{3}L_{3}^{3}{\rm{(bot\ right)}}\leftrightarrow{\rm{(bot\ right)}}BCL_{1}R_{1}^{n-1}R_{2}R_{3}L_{3}L_{2}R_{3}L_{3},\end{equation}
\begin{equation}\label{8.9.17}UCL_{1}R_{1}^{n}R_{2}R_{3}L_{3}^{3}{\rm{(top\ left)}}\leftrightarrow{\rm{(top\ left)}}BCL_{1}R_{1}^{n-1}R_{2}R_{3}L_{3}L_{2}R_{3}L_{3}.\end{equation}

  No preimage of (\ref{8.9.14}) to (\ref{8.9.17}) intersects any set of the form $\xi _{m,0}(\partial 'U(\gamma _{3},\gamma _{4}))$, for similar reasons to (\ref{8.9.6}), that is, for similar reasons to (\ref{8.9.7}) and (\ref{8.9.8}). 

Now we consider preimages of (\ref{8.9.3}). The first preimages are
\begin{equation}\label{8.9.9}L_{3}L_{2}R_{3}L_{3}^{2}L_{2}R_{3}{\rm{(bot)}}\leftrightarrow {\rm{(bot\  right)}}BCL_{1}R_{2}R_{3}L_{3}L_{2}R_{3},\end{equation}
and 
\begin{equation}\label{8.9.10}R_{3}L_{2}R_{3}L_{3}^{2}L_{2}R_{3}{\rm{(top)}}\leftrightarrow {\rm{(top\  left)}}UCL_{1}R_{2}R_{3}L_{3}L_{2}R_{3}.\end{equation}
The track in (\ref{8.9.9}) does not intersect any set $\xi _{m,0}(\partial 'U(\gamma _{3},\gamma _{4}))$ for $6\leq m$. We see this as follows. The longest dashed boundary $\xi _{m,0}(\partial 'U(\gamma _{3},\gamma _{4}))$ which spans both $D(BC)$ and $D(L_{3})$ passes on the outside of the track of (\ref{8.9.9}), and apart from this the there are no other dashed boundaries which intersect either $BCL_{1}R_{2}R_{3}L_{3}L_{2}R_{3}$ or $D(L_{3}L_{2}R_{3}L_{3}^{2}L_{2}R_{3})$ which are not preimages of others. The same is true if we adjoin any $(u,u')$, by the rules of \ref{8.16}. For (\ref{8.9.10}): if we prefix by $BCL_{1}R_{2}(R_{3}L_{2})^{n}$ or by $L_{3}L_{2}(R_{3}L_{2})^{n}$ respectively then there are no intersections with any $\xi _{m,0}(\partial 'U(\beta _{5},\beta _{6}))$ for $9+2n\leq m$ or $8+2n\leq m$ respectively, by inspection. The same is true for other choices of prefix, by the rules governing the sets $\partial 'U(\gamma _{3},\gamma _{4})$ in \ref{8.16}.

The case of (\ref{8.9.4}) is similar. All preimages  are further preimages of
$$L_{1}R_{1}^{n}R_{2}R_{3}L_{3}^{2}L_{2}R_{3}{\rm{(top)}}\leftrightarrow {\rm{(bot)}}L_{1}R_{1}^{n+1}R_{2}R_{3}L_{3}L_{2}R_{3}.$$
Prefixing by $BC$, we have
\begin{equation}\label{8.9.11}BCL_{1}R_{1}^{n}R_{2}R_{3}L_{3}^{2}L_{2}R_{3}{\rm{(bot\ left)}}\leftrightarrow {\rm{(bot\  right)}}BCL_{1}R_{1}^{n+1}R_{2}R_{3}L_{3}L_{2}R_{3}.\end{equation}
We see that there are no intersections with $\xi _{m,0}(\partial 'U(\gamma _{3},\gamma _{4}))$ for $6+n\leq m$, because any set $\xi _{m,0}(\partial 'U(\gamma _{3},\gamma _{4}))$ with $w(\gamma _{3})$ starting with $BCL_{1}R_{1}^{n}R_{2}R_{3}L_{3}^{2}$ and with $\gamma _{3}<\gamma _{4}$ has $w(\gamma _{4})$ starting with $BCL_{1}R_{1}^{n+1}R_{2}BC$, which is shorter than the track of (\ref{8.9.11}).
Then there are no intersections by attaching prefixes, by the rules governing the sets $\partial 'U(\gamma _{3},\gamma _{4})$. The same is true if we attach $UC$ and a prefix.

Now we consider what happens when we drop the assumption that successive components of $s^{-k}(E_{m-r+1,m-r})$ avoid $\beta _{m-r+k+1}$. So suppose a component of $s^{-k}(E_{m-r+1,m-r})$ intersects $\beta _{m-r+k+1}$. We point out that this can happen, and, indeed, it is possible for a component of $s^{-k}(E_{m-r+1,m-r})$ to intersect the initial segment of $\beta _{m-r+k+1}$. Whenever there is intersection with an initial segment, there is an intersection between a component of $s^{-i}(E_{m-r+1,m-r})$ and $\beta _{2/7}$ for some $i\leq k$, with $i<k$ unless the initial segments of $\beta _{m-r+k+1}$ and $\beta _{2/7}$ coincide, and taking further preimages of this, as necessary,  produces the intersection between $s^{-k}(E_{m-r+1,m-r})$ and $\beta _{m-r+k+1}$. In the same way, intersections between $s^{-k}(E_{m-r+1,m-r})$ and any segment of $\beta _{m-r+k+1}$ arise as preimages of an intersection between a path $\gamma $ whose first $S^{1}$-crossing is at $e^{2\pi i(2/7)}$ and a component $E'$ $s^{-i}(E_{m-r+1,m-r})$ for some $i\leq k$. We consider $(\sigma _{\gamma }\circ s)^{-1}(E')$. This has two components, both close to $s^{-1}(\beta _{2/7})$. We then consider the backward orbits as before, as in (\ref{8.9.1}) to (\ref{8.9.10}). The backward orbits are on either side of the backward orbits of $s^{-1}(\beta _{2/7})$ and   we see that preimages which come close to $\xi _{m,0}(\partial 'U(\gamma _{3},\gamma _{4}))$ -- given by (\ref{8.9.5}), (\ref{8.9.13}), (\ref{8.9.14}) and (\ref{8.9.15}) as before --- either inside, or outside, $\xi _{m,0}(\partial 'U(\gamma _{3},\gamma _{4}))$, and intersections are avoided.

This completes the proof that the support of $\xi _{m,1}\circ \xi _{m,0}^{-1}$ intersects exactly those sets $\xi _{m,0}(\partial 'U(\gamma _{3},\gamma _{4}))$ that it is required to. We still have to show that the support of $\xi _{m}\circ \xi _{m,1}^{-1}$ does not intersect $\xi _{m,1}(\partial 'U(\gamma _{1},\gamma _{2}))$.  It suffices to show that $E_{m,m-r}$ does not intersect any set $\xi _{m,1}(\partial 'U(\gamma _{3},\gamma _{4}))$ for any pair $(\gamma _{3},\gamma _{4})$ with $w_{1}'(\gamma _{3})=w_{1}'(\beta _{k_{1},1})$, for  $m-r>k_{1}$. As in the previous case, we consider intersections with $s^{1-r}(E_{m-r+1,m-r})$. As before, it is certainly true that  $s^{1-r}(E_{m-r+1,m-r})$ does not intersect $\xi _{m,1}(\partial 'U(\gamma _{3},\gamma _{4}))$ . So we need to show that $s^{-r}(\beta _{m-r})$ does not intersect $\xi _{m,1}(\partial 'U(\beta _{5},\beta _{6}))$. So we use the same analysis as in (\ref{8.9.1}) to (\ref{8.9.4}). Intersections with backward orbits of sets as in (\ref{8.9.1}) and (\ref{8.9.2}) are now avoided, precisely because of the definitions of sets $U(\beta _{5},\beta _{6})$ for adjacent pairs $(\gamma _{3},\gamma _{4})$ for which $w_{1}'(\gamma _{3})=w_{1}'(\beta _{1})$. Intersections with backward orbits of sets as in (\ref{8.9.3}) and (\ref{8.9.4}) are avoided as before. So there are no intersections, as required.\Box

\section{Unsatisfying the Inadmissibility Criterion of 
\ref{4.4}}\label{8.21}

To complete the proof of the Main Theorem \ref{7.8}, we need to show 
that, for all $m$ and $p$, the paths of $R_{m,p}$ and $R_{m,p}'$ do not satisfy the 
Inadmissibility 
Criterion of \ref{4.4}, so that they lie in $D'$. We have not defined the paths in $R_{m,p}$ and $R_{m,p}'$ at all precisely, but the only properties we shall use is that the the sets $U(\beta _{1},\beta _{2})$, for adjacent pairs $(\beta _{1},\beta _{2})$ in $R_{m,p}$, lie in $U^{p}$, and the first $S^{1}$-crossing of any path in $R_{m,p}$ is at $e^{2\pi ix}$ for $x\in [q_{p},q_{p+1})$. 
First we consider the case of $R_{m,p}$. We need to show the 
nonexistence of $(\alpha, Q)$ as in the 
Inadmissibility Criterion. As pointed out in \ref{4.4}, this means 
that we have to analyse the lamination $L(\beta )$, for $\beta \in 
R_{m,p}$, and show that period $2$ and period 
$3$ leaves in $L(\beta )$ cannot combine to form a set $Q$. So let 
$\beta \in R_{m,p}$. Then the only intersections between 
$\beta $ and $s^{-j}\beta $ for $j\leq 3$ are between the first 
segment on $\beta $ --- that is, the segment nearest to 
$v_{2}=\infty $ --- and the diagonal segment on $s^{-1}(\beta )$. 
This 
is because the first intersection of $\beta $ with $S^{1}$ is at 
$e^{2\pi ix}$ for $x\in [{2\over 7},{1\over 3})$, and all subsequent 
intersections are such that if $I\in S^{1}$ is an interval with 
endpoints in common with a component of $\beta \cap \{ z:\vert z\vert 
>1\} $, then $\beta\cap s^{-j}I=\emptyset $ for $1\leq j\leq 3$. So 
the only 
leaves of $L(\beta )$ of period $\leq 3$ are the same as those for 
$L_{3/7}$ or $L_{3/7}\cup L_{6/7}^{-1}$. So the Inadmissibility 
Criterion is not satisfied, as required. 

 The 
proof that the paths in $R_{m,p}'$ lie in the admissible domain $D'$ is also 
indirect, and uses induction. To start the induction, $\beta 
_{5/7}\subset D'$. It suffices to 
show that if $\beta _{1}'\in R_{m+1,p}'$ is adjacent to $\beta 
_{0}'\in 
R_{m,p}'$ and $\beta _{0}'\subset D'$, then $\beta _{1}'\subset D'$ 
also. So suppose given such an adjacent pair $(\beta _{0}',\beta 
_{1}')$ paired with adjacent pair $(\beta _{0},\beta _{1})$ in 
$R_{m,p}$, and suppose that the matching is given by $(\psi 
,\alpha )$, that is,  $\alpha \in \pi _{1}(\overline{\mathbb 
C}\setminus Z_{m+1})$, and $[\psi ]\in {\rm{Mod}}(\overline{\mathbb 
C},Y_{m+2})$ with
$$(s,Y_{m+2})\simeq _{\psi }(\sigma _{\alpha }\circ s,Y_{m+2}),$$
$$\alpha *\psi (\beta _{i}')=\beta _{i}{\rm{\ rel\ }}Y_{m+1},\ \ 
i=0,\ 1.$$
Then since $\beta _{0}$ and $\beta _{0}'\subset D'$, the element 
$\gamma $ of $\pi _{1}(B_{3,m+1},s)$ with $\rho (\gamma )=\alpha $ and  
$\Phi _{2}(\gamma )=[\psi ^{-1}]$ must actually be in $\pi 
_{1}(V_{3,m+1},s)$ and preserves $D'$, so that $\beta _{1}'\subset 
D'$.\Box
\chapter{Open Questions}\label{9}

\section{Continuity of the Model}\label{9.1}

We end this paper by considering two lines of enquiry which arise. I have given little thought to the first of these, and have certainly made no progress myself, although others have made interesting advances in related research. Some progress has been made on the second line of enquiry

The fundamental domains that have been constructed for $V_{3,m}$ give 
rise to an inverse limit $\lim _{m\to \infty }M_{m}$, which provides a 
model for the closure in 
$V_{3}$ of the union of all type II and type III hyperbolic 
 components within $V_{3}$. This union is not dense in $V_{3}$, 
because 
 it omits all type IV hyperbolic components. The hope is that a 
 quotient of it is a 
 model for the complement of the type IV hyperbolic components.
 
 The definition of the inverse limit uses the fact that the paths of 
 $R_{m}(.)$ are defined up to isotopy preserving $Z_{n}(s_{.})$ for 
 all $n$. The set $M_{m}$ will be a union of sets 
 $$M_{m}(a_{0}),\ M_{m}(\overline{a_{0}}),\ M_{m}(a_{1},-),\ 
 M_{m}(a_{1},+),$$
with some identifications on boundaries, but otherwise these sets are 
disjoint. We can form sets $U(\beta _{1},\beta _{2})$ for $(\beta 
_{1},\beta 
 _{2})$ an adjacent pair in any of the sets $R_{m}(a_{0})$, 
$R_{m}(\overline{a_{0}})$, $R_{m}(a_{1},-)$, exactly as for 
$R_{m}(a_{1},+)$, except that there is a choice of which of the two 
pairs in a matching pair of adjacent pairs to take. This does not 
matter, because in these cases, if $((\beta _{1},\beta _{2}),(\beta 
_{1}',\beta _{2}'))$ is a matching pair of adjacent pairs, then the 
boundary of the region $U(\beta _{1},\beta _{2})$ contains both 
$\overline{\beta _{1}}*\beta _{2}$ and $\overline{\beta _{1}'}*\beta 
_{2}'$, up to homotopy. This region can equally well be called 
$U(\beta _{1}',\beta _{2}')$. Then similarly to \ref{8.1}, we can 
define $V(\beta _{1},\beta _{2},n)$ to be the union of all sets 
$U(\beta _{3},\beta _{4})$ with $(\beta _{3},\beta _{4})$ an adjacent 
pair in $R_{n}(.)$ between $\beta _{1}$ and $\beta _{2}$, where 
$(\beta _{1},\beta _{2})$ is an adjacent pair in $R_{m}(.)$ for some 
$m\leq n$. Then, using the notation of \ref{4.5} for $\beta _{1/7,c_{1}}$ and so on, we define
 
$$M_{m}(a_{0})=V(\beta _{1/7,c_{1}},\beta _{2/7,s(v_{1})},m)/\sim ,$$ 
$$M_{m}(\overline{a_{0}})=V(\beta _{6/7,c_{1}},\beta 
_{5/7,s(v_{1})},m)/\sim ,$$
$$M_{m}(a_{1},-)=V(\beta _{1/7},\beta _{6/7},m)/\sim .$$
In each case, the relation $\sim $ is simply the lamination 
equivalence relation for $L_{q}$, $q={1\over 7}$, ${6\over 7}$, 
${3\over 7}$. We also define
$$M_{m}(a_{1},+)=\cup \{ (U(\beta _{1},\beta _{2})\times \{ (\beta _{1},\beta _{2})\} )/\sim :(\beta _{1},\beta _{2}){\rm{\ adjacent\ in\ }}\cup _{p\geq 0}R_{m,p}\} ,$$
where, in this case, $\sim $ is still to be determined completely, but for all pairs except  when $\beta _{1}=\beta _{2/7}$, the relation $\sim $ identifies the solid and dashed boundary of $U(\beta _{1},\beta _{2})$. In the case of $\beta _{1}=\beta _{2/7}$, we only identify parts of the solid and dashed boundary nearer the endpoint of $\beta _{2}$.

Now we need to determine how to glue together the different sets $M_{m}(a_{0})$ and $M_{m}(\overline{a_{0}})$ and $M_{m}(a_{1},\pm )$.
Parts 
of the boundaries in the respective cases are 
$$\beta _{1/7,c_{1}}\cup \beta _{2/7,s(v_{1})}
\cup \beta _{4/7,c_{1}}\cup \beta _{4/7,s(v_{1})},$$
$$\beta 
_{6/7,c_{1}}\cup \beta 
_{5/7,s(v_{1})}\cup \beta _{3/7,c_{1}}\cup \beta _{3/7,s(v_{1})},$$
$$\beta _{1/7}\cup \beta _{6/7},$$
and parts of $\beta _{2/7}$ and $\psi _{m,2/7}(\beta _{5/7})$ which form parts of the solid and dashed boundary of $U(\beta _{2/7},\beta _{2})$.
These should be identified.  We identify  the free part of $\beta _{2/7}$ --- that is, the part which forms part of the solid boundary of $U(\beta _{2/7},\beta _{2})$, and which has not been identified with dashed boundary of $U(\beta _{2/7},\beta _{2})$ --- with 
part of $\beta _{4/7,s(v_{1})}$. Similarly, we identify part of $\psi _{m,2/7}(\beta 
_{5/7})$ with part of $\beta _{3/7,s(v_{1})}$. These identifications 
are not so important. There is a natural map from  $M_{m}(.)$ to 
$M_{m+1}(.)$ 
defined as follows. Let $((\beta _{1},\beta _{2}),(\beta _{1}',\beta 
_{2}'))$ be a matching pair of adjacent pairs in $R_{m}(.)$. 
In the case of $a_{0}$, $\overline{a_{0}}$ or  
$(a_{1},-)$, the sets $U(\beta _{1},\beta_{2})$ and $V(\beta 
_{1},\beta _{2},m+1)$ have the same boundary and we choose a map from 
$M_{m}(.)$ to $M_{m+1}(.)$ which is the identity on this boundary. In 
the case of $(a_{1},+)$ we take $(\beta _{1},\beta _{2})$ to be an 
adjacent pair in $R_{m,p}$, for some $p$, matched with $(\beta 
_{1}',\beta _{2}')$ from $R_{m,p}'$. Then we map   $V(\beta 
_{1},\beta 
_{2},m+1)$ to $U(\beta _{1},\beta _{2})$, by taking the identity on 
$\overline{\beta _{1}}*\beta 
_{2}$ and $\xi _{m}^{-1}$ on $\psi 
_{m+1}(\overline{\beta _{1}'}*\beta _{2}')$. Then the inverse 
limits 
$$\lim _{m\to +\infty }M_{m}(a_{0}),\ \lim _{m\to +\infty 
}M_{m}(\overline{a_{0}}),\ \lim _{m\to +\infty }M_{m}(a_{1},-),$$
$$\lim _{m\to +\infty }M_{m}(a_{1},+)$$
are well-defined, with neighbourhood base given by sets $V(\beta 
_{1},\beta _{2},m)$. The analogue for this inverse limit space in the 
case of quadratic polynomials is simply the complement in the 
complex 
plane of the open unit disc. The combinatorial model for the 
Mandelbrot set complement and boundary is a quotient space of this, 
and the combinatorial model for the Mandelbrot set itself is a 
quotient of the closed unit disc. The big question is then the 
following.

\medskip

{\em{Question 1.}} Is  a quotient of $\lim _{m\to \infty 
}M_{m}(.)$ homeomorphic 
    to the closure of the union of type II and type III hyperbolic 
    components in $V_{3}(.)$?

\medskip

This question has varying degrees of difficulty. The corresponding 
question for the Mandelbrot set is answered relatively easily in one 
direction. It is known that the combinatorial model for 
the Mandelbrot set is a continuous image of the true Mandelbrot set, 
because there is a continous monotone map on the union of all type IV 
hyperbolic components and parameter rays of rational argument, with 
dense image 
in the combinatorial Mandelbrot set and this map has 
a unique continuous monotone extension. The difficulty lies in 
showing that 
this map is injective, and hence a 
homeomorphism. The parameter rays in the complement of the Mandelbrot 
set are the images of radial lines under the uniformising map from 
the exterior of the unit disc in the complex plane, to the exterior 
of the Mandelbrot set. The big Mandelbrot set question is then equivalent to the 
question of whether limits along rays exist. These parameter rays do 
have an analogue in $V_{3}(a_{0})$, $V_{3}(\overline{a_{0}})$. Type 
III components 
in $V_{3}(a_{0})$ are connected together in a particularly simple 
way, mimicking connections between between Fatou components of 
$h_{a_{0}}$, excluding those Fatou components in the ``forbidden 
limb''.  So the 
analogue, in $V_{3}(a_{0})$, of the parameter rays in the complement 
of the 
Mandelbrot set, is a collection of paths in the union of type III 
hyperbolic components and connecting points between two hyperbolic 
components, which extend through infinitely many hyperbolic 
components, with no backtracking, 
to limit points in the closure of the union. So these are  
limits of rays $\omega _{a}$, where $\omega _{a}$ are the paths first 
mentioned in \ref{4.5}, and defined there in terms of image under 
$\rho (\omega _{a},s_{1/7})$. The result about parameter rays of 
rational argument having limits then has a possible analogue given by 
the following question.

\medskip

{\em{Question 2.}} Does the 1-1 correspondence 
${\rm{endpoint}}(\rho (\omega _{a},s_{1/7}))\mapsto {\rm{endpoint}}(\omega _{a})$
extend continuously to  paths with periodic endpoints?

\medskip

The analogue of this for $V_{2}$ is proved in \cite{Asp-Yam}. I think the analogue for  $V_{3}(a_{0})$ and $V_{3}(\overline{a_{0}})$ should be reasonably straightforward.  The question for $V(a_{1},-)$ looks harder. There is an 
analogue for $V(a_{1},+)$ and that looks harder still. 

Related to continuity questions, but probably considerably easier, 
there are questions about type IV hyperbolic 
components. There has been no attempt in this paper to describe type 
IV 
hyperbolic components, because the method uses the results of 
\cite{R3}, and that theory was not developed in the type IV case 
because of considerable technical difficulties. But it may well be 
possible to describe all type IV hyperbolic components by using the type 
III ones. In fact, this was the technique used in \cite{R1}, where 
the idea used was that a parabolic parameter value on the boundary of 
a type IV 
hyperbolic component is approximated by a sequence of type III 
hyperbolic centres. One can ask if the corresponding sequence in the 
combinatorial model converges. If so, one would surely have a 
description of the type IV centres up to Thurston equivalence. 
More generally, one can ask:

\medskip

{\em{Question 3.}} Can a complete description be given of the 
    type IV hyperbolic components in $V_{3}$, and of their centres up 
    to Thurston equivalence, by using the description of the type III 
    hyperbolic components?

\medskip

Question 3 is 
related to giving a complete description of all type IV hyperbolic 
components. As indicated in \ref{2.8}, it is known that all type IV hyperbolic components in 
$V_{3}(a)$ for  $a=a_{0}$ or $\overline{a_{0}}$, are matings. I assume that all type IV 
hyperbolic components are matings in
 $V(a_{1},-)$ also. A complete proof would need a a generalisation 
of  
the main result of  \cite{R3}, and of the present paper, 
in order to deal with type IV centres, 
or some circumnavigation, which might well be possible, and certainly 
desirable.   In order 
to describe completely all type IV hyperbolic components in 
$V_{3}(a_{1},+)$, one would probably first have to answer the 
following purely combinatorial question, which should not be too hard.

\medskip

{\em{Question 4.}} What are all possible limits of paths in 
    $\cup _{m,p}R_{m,p}$ such that the limits have periodic endpoints? 
    Is the number of paths with endpoint of period $n$ 
precisely twice the 
     number of type IV hyperbolic components of period $n$ in 
    $V_{3}(a_{1})$?
 
\medskip

If the 
answer to question 2 is ``yes'' for periodic rays, then the answer to  question 5 is also 
``yes'', by the same argument as for the Mandelbrot set. The approximate analogue to this question for the Mandelbrot set is the known fact that the Mandelbrot set maps continuously onto its combinatorial model, by a map which sends hyperbolic components to their analogues in the model. (This known fact is an analogue of the result that a continuous monotone map from the rationals to a subset of $\mathbb R$ with dense image has a unique continuous extension to a homeomorphism from $\mathbb R$ to $\mathbb R$.)

\medskip 

    {\em{Question 5.}}  Does the map 
    $\omega _{a}\mapsto \rho (\omega _{a},s_{1/7})$ 
    extend continuously to a map of the closure of the union all 
    type III components in $V(a_{0})$ onto a quotient of 
    $\lim _{m\to \infty }M_{m}(a_{0})$?
  
\medskip

A positive answer to question 2 in the case of $V(a_{1},-)$ would not 
so easily yield a positive answer to question 5 in the case of 
$V(a_{1},-)$, because there is no clear definition of parameter rays 
in  $V(a_{1},-)$. Interesting work has been done recently by Timorin \cite{Tim2} on the boundaries of type III components, away from intersection points with boundaries of type II components. Away from such intersection points, he has described all boundary points up to topological conjugacy. His topolological models are described a little differently, but certainly related to those used in \cite{R1}.  

 A Markov partition which has come to be known as the Yoccoz puzzle 
has been an outstanding tool towards showing that the Mandelbrot set 
is 
homeomorphic to its combinatorial model, although the problem is not 
usually phrased in this way.  This partition 
is a Markov partition for all polynomials in a {\em{limb}} of the 
Mandelbrot set. It pulls back to give successively finer Markov 
partitions on nonrenomalisable polynomials in the limb, and is a 
generating partition for all the nonrenormlisable polynomials. The 
refined partitions are not topologically or combinatorially  the same 
for all 
polynomials, but are controlled by the {\em{critical puzzle pieces}}, 
which occur at infinitely many levels. A partition of parameter space 
is obtained by partitioning according to the combinatorial nature of 
the refined partition up to level $n$, and in particular, of the 
critical puzzle pieces up to level $n$. So a Markov partiton common 
to a set of polynomials gives rise to a partition of parameter space 
which Yoccoz showed to be a neighbourhood base at each 
nonrenormalisable point (\cite{H, Ro3}).  Numerous refinements 
of this result have 
been 
obtained over the past 15 years or so e.g. \cite{L3}. There have 
also, as already 
mentioned, been analogous results for cubic polynomials, in the 
complement of the connectedness locus, and in nonArchimedean 
dyanmics, 
and in other cases (\cite{Ro1, Ro2, Pe}).
Recently, Lyubich and Kahn and
collaborators developed a Quasi-additivity Law \cite{L1} to extend the methods to unicritical 
polynomials  (\cite{L2, K-L3}), 
a result which was, until now, elusive, and have also made other extensions (\cite{Ka, K-L1, K-L2}) which appear to be very important.

Markov partitions which persist arise in each of the four parts of 
$V_{3}$ which we have studied, $V_{3}(a_{0})$, 
$V_{3}(\overline{a_{0}})$,  and $V(a_{1},\pm)$, associated to the 
polynomials $h_{a}$ for $a=a_{0}$, $\overline{a_{0}}$, $a_{1}$. As already mentioned in \ref{2.8}, a Yoccoz puzzle for $V_{2}$, introduced nonrigorously in Luo's thesis \cite{Luo} was used by Aspenberg-Yampolsky \cite{Asp-Yam} and Timorin \cite{Tim} to prove essentially the same  results about nonrenormalisable points that were proved by Yoccoz for the Mandelbrot set. It is likely that $V(a_{0})$ and $V(\overline{a_{0}})$ can be treated in the same way. Recent discussion suggests that there is a chance that $V_{3}(a_{1},-)$ will also be eligible to such treatment. However, even  
in the simpler cases, the Markov partitions sketched out here are of a somewhat different nature to the Yoccoz puzzle. So if the answers to the following questions are essentially contained in the work of \cite{Asp-Yam, Tim}, this needs checking, to say the least. For concreteness we consider 
the case of $V_{3}(a_{0})$.  An obvious question in the reverse 
direction of question 5, and a generalisation of question 2, is:

\medskip

{\em{Question 6.}} For $a=a_{0}$ or $\overline{a_{0}}$,  does the map 
    $\rho (\omega _{a},s_{1/7})\mapsto \omega _{a} $ 
    extend continuously to a map of $\lim _{m\to \infty 
}M_{m}(a)$ onto 
    the closure of the union all 
    type III components in $V(a)$? What about the case of $\lim _{m\to \infty }M_{m}(a_{1},-)$ and $V(a_{1},-)$?

\medskip

    Intriguingly, a rather similar-looking question to question 6 
about the parameter 
space of cubic polynomials has been answered by Kiwi and others. The 
techniques used in that case appear to be unavailable here, because 
in the cubic polynomial case there is an unbounded Fatou component in 
which one can draw rays of rational argument and show that these 
extend continuously. Rays of rational argument figure in the Yoccoz 
puzzle too of course. So perhaps something can be done in the case of $a=a_{0}$ or $\overline{a_{0}}$, and even in the case of $(a_{1},-)$.

The corresponding question for  $V(a_{1},+)$ is harder to state.  We
use the sets $V(\beta _{1},\beta _{2},n)$ 
of \ref{8.1} and \ref{8.17} for adjacent pairs $(\beta _{1},\beta _{2})$ in $R_{m,0}$. These sets are not all disjoint, but sets $V(\beta _{1},\beta _{2},n)$ with the same first $S^{1}$-crossing are disjoint. We have developed the theory only in the case of $R_{m,0}$, but a similar procedure can be carried out in $R_{m,p}$ for any $p\in \mathbb Z$ with $p\geq 0$.

One can then ask:

\medskip

{\em{Question 7.}} 
Let $(\beta _{1},\beta _{2})$ be adjacent in $R_{m,0}$, or more generally, in $\cup _{p}R_{m,p}$. Does the natural map from 
    $$V(\beta _{1},\beta _{2},m)\cap 
    (\cup _{n}Z_{n}(s))$$
     to centers of hyperbolic 
    components in $V_{3}$, given by mapping endpoints of paths $\rho 
    (\omega _{a})$ to endpoints of paths of $\omega _{a}$,
    extend continuously to \\ $\lim _{n\to \infty }V(\beta _{1},\beta _{2},n)$?

\medskip
    
   {\em{Question 8.}} If the answer to Question 6 is ``yes'', does 
    the diameter image of a set $\lim _{n\to \infty }V(\beta _{1},\beta _{2},n)$ tend to $0$ 
uniformly with $n$? (Remember that $\beta _{1}$ and $\beta _{2}$ are assumed adjacent in $\cup _{p}R_{m,p}$.)
    
\medskip

 These questions are of course, intertwined and there are 
partial implications in both directions. But it does not seem 
worth spelling out implications without more idea on how to tackle 
such questions. 

The famous quadratic polynomial analogue of the following, has not been 
answered but would follow from homeomorphism between the Mandelbrot 
set and its combinatorial model:

\medskip

{\em{Question 9.}}  Is the complement in $V_{3}$ of the 
closure of the union of the type III 
hyperbolic components exactly equal to the union of type IV 
hyperbolic components?

\medskip

    \section{Some nontrivial equivalences between 
captures}\label{9.2}

 It will have been noticed that the part of the  fundamental domain corresponding to $V_{3}(a_{1},+)$ has 
relatively 
few capture paths in it. In fact the proportion of such paths, for 
each 
$m$, is exponentially small (although growing exponentially with 
$m$).  I suspect that it is impossible to choose a fundamental domain with a 
significantly higher proportion of capture paths, although I cannot 
be 
sure.  A rather elementary 
point is that captures $\sigma _{r}\circ s_{3/7}$ and $\sigma 
_{1-r}\circ s_{3/7}$ are usually not Thurston equivalent, even when 
$\beta _{r}$ and $\beta _{1-r}$ end at the same point of $Z_{m}(s)$. 
An example which was mentioned in \ref{7.5} is given by $r={65\over 224}$, so that $\beta _{r}=\beta 
(L_{3}L_{2}R_{3}L_{3}L_{2}C)$. Then $\sigma _{1-r}\circ s_{3/7}$ is 
Thurston 
equivalent to $\sigma _{\beta '}\circ s_{3/7}$, where $\beta '=\beta 
(L_{3}^{4}L_{2}C)$, and hence not Thurston equivalent to $\sigma 
_{r}\circ s_{3/7}$, because it can only be equivalent to one point in 
the fundamental domain.

It is also clear that Thurston equivalences between 
different captures can be quite complicated. Here is an example. We claim that the captures 
$\sigma _{r}\circ s_{3/7}$ and $\sigma _{q}\circ s_{3/7}$ are 
Thurston equivalent, where $r={271\over 896}$ and $q={635\over 896}$. 
Note that $r+q\neq 1$. In fact, $1-r={625\over 896}$ and $\sigma 
_{1-q}\circ s$ and $\sigma _{1-r}\circ s$ are also Thurston 
equivalent.  We have $e^{\pm 2\pi ir}\in \partial 
D(L_{3}(L_{2}R_{3})^{2}L_{3}L_{2}C)$ and $e^{\pm 2\pi iq}\in \partial 
D(L_{3}L_{2}R_{3}L_{3}^{3}L_{2}C)$. Write 
$w=L_{3}(L_{2}R_{3})^{2}L_{3}L_{2}C$
 and $w'=L_{3}L_{2}R_{3}L_{3}^{3}L_{2}C$. Now $w_{1}(w',0)=w'$ and $\psi _{7,2/7}$ maps $e^{\pm 2\pi iq}$ to  $e^{\pm 2\pi ir}$. Although $\psi _{7,2/7}$ does not map $\beta _{q}$ to $\beta _{r}$, $\psi _{7,2/7}(\beta _{q})$ and $\beta _{r}$ are homotopic under a homotopy which fixes endpoints and preserves 
the forward orbit of the common second endpoint: because for any 
suffix $u$ of 
$w$, $D(u)$ is not between $D(w)$ and $D(C)$. So $\sigma _{q}\circ s_{3/7}$ and $\sigma _{r}\circ s_{3/7}$ are Thurston equivalent.

A question that remains is: how to identify all captures, up to 
Thurston equivalence, maps in the fundamental domain. Vertices of the fundamental domain in the three ``easy'' parts of the fundamental domain are represented by captures in any case. Easy identifications were described in \ref{2.8}. So we already know how the following captures are represented in the fundamental domain:
$$\sigma _{r}\circ s_{1/7}{\rm{\ for\ }}r\in (2/7,8/7),$$
$$\sigma _{r}\circ s_{6/7}{\rm{\ for\ }}r\in (-1/7,5/7),$$
$$\sigma _{r}\circ s_{3/7}{\rm{\ for\ }}r\in (-1/7,1/7).$$
 So the case to consider is the ``hard'' case:
 $$\sigma _{r}\circ s_{3/7}{\rm{\ for\ }}r\in (2/7,1/3)\cup (5/7,2/3).$$
 This can be done 
on a case-by-case basis, simply because, given a fundamental domain 
$F$ 
in the unit disc $D$,  for a free group of hyperbolic isometries with 
all 
vertices of $F$ on $\partial D$, one can determine, for each vertex 
of a translate of $F$, the pair of vertices of $F$ to whose orbit it 
belongs.  Now $\beta _{r}$ 
lifts to the  set $D'$ in the unit disc. Identifying $D'$ with the 
universal cover of $V_{3,m}$, as the Resident's View says we can do, 
the lift of $\beta $ to $D'$ has second endpoint at a lift of a 
puncture of $V_{3,m}$, and therefore the lift of $\beta $ is contained 
in finitely many translates of the fundamental domain. Now $\sigma 
_{r}\circ s_{3/7}$ must be Thurston equivalent to one of the 
maps represented in the fundamental domain. A priori, this map could be in $V_{3,m}(a_{0})$ or $V_{3,m}(\overline{a_{0}})$, or $V_{3,m}(a_{1},\pm )$. However, it appears to be the case  that the map is in $V_{3,m}(a_{1},+)$, and is therefore of the form
$$\sigma _{\zeta }\circ s_{3/7}$$
for  a unique $\zeta =\zeta (r)\in \cup _{p}R_{m,p}$. It also seems to be true   that if $r\in (q_{p},q_{p+1}]\cup [1-q_{p+1},1-q_{p})$ then $\zeta (r)\in R_{m,p}$. It was pointed out in \ref{2.8} that the set of $r\in (2/7,1/3)\cup (2/3,5/7)$ giving preperiod $m$ captures $\sigma _{r}\circ s_{3/7}$, after reducing to exclude obvious Thurston equivalences, has cardinality asymptotic to $\frac{1}{21}2^{m}$, just half the asymptotic number of preperiod $m$ type III hyperbolic components in $V_{3,m}(a_{1},+)$. 

In \ref{3.4} we mentioned Jonathan Woolf's simple recurrence formula for the mean image size 
$$n_{2}\left( 1-\left( 1-\frac{1}{n_{2}}\right) ^{n_{1}}\right) $$
of a map from a set $A$ of $n_{1}$ elements to a set $B$ of $n_{2}$ elements. We are, of course, interested in the case when $n_{1}/n_{2}\to 1/2$. Much more detailed information can be obtained by using the saddle point method. It can be shown, for example, that, if $n_{2}\to \infty $ and $n_{1}/n_{2}\to a$, the  ratio of the image size to $n_{2}$ tends to $1-e^{-a}$ with probability one, and a limiting normal distribution of image size can be obtained.
It can also be shown that the  the mean of the maximal inverse image size for maps $f$ from $A$ to $B$:
$${\rm{Max}}(\{ \# (f^{-1}(b)):b\in B\} ),$$
is at least $c\log n_{2}/\log \log n_{2}$ for any $c<1$, with probability tending to one as $n_{2}\to \infty $, if $n_{1}/n_{2}\to a$.
This raises the following questions, concerning the specific map from the set  of capture maps to maps in the fundamental domain:

 \begin{itemize}
    \item Can the proportion of preperiod $m$ captures, up to 
Thurston 
equivalence, in $V_{3,m}(a_{1},+)$, be computed? If so is it
$c2^{m}(1+o(1))$ for some $c>0$ which can be determined?
\item For each integer $N$, are there $N$ different rational 
numbers $r_{j}\in ({2\over 7},{1\over 3})\cup ({2\over 3},{5\over 
7})$ 
such that the captures $\sigma _{r_{j}}\circ s_{3/7}$ are all 
Thurston 
equivalent?
\item Are there computable asymptotes in $m$ for the average 
number of 
 captures $\sigma _{r}\circ s_{3/7}$ 
which are Thurston-equivalent to a preperiod $m$ capture $\sigma 
_{q}\circ 
s_{3/7}$, 
for $q\in ({2\over 7},{1\over 3})\cup ({2\over 3},{5\over 7})$?
\end{itemize}

I believe that:
\begin{itemize}
\item  there is a constant $c$ as in the first question, and it can probably be estimated, at least;
\item the answer to the second question is ``yes'': see \cite{R5}.

\item qualifying this ``yes'', the {\em{average }} number of captures in a Thurston equivalence class is probably boundedly finite with probability tending to one as $m\to \infty $.
\end{itemize}

\backmatter



\begin{thebibliography}{99}

\bibitem{Asp-Yam} Aspenberg, M. and Yampolsky, M., Mating non-renormalisable quadratic polynomials, Comm. Math. Phys. 287 (2009), 1-40.

\bibitem{L2} Avila, A., Kahn, J., Lyubich, M. and Shen, W., 
Combinatorial Rigidity for Unicrtical Polynomials, http://front.math.ucdavis.edu/math.DS/0507240, to appear in Ann. of Math..

\bibitem{Bar} Bartholdi, L. and Nekrashevych, V., Thurston equivalence of topological polynomials. Acta Math., 197 (2006), 1-51.

\bibitem{Bow1} Bowditch, B., Geometric models for hyperbolic 3-manifolds. http://www.warwick.ac.uk/~masgak/preprints.html

\bibitem{Bow2} Bowditch, B., End invariants of hyperbolic 3-manifolds. 
http://www.warwick.ac.uk/~masgak/preprints.html

\bibitem{Bow3} Bowditch, B., Length bounds on curves arising from tight geodesics. Geom. Funct. Anal.17 (2007), 1001-1042.

\bibitem{BBES} Brock, J., Bromberg, K., Evans, R., Souto, J., In preparation.

\bibitem{B-H1} Branner, B. and Hubbard, J.H., The iteration of
cubic polynomials.  Part I: The global topology of parameter space. 
Acta Math., 160 (1988), 143-206.

\bibitem{B-H2} Branner, B. and Hubbard, J.H., The iteration of cubic 
polynomials. Part II: Patterns and parapatterns. Acta Math., 
169 (1992), 229-325. 

\bibitem{B-C-M} Brock, J., Canary, R. and Minsky, Y., The classification of Kleinian surface groups II: The Ending Laminations Conjecture. Preprint 2004. http://www.math.yale.edu/users/yair/research
    
\bibitem{D-H1} Douady, A. and Hubbard, J.H., Etudes 
dynamiques des polyn\^omes complexes, avec la 
collaboration de P. Lavaurs, Tan Lei, P. 
Sentenac. Parts I and II, Publications 
Math\'ematiques d'Orsay, 1985.

\bibitem{D-H2} Douady, A. and Hubbard, J.H., A proof 
of Thurston's topological characterization of 
rational functions. Acta Math., 171 (1993), 263-297.

\bibitem{AE} Epstein, A.L., Bounded hyperbolic components of 
quadratic 
rational maps. Ergod. Th. Dyn. Sys., 20 (2000), 727-748.

\bibitem{H} Hubbard, J.H., Local connectivity of Julia sets and 
bifurcation loci: three theorems of J.C. Yoccoz, in Topological 
methods 
in modern mathematics (SUNY at Stony Brook, NY 1991), pp 467-511. 
Publish or Perish 
Houston TX 1993.

\bibitem{L1} Kahn, J. and Lyubich, M., The Quasi-Additivity Law in 
Conformal Geometry. http://front.math.ucdavis.edu/math.DS/0505191, to appear in Ann. of Math..

\bibitem{Ka} Kahn, J., A priori bounds for some infinitely renormalizable quadratics: I. Bounded Primitive Combinatorics, http://front.math.ucdavis.edu/math.DS/0609045

\bibitem{K-L3} Kahn, J. and Lyubich, M., Local connectivity of Julia sets for unicritical polynomials, http://front.math.ucdavis.edu/math.DS/0505194, to appear in Ann. of Math..

\bibitem{K-L1} Kahn, J. and Lyubich, M., A priori bounds for some infinitely renormalizable quadratics: II. Decorations, http://front.math.ucdavis.edu/math.DS/0609046, to appear in Ann. Math. Ec. Norm. Sup..

\bibitem{K-L2} Kahn, J. and Lyubich, M., A priori bounds for some infinitely renormalizable quadratics: III. Molecules, http://front.math.ucdavis.edu/0712.2444.


\bibitem{K2} Kiwi, J.,ational Laminations of Complex Polynomials, 
in Laminations and Foliations in Dynamics, Geometry and Topology  
(SUNY 
at Stony Brook, 1998, eds M. Lyubich et al.) Contemporary Math. 269 
(2001), 111-154.

\bibitem{K3} Kiwi, J.,  Real Laminations and the topological dynamics 
of complex polynomials, Advances in Math., 184 (2004), 207-267.

\bibitem{K4} Kiwi, J., Combinatorial continuity in complex polynomial 
dynamics. Proc. London Math. Soc., 91 (2005), 215-248.

\bibitem{K1} Kiwi, J., Puiseux series, polynomial dynamics and 
iteration of complex cubic polynomials. To appear in Annales de l'Institut Fourier.

\bibitem{Luo}  Luo, J., Combinatorics 
and Holomorphic Dynamics: Captures, Matings and Newton's Method. 
Thesis, 
Cornell University, 1995.

\bibitem{L3} Lyubich, M., Dynamics of quadratic polynomials III: 
Parapuzzles and SBR measures, in Geometrie complexe et systemes 
dynamiques (Orsay 1995), pp 173-200. Ast\'erisque 261 (2000).


\bibitem{M1} Milnor, J., Geometry and dynamics of quadratic rational 
maps. Exper. Math., 2 (1993), 37-83.

\bibitem{M2} Milnor, J., Rational maps with two critical points. 
Exper. Math., 9 (2000), 481-522.

\bibitem{Min} Minsky, Y., The classification of Kleinian surface groups I. http://arxiv.org/pdf/math/0302208, to appear in Ann. of Math.. 

\bibitem{Pe} Petersen, C.L., Puzzles and Siegel Disks,  Progress 
in 
holomorphic dynamics, 50-83 Pitman Res. Notes Math. Ser., 387, 
Longman Harlow 1998.

\bibitem{R4} Rees, M., Components of degree two 
hyperbolic rational maps. Invent. Math., 100
(1990), 357-382.

\bibitem{R1} Rees, M., A Partial Description of the 
Parameter space of Rational Maps of Degree Two: 
Part 1. Acta Math., 168 (1992), 11-87.

\bibitem{R2} Rees, M., A Partial Description of the 
Parameter space of Rational Maps of Degree Two: 
Part 2.  Proc. Lond. Math. Soc., 70 (1995), 
644-690.

\bibitem{R3} Rees, M., Views of Parameter 
Space, Topographer and Resident.    
Ast\'erisque, 288 (2003).

\bibitem{R6} Rees, M., The Ending Laminations Theorem direct from Teichm\"uller geodesics.  http://arxiv.org/pdf/math/0404007

\bibitem{R5} Rees, M., Multiple equivalent matings with the aeroplane polynomial, to appear in Ergod. Th. and Dynam. Sys..

\bibitem{Ro1} Roesch, P., On local connectivity for the Julia set of rational maps: Newton's famous example. Annals of Mathematics 168 (2008), 1-48.

\bibitem{Ro2} Roesch, P., Puzzles de Yoccoz pour les application a 
allure rationelle. Enseign. Math., 45 (1999), 133-168.

\bibitem{Ro3} Roesch, P., Holomorphic motions and puzzles (following 
M. Shishikura) in The Mandelbrot Set, Theme and Variations, 
pp117-132. LMS Lecture Note Series 274, Cambridge U.P. 2000.

\bibitem{Ro4} Roesch, P., Hyperbolic components of polynomials with a fixed critical point of maximal order. Ann. Sci. Ec. Norm. Sup. 40 (2007), 1-53.

\bibitem{Sti} Stimson, J., Degree two rational maps with a 
periodic 
critical point. Thesis, University of Liverpool, July 1993.

\bibitem{TL} Tan Lei: Matings of Quadratic Polynomials. 
Ergod. Th. and Dynam. Sys., 12 (1992), 589-620.

\bibitem{T} Thurston, W.P., On the Combinatorics of 
Iterated Rational Maps. Preprint, Princteton 
University and I.A.S., 1985. To appear in the J.H. Hubbard 60th birthday 
volume, Complex dynamics:
families and friends, AK Peters, 2008, 3-109. Editorial notes by D. Schleicher and N. Selinger. 

\bibitem{Tim} Timorin, V., External boundary of $M_{2}$, Proceedings of conference dedicated to Milnor's 75th birthday, Proceedings of the Fields Institute 53 (2006), 225-267.


\bibitem{Tim2} Timorin, V., Topological Regluing of Holomorphic Functions.  Invent. Math. 179 (2010), 461-506.

\bibitem{W} Wittner, B., On the bifurcation loci of
rational maps of degree two.  Thesis, Cornell University, 1988.

\end{thebibliography}


\end{document}